\documentclass[11pt]{memo-l}
\usepackage{amsfonts, amsmath, amssymb, latexsym, pictex, mymacros, eucal}

\numberwithin{equation}{chapter}
  \def \F {\mathbb F}
  
\def \Z {\mathbb Z}
\newcommand{\m}{\medskip}
\newcommand{\bi}{\bigskip}
\newcommand{\n}{\noindent}
\newcommand{\ad}{\hbox{\rm ad\,}}
\newcommand{\al}{\alpha}
\newcommand{\Ann}{\hbox{\rm Ann}}

\newcommand{\be} {\beta}

\newcommand{\csp}{\mathfrak{csp}_{2m}}
\newcommand{\Der}{\hbox {\rm Der}}
\newcommand{\di}{{\tt div}}
 \newcommand{\En}{\hbox {\rm End}}
 \newcommand{\eqdef}{:=}

 \newcommand{\e}{\varepsilon}
 
 \newcommand{\g}{\mathfrak g}
 
 \newcommand{\h}{\mathfrak h}
 \newcommand{\hf}{\frac{1}{2}}
 \newcommand{\Hom}{\hbox {\rm Hom}}
  \newcommand{\id}{\hbox {\rm id}}
 
 \newcommand{\lam}{\lambda}
 \newcommand{\la}{\langle}
 \newcommand{\M} {\mathcal{M}}
 \newcommand{\one}{1}  
 \newcommand{\ot}{\otimes}
 \newcommand{\ov}{\overline} 
 \newcommand{\ra}{\rangle}
  \newcommand{\rphi} {\buildrel{\varphi_i} \over  \rightarrow}

 \newcommand{\spa} {\hbox {\rm  span}}
 \newcommand{\spm}{\mathfrak{sp}_{2m}}
 \newcommand{\tf}{\mathfrak t}
  
 \newcommand{\W}{Witt}
 \newcommand{\un} {\underline}

\newtheorem{Thm}[equation] {Theorem}
\newtheorem{Lem}[equation] {Lemma}
\newtheorem{Cor}[equation] {Corollary}
\newtheorem{Def}[equation] {Definition}
\newtheorem{Pro}[equation] {Proposition}

\newtheorem{pgraph}[equation] {}

\newtheorem{Rem}[equation] {Remark}
\newtheorem{Rems}[equation] {Remarks}

\newcommand{\pf} {\medskip \noindent \textbf{Proof.\  \  }}

\begin{document}
\frontmatter
\title{The Recognition Theorem \\
for Graded Lie Algebras \\
 in Prime Characteristic }
 
\author{Georgia Benkart$^*$}
\address{Department of Mathematics, University of Wisconsin-Madison, Madison, Wisconsin 53706  USA}

\email{benkart@math.wisc.edu}
\thanks{{$^*$}The first author gratefully acknowledges support from  National Science Foundation  grant \#{}DMS--0245082  and from the E.B. Van Vleck Professorship at the University of Wisconsin-Madison.}

\author{Thomas Gregory}
\address{Department of Mathematics, The Ohio State University at Mansfield,
Mansfield, Ohio 44906  USA}
\email{tgregory@math.ohio-state.edu}

\author{Alexander Premet}
\address{School of Mathematics, The University of Manchester, \break
Manchester M13 9PL, UK}
\email{sashap@maths.man.ac.uk} 

\date{\today}

\begin{abstract}
The ``Recognition Theorem''  for graded Lie algebras is an essential ingredient in
the classification of finite-dimensional simple Lie algebras over an algebraically closed field of
characteristic $p > 3$.  The main goal of this monograph is to present the first complete
proof of this fundamental result.
\end{abstract}

\maketitle 

\tableofcontents

%-----------------------------------------------------------------------
% Beginning of intro.tex - 8-19-05 version
%-----------------------------------------------------------------------
%
% AMS-LaTeX 1.2 sample file for a monograph, based on amsbook.cls.
% This is a data file input by chapter.tex.
%%%%%%%%%%%%%%%%%%%%%%%%%%%%%%%%%%%%%%%%%%%%%%%%%%%%%%%%%%%%%%%%%%%%
 
\chapter*{Introduction}

The focus of this work is the following Main Theorem,  often
referred  to as the ``Recognition Theorem,''  \index{Recognition Theorem}  because of its extensive
use in recognizing  certain graded
Lie algebras from their null components.

\bi \begin{Thm} \  Let $\g = \bigoplus_{j= -q}^r \g_j$ be a
finite-dimensional graded Lie algebra over an algebraically closed
field $\mathbb F$ of characteristic $p > 3.$ Assume that:
\begin{itemize}
\item[{\rm (a)}]\  $\g_0$ is a direct sum of ideals, each of which is abelian, 
a classical simple Lie algebra,  or   one of the Lie algebras
$\mathfrak{gl}_n$, $\mathfrak{sl}_n$, or
$\mathfrak{pgl}_n$  with $p\mid n$;

\item[{\rm (b)}]\  $\g_{-1}$ is an irreducible $\g_0$-module;

\item[{\rm (c)}]\   If  $x \in \bigoplus_{j \geq 0} \g_j$ and $[x,\g_{-1}] = 0$,  then
$x = 0$;

\item[{\rm (d)}]\  If  $x \in \bigoplus_{j \geq 0} \g_{-j}$  and $[x,\g_{1}] = 0$,  then
$x = 0$.
\end{itemize}

\noindent Then $\g$ is isomorphic as a graded Lie algebra to one of
the following:

\begin{itemize}

\item[{\rm (1)}]\  a classical simple Lie algebra with a standard grading;

\item[{\rm (2)}]\  $\mathfrak{pgl}_{m}$ for some $m$ such that $p \mid m$  with a standard grading;

 \item[{\rm (3)}]\ a Cartan type Lie algebra with the natural grading or its reverse;

\item[{\rm (4)}]\  a Melikyan\footnote{transliterated as Melikian in many references
such as  \cite{St},  for example} algebra (in characteristic 5)  with either the natural grading or its
reverse.  \end{itemize}   \end{Thm}

The classical simple Lie algebras in this theorem are the algebras obtained by
reduction modulo $p$ of the finite-dimensional complex simple Lie algebras,  as 
described in \cite[Sec.~10]{S}
(see also Section \ref{sec:2.2}).  Thus, they are of type 
A$_{n-1}$, $p \not  |\, n$, \ B$_n$, C$_n$, D$_n$, E$_6$, E$_7$, E$_8$, F$_4$,
G$_2$, or they are isomorphic to $\mathfrak{psl}_n$ where $p\mid n$.  

The Recognition Theorem  is an essential ingredient in the classification of
the finite-dimensional simple Lie algebras over algebraically
closed fields of characteristic $p
> 3$. In a sense, the whole classification theory is built around
this theorem, as  the theory  aims to show that any finite-dimensional simple Lie
algebra $L$ admits a filtration $L=L_{-q}\supset\ldots\supset
L_0\supset\cdots\supset L_r\supset L_{r+1}=0$ such that the
corresponding graded Lie algebra
$\g\,=\,\bigoplus_{i=-q}^r\,\g_i$, where $\g_i = L_i/L_{i+1}$,  satisfies  
conditions (a)-(d) above. The Recognition Theorem is used several
times throughout the classification; its first application results
in a complete list of the simple Lie algebras of absolute toral
rank two, and its last application yields a crucial characterization
of the Melikyan algebras, thereby completing the classification.
The main goal of this monograph is to present the first complete
proof of this fundamental theorem.

V.G.~Kac first undertook to prove the Recognition Theorem in
\cite{K2}.  This pioneering work was ahead of its time.  In 1970, 
very little was known about rational representations of simple
algebraic groups in prime characteristic, and  the Melikyan algebras
were discovered only in the 1980s.  Despite that, Kac made
many deep and important observations towards the proof
of the theorem  in \cite{K2}.  Most of them are incorporated in Chapters~3 and 4
of this monograph.

Historical accounts  of the classification of simple Lie algebras
of characteristic $p > 0$ can be found  for example in \cite{M},
\cite{Wi3}, and \cite{B}.   Investigation of the
finite-dimensional simple Lie algebras over algebraically closed
fields of positive characteristic began in the 1930s  in the work
of Jacobson, Witt,  and Zassenhaus. During the next quarter
century, many examples of such Lie algebras were discovered.   In
\cite{S}, written in 1967, Seligman spoke of a ``rather awkward
array of simple modular Lie algebras which would be totally
unexpected to one acquainted only with the non-modular case.''
Seligman's book contained a characterization of the {\it
classical} Lie algebras of characteristic $p > 3$; that is, those
obtained from $\Z$-forms of the finite-dimensional simple Lie
algebras over $\mathbb C$ by reducing modulo $p$. (See Section
\ref{sec:2.2}.) It was about the same time that Kostrikin and
\v Safarevi\v c \cite{KS} observed a similarity between the known
nonclassical simple Lie algebras of prime characteristic and the
four families $W$, $S$, $H$, $K$ (Witt, special, Hamiltonian,
contact)  of infinite-dimensional complex Lie algebras arising in
Cartan's work on Lie pseudogroups. They called their analogues
``Lie algebras of Cartan type'' and formulated a conjecture which
shaped research on the subject during the next thirty years.  The
\emph{Kostrikin-\v Safarevi\v c Conjecture of 1966} states
\begin{itemize}
\item[]  \emph{Over an algebraically closed field
of characteristic $p > 5$,  a finite-dimensional \, {\tt restricted}
simple Lie algebra is classical or of Cartan type.}
\end{itemize}
\medskip
In 1984, Block and Wilson \cite{BW} succeeded
in proving this conjecture for algebraically
closed fields of characteristic
 $p > 7.$

   If the notion of a Cartan type Lie algebra is expanded to include both the restricted
   and nonrestricted ones  as well as  their filtered deformations
   (determined later by Kac \cite{K3}, Wilson \cite{Wi2}, and Skryabin \cite{Sk1}), 
   then one can formulate the
{\it Generalized Kostrikin-\v Safarevi\v c
Conjecture} by simply erasing the word {\tt restricted}  in the statement
above.

   The Generalized Kostrikin-\v Safarevi\v c
Conjecture is now a theorem for $p > 7$.  First announced by
Strade and Wilson  \cite{SW} in 1991, its proof is spread over a
number of papers. As mentioned above, the proof depends in a
critical way on the above Recognition Theorem (hence on our monograph). We refer the
interested reader to \cite{St} for
 a comprehensive
 exposition of the classification.

 Recent work of Strade and the third author  \cite{PS1}-\cite{PS5}
 has made significant progress on
 the problem of
 classifying  the finite-dimensional simple Lie
algebras over fields of  characteristic 7 and of characteristic 5,
where the Generalized Kostrikin-\v Safarevi\v c Conjecture is
known to fail because of the  Melikyan algebras. (See Section
2.45.)      The  Classification Theorem announced in \cite[p.~7]{St} (and also in \cite{P2}), 
asserts  ``{\it Every finite dimensional simple Lie algebra over an
algebraically closed field of characteristic $p>3$ is of classical, Cartan
or Melikyan type.}''        The Recognition Theorem  
plays a vital role in this extension of the classification theory to $p = 5$ and 7.     
  In characteristics 2 and 3, many more simple Lie algebras
which are neither classical nor Cartan type are known (characteristic 3 examples
can be found in \cite[Sec.~4.4] {St}). The papers \cite{KO}, 
\cite{BKK}, \cite{BGK}, and \cite{GK} prove recognition theorems for graded Lie
algebras of characteristic 3 under various assumptions on the
gradation spaces. One of the main challenges in characteristics
$2$ and $3$ will be to remove such extra assumptions and determine
{\it all} finite-dimensional graded Lie algebras $\g$ satisfying
conditions (b), (c), (d) above with the graded component
$\g_0$ isomorphic to the Lie algebra of a reductive group. Once
this is accomplished,  one might be able to formulate a plausible
analogue of the Generalized Kostrikin-\v Safarevi\v c Conjecture for
$p=2$ and $3$ and to begin the classification work in a systematic way.

Our monograph  consists of four chapters.  In the first,  we establish
general properties of  graded Lie algebras and use them to show
that in a finite-dimensional  graded Lie algebra
satisfying conditions (a)-(d)  of the Main Theorem, the
representation of the commutator ideal $\g_0^{(1)}$ of the null component $\g_{0}$
on $\g_{-1}$ must be restricted.    In Chapter 2, we
gather useful information about known graded Lie algebras,  both classical and
nonclassical.    Chapter 3 deals
with the case in which   $\g_{-1}$
and $\g_{1}$ are dual $\g_0$-modules, the so-called contragredient case,
which leads to the
classical or  Melikyan Lie algebras.
Chapter 4 treats the noncontragredient case, and there the
graded Lie algebras are shown to be of Cartan type.

\bigskip  
\begin{center}{\textbf {Acknowledgments}} \end{center}  \m

Work on this monograph began under the sponsorship
of  National \break Science Foundation US-FSU Cooperative Program
Grant \  $\#$DMS-9115984, which funded the visits of Alexei Kostrikin,
Michael Kuznetsov, and  Alexander Premet to the University of Wisconsin-Madison.
We are grateful for the Foundation's support and for the enthusiasm 
that  Alexei Kostrikin showed in urging
us to undertake this project.   We regret that he did not live to see its
completion.       We  also wish to express our  sincere
thanks to  Richard Block, Michael Kuznetsov, Hayk Melikyan,
Yuri Razmyslov, George Seligman, Serge Skryabin, Helmut Strade, and
Robert Wilson for their interest and encouragement.

\setcounter{page}{1}

\mainmatter 
%-----------------------------------------------------------------------
% Beginning of chap1.tex - 9-5-05 version
%-----------------------------------------------------------------------
%
% AMS-LaTeX 1.2 sample file for a monograph, based on amsbook.cls.
% This is a data file input by chapter.tex.
%%%%%%%%%%%%%%%%%%%%%%%%%%%%%%%%%%%%%%%%%%%%%%%%%%%%%%%%%%%%%%%%%%%%%%%%%%%%%%%%%%%%%%%%%%%%%%%%%%%%%%%%%%%%%%%%%%%%%%%
%%%%%%%%%%%%%%%%%%%%%%%%%%%%%%%%%%%%%%%%%%%%%%%%%%%%%%%%%%%%%%%%%%%%%%%%

\setcounter{page}{1} 
  \chapter{Graded Lie Algebras}  
\bi

\section {{\ Introduction} \label{sec:1.1}}  

\m

 In this chapter,  we develop results about general graded Lie
algebras.  Later (starting in  Section \ref{sec:1.8}) and in
subsequent chapters, we specialize to modular graded Lie algebras
satisfying the hypotheses of the Recognition Theorem.

\m To begin, our focus is on Lie algebras over an arbitrary field
${\F}$ having an integer grading,

$$\g = \bigoplus_{i=-q}^{r} \g_i,$$

\n where $[\g_i,\g_j] \subseteq \g_{i +j}$ if $-q \leq i+j \leq r$
and $[\g_i,\g_j] = 0$ otherwise.  Then $\g_0$ is a subalgebra of
$\g,$ and each subspace $\g_j$ is a $\g_0$-module under the
adjoint action.  The spaces

$$\g_{\leq 0} \eqdef \g_- \oplus \g_0 \quad  \text{and} \quad  \g_{\geq
0} \eqdef \g_0 \oplus \g_+$$

\n are also subalgebras of $\g$, where

$$\g_{-} \eqdef \,\bigoplus_{i = 1}^q \g_{-i}\quad  \text{and} \quad
\g_{+} \eqdef\, \bigoplus_{j = 1}^r \g_j$$

\n are nilpotent ideals of $\g_{\leq 0}$ and $\g_{\geq 0}$
respectively. If $\g_{-q}$ and $\g_r$ are nonzero, then $q$ is
said to be the {\em depth} and $r$ the {\em height} of $\g$.  We
assume that $q,r \geq 1$ and $q$ is finite, but in this section
and the next allow the
possibility that the height $r$ is infinite.  The following
conditions play a key role in this investigation:

\begin{equation}\label{eq:1.2} \g_{-1}~~ \hbox{\rm is~an~irreducible~}~\g_0
\hbox{\rm -module}.  \end{equation}
\begin{equation}\label{eq:1.3} \hbox{\rm If}\ \ x
\in \g_{\geq 0} \ \ \hbox{\rm and } \ \ [x,\g_{-1}] = 0, \ \ \hbox{\rm
then }\ \ x = 0.  \end{equation}

\n Property \eqref{eq:1.2} is termed \ {\em irreducibility} \ and
\eqref{eq:1.3} is called \ {\em transitivity}.  When we say an algebra is
irreducible and transitive, we mean that both \eqref{eq:1.2} and
\eqref{eq:1.3} hold. On occasion we refer to algebras satisfying
the following constraint as being \ $1$-{\em transitive}, or
having \ $1$-{\em transitivity}:

\begin{equation}\label{eq:1.4} \hbox{\rm If} \ \ x \in
\g_{\leq 0} \ \ \hbox{\rm and }\ \ [x,\g_1] = 0, \ \ \hbox{\rm then }
\ \ x = 0.  \end{equation} 

\m
\section  {\ The Weisfeiler radical \label{sec:1.2}} 

 \m Every graded
Lie algebra $\g$ has a radical, which was first introduced by
Weisfeiler [W], and which is constructed as follows: \quad Set
$\M^0(\g) = 0$, and for $i \geq 0$ define $\M^{i+1}(\g)$
inductively by

\begin{equation}\label{eq:1.6} \M^{i+1}(\g) = \{x \in \g_{-} \mid [x,\g_+]
\subseteq \M^i(\g) \}.  \end{equation}

\n Then

\begin{equation}\label{eq:1.7}\begin{gathered} 0 = \M^0(\g) \subseteq \M^1(\g)
\subseteq \dots \subseteq \M^{i-1}(\g) \subseteq \M^i(\g)
\subseteq
\dots \ \ \text{\rm and} \\
\M(\g) \eqdef \,\bigcup_i \, \M^i(\g) \end{gathered}
\end{equation}

\n is called the \ {\em Weisfeiler radical} \ of $\g$.  By its
definition, $\M(\g)$ is a subspace of $\g_-$ invariant under
bracketing by $\g_+$, and for $j = 0,1, \dots, q,$

\begin{eqnarray}\label{eq:1.8} [[\M^i(\g),\g_{-j}],\g_+] &\subseteq&
[[\M^i(\g),\g_+],\g_{-j}] + [\M^i(\g) ,[\g_{-j},\g_+]] \\
&\subseteq& [\M^{i-1}(\g),\g_{-j}] + [\M^i(\g),\sum_{\ell \geq
1}\g_{-j+\ell}]. \nonumber
\end{eqnarray}

\n Now when $j = 0$ and $i = 1$, the right side of \eqref{eq:1.8}
is zero, which implies that $[\M^1(\g),\g_0]$ $\subseteq
\M^1(\g)$.  We may assume that $[\M^{i-1}(\g),\g_0]$ $\subseteq
\M^{i-1}(\g)$.  Then \eqref{eq:1.8} shows that $[\M^i(\g),\g_0]
\subseteq \M^i(\g)$. Suppose we know that

\begin{equation}\label{eq:1.9}[\M^i(\g),\g_{-k}] \subseteq
\M^{i+k}(\g) \end{equation}

\n for all 0 $\leq$ $k < j$.  Then by \eqref{eq:1.8} and induction
on $i$ we have,

$$[[\M^i(\g),\g_{-j}],\g_+]\subseteq [\M^{i-1}(\g),\g_{-j}]
+ [\M^i(\g),\sum_{\ell \geq 1}\g_{-j+\ell}] \subseteq
\M^{i+j-1}(\g).$$  Consequently, $[\M^i(\g),\g_{-j}]$ $\subseteq
\M^{i+j}(\g)$ for $j = 0,$ $1, \dots $.  Thus $\M(\g)$ is an ideal
of $\g$, and it exhibits the following characteristics enjoyed by
a ``radical'':
\bi 

\begin{Pro} \label{Pro:1.10} \
\begin{enumerate}
\item[{\rm (i)}]  $\M(\g)$ is a graded ideal of $\g$ contained in
$\g_-$.

\item[{\rm (ii)}]  Suppose that $\g$ is irreducible \eqref{eq:1.2}  and
transitive \eqref{eq:1.3}, and let $J$ be an ideal of $\g$ contained in $\g_-$.
Then $J$ $\subseteq \M(\g)$ $\subseteq \sum_{i \geq 2} \g_{-i}$.
Thus, $\M(\g)$ is the sum of all the ideals of $\g$ contained in
$\g_-$.  Moreover, $\g/\M(\g)$ is irreducible and transitive, and
$\M(\g/\M(\g)) = 0$.

\item[{\rm (iii)}] If $\g$ is irreducible and transitive,  and
$[\g_{-i},$ $\g_{-1}]$ $= \g_{-(i+1)}$ for all $i \geq 1$, then
for any ideal $J$ of $\g$ either $J \subseteq \M(\g)$ or $J
\supseteq \g_{-}$. 

\item[{\rm (iv)}] If $\g$ is 1-transitive \eqref{eq:1.4}, then $\M(\g) = 0$.
\end{enumerate} \end{Pro}

\pf  (i) We know from the calculations above that $\M(\g)$ is
an ideal of $\g$.  If $x = x_{-q}+ \cdots + x_{-1} \in \M(\g)$, then
one can see by bracketing with
homogeneous elements of positive degree
that the homogeneous components $x_{-i} \in \g_{-i}$ of $x$ belong
to $\M(\g)$ also, and the quotient algebra $\ov \g \eqdef
\g/\M(\g)$ is graded.

\m (ii) Suppose that $\g$ is irreducible and transitive.  If $J$
is any ideal of $\g$ contained in $\g_-$ and if $J \cap \g_{-1}
\neq 0,$ then $J \cap \g_{-1} = \g_{-1} $ by irreducibility, since $J \cap \g_{-1}$ is a $\g_0$-submodule of
$\g_{-1}$.  But if $J \supseteq \g_{-1},$ then $J \supseteq
[\g_{-1},\g_1] \neq 0$ by transitivity,  to
contradict $J \subseteq \g_-$.  Hence, when $\g$ is irreducible
and transitive, every ideal $J$ of $\g$ contained in $\g_-$ has
trivial intersection with $\g_{-1}$.  In particular, since
$\M(\g)\subseteq \g_- $ and $\M(\g)$ is graded, we must have
$\M(\g) \subseteq \sum_{i \geq 2} \g_{-i}.$ It follows that $
\g/\M(\g)$ will inherit the properties of irreducibility
and transitivity.

Now suppose that $J$ is an ideal of $\g$ contained in $\g_-$ and
let $x_{-q} + \cdots + x_{-1} \in J$, where $x_{-i} \in \g_{-i}$ for all $i$.  Then it must be that
$x_{-i} \in \M^i(\g)$ for $i=1,\dots,q$.  Consequently $J
\subseteq \M(\g) \subseteq \sum_{i \geq 2}\, \g_{-i}$, and
$\M(\g)$ is the sum of all the ideals of $\g$ contained in $\g_-$.
Any ideal of $\ov \g = \g/\M(\g)$ contained in $\ov \g_-$ has the
form $\ov K$ where $K$ is an ideal of $\g$ contained in $\g_-$.
But then $K \subseteq \M( \g)$ so $\ov K = 0$. Hence, $\M(\ov \g)
= 0$.

\m (iii) Under the hypotheses in (iii),  suppose that $J$ is an
ideal of $\g$, and define the following subspaces of $\g$:

\begin{equation}\label{eq:1.11} Y_j \eqdef \left \{ y \in \g_{j} \ \Bigg | \ y + z \in J ~~
\hbox{ for~ some~~} z \in \sum_{i \leq j-1} \g_{i}\right \}.
\end{equation}

\n Then $Y_j$ is a $\g_0$-submodule of $\g_{j}$ for each $j$.
Now either $J \subseteq \g_{-}$ (and hence $J \subseteq \M(\g)$ by (ii))  or there exists an
element $x = x_{-q} + \dots + x_k \in J$ with $x_k \neq 0$ and $k
\geq 0$.  Since in the second case $(\ad \g_{-1})^{k+1}(x) \neq
0$ by the transitivity of $\g$, we see that $J \cap
\g_{-} \neq 0$ and $Y_{-1} \neq 0$. Hence by irreducibility,  $Y_{-1}$ 
must equal $\g_{-1}$. Assume we have
shown that $Y_{-i} = \g_{-i}$.  Then $Y_{-(i+1)} \supseteq
[Y_{-i}, \g_{-1}] = [\g_{-i}, \g_{-1}] = \g_{-(i+1)}$ so that
$Y_{-(i+1)} = \g_{-(i+1)}$.  In particular, $J \supseteq \g_{-q}$.
Suppose that $J \supseteq \g_{-t}$ for all $t$ with $i+1 \leq t \leq q$.
Since $Y_{-i} = \g_{-i}$, it follows that $J \supseteq
\g_{-i}$.  Consequently $J \supseteq \g_{-}$, as asserted. 

\m (iv)  Suppose $\M^k(\g) = 0$ for $0 \leq k < i$.  Then $[\M^i(\g), \g_1] \subseteq \M^{i-1}(\g) = 0$.
Thus, if $\g$ is 1-transitive,
we must have  
$\M^i(\g) = 0$ also.    Therefore,  1-transitivity implies   $\M(\g) = \bigcup \M^i(\g) =  0$.   \qed
 
\m \begin{Pro} \label{Pro:1.13}  \ Assume $\g =
\bigoplus_{i=-q}^r \g_i$ is an irreducible, transitive graded Lie
algebra such that $\g_{-i} = [\g_{-1},\g_{-i+1}]$ for all $i \geq
2$.  If $\M(\g) = 0$, then $\g$ contains no abelian ideals.
\end{Pro}

\pf By Proposition \ref{Pro:1.10}\,(iii),  any nonzero ideal $J$ of
$\g$ must contain $\g_{- }$  and hence must contain $\g_{-1}
\oplus [\g_{-1},\g_1]$.  By transitivity \eqref{eq:1.3}, $[J,J]
\supseteq [\g_{-1},[\g_{-1},\g_1]] \neq 0$. Consequently, $J$
cannot be abelian. \ \qed
\m  \smallskip

 \begin{Cor} \label{Cor:1.14} \  If $\g =
\bigoplus_{i=-q}^r \g_i$ is a irreducible, transitive graded Lie
algebra such that $\g_{-i} = [\g_{-1},\g_{-i+1}]$ for all $i \geq
2$, then $\g/\M(\g)$ is semisimple.  \end{Cor}

 \pf This is an immediate consequence of Proposition \ref{Pro:1.13}
and the fact that  $\M(\ov \g) = 0$ for $\ov \g = \g/\M(\g)$  by (iii) of Proposition \ref{Pro:1.10}. \ \qed

\m

\begin{Lem} \label{Lem:1.15} \  Assume $\g =
\bigoplus_{i=-q}^r \g_i$  is a graded Lie algebra such that $\g_+$ is generated by $\g_1$ and
$\M(\g) = 0$.   If  $x \, \in \, \g_{-}$
and $[x,\g_{1}] = 0$, then $x = 0$.  \end{Lem}

\pf   The hypotheses imply that $[x,\g_+] = 0$, so that $x \in \M^1(\g) = 0$.   \qed

\bi
\section  {\ The  minimal ideal \
$\boldsymbol{\mathcal  I}$ \label{sec:1.3} }

In this section we show that graded Lie algebras satisfying certain conditions must contain
a unique minimal ideal $\mathcal  I$  which is graded, and we derive some properties
of $\mathcal  I$.

\bi \begin{Pro} \label{Pro:1.17} \  {\rm (Compare [W], Prop.
1.61, Cor.  1.62, Cor.  1.66.)} \ Assume $\g = \bigoplus_{i=-q}^r
\g_i$ is an irreducible, transitive graded Lie algebra such that
$\g_{-i} = [\g_{-1},\g_{-i+1}]$ for all $i \geq 2$.  If
$\M(\g) = 0$, then $\g$  has a
unique minimal ideal ${\mathcal  I}$ which is graded and contains
$\g_{-}$.
\end{Pro}

\pf  Since $\M(\g) = 0$, every nonzero ideal must contain
$\g_{-}$ by Proposition \ref{Pro:1.10}\,(iii).    Therefore the intersection
$\mathcal  I$ of all the ideals is the unique minimal ideal of $\g$, and
it contains $\g_{-}$.   
 Let ${\mathcal  I}_j \eqdef {\mathcal  I}
\cap  \g_j$ for each $j$, and observe that ${\mathcal  I}_{-i} =
 \g_{-i}$ for all $i \geq 1$.  The space ${\mathcal I}' =
\bigoplus_j {\mathcal  I}_j$ is a nonzero ideal of $\g$, since $[\g_i, {\mathcal I}_j] 
\subseteq  \g_{ i+j} \cap {\mathcal  I} =
{\mathcal  I}_{i+j}$.  By the minimality of ${\mathcal  I}$, we
conclude that ${\mathcal  I} = \bigoplus_j {\mathcal  I}_j$.  \ \qed

\m
\begin{Cor}\label{Cor:1.18}  \ Assume $\g = \bigoplus_{i=-q}^r
\g_i$ is an irreducible, transitive graded Lie algebra such that
$\g_{-i} = [\g_{-1},\g_{-i+1}]$ for all $i \geq 2$.  Then $\ov \g = \g/\M(\g)$ has
a unique minimal ideal which is graded and contains $\ov \g_{-}$.
\end{Cor}
 
\m {\it  {F}rom now on, we assume that the height $r$ of the graded Lie algebras is  finite. } 

\m   Assume $\g$ is an
irreducible, transitive graded Lie algebra with $\M(\g)=
0$ and with $\g_{-}$ generated by $\g_{-1}$, and let $\mathcal  I$ be the unique minimal
ideal of $\g$ from Proposition \ref{Pro:1.17}.    Then 
 
\begin{equation}\label{eq:1.19} {\mathcal  R} \eqdef \sum_{k \ge 1}(\ad \,
\g_{-1})^{k}\g_{r} \subseteq {\mathcal  I} \end{equation}

\noindent   is stable under the action of both $\ad
\, \g_{-1}$ and $\ad \, \g_{0}$.  By the Jacobi identity,
$[\g_{+},{\mathcal  R}] \subseteq {\mathcal  R} + \g_{r}$. As
$\g_{-}$ is generated by $\g_{-1}$ by assumption,  ${\mathcal  R} +
\g_{r}$ is a graded ideal of $\g$.    It follows that ${\mathcal 
I} \subseteq {\mathcal R} + \g_{r}$.   Moreover,

$${\mathcal I} = {\mathcal I} \cap ({\mathcal R} + \g_{r}) =
{\mathcal R} + ({\mathcal I} \cap \g_{r}) \subseteq {\mathcal I}$$

\noindent so ${\mathcal I} = {\mathcal R} + ({\mathcal I} \cap
\g_{r})$.  This establishes the following result: \bi

\begin{Lem} \label{Lem:1.20} \  Any transitive, irreducible graded Lie algebra
$\g$ with $\M(\g)$ = $0$ and $\g_-$ generated by $\g_{-1}$ contains
a unique minimal graded ideal ${\mathcal I}$ such that

$${\mathcal I} =
\sum_{k \ge 1}(\ad\g_{-1})^{k}\g_{r} + ({\mathcal I} \cap
\g_{r}).$$

\noindent Thus, ${\mathcal I} = \bigoplus_{i=-q}^s {\mathcal I}_i$
where ${\mathcal I}_s \neq 0$, and  $s = r$ if ${\mathcal I} \cap
\g_{r} \ne 0$, and $s = r - 1$ otherwise.  Moreover,

$${\mathcal I}_{i} = \begin{cases} (\ad \, \g_{-1})^{r - i}\g_{r} & \qquad
\hbox{\rm if } \  0 \leq i \le r-1,  \\
  \g_{i} & \qquad \hbox{\rm if} \  -q \le i \le -1.  \end{cases} $$

\end{Lem}

 \m 

\section  {\ The graded algebras \
$\boldsymbol{{\mathcal B}(V_{-t})}$ and $\boldsymbol{{\mathcal B}(V_t)}$ \label{sec:1.4}}

 \m There is a general procedure described by Benkart and
Gregory in \cite[Sec.~3]{BG}  for producing quotients of subalgebras in a graded Lie
algebra $\g = \bigoplus_{i=-q}^r \g_i$.  It has been applied
subsequently in several other settings; for example, it was used
by Skryabin in \cite{Sk2} to derive information about the solvable
radical of $\g_0$.

\m The procedure begins with a subalgebra $F_0$ of $\g_0$ and an
$F_0$-submodule $F_{-1}$ of $\g_{-t}$ for some $t$ with $1 \leq t
\leq q$.  For $i > 0$ define:

\begin{eqnarray} && F_{-i} = [F_{-1}, F_{-i+1}],\label{eq:1.22} \\
&& F_i =  \{ y \in \g_{it} \mid [y, F_{-1}] \subseteq F_{i-1}\}.
\label{eq:1.23}
\end{eqnarray}

\n We claim that $F = \bigoplus_j F_j$ is a subalgebra of $\g$. To
see this,  note that $[F_{-1}, F_j] \subseteq F_{j-1}$ for all $j$.
Assume for $1 \leq i < k$ that $[F_{-i}, F_j] \subseteq F_{j-i}$
for all $j$. Then
\begin{eqnarray*}[F_{-k}, F_j] = [[F_{-1}, F_{-k+1}], F_j] &\subseteq&
[[F_{-1}, F_j], F_{-k+1}] + [F_{-1},[F_{ -k+1}, F_j]] \\
&\subseteq& F_{j- k}\end{eqnarray*} for all $j$.  We show that
$[F_i, F_j] \subseteq F_{i+j} $ for $i ,j \geq 0$ by induction on
$i+j$.  If $i + j = 0$, the result is clear, so we may assume that
$i+j > 0$ and that the result holds for values $< i+ j$.  Then
$[F_{-1}, [F_i, F_j]] \subseteq [F_{i-1}, F_j ] + [ F_i,
F_{j-1}]$.   Now if $i-1$ or $j-1$ is negative, we may use what has
been established previously.  Otherwise, we may apply the
inductive hypothesis.  In either event, the sum on the right is in
$F_{i+j-1}$, and $[F_i,F_j] \subseteq F_{i+j}$ by definition.  To
summarize, we have shown:

\bi
\begin{Pro} \label{Pro:1.24} \ Let $\g$ be an arbitrary graded Lie
algebra.  Let $F_0 $ be a subalgebra of $\g_0$, and for $t > 0$
let $F_{-1}$ be an  $F_0$-submodule of $\g_{-t}$.  If $F_{-i}$ and
$F_i$ are defined by \eqref{eq:1.22} and \eqref{eq:1.23}, then $F
= \bigoplus_j F_j$ is a graded subalgebra of $\g$.  \end{Pro}\m

\begin{Pro} \label{Pro:1.25} \ Let $F = \bigoplus_j
F_j$ be any graded Lie algebra satisfying \eqref{eq:1.22},  and assume
that $F_{-1}$ is a nontrivial $F_0$-module under the adjoint action.  Set
$A_{-i} = 0$ for all $i > 0$.  For $i \geq 0,$ define $A_i = \{ y
\in F_i \mid [y, F_{-1}] \subseteq A_{i-1}\}$.  Then $A =
\bigoplus_j A_j$ is a graded ideal of $F$.  Moreover the graded
Lie algebra $\mathcal B = \bigoplus_i {\mathcal B}_i = F / A$
satisfies \eqref{eq:1.3} (transitivity) and \eqref{eq:1.22}.
\end{Pro}

\pf  Assume $F$ is any graded Lie algebra satisfying the assumptions  and $A_j$
 is as in the statement of the proposition.   Then $[F_{-1},A] \subseteq A$, and 
 it follows by induction and the Jacobi identity that $[F_{-i}, A]
 \subseteq A$ for all $i > 0$.  We suppose that $i, j \geq 0$ and prove
 that $[F_i,A_j] \subseteq A$ for all such $i$ and $j$ by induction on
 $i + j$.  The result is clear if $i + j = 0$, for the annhilator of
 $F_{-1} $ is an ideal of $F_0$.  Now $[F_{-1}, [F_i, A_j]] \subseteq [
 F_{i-1}, A_j] + [F_i,A_{j-1 }] \subseteq A_{j+i-1}$ by induction.
 Thus, $[F_i,A_j] \subseteq A_{i+j}$, and $A$ is an ideal as desired.

That \eqref{eq:1.22} holds in ${\mathcal B} = F/A$ is obvious
since $A_j = 0$ for $j < 0$. Suppose for $\ov y \in {\mathcal 
B}_j$ with $j \geq 0$ that $[\ov y, {\mathcal B}_{-1}] = 0$.  If
$\ov y = y + A_j$, then $ [y, F_{-1}] \subseteq A_{j-1}$, and that
implies $y \in A_j$ so that $\ov y = 0$.  Thus, ${\mathcal B}$ is
transitive.  \ \qed

\bi We note that in this construction,  $A_0 =
\Ann_{F_0}(F_{-1})$, the annihilator of $F_{-1}$ in $F_0$.

\bi By symmetry we have:

\bi \begin{Pro} \label{Pro:1.26}  \ Let $\g$ be an arbitrary
graded Lie algebra.  Let $F_0$ be a subalgebra of $\g_0$, and let
$F_1$ be an $F_0$-submodule of $\g_t$ for some $t \geq 1$.  For $i
> 0$, define $F_i$ and $F_{-i}$ as follows:

\begin{equation}\label{eq:1.27} F_{i} = [F_1, F_{i-1}], \end{equation}
\begin{equation}\label{eq:1.28}
 F_{-i} = \{ y \in \g_{-it} \mid [y, F_1] \subseteq F_{-i+1}\}.
\end{equation}

\m \n Then $F= \bigoplus_j F_j$ is a graded subalgebra of $\g$.
\end{Pro} \m

\begin{Pro} \label{Pro:1.29}  \ Let $F = \bigoplus_j F_j$ be any graded
Lie algebra satisfying \eqref{eq:1.27}, and assume that $F_1$ is
a nontrivial $F_0$-module under the adjoint action.  Set $A_{i} = 0$ for all
$i > 0$. For $i \geq 0,$ define $A_{-i} = \{ y \in F_{-i} \mid [y,
F_1] \subseteq A_{-i+1}\}$.  Then $A = \bigoplus_j A_j$ is a
graded ideal of $F$. Moreover in the graded Lie algebra ${\mathcal B} =
\bigoplus_i {\mathcal B}_i = F / A$,  conditions \eqref{eq:1.4}
(1-transitivity) and \eqref{eq:1.27}  hold.
\end{Pro}

 \bi
  There are several  important instances of this construction that
will play a crucial role in the remainder of the work.  For the
graded Lie algebra $\g =\bigoplus_{i=-q}^r \g_i$, and for $1 \leq
t \leq q$, the ingredients in this special case are the following:

\begin{equation}\label{eq:1.30} F_0 = \g_0 \ \hbox{\rm and } \ F_{-1} = V_{-t}, ~~ \hbox{\rm
 a nontrivial}~~\g_0\hbox{\rm -submodule of}~~\g_{-t}, \end{equation}
\begin{equation}\label{eq:1.31}F_{-i} =
 [F_{-1}, F_{-i+1}] \ \hbox{ and } \ F_i = \{ y \in \g_{it} \mid [y,
 F_{-1}] \subseteq F_{i-1}\} = \g_{it} \ \ (i \geq 1), \end{equation}
\begin{equation}\label{eq:1.32} A_0 = \Ann_{\g_0}(F_{-1}) =
\{x \in \g_0 \mid [x,F_{-1}] = 0 \},
\end{equation}
\begin{equation}\label{eq:1.33}  A_i = \{ x \in \g_{it} \mid [x,F_{-1}] \subseteq
A_{i-1}\} \quad \hbox{ for} \ \ i \geq 1,
\end{equation}
\begin{equation}\label{eq:1.34} F = F(V_{-t}) = \bigoplus_j F_j ~~
\hbox{ \rm and }~~ A = A(V_{-t}) = \bigoplus_{i \geq 0} A_i,
\end{equation}
\begin{equation}\label{eq:1.35}
{\mathcal B} = {\mathcal B}(V_{-t}) = F(V_{-t})/A(V_{-t}).
\end{equation}

 \m
Usually we take $V_{-t}$ to be an {\em irreducible}
$\g_0$-submodule of $\g_{-t}$, although this is not required for
the construction.  We try to make choices of $t$ and submodules
$V_{-t}$ that enable us to show that the homogeneous components
${\mathcal B}_0$ and ${\mathcal B}_1$ of ${\mathcal B}$ are both
nonzero, for then the graded algebra ${\mathcal B}$ has at least
three homogeneous components, and  we can
deduce much more information about the structure of the 
graded algebra  ${\mathcal B}$.  When $t > \lfloor q/2\rfloor$, the
algebra ${\mathcal B}$ has depth one.  This proves very helpful as
it allows us to transfer results about depth-one algebras to
general graded Lie algebras.

\m The analogous construction starting from a $\g_0$-submodule
$V_t$ of $\g_t$, for $1 \leq t \leq r$, has the following
constituents: \smallskip

\begin{equation}\label{eq:1.36} F_0 = \g_0 \ \hbox{\rm and }\ F_{1} =
V_t, ~~ \hbox{a nontrivial}~~\g_0 \hbox{\rm -submodule~of} \ \g_{t},
\end{equation}
\begin{equation}\label{eq:1.37} F_{i} = [F_{1}, F_{i-1}]\ \
\hbox{\rm and } \ \ F_{-i} = \{ y \in \g_{-it} \mid [y, F_{1}]
\subseteq F_{-i+1}\}= \g_{-it} \ \ (i \geq 1), \end{equation}
\begin{equation}\label{eq:1.38} A_0 =
\Ann_{\g_0}(F_{1}) = \{x \in \g_0 \mid [x,F_1] = 0 \},
\end{equation}
\begin{equation}\label{eq:1.39} A_{-i}= \{ x \in \g_{-it} \mid [x,F_{1}] \subseteq A_{-i+1}\} \quad
\hbox{ for}\ \ i \geq 1, \end{equation}
\begin{equation}\label{eq:1.40} F = F(V_t) = \bigoplus_j F_j ~~ \hbox{ and }~~ A = A(V_t) =
\bigoplus_{i \geq 0} A_{-i},  \end{equation}
\begin{equation}\label{eq:1.41}
{\mathcal B} = {\mathcal B}(V_t) = F(V_t)/A(V_t).  \end{equation}

\m
\begin{Rems} \label{subsec:1.42} 
{\rm When $t > \lfloor r/2 \rfloor$, the algebra ${\mathcal 
B}(V_t)$ constructed in \eqref{eq:1.36}--\eqref{eq:1.41} has
height one and is $1$-transitive by Proposition \ref{Pro:1.29}.
By assigning ${\mathcal B}(V_t)$ the opposite grading where
${\mathcal B}(V_t)_{-i}$ is set equal to ${\mathcal B}(V_t)_i$, we
obtain a depth-one, transitive algebra, which is irreducible if
and only if the module $V_t$ is irreducible. When speaking about
${\mathcal B}(V_t)$ in subsequent sections, we always assume it is
graded with the opposite grading and hence always regard it as a
depth-one algebra.} \end{Rems}

\m
 \section{\ The local
subalgebra  \label{sec:1.5} } 

 \m There is yet another subalgebra of a graded Lie algebra
$\g =\bigoplus_{i=-q}^r \g_i$ which plays a prominent role.  This
is the subalgebra $\widehat \g$ of $\g$   generated by
the \ {\em local part} \ $\g_{-1} \oplus \g_0 \oplus \g_1$.  Thus,
\begin{equation}\label{eq:1.44} \widehat \g = \bigoplus_{i=-q'}^{r'}
\widehat \g_{i}
\end{equation}

\n where $\widehat \g_i = \g_i$ for $i = -1,0,1,$ \ $\widehat
\g_{-j} = [\widehat \g_{-j + 1}, \widehat \g_{-1}]$, \ and \
$\widehat \g_{i} = [\widehat \g_{i-1}, \widehat \g_{1}]$ \
for all $i,j \geq 2$.

\bi
\begin{Pro} \label{Pro:1.45} \ Let $\widehat \g$ be the Lie algebra
generated by the local part $\g_{-1} \oplus \g_0 \oplus \g_1$ of
the graded algebra $\g$.  Set ${\mathcal A}^- =
\bigoplus_j{\mathcal A}^-_j$ where ${\mathcal A}^-_{-i} = 0$ for
$i > 0$ and ${\mathcal A}^-_i = \{y \in \widehat \g_i \mid
[y,\widehat \g_{-1}] \subseteq {\mathcal A}^-_{i-1}\}$ for all
$i \geq 0$.  Similarly, let ${\mathcal A}^+ = \bigoplus_j{\mathcal 
A}^+_j$ where ${\mathcal A}^+_i = 0$ for $i > 0$ and ${\mathcal 
A}^+_{-i} = \{y \in \widehat \g_{-i} \mid [y,\widehat \g_1]
\subseteq {\mathcal A}^+_{-i+1}\}$ for all $i \geq 0$.  Then
$\widehat \g/{\mathcal A}^-$ is transitive and $\widehat
\g/{\mathcal A}^+$ is  1-transitive.  Thus, $\widehat \g$ is
transitive \eqref{eq:1.3} (respectively, 1-transitive \eqref{eq:1.4}) if and only
if ${\mathcal A}^- = 0$ (respectively, ${\mathcal A}^+ = 0$).
\end{Pro}

\pf The statements about the transitivity  of
$\widehat \g/{\mathcal A}^-$ and the $1$-transitivity
of $\widehat \g/{\mathcal A}^+$ follow directly
from Propositions \ref{Pro:1.25} and \ref{Pro:1.29}.  Now if $y
\in \widehat \g_i$ for some $i \geq 0$, then $y \in {\mathcal 
A}^-_i$ if and only if $(\ad \widehat \g_{-1})^{i+1}(y) = 0$.
If $\widehat \g$ is transitive, then $(\ad \, \widehat
\g_{-1})^{i+1}(y) = 0$ gives $y = 0$. Thus, transitivity
\eqref{eq:1.3} implies ${\mathcal A}^-_i = 0 $ for all $i \geq 0$.
Conversely, when ${\mathcal A}^-_i = 0$  for all $i \geq 0$,
then $\widehat \g = \widehat \g/{\mathcal A}^-$ is transitive.
Similar arguments yield the other assertions.  \ \qed
 
\m
\section  {\ General
properties of  graded Lie algebras  \label{sec:1.6}}  

\m In this section we
will derive some very general properties of graded Lie algebras
$\g = \bigoplus_{i=-q}^r \g_i$ which satisfy

\begin{eqnarray}\label{sec:1.6conds}
&&\hspace{-.25 truein} \ \ \hbox{\rm (i)} \quad \g \ \hbox{\rm  is irreducible \eqref{eq:1.2} and transitive \eqref{eq:1.3};} \\
&&\hspace{-.25 truein} \ \hbox{\rm (ii)} \quad \M(\g) = 0; \nonumber \\
&&\hspace{-.25 truein} \hbox{\rm (iii)} \quad   \g_{-} \ \hbox{\rm is
generated by} \  \g_{-1}; \  \hbox{\rm and} \nonumber \\
&&\hspace{-.25 truein}\, \hbox{\rm (iv)} \quad   \g \  \hbox{\rm is finite-dimensional over a field
$\F$ of characteristic $\neq 2, 3$}. \nonumber \end{eqnarray}
 
\m Our first objective will be to show that transitivity
holds in all of $\g \setminus \g_{-q}$.  Hence we
prove for any $x \in \g \setminus \g_{-q}$ that  $[x,\g_{-1}] \neq  0$.  We  may assume $x = x_{-q} + \cdots + x_j$, where $x_i \in \g_i$ for all $i$, 
and $x_j \neq 0$.   It follows that $[x_j, \g_{-1}] = 0$.   As $\g$ is assumed to be transitive, we
may suppose $-q+1 \leq j \leq -1$.    Then  since $\g_{-}$ is
generated by $\g_{-1}$,

$$J \eqdef \sum_{k,\ell \ge 0}(\ad \, \g_{+})^{k}(\ad \,
\g_{0})^{\ell}x_j$$

\noindent would be a graded ideal of $\g$.  But $j \geq  -q +1$, so
$J$ cannot contain the unique minimal graded ideal $\mathcal I$,
to contradict Lemma \ref{Lem:1.20}.     Thus, we have demonstrated
the first part of the next result.

\bi \begin{Lem}\label{Lem:1.51} \ Suppose $\g = \bigoplus_{j=-q}^r \g_j$ is a graded Lie
algebra satisfying (i)-(iii) of \eqref{sec:1.6conds}.  Then the following hold:
\begin{enumerate}
\item[{\rm(a)}] If $x \, \in \, \g \setminus \g_{-q}$, then $[x, \g_{-1}] \neq 0$. 
\item[{\rm(b)}]
$\g_{-q}$ is an irreducible $\g_0$-module. \item[{\rm(c)}]  If $K$
is an ideal of $\g$ and  $K \cap (\g_{-1} \oplus \g_0 \oplus \g_1)
\neq 0$, then $K  \supseteq \g_{-}$.  \end{enumerate}
\end{Lem}

\pf   For (b),
note that similar to the argument above, if $M$
is any proper $\g_0$-submodule of $\g_{-q}$, then

$$J \eqdef
\sum_{k,\ell \ge 0}(\ad \, \g_{+})^{k}(\ad \, \g_{0})^{\ell}M$$

\n is a graded ideal of $\g$, which fails to contain $\g_{-q}$,
hence fails to contain the unique minimal graded ideal  $\mathcal I$. To establish (c), observe
that when $K$ is an ideal such that $K \,\cap\, (\g_{-1} \oplus \g_0
\oplus \g_1) \neq 0$, then by transitivity \eqref{eq:1.3}, \ $K\,
\cap\, \g_{-1} \neq 0$.  By irreducibility \eqref{eq:1.2}, $\,K \cap\,\g_{-1} = \g_{-1}$,
and since $\g_{-}$ is generated by $\g_{-1}$,
$K$ must contain $\g_{-}$.  \ \qed

\bi Suppose now that $V_{-t}$ is any irreducible
$\g_{0}$-submodule of $\g_{-t}$, where $q/2  < t \leq q$.
(For example, $V_{-q} = \g_{-q}$ is a choice we have in mind in
view of Lemma \ref{Lem:1.51}\,(b).) As in
\eqref{eq:1.30}--\eqref{eq:1.35},

\begin{eqnarray} && F = F(V_{-t}) = \bigoplus_{i\ge-1}
{F_i}, \ \ \hbox{\rm where} \ \ F_{-1} = V_{-t},   \\
&& \hskip .3 truein \ \   F_{-j} = [F_{-1},F_{-j+1}],
 \ \ \hbox{\rm and} \ \ F_{i} = \g_{it} \ \ \hbox{\rm
for} \ \ i \geq 0, \nonumber \\
&& A = A(V_{-t}) = \bigoplus_{i\ge-1}{A_i},  \
\ \hbox{\rm where} \ \ A_{-1} = 0 \ \ \hbox{\rm and} \label{eq:1.52} \\
&& \hskip .3 truein \ \ A_i = \{ y \in F_i \mid [y, F_{-1}]
\subseteq
A_{i-1}\} \ \ \hbox{\rm for} \ \ i \geq 0, \nonumber \\
&& {\mathcal B} = {\mathcal B}(V_{-t}) = \bigoplus_i {\mathcal 
B}_i = F / A.  \nonumber  \end{eqnarray}

\n The algebra ${\mathcal B}$ satisfies \eqref{eq:1.2}
(irreducibility), \eqref{eq:1.3} (transitivity), and
\eqref{eq:1.22}.  Note that $ A_{j} = \{x \in \g_{jt} \mid (\ad
V_{-t})^{j+1}x = 0 \}$, and $\mathcal B_{j} \ne 0
\Longleftrightarrow A_{j} \ne F_{j} \Longleftrightarrow (\ad
V_{-t})^{j+1}\g_{tj} \ne 0.$ In particular,

\begin{equation}\label{eq:1.53} {\mathcal B}_{1} \neq 0 \Longleftrightarrow (\ad
V_{-t})^{2}\g_{t} \ne 0.  \end{equation}

\m Next we will set  $F = F(\g_{-q})$, $A = A(\g_{-q})$,
and  ${\mathcal B} = {\mathcal B}(\g_{-q}) = F/A = \bigoplus_i
{\mathcal B}_i$, and we will argue that ${\mathcal B}_{1} \neq 0$ when $q \leq r$.
But to accomplish this, we require
some additional information about the action of $\ad \, \g_{-q}$
on the various spaces $\g_i$.  Clearly, $\ad \, \g_{-q}$ acts
nilpotently on $\g$.

\bi \begin{Lem}\label{Lem:1.54}\  Suppose $\g = \bigoplus_{j=-q}^r \g_j$ is a graded Lie
algebra satisfying (i)-(iii) of  \eqref{sec:1.6conds}.  Let $\mathcal I = \bigoplus_{j=-q}^s \mathcal I_j$
be the unique minimal graded ideal of $\g$ as in Lemma  \ref{Lem:1.20}, and assume
$\ell \geq 0$ satisfies $\ell q-q \le s$.
 Then $(\ad \, \g_{-q})^{\ell}{\mathcal I} = 0$ if
and only if $(\ad \, \g_{-q})^{\ell}{\mathcal I}_{j} = 0$ for some
$j$ such that $\ell q-q \le j \le s$.  \end{Lem}

\pf  What is being asserted in this rather technical lemma is
that if $(\ad \, \g_{-q})^{\ell}{\mathcal I}_{j}$ $= 0$ for some
$j$ $\geq \ell q - q$ (i.e., $j$ is large enough so that $(\ad \,
\g_{-q})^{\ell}{\mathcal I}_{j}$ $= 0$ is not an automatic
consequence of having $-\ell q +j$ $< -q$), then $(\ad \,
\g_{-q})^{\ell}$ annihilates ${\mathcal I}$.  To see this, suppose
that

$$(\ad \, \g_{-q})^{\ell}{\mathcal I}_{j} = 0$$

\noindent for some $j$ such that $\ell q-q$ $\leq j$ $\leq s$.
Then it is necessary to show that $(\ad \,
\g_{-q})^{\ell}{\mathcal I}$ $= 0$, which we will do by proving
that $(\ad \, \g_{-q})^{\ell}{\mathcal I}_{k}$ $= 0$ for all $k$
with $ -q$ $\le k$ $\le s$.  Indeed if $-q$ $\le$  $k$ $< \ell
q-q$, then $(\ad \, \g_{-q})^{\ell}{\mathcal I}_{k}$ $\subseteq
\g_{k-q\ell}$ $= 0$.  To prove the result for $\ell q-q$ $\le k$
$\le s$, we distinguish two cases: (1) $\ell q-q$ $\le k$ $< j$
and (2) $j$ $< k$ $\le s$.

In case (1), we have by assumption and Lemma \ref{Lem:1.20} that

$$0 =
(\ad \, \g_{-1})^{j-k}(\ad \, \g_{-q})^{\ell}{\mathcal I}_{j} =
(\ad \, \g_{-q})^{\ell}(\ad \, \g_{-1})^{j-k}{\mathcal I}_{j} =
(\ad \, \g_{-q})^{\ell}{\mathcal I}_{k}.$$

\noindent In the second case,

$$(\ad \, \g_{-1})^{k-j}(\ad \, \g_{-q})^{\ell}{\mathcal I}_{k} = (\ad
\, \g_{-q})^{\ell}(\ad \, \g_{-1})^{k-j}{\mathcal I}_{k} = (\ad \,
\g_{-q})^{\ell}{\mathcal I}_{j} = 0$$

\noindent must hold so that by Lemma \ref{Lem:1.51} (a), $(\ad \,
\g_{-q})^{\ell}{\mathcal I}_{k} = 0$.  The other direction is
obvious.  \ \qed

\bi  In view of \eqref{eq:1.53}, showing ${\mathcal B}_1 \ne 0$
for ${\mathcal B} = {\mathcal B}(\g_{-q})$ is equivalent to
showing that $(\ad \, \g_{-q})^{2}\g_{q} \ne 0$.  As a starting
point we have

\bi
\begin{Lem} \label{Lem:1.55} \, Under the assumptions of   \eqref{sec:1.6conds}, \
$(\ad \, \g_{-q}){\mathcal I}_{j} \ne 0$ for $0 \le j \le s$,
where $\mathcal I = \bigoplus_{j=-q}^s \mathcal I_j$
is the unique minimal graded ideal of $\g$ as in Lemma  \ref{Lem:1.20}.  In particular, since
$\g$ is irreducible, $[\g_{-q},{\mathcal I}_{q-1}] = \g_{-1}$
whenever $s \ge q-1$.
\end{Lem}

\pf If $(\ad \g_{-q}){\mathcal I}_{j} = 0$ for some $j$,  \ $0
\leq j\leq s$, then by Lemma \ref{Lem:1.54} (with $\ell=1$), we
have $(\ad \g_{-q}){\mathcal I} = 0$.  Since $\g_{-}$ is 
assumed to be generated by $\g_{-1}$, and since
$[\g_{-1},\g_{+}] \subseteq {\mathcal I}$,  it follows that

$$ \sum_{i \ge 1}(\ad \,  \g_{+})^{i}{\mathcal I}_{-q}$$

\noindent is an ideal of $\g$ not containing ${\mathcal I}_{-q} = \g_{-q}$,
and hence properly contained in ${\mathcal I}$, to contradict the
minimality of ${\mathcal I}$.  \qed

\bi
\begin{Rem} \label{Rem:1.56} {\rm  To complete the proof that
${\mathcal B}_1 \neq 0$ for
${\mathcal B}= {\mathcal B}(\g_{-q})$ when $q \leq r$, we use the notion of a
weakly closed set.  Recall that a subset $\mathcal E$ of the
endomorphism algebra $\En(W)$ on a vector space $W$ over $\F$ is \
{\em weakly closed} \ if for every ordered pair $(x,y) \in
\mathcal E^2$, there is a scalar $f(x,y) \in {\F}$ such that $xy +
f(x,y)yx \in \mathcal E$.  If $\mathcal E$ is nil (that is, every
$x\in\mathcal E$ is a nilpotent transformation), then the
associative subalgebra of $\En(W)$ generated by $\mathcal E$ is
nilpotent (see \cite[p.~33]{J}).} \end{Rem}
\bigskip

\begin{Lem} \label{Lem:1.57} \ Let $U$ and $V$ be $\g_{0}$-submodules of
$\g$ such that $[U, V] \subseteq \g_{0}$ and $[U, [U, V]] = 0$.
Let $M$ be a $\g_{0}$-submodule of $\g$. Then the following are
true:
\begin{enumerate}
\item[{\rm(i)}] $\ad_M[U,V] \eqdef \left \{\ad [u, v] \big |_{M} \, \big | \, u
\, \in U, v \, \in V \right \}$ is a weakly closed subset of $\En(M)$.
\item[{\rm(ii)}] If $\ad_M [U,V]$ consists of nilpotent
transformations, then the associative subalgebra of $\En(M)$
generated by $\ad_M[U,V]$ acts nilpotently on $M$.  If
in addition  $M$ is an irreducible $\g_{0}$-module, then
$[[U,V],M] = 0$.  \end{enumerate}
\end{Lem}

\pf  Let $u_{1}$ and $u_{2}$ be any elements of $U$, and let
$v_{1}$ and $v_{2}$ be any elements of $V$.  Then we have

\begin{eqnarray*} \left [\ad \, [u_{1}, v_{1}],\ad \, [u_{2}, v_{2}]\right] &=& \ad \,
[[u_{1}, v_{1}], [u_{2}, v_{2}]] \\
&=& \ad \, [[[u_{1}, v_{1}], u_{2}], v_{2}] +\ad \, [u_{2},
[[u_{1},
v_{1}], v_{2}]] \\
&=& \ad \, [u_{2}, [[u_{1}, v_{1}], v_{2}]].  \end{eqnarray*}

\noindent Since $[u_{1}, v_{1}] \, \in \, \g_{0}$ and $V$ is
assumed to be a $\g_{0}$-submodule of $\g$, it follows that
$[[u_{1}, v_{1}], v_{2}] \, \in \, V$, so $\ad_M[U,V]$ is weakly
closed. Thus if  $\ad_M[U,V]$ consists of nilpotent
transformations, then the associative subalgebra of $\En(M)$
generated by $\ad_M[U,V]$ is nilpotent.  But then
$\Ann_M\big([U,V]\big) \neq 0$.  Since $\Ann_M\big([U,V]\big)$ is
a $\g_0$ submodule of $M$, it must equal $M$ whenever $M$ is
irreducible, implying $[[U,V],M] = 0$ in this case. \ \qed

\bi In the next sequence of results, as in Lemma \ref{Lem:1.57},
we adopt the notation $\ad_V(U)$ (respectively $\ad_{V}u$) to
indicate the restriction of $\ad U$ (respectively $\ad u$) to the
subspace $V$.

\m We now have the requisites for showing that ${\mathcal B}_1 \ne
0$ when $q \leq r$ for ${\mathcal B} = {\mathcal B}(\g_{-q})$, or
equivalently, that $(\ad \, \g_{-q})^{2}\g_{q} \ne 0$ when $\g$ is finite-dimensional.  Suppose to
the contrary that $(\ad \, \g_{-q})^{2}\g_{q} = 0$.  If $q \leq s$,  Lemma
\ref{Lem:1.54} (with $\ell = 2$)  implies that $(\ad \, \g_{-q})^{2}{\mathcal I} =
0$.  If instead $s < q$, then because $r - 1 \le s < q \le r$, we 
must have $s = r - 1 = q-1$, so again we have $(\ad \,
\g_{-q})^{2}{\mathcal I} \subseteq \g_{-q-1} \oplus \g_{-q-2}
\oplus \cdots = 0$.  In any event, for all $y \in \g_{-q}$ and $z
\in \g_{q}$,

$$0 = (\ad y)^{2}(\ad z)^{2}\g_{-1} = 2(\ad \,
[y,z])^{2}\g_{-1},$$

\noindent because $\g_{-1}={\mathcal I}_{-1}$ and
$[\g_{-q},\g_{-1}] = 0 =(\ad \g_{-q})^{2}{\mathcal I}$.  Since by Lemma \ref{Lem:1.57} (i),
$$Y\eqdef \left \{\ad_{\g_{-1}}[y,z]\mid y\in \g_{-q}, \,z \in \g_{q}
\right \}$$ is a weakly closed set, and since it consists of nilpotent transformations of
$\g_{-1}$,  the associative subalgebra of $\En(\g_{-1})$
generated by $Y$ acts nilpotently on $\g_{-1}$.   As $\g_{-1}$
is $\g_{0}$-irreducible, Lemma \ref{Lem:1.57} (ii) implies that
$[[\g_{-q},$ $\g_{q}],$ $\g_{-1}]$ $= 0$.  By transitivity
\eqref{eq:1.3}, $[\g_{-q},$ $\g_{q}]$ $= 0$.  Thus,
$[\g_{-q}, \mathcal I_{q-1}] \subseteq [\g_{-q}, [\g_{-1}, \g_q]] = 0$, to contradict Lemma
\ref{Lem:1.55}.  Thus, we have established

\bi \begin{Lem} \label{Lem:1.58}\ Suppose $\g$
is a graded Lie algebra satisfying
the conditions of \eqref{sec:1.6conds}.  If  $q \leq r$,   then  $(\ad \g_{-q})^{2}\g_{q} \neq 0$.
Consequently, ${\mathcal B}_1 \neq 0$ for ${\mathcal B} =
{\mathcal B}(\g_{-q}) = \bigoplus_i {\mathcal B}_i$.  \end{Lem}
 \bigskip

Our next goal is to prove when $q \geq 2$ that ${\mathcal B}'_1$ is nonzero
for ${\mathcal B}' = {\mathcal B}(V_{-q+1})$,  where $V_{-q+1}$ is any
irreducible $\g_{0}$-submodule of $\g_{-q+1}$.  This will require
some results showing that annihilators in $\g_{0}$ of pairs of
spaces $\g_j$ and $\g_{j+1}$ intersect trivially.  In fact, we
prove a slightly more general result.  For $-q \leq j \leq r - 1$,
let $V_{j+1}$ be any nonzero $\g_{0}$-submodule of $\g_{j+1}$, and
set $Q = \Ann_{\g_{0}}({\mathcal I}_{j}) \cap \Ann_{\g_{0}}(V_{j+1})$.
Under the assumption $Q \ne 0$, transitivity \eqref{eq:1.3} and
irreducibility \eqref{eq:1.2} force $[Q,\g_{-1}] = \g_{-1}$.  Then

$$[\g_{-1},V_{j+1}] = [[Q,\g_{-1}],V_{j+1}] =
[Q,[\g_{-1},V_{j+1}]] \subseteq [Q,{\mathcal I}_{j}] = 0$$

\noindent to contradict transitivity (we have used the inclusion
$\g_{-1}\subset {\mathcal I})$. Consequently, we have

 \bi
\begin{Lem} \label{Lem:1.59} \ Assume  $q \geq 2$ and the hypotheses of    \eqref{sec:1.6conds}
hold.  If  $V_{j+1} \neq 0$ is a $\g_{0}$-submodule of
$\g_{j+1}$ for some $j$ such that  $-q \le j \le r - 1$, then
$\Ann_{\g_{0}}({\mathcal I}_{j}) \cap \Ann_{\g_{0}}(V_{j+1}) = 0$.
\end{Lem} 

\smallskip
\begin{Lem} \label{Lem:1.60}  \ Under the assumptions of
 \eqref{sec:1.6conds}, if $q \geq 2$ and  $V_{-q+1} \neq 0$ is a $\g_{0}$-submodule
of $\g_{-q+1}$,   then $\Ann_{\g_{j}}(\g_{-q}) \cap
\Ann_{\g_{j}}(V_{-q+1}) = 0$ for all  $j$ such that $0 \leq j \leq
r-1$.   \end{Lem}

\pf  Suppose $N\eqdef \Ann_{\g_{j}}(\g_{-q})\cap
 \Ann_{\g_{j}}(V_{-q+1})$.  Then 
 
 $$[\g_{-q}, [\g_{-1},N]]=
 [\g_{-1},[\g_{-q},N]]= 0.$$  
 
 \noindent In addition,  we have that $[V_{-q+1},$ $[\g_{-1},$ $N]]$ $\subseteq$
 $[\g_{-q},$ $N]$ $+ [\g_{-1},$ $[V_{-q+1},$ $N]]$ $= 0$.   Hence  $[\g_{-1},$ $N]$ $\subseteq$
 $\Ann_{\g_{j-1}}(\g_{-q}) \cap
 \Ann_{\g_{j-1}}(V_{-q+1})$.  Repeating this argument $j - 1$ times, we
 determine that $(\ad \g_{-1})^{j}N\subseteq\Ann_{\g_{0}}(\g_{-q})
 \cap\Ann_{\g_{0}}(V _{-q+1})$, which is zero by Lemma \ref{Lem:1.59}.  By
 transitivity \eqref{eq:1.3}, $N$ $= 0$ as claimed. \ \qed

 \bi In order to show that ${\mathcal B}'_1$ $\neq 0$ where ${\mathcal 
 B}'$ $= {\mathcal B}(V_{-q+1})$ for some irreducible
 $\g_{0}$-submodule $V_{-q+1}$ of $\g_{-q+1}$, it will suffice to argue
 that $(\ad \, V_{-q+1})^{2}\g_{q-1}$ $\ne 0$ (compare \eqref{eq:1.53}).  At this
 stage we can prove

\bi \begin{Lem} \label{Lem:1.61}\ With assumptions as in
\eqref{sec:1.6conds}, \  if $q \geq 2$, then  \break  $\Ann_{\g_{q-1}}(V_{-q+1}) = 0$ for any
irreducible $\g_{0}$-submodule $V_{-q+1}$ of $\g_{-q+1}$.
\end{Lem}

\pf  Set $N \eqdef \Ann_{\g_{q-1}}(V_{-q+1})$, and assume that $N \ne
 0$.  Then, by Lemma \ref{Lem:1.60}, $[\g_{-q},N] \ne 0$, and so by
 irreducibility \eqref{eq:1.2}, $[\g_{-q},N] = \g_{-1}$.  But then

$$[\g_{-1},V_{-q+1}] = [[\g_{-q},N],V_{-q+1}] = [\g_{-q},[N,V_{-q+1}]] =
0$$

\noindent to contradict Lemma \ref{Lem:1.51} (a). \  \qed \bi

    From Lemmas \ref{Lem:1.51} (b) and \ref{Lem:1.55},
    we know that $\g_{-q}$ is an irreducible
$\g_{0}$-module and that $\g_{-q} = [\g_{-q}, \g_{0}]$.  An
analogous statement holds with $\g_i$ in place of $\g_0$ under
some mild assumptions.

\bi
\begin{Lem} \label{Lem:1.62} \ Suppose $\g$ is a graded Lie algebra satisfying the
conditions of   \eqref{sec:1.6conds}.  If $2 \leq q \leq r$, then $\g_{-q+i}
= [\g_{-q}, \g_{i}]$ for all $i$ such that $0 \le i \le q - 1$.
\end{Lem}

\pf  We may assume that $i>0$.
When $q \leq r$, then by Lemma \ref{Lem:1.55} $[\g_{-q}, \g_{q-1}]
=
 \g_{-1}$.  Since  $\g_{-}$ is generated by $\g_{-1}$, we have
\begin{eqnarray*} \g_{-q+i} &=& (\ad \, \g_{-1})^{q-i-1}(\g_{-1})\\
&=& (\ad \, \g_{-1})^{q-i-1}[\g_{-q}, \g_{q-1}]\\ &=& [\g_{-q}, (\ad
\, \g_{-1})^{q-i-1}\g_{q-1}]  \subseteq [\g_{-q}, \g_{i}]. \quad
\square
\end{eqnarray*}

\m We now have the tools for showing when $3 < q \leq r$ that
${\mathcal B}'_1$ is nonzero for ${\mathcal B}' = {\mathcal 
B}(V_{-q+1})$ and $V_{-q+1}$ any irreducible $\g_{0}$-submodule of
$\g_{-q+1}$.

\bi\begin{Lem} \label{Lem:1.63} \  Assume $\g$ is a graded Lie
algebra satisfying the conditions of   \eqref{sec:1.6conds}, and
suppose further that $3 < q \leq r$.  Let $V_{-q+1}$ be an
irreducible $\g_0$-submodule of $\g_{-q+1}$. Then ${\mathcal 
B}'_{1} \neq 0$ for ${\mathcal B}' = {\mathcal B}(V_{-q+1}) =
\bigoplus_i {\mathcal B}'_i$, (equivalently
$[V_{-q+1},[V_{-q+1},\,\g_{q-1}]] \ne 0$).  \end{Lem}

\pf   Suppose to the contrary that
$[V_{-q+1},[V_{-q+1},\,\g_{q-1}]] = 0.$ Then it follows from that and Lemma
\ref{Lem:1.57} (i) that $Z \eqdef \{\ad_{\g_{-q}}\,[y,z] \mid y \,
\in V_{-q+1}, \,z \, \in \g_{q-1} \}$ is a weakly closed subset of
$\En(\g_{-q})$.   First we suppose that

\begin{equation}\label{eq:1.64} (\ad \, y)^{3}(\ad \, z)^{3}\g_{-q} = 0
\end{equation}

\noindent for all $y \, \in V_{-q+1}$ and all $z \, \in \g_{q-1}$.
Since $[\g_{-q+1},\g_{-q}]=0$, our initial assumption now implies
that $0 = (\ad \, y)^{3}(\ad \, z)^{3}\g_{-q} = 6(\ad \,
[y,z])^{3}\g_{-q}$.  Then $Z $ consists of nilpotent
transformations, and the associative subalgebra of $\En(\g_{-q})$
generated by $Z$ acts nilpotently on $\g_{-q}$.  As $\g_{-q}$ is
$\g_{0}$-irreducible, it follows from Lemma \ref{Lem:1.57} (ii)
that $[\g_{-q},[V_{-q+1},\g_{q-1}]] = 0$.  However then,

$$0 = [\g_{-q},[V_{-q+1},\g_{q-1}]] = [V_{-q+1},[\g_{-q},\g_{q-1}]] =
[V_{-q+1},\g_{-1}]$$

\noindent to contradict Lemma \ref{Lem:1.51} (a).  Thus,  
assumption (\ref{eq:1.64}) has led to a contradiction, so we have
from Lemma \ref{Lem:1.51}(b) that

$$\g_{-q} = (\ad \, V_{-q+1})^{3}(\ad \, \g_{q-1})^{3}\g_{-q}.$$
Note that $V_{-q+1}$ commutes with  $\g_{-q+2}$ because $q>3$. By
Lemma \ref{Lem:1.62},
\begin{eqnarray*} 0\,\ne\,\g_{-q+1} &=& [\g_{1},\g_{-q}] = [\g_{1},(\ad \,
V_{-q+1})^{3}(\ad \, \g_{q-1})^{3}\g_{-q}] \\
&\subseteq& (\ad \, V_{-q+1})^{2}\g_{q-1}
\end{eqnarray*}
contrary to our initial assumption. Thus
$[V_{-q+1},[V_{-q+1},\g_{q-1}]] \neq 0$, as asserted. \ \qed

\m
\section {\ Restricted Lie algebras \label{sec:1.7}}  

\m Recall that a Lie algebra $L$ of prime characteristic is  {\it restricted}  (see 
\cite[Chap.~V, Sec.~7]{J} or \cite[Sec.~2.1]{SF})
if it has a mapping $[p]: L
\rightarrow L$,\  $x \mapsto x^{[p]}$,  satisfying

\begin{eqnarray}\label{eq:1.66} \ad x^{[p]} &=& (\ad x)^p\nonumber\\
(\xi x)^{[p]} &=& \xi^p x^{[p]}\\
(x+y)^{[p]} &=& x^{[p]} + y^{[p]} + \sum_{i=1}^{p-1} s_i(x,y),
\nonumber
\end{eqnarray}

\n for all $x,y \in L$, $\xi \in {\F}$, where $s_i(x,y)$ is given
by

$$\big(\ad(tx+y)\big)^{p-1}(x) = \sum_{i=1}^{p-1}t^{i-1}is_i(x,y)$$
in ${\F}[t] \ot_{\F} L$ (tensor symbols are omitted in writing
elements).  The summands  $s_i(x,y)$ are commutators in  $x,y$.

\m
The inspiration for the notion of the $[p]$-mapping comes from  the ordinary $p$th power in an
associative algebra.  Indeed, elements $X$ and $Y$ of an associative algebra satisfy
\begin{equation}\label{eq:ppower} (tX+Y)^{p} = t^pX^{p} + Y^{p} + \sum_{i=1}^{p-1} t^i s_{i}(X,Y)
\end{equation}
where $t$ is an indeterminate as above, and 
$s_i(X,Y)$ is a polynomial in $X,Y$ of total degree $p$.   Differentiating both sides with
respect to $t$ gives
\begin{equation*}  \sum_{j=0}^{p-1} (tX+Y)^{p-1-j}X (tX+Y)^j
= \sum_{i=1}^{p-1} t^{i-1} is_{i}(X,Y). \end{equation*}
Now because the relation $\displaystyle{{p-1 \choose j} = (-1)^j}$ holds in characteristic $p$,
and because $\ad(tX+Y)$ is the difference of the left multiplication and right multiplication operators
of $tX+Y$ which commute,  we have
\begin{eqnarray*} 
\big(\ad(tX+Y)\big)^{p-1}(X) &=& \sum_{j=0}^{p-1}{p-1 \choose j}  (-1)^j
(tX+Y)^{p-1-j}X (tX+Y)^j \\  &=& \sum_{j=0}^{p-1} (tX+Y)^{p-1-j}X (tX+Y)^j \\
&=& \sum_{i=1}^{p-1} t^{i-1} is_{i}(X,Y),  \end{eqnarray*}   
which is the analogue of the above equation in the associative setting.  
\m

A representation $\varrho: L \rightarrow \mathfrak {gl}(V)$ of a restricted
Lie algebra $(L, [p])$ is said to be  {\it restricted} if $\varrho(x^{[p]}) = \varrho(x)^p$
for all $x \in L$.    In particular, the adjoint representation of $L$ is
restricted by  \eqref{eq:1.66}.      The next result is apparent for
any restricted representation, but in fact  it holds
for {\em any} representation of a restricted Lie algebra.   
\m

\begin{Pro} \label{Pro:1.67} \ {\rm (Compare \cite{WK}.)}
 Let  $\varrho: L \rightarrow \mathfrak {gl}(V)$ be any
 representation of a restricted
Lie algebra $(L, [p])$.  Then 
$$ \big(\varrho(x+y)\big)^{p} -\varrho \big((x+y)^{[p]}\big) = \big(\varrho (x)\big)^{p} -\varrho (x^{[p]})
+ \big(\varrho (y)\big)^{p} -\varrho(y^{[p]}) $$ \n for all $x,y \in L$.  \end{Pro}

\pf  For elements $x$ and $y$ of a restricted Lie algebra $L$, we
have from \eqref{eq:ppower} with $X = \varrho(x)$, $Y = \varrho(y)$,  and $t = 1$, 
\begin{eqnarray}\label{eq:1.68} \big(\varrho(x+y)\big)^{p}
&=& \big(\varrho(
x)\big)^{p} + \big(\varrho(y)\big)^{p} + \sum_{k=1}^{p-1}s_{i}\big(\varrho(x), \varrho(y) \big)\\
&=& \big(\varrho(x)\big)^{p} + \big(\varrho(y)\big)^{p} + \varrho \left
(\sum_{k=1}^{p-1}s_{i}(x,y) \right).  \nonumber  \end{eqnarray}
\noindent But from \eqref{eq:1.66} it follows that
\begin{equation}\label{eq:1.69}\varrho \big((x+y)^{[p]}\big) = \varrho( x^{[p]}) + \varrho(
y^{[p]}) + \varrho \left (\sum_{k=1}^{p-1}s_{i}(x,y)\right)
\end{equation}
 \medskip
  \noindent
Subtracting \eqref{eq:1.69} from \eqref{eq:1.68}, we obtain
the desired relation.  \  \qed 

\bi 
\begin{Rem} \label{Rem:1.59}  {\rm Assume $\F$ is
an algebraically closed field of
characteristic $p > 0$.  For any $x$ in a restricted Lie algebra $L$ over $\F$, 
the element  $x^p-x^{[p]}$ is central in the
universal enveloping algebra $\mathfrak{U}(L)$  of $L$.  This can be readily
seen from the fact that $\ad x$,  when applied to $\mathfrak{U}(L)$,
is the difference of the left and right multiplication operators of $x$ on $\mathfrak{U}(L)$
and so  
\begin{eqnarray*} 0 &=& (\ad x)^p -\ad(x^{[p]})  = (L_x - R_x)^p -\ad(x^{[p]})\\
&=& (L_x)^p - (R_x)^p -\ad(x^{[p]})\\
&=& L_{x^p}-R_{x^p}-\ad(x^{[p]})  =  \ad \big( x^p-x^{[p]}\big).  \end{eqnarray*}
\noindent    This result leads to two important consequences.
The first is that if $L$ is finite-dimensional, then every irreducible  $L$-module
is finite-dimensional (see for example, \cite[Thm.~2.4]{SF}).  
The second is that $x^p-x^{[p]}$ must act as a scalar on any irreducible
representation $\varrho:L \rightarrow \mathfrak{gl}(V)$ of $L$ (hence of $\mathfrak{U}(L)$)
by Schur's Lemma.  
Let $\chi: L \rightarrow \F$ be such that
\begin{equation}\label{eq:centchar} \varrho(x)^p - \varrho(x^{[p]}) = \varrho\big(x^p - x^{[p]}\big) = \chi(x)^p\, \id \end{equation}
for all $x \in L$.   Then it follows from Proposition \ref{Pro:1.67}
that $\chi$ is linear.  Often $\chi$ is referred to
as the  $p$-character of $\varrho$.    The representation $\varrho$ is
restricted precisely when $\chi$ is identically zero. }   
    \end{Rem}
\bi

Subsequent sections of this
chapter  will be devoted to showing that under certain 
hypotheses on a graded Lie algebra $\g = \bigoplus_{j=-q}^r \g_j$, the adjoint representation of
$\g_0^{(1)}:=[\g_{0},\g_{0}]$ on each homogeneous component
$\g_{-j}$ for $j \geq 1$  is restricted.     Although a homogeneous
component $\g_{-j}$ may not be an irreducible $\g_0^{(1)}$-module, nonetheless
it will be shown to have a $p$-character.  
Our plan of attack will be
to demonstrate that these $p$-characters are zero. 
This will be a consequence of the main theorem on restrictedness
(Theorem \ref{Thm:1.48}), which we discuss next,  and of Lemma \ref{Lem:1.71} below. 
 
\bi \section {\ The main
 theorem on restrictedness (Theorem \ref{Thm:1.48}) \label{sec:1.8}}

  \m With these preliminaries at hand,
 we turn our attention now to this chapter's main goal, which will be
 to prove

 \bi \begin{Thm} \label{Thm:1.48} \ Suppose $\g = \bigoplus_{i=-q}^r
 \g_i$ is a finite-dimensional irreducible, transitive graded Lie
 algebra over an algebraically closed field ${\F}$ of
 characteristic $p > 3$ such that
 $[\g_{1},[\g_{-1},\g_1]]\ne 0$ and
 $\g_0$ is a direct
 sum of ideals, each of which is one of the following:
 \begin{enumerate}
 \item[{\rm(a)}] abelian,
 \item[{\rm(b)}] a classical simple Lie
 algebra (including, possibly, $\mathfrak{psl}_n$ with $p\mid n$),
  \item[{\rm(c)}] a Lie algebra isomorphic to $\mathfrak{sl}_n$,
  $\mathfrak{gl}_n$,  or $\mathfrak{pgl}_n$ with $p\mid n$.
  \end{enumerate}
  Then the representation of
  $\g_0^{(1)}:=[\g_{0},\g_{0}]$ on $\g_{-1}$ is restricted. \end{Thm}
\m

 \section {\ Remarks on restrictedness  \label{sec:1.9} }

  \m In Kac's version of the Recognition Theorem in \cite{K2}, the representation of
  $\g_0$ on $\g_{-1}$ was {\em assumed} to be restricted.  This
  hypothesis was shown to be unnecessary  in the depth-one case in \cite{G2} (see also \cite{G1} and  \cite{Ku1}),
  and in  \cite{BG},  which treated the arbitrary depth case  for $p>5$.
  Below we provide a new proof
  for $p>3$ based on the ideas from \cite{BG} and \cite{G2}.    Our proof will
 complete the original argument in \cite{BG}.    As shown in \cite{BKK}, the
  assumption
  on $p$ in the statement of Theorem \ref{Thm:1.48} cannot be relaxed even in
  the depth-one case. Here is a brief outline of how the argument will proceed.

 \m As $\g_0$ is
classical reductive, the sum $\g_0^\odot$ of the nonabelian ideals  is a
restricted Lie algebra under a certain natural $[p]$-mapping.  For
every $x \in \g_0^\odot$, $(\ad x)^p - \ad x^{[p]}$ must act as a
scalar, say $\chi(x)^p$,  on the irreducible $\g_0$-module
$\g_{-1}$ (compare Remark \ref{Rem:1.59}).  The representation will be restricted
if the $p$-character $\chi(x) = 0$ for all $x$.  First we prove that $\chi(x) = 0$
for all $x \in [\g_0^\odot,\g_0^\odot] = [\g_0,\g_0] = \g_0^{(1)}$ in the special case
that the depth $q$ of the graded Lie algebra equals one.  Next we
consider the depth-one Lie algebra $\bigoplus_{j\ge-1} \g_{jq}$
when $q \leq r$ or the algebra $\bigoplus_{j\le 1} \g_{jr}$ when
$q > r$.   We also need to work with the depth-one algebra
$\bigoplus_{j\ge-1} \g_{j(q-1)}$  for $q\le r$ or $\bigoplus_{j\le
1} \g_{j(r-1)}$ for $q > r$. Applying the depth-one result to the
quotients ${\mathcal B}(\g_{-q})$ and ${\mathcal B}(\g_{-q+1})$ of
the subalgebras $\bigoplus_{j\ge-1} \g_{jq}$ and
$\bigoplus_{j\ge-1} \g_{j(q-1)}$, which come from Proposition
\ref{Pro:1.25}, (or to their counterparts when $q > r$), and using
Lemma \ref{Lem:1.51} (a), we conclude $q\chi(x) = 0$ and
$(q-1)\chi(x) = 0$ for all $x$.  This forces $\chi(x) = 0$ for all
$x \in \g_0^{(1)}$, so that the representation is restricted as
asserted.  

\m \section {The action
of $\boldsymbol{\g_0}$ on $\boldsymbol{\g_{-j}}$  \label{sec:1.10}}

 \m Henceforth in this chapter,  we
assume that $\g$ is a graded Lie algebra satisfying the hypotheses
of Theorem \ref{Thm:1.48}.  Since by assumption, $\g_{-1}$ is an
irreducible $\g_0$-module, any central element of $\g_0$ must act
as a scalar on $\g_{-1}$ by Schur's Lemma.  By transitivity
\eqref{eq:1.3}, the center can be at most one-dimensional.
Consequently, there is at most one summand of type (a), and when
such a summand exists, it is the one-dimensional center of $\g_0$.
In this case, no summands isomorphic to $\mathfrak{sl}_n$ or
$\mathfrak{gl}_n$ with $p\mid n$ occur.  Similarly there is at
most one ideal isomorphic to $\mathfrak{sl}_n$ or
$\mathfrak{gl}_n$ with $p\mid n$, and if one occurs, there are no
summands of the abelian type (a).

\m Let $\tf$  be the sum of the maximal toral subalgebras of the summands of $\g_0$.
 Thus, $\tf$ contains the abelian summand if one occurs, and
$\tf$ is a maximal toral subalgebra, hence a Cartan subalgebra, of $\g_0$.  The algebra $\g_0$ has a
root space decomposition, $\g_0 =  \tf \oplus
\bigoplus_{\alpha\in{\tf}^* \setminus \{0\}} \, \g_0^{\alpha}$ relative to $\tf$, where
$$\g_0^{\alpha} = \{x \in \g_0 \mid [t,x] = \alpha(t)x \,\,\mbox{ for
all }\,\, t \in \tf \}.$$  By  \cite[Lem.~II.3.2]{S}, the root spaces of
the summands, hence of $\g_0$, are one-dimensional. This result
will be important for us; it holds under our assumption on $p$ but
may fail in characteristic $2$ or $3$ (this can be seen, for
example, by examining the  roots of the Lie algebra
$\mathfrak{psl}_3$ in characteristic $3$).

\m It follows from \cite[Lem.~II.3.2]{S}  that the set of roots of
$\g_0$ relative to $\tf$ is obtained by  mod $p$ reduction from a
genuine (reduced) root system $\Phi$ in the sense of \cite[Chap.~III]{H3}  or
\cite[Chap.~IV--VI]{Bou1}  (see Chapter 2 below for more detail).  Let
$\Delta$ be a basis of simple roots in $\Phi$. We say that a root
$\alpha$ of $\g_0$ relative to $\tf$ is {\em positive} (which we write
$\alpha > 0$) if $\al$ corresponds to an element from the positive
subsystem of $\Phi$ relative to $\Delta$.  A root is negative
($\alpha<0$) if $-\alpha>0$. We suppose that the vector
$e_{\alpha}$ spans the root space $\g_0^{\alpha}$.  For $\alpha
> 0$, we assume that the vectors $e_\alpha, e_{-\alpha}$, and
$h_{\alpha} := [e_{\alpha},e_{-\alpha}]$ determine a canonical
basis for a copy of $\mathfrak{sl}_2$.  Thus,
$$[h_\alpha, e_{\alpha}] = 2 e_{\alpha}, \qquad [h_\alpha,e_{-\alpha}] = -2 e_{-\alpha}.$$

\m

In classical algebras, there is a `standard' $p$th power map
satisfying $e_{\alpha}^{[p]} = 0$ for all roots $\alpha$, and
$h_{\alpha}^{[p]} = h_{\alpha}$ for all $\alpha > 0$.  (See
Chapter 2 below.)    Thus, when $\g$ satisfies the hypotheses of
Theorem \ref{Thm:1.48}, the standard $p${th} power map on the
various classical summands of $\g_0$ can be extended to the sum $\g_0^\odot$ of
the nonabelian ideals giving it the structure of a restricted Lie
algebra, which we assume is fixed from now on.

\m As in Remark \ref{Rem:1.59},   there exists a linear map
$\chi: \g_0^\odot \rightarrow {\F}$  such that for all $x \in \g_0^\odot$,

\begin{equation}\label{eq:1.70} (\ad x)^p - \ad x^{[p]} = \chi(x)^p \id
\end{equation}
\n on $\g_{-1}$.  Now when $\g_{-j} = [\g_{-j+1}, \g_{-1}]$ for all
$j=1,\dots,q$, we may assume by induction that $(\ad x)^p - \ad x^{[p]} =
(j-1)\chi(x)\id$ on $\g_{-j+1}$.  Then since $(\ad x)^p$ and $\ad
x^{[p]}$ are derivations, it follows from applying $(\ad x)^p -
\ad x^{[p]}$ to $\g_{-j} = [\g_{-j+1}, \g_{-1}]$ that $(\ad x)^p -
\ad x^{[p]} = j\chi(x)^p \id$ on $\g_{-j}$.  To summarize we have:
\bi
\begin{Lem} \label{Lem:1.71}\ {\rm (Compare \cite[Lem.~1]{BG}.)}
 Assume
$\g$ is a graded Lie algebra over an algebraically closed field
${\F}$ of characteristic $p > 3$ satisfying the conditions of
Theorem \ref{Thm:1.48}.  Let $\g_0^\odot$ be the sum of the nonabelian
ideals of $\g_0$. 
\begin{enumerate}
\item[{\rm (i)}]  Then there is a linear functional $\chi$ on
$\g_0^\odot$ such that $(\ad x)^p - \ad x^{[p]} = \chi(x)^p \id$ on
$\g_{-1}$ for all $x \in \g_0^\odot$.  \item[{\rm (ii)}] If $\g_{-j} =
[\g_{-j+1}, \g_{-1}]$ for $j=1,\dots,q$, then $(\ad x)^p - \ad
x^{[p]} = j\chi(x)^p \id$ on $\g_{-j}$ for all $x \in \g_0^\odot$.
\end{enumerate} \end{Lem}

\bi  When $\g$ is a graded Lie algebra satisfying the conditions
of Theorem \ref{Thm:1.48}, the maximal toral subalgebra $\tf$ of $\g_0$ is
abelian, and so any finite-dimensional $\g_0$-module $M$
decomposes into \ {\em weight spaces} (common generalized
eigenspaces) relative to $\tf$.  Thus, if $\varphi: \g_0
\rightarrow \mathfrak{gl}(M)$ is the corresponding representation
and $d = \dim M$, then
\ $M = \bigoplus_{\lam \in \tf^*} M^\lam$ where
$$M^{\lam} =
\{v \in M \mid \big(\varphi(t)-\lam(t)\id \big)^d\,v = 0\,\,
\mbox{ for all }\,\, t \in \tf\}.$$   The functional $\lam \in
\tf^{*}$ is a {\em weight} of $\tf$ on $M$ if $M^{\lam}\neq 0$.  If
$v \in M^\lam$ and $\big(\varphi(t)-\lam(t)\id \big)^d v = 0$, then
$\big(\varphi(t)-(\lam+\al)(t)\id\big)^d x.v = 0$ for all root
vectors $x \in \g_0^{\al}$ so that

\begin{equation}\label{eq:1.72} \g_0^{\al\,}.M^{\lam} \subseteq M^{\lam+\al}.
\end{equation}

\n Moreover, the sum $M' =
\bigoplus_{\lam \in \h^*} (M^\lam)'$ of the genuine eigenspaces 
 $(M^\lam)'= \{v \in M \mid
\big(\varphi(t)-\lam(t)\id\big)v = 0$ for all $\,t \in \tf\}$ is a
nonzero $\g_0$-submodule of $M$.  As a consequence, we have the
following well-known result:
\bi

\begin{Lem} \label{Lem:1.73} \ On any (finite-dimensional)
irreducible $\g_0$-module, the maximal toral subalgebra $\tf$ acts
semisimply.  \end{Lem} \m

\section {The depth-one
case of Theorem 1.63   \label{sec:1.11}}

\m We begin tackling Theorem \ref{Thm:1.48} by first proving it
when the depth $q=1$.  Initially we show under the hypotheses of
that theorem and the condition $q=1$, that the $\g_0$-module
$\g_{-1}$ possesses a common 
eigenvector for the adjoint action of the solvable subalgebra  $\mathfrak b^+ = \tf \oplus \bigoplus_{\alpha > 0} \g_0^\alpha$ 
of $\g_0$  (and by symmetry,  a common eigenvector for the adjoint action of
$\mathfrak b^- = \tf \oplus \bigoplus_{\alpha > 0} \g_0^{-\alpha}$).
In the process of establishing this, we obtain valuable
information about how root vectors of $\g_0$ act on $\g_{-1}$.

\m  From Lemma \ref{Lem:1.73}, we know that the maximal toral subalgebra
$\mathfrak t$ acts semisimply on $\g_{-1}$, say $\g_{-1} =
\bigoplus_{\lam\in {\mathfrak t}^*}\, \g_{-1}^{\lam}$ is its eigenspace
(weight space) decomposition.  We want to prove that for some
$\lam$, there exists a nonzero $v \in \g_{-1}^{\lam}$ such that
$[x,v] = 0$ for all $x \in \g_{0}^{\al}$, \  $\al > 0$.  Such
a vector $v$ will generate $\g_{-1}$ by irreducibility
\eqref{eq:1.2} and will be a common eigenvector for $\ad \mathfrak b^+$.
In particular, we show
\bi

\begin{Thm} \label{Thm:1.75} \ Under the assumptions of
Theorem \ref{Thm:1.48} and the hypothesis that the depth $q = 1$,
the $\g_0$-module $\g_{-1}$ contains a common eigenvector
for $\mathfrak b^+ = \tf \oplus \bigoplus_{\alpha > 0} \g_0^\alpha$ 
(and
by symmetry, for $\mathfrak b^- = \tf \oplus \bigoplus_{\alpha > 0} \g_0^{-\alpha}$).
\end{Thm}

\pf  Recall from Lemma \ref{Lem:1.71} (i) that there is a
linear functional $\chi$ on the sum $\g_0^\odot$ of the nonabelian
ideals such that (see \eqref{eq:1.70}) $(\ad x)^p - \ad x^{[p]} =
\chi(x)^p \id$ holds on $\g_{-1}$ for all $x \in \g_0^\odot$.  If $x \in
\g_0^{\al}$ is a root vector of $\g_0$, then $x^{[p]} = 0$, and
$(\ad x)^p = \chi(x)^p \id$ on $\g_{-1}$. If $\chi(e_{\alpha}) = 0$
for all $\al > 0$, then the subspace
$\mathfrak{n}^+:=\,\bigoplus_{\al > 0}\, \g_0^{\al}$ is a Lie
algebra acting as nilpotent transformations on $\g_{-1}$, so it
possesses a common eigenvector (compare Remark \ref{Rem:1.56}).
As the space of those eigenvectors is invariant under $\ad \tf$,
there will be a common eigenvector for $\ad \mathfrak b^+$ in $\g_{-1}$.  We may assume
then that $\chi(e_{\alpha}) \neq 0$ for some $\al >0$.

Let $M$ be any irreducible $\g_{0}$-submodule of $\g_{1}$. Then
the product $P:=[\g_{-1},M]$ is a nonzero ideal of $\g_{0}$ by
transitivity \eqref{eq:1.3}. Recall that all root spaces of $\g_0$
relative to $\tf$ are one-dimensional. Hence if $P\not\subseteq\tf$,
then $P$ contains a root vector, $e_{\gamma}$ say. In this case
$e_{-\gamma}\,\in[\g_0^{-\gamma},[\g_0^{-\gamma},e_{\gamma}]]$ is
in $P$ as well. So it can be assumed further that $\gamma>0$. As
the Lie algebra $\mathfrak{n}^+$  is nilpotent, the
$\text{ad}\,\mathfrak{n}^+$-module generated by $e_\gamma$
contains a root vector $e_{\delta}$ (with $\delta >0$)
such that $[\mathfrak{n}^+,e_{\delta}]=0$. If $P\subseteq
\tf$, then each root of $\g_0$ relative to $\tf$ vanishes on $P$.
Since $\tf$ is abelian, this shows that $P$ is contained
in the center $\mathfrak{Z}(\g_0)$ of $\g_0$. Since
$\dim\,\mathfrak{Z}(\g_0)\le 1$, it must be that $P=\mathfrak{Z}(\g_0)$ in this case.

If $P\not\subseteq \tf$,  we set $e=e_{\delta}$. If $P\subseteq
\tf$,  we take for $e$ any nonzero element in $\mathfrak{Z}(\g_0)=P$.  In both
cases, we let $\beta$ denote the weight of $e$ (so that
$\beta\in\{0,\delta\}$). It follows from Lemma
\ref{Lem:1.73} that $e \in \sum_{\lam\in\tf^*}\, [\g_{-1}^{\lam},
M^{\be-\lam}]$.  For some $\lam$, we must have $[\g_{-1}^{\lam},
M^{\be-\lam}] = {\F} e$ if $\be$ is a root, because the root space is
one-dimensional.   Then $e=[u,v]$ for some $u \in \g_{-1}^{\lam}$,
$\,v \in M^{\be-\lam}$.  If $\be = 0$,
we may assume $e$ has that form.    By transitivity, $[\g_{-1}, e] \neq 0$,
and since $\g_{-1}$ is spanned by weight vectors, there exists $w
\in \g_{-1}^\nu$ for some $\nu$ such that $[w,e] \neq 0$.

Since the depth $q=1$, we have  $[(\ad e_\alpha)^k(w), u]=0$ for
all $k\ge 0$. For $k=0,1,\dots,p-1$ we set $x_k=[(\ad
e_\alpha)^k(w),v]$ and observe that $x_k\in\g_0^{\eta+k\alpha}$
where $\eta=\beta-\lam+\nu$. Since $e$ commutes with $e_\alpha$
and $(\ad e_\alpha)^p\big([w,e]\big) = \chi(e_\alpha)^p [w,e]$, a
nonzero multiple of $[w,e]$, it must be that
$$[u, x_k]= [(\ad e_\alpha)^k(w),[u,v]]=
(\ad e_\alpha)^k \big([w,e]\big)\ne 0$$ for all $k\ge 0$.  But
then $\g_0^{\eta+k\alpha}\ne 0$ for $k=0,1,\dots, p-1$.  If
$\eta$ is not a multiple of $\al$, this is impossible, since
root strings in classical Lie algebras have length at most 4. If
$\eta$ is a multiple of $\al$, then all multiples $k\al$ for
$k=1, \dots, p-1$ are roots, which again cannot happen for $p > 3$
(both statements follow from the Mills-Seligman axioms which
characterize direct sums of classical simple Lie algebras; see
Theorem \ref{Thm:2.6} of Chapter 2 for more details). We have reached a contradiction, so we
are forced to draw the following two conclusions:  Under the
hypotheses of Theorem \ref{Thm:1.48} and the assumption $q=1$,
\begin{eqnarray*}
&&\chi(x) =  0 \ \ \hbox{\rm for all} \ \ x \in \g_0^{\al}\ \
\hbox{\rm and all}
\ \al > 0; \\
&& \hbox{\rm the} \ \g_0\hbox{\rm-module} \ \ \g_{-1} \ \
\hbox{\rm has a common eigenvector for $\ad \mathfrak b^+$} \end{eqnarray*}

\noindent By symmetrical arguments,
\begin{eqnarray*}
&&\chi(x) =  0 \ \ \hbox{\rm for all} \ \ x \in \g_0^{-\al}\ \
\hbox{\rm and all}
\ \al > 0; \\
&& \hbox{\rm the} \ \g_0\hbox{\rm-module} \ \ \g_{-1} \ \
\hbox{\rm has a common eigenvector for $\ad \mathfrak b^-$.} \quad  \square \end{eqnarray*}

\m

\section  {Proof of
Theorem 1.63 in the depth-one case  \label{sec:1.12}}

 \m
Recall,  as in the proof of Theorem \ref{Thm:1.75},  that when $\g$ is
a graded Lie algebra over an algebraically closed field ${\F}$ of
characteristic $p > 3$ satisfying the hypotheses of Theorem
\ref{Thm:1.48},  there is a linear functional $\chi$ on the
sum $\g_0^\odot$ of the nonabelian ideals such that  $(\ad x)^p - \ad
x^{[p]} = \chi(x)^p \id$ holds on $\g_{-1}$ for all $x \in \g_0^\odot$.
{F}rom our proof of Theorem \ref{Thm:1.75} we know that
$\chi(x) = 0$ for all root vectors $x$ of $\g_0$ when $q = 1$ .

Suppose that $e = e_{\alpha}\in \g_0^{\al}$ and $f = e_{-\alpha}
\in \g_{0}^{-\al}$ are such that they along with $h = [e,f] \in
\tf$ determine a canonical basis of $\mathfrak{sl}_2$ as before.
Then since $\exp(\ad e)$ is an automorphism of $\g_0$ (see Chapter
2 to follow for the details), the element $f' = \exp(\ad e)(f) = f
+ h -e $ is a root vector relative to the maximal toral subalgebra $\tf' =
\exp(\ad e)(\tf)$.  As above, \ $\chi(f') = 0$.  Then from
Proposition \ref{Pro:1.67} we deduce that on $\g_{-1}$,
 
\begin{eqnarray*} 0 &=& (\ad f')^{p} -\ad (f')^{[p]} \\
&=& \big (\ad(f + h - e)\big)^{p} -\ad(f + h -
e)^{[p]} \\
&=& (\ad f)^{p} -\ad f^{[p]} + (\ad h)^{p} -\ad h^{[p]}
- (\ad e)^{p} + \ad e^{[p]}  \\
&=& (\ad h)^{p} -\ad h ^{[p]}.  \end{eqnarray*}

\n Thus, $\chi(h) = 0$.  Since the elements $e = e_\al,f = e_{-\al},h = h_\al$, as $\al$
ranges over all the roots of $\g_0$,  span $[\g_0^\odot,\g_0^\odot] =
 \g_0^{(1)}$, we may use the semilinearity of the $[p]$-mapping in
\eqref{eq:1.66} to see that the representation of $\g_0^{(1)} $ on
$\g_{-1}$ is restricted.  Consequently, we have proven
\m

\begin{Lem} \label{Lem:1.77} \ If $\g = \bigoplus_{i=-1}^r
\g_i$ is a depth-one graded Lie algebra over an algebraically
closed field of characteristic $p > 3$ satisfying the
assumptions of Theorem \ref{Thm:1.48}, then the representation of
$\g_0^{(1)}$ on $\g_{-1}$ is restricted.  \end{Lem}
\medskip

\begin{Rem} \label{irr}
{\rm If $M$ is an irreducible $\g_0$-module, which is restricted
as a $\g_0^{(1)}$-module, then $M$ is already irreducible over
$\g_0^{(1)}$. Indeed, let $\mathfrak{U}_{[p]}(\mathfrak n^-)$ denote the
restricted enveloping algebra of $\mathfrak{n}^- :=  \bigoplus_{\alpha > 0} \g_0^{-\alpha}$, and let $J$ be
the augmentation ideal of $\mathfrak{U}_{[p]}(\mathfrak n^-)$. It is well-known
(and easy to see) that $J$ is the Jacobson radical of
$\mathfrak{U}_{[p]}(\mathfrak n^-)$. Clearly, the subspace $JM$ is
$\tf$-stable.    By Nakayama's lemma, $J.M\ne M$. Let $v_0$ be a
common eigenvector for $\mathfrak b^+ = \tf \oplus \mathfrak n^+$ in $M$, and let $\lambda\in\tf^*$ be the
weight of $v_0$. Since $M= \mathfrak{U}^{[p]}(\mathfrak{n}^-)v_0$ by the
irreducibility of $M$, we have $v_0\not\in JM$. It follows that
$JM$ has codimension one in $M$. Let $M_0=\{v\in M \mid
\mathfrak{n}^+\,v=0\}.$ The subspace $M_0$ is $\tf$-stable, hence
decomposes into weight spaces for $\tf$. If $\mu\ne \lambda$, then
$M_0^\mu\subseteq JM$, and hence $\mathfrak{U}_{[p]}(\mathfrak n^-)M_0^\mu$ is
a $\g_0$-submodule of $M$ contained in $JM$. Therefore,
$M_0=M_0^\lambda$. If $\dim M_0^\lambda>1$, then
$\mathfrak{U}_{[p]}(\mathfrak n^-)(M_0^\lambda \cap JM)$ is a nonzero
$\g_0$-submodule of $M$ contained in $JM$. Since this contradicts
the irreducibility of $M$, we derive that $M_0=\F\,v_0$. 
As $\mathfrak n^+ \subset \g_0^{(1)}$, every  $\g_0^{(1)}$-submodule $N$ of $M$ 
must contain $v_0$.  But then
$N \supseteq \mathfrak{U}_{[p]}(\mathfrak n^-)v_0 = M$, to show $M$ is irreducible
over $\g_0^{(1)}$.}
\end{Rem} \m
\section {Quotients of $\boldsymbol{\g_0}$\label{sec:1.13} }
 \m

In demonstrating Theorem \ref{Thm:1.48} for arbitrary depths,  we will require  some
knowledge about quotients of $\g_0$.   

\bi \begin{Lem} \label{Lem:1.79} \ Suppose $\g =
\bigoplus_{i=-q}^r\, \g_i$ is a graded Lie algebra satisfying the
hypotheses of Theorem \ref{Thm:1.48},  and let $A_0$ be a
nonzero ideal of
$\g_0$. Then
the following are true:
\begin{itemize}
\item[{\rm (i})]  If  $\dim\,\mathfrak{Z}(\g_0/A_0)\le 1$, then
the Lie algebra
$\g_0/A_0$ is isomorphic to a direct sum of ideals of type (a),
(b), (c);

\item[{\rm (ii)}]
Any solvable ideal of $A_0$ is contained in ${\mathfrak Z}(\g_0)$.
\end{itemize}
\end{Lem}

\pf (i)   We know that $\g_0 = \bigoplus_{i=1}^\ell \g_0^{[i]}$, where each
$\g_0^{[i]}$ is of type (a), (b), or (c).   Let
$\pi_i$ denote the canonical projection from $\g_0$ onto $\g_0^{[i]}$ for
$1\le i\le \ell$. We may assume after possibly renumbering that
for some $k\ge 0$ we have
$\pi_i(A_0)\not\subseteq{\mathfrak Z}(\g_0^{[i]})$ if $i\le k$ and
$\pi_i(A_0)\subseteq {\mathfrak Z}(\g_0^{[i]})$
if $k < i\le \ell$. Since ${\mathfrak Z}(\g_0^{[i]})\subseteq {\mathfrak Z}(\g_0)$,  and
any noncentral ideal of $\g_0^{[i]}$ contains $(\g_0^{[i]})^{(1)}$, it must be that
\begin{equation}\label{eq:1.785} (\g_0^{[1]})^{(1)}\oplus \cdots \oplus (\g_0^{[k]})^{(1)}\subseteq A_0\subseteq
\g_0^{[1]}\oplus \cdots \oplus \g_0^{[k]} \oplus \mathfrak Z(\g_0^{[k+1]}) \oplus
\cdots \oplus  \mathfrak Z(\g_0^{[\ell]}) . \end{equation}
Set $K : = (\g_0^{[1]})^{(1)}\oplus \cdots \oplus (\g_0^{[k]})^{(1)}$. 
For each $i$, let  $\eta_i: \g_0^{[i]} \rightarrow \g_0^{[i]}/ (\g_0^{[i]})^{(1)}$ be
 the canonical map, and set $\psi = \sum_{i=1}^k \eta_i \circ \pi_i + \sum_{t=k+1}^\ell \pi_t$. 
 Then 
 $$\psi: \g_0 \rightarrow Q := \left(\bigoplus_{i=1}^k  \g_0^{[i]}/ (\g_0^{[i]})^{(1)}\right)
 \oplus \left(\bigoplus_{t=k+1}^\ell \g_0^{[t]}\right),$$ 
 and the kernel of $\psi$ is $K$, so there is an induced isomorphism $\ov \psi: \g_0/K \rightarrow Q$.
 Now if $S := \bigoplus_{i=1}^k  \g_0^{[i]}/ (\g_0^{[i]})^{(1)}$,  then 
 $$R: = \ov \psi (A_0/K) \subseteq S \oplus \mathfrak Z(\g_0^{[k+1]}) \oplus
\cdots \oplus  \mathfrak Z(\g_0^{[\ell]}) = \mathfrak Z(Q).$$ 
Thus, 

$$(Q/R)/\big(\mathfrak Z(Q)/R\big) \cong Q/\mathfrak Z(Q) \cong \bigoplus_{t=k+1}^\ell \g_0^{[t]}/ \mathfrak Z(\g_0^{[t]}),$$  which
is semisimple and is the direct sum of ideals of type (b)  or (c) (only $\mathfrak{pgl}_n$
with $p \mid n$ is possible in case (c)).  
As $\mathfrak Z(Q)/R \subseteq \mathfrak Z(Q/R)$, we must have equality,
$\mathfrak Z(Q)/R  = \mathfrak Z(Q/R)$.

Now \, $\g_0/A_0 \cong (\g_0/K)/(A_0/K) \cong Q/R$, \, so if
\,$\dim \mathfrak Z(\g_0/A_0) \leq 1$,\, then \break $\dim \mathfrak Z(Q)/R \leq 1$.  
In particular, if $R \subseteq S$, then $\g_0/A_0 \cong Q/R \cong S/R \oplus \bigoplus_{t=k+1}^\ell \g_0^{[t]}$,
which is the sum of ideals of type (a), (b), or (c).   If  $R \not \subseteq S$, then we can write
$R = (R \cap S) \oplus \F (s+z)$  where  $s \in S$,  $z \in \mathfrak Z(\g_0^{[j]}) = \mathfrak Z(\g_0)$
(see the comments in Section \ref{sec:1.10}) for some $j$,  and $z \neq 0$.    Moreover,
$Q = (R \cap S) \oplus \F (s+z) \oplus \F s' \oplus  \bigoplus_{t=k+1}^\ell \g_0^{[t]}$ for
some nonzero $s' \in S$.
Thus, $\g_0/A_0 \cong Q/R \cong \F s'   \oplus  \bigoplus_{t=k+1}^\ell \g_0^{[t]}$,
which is a direct sum of ideals of type (a), (b), or (c).     
\m

\noindent
(ii) With notation as in the proof of (i), observe first  that $A_0^{(1)} = K =  \bigoplus_{i=1}^k (\g_0^{[i]})^{(1)}$. 
Now let $J$ be a solvable ideal
of $A_0$. Then $[J,A_0]$ is
a solvable ideal of $A_0^{(1)}$. Since each nonzero $(\g_0^{[i]})^{(1)}$
is nonsolvable,  and any proper ideal of $(\g_0^{[i]})^{(1)}$ lies in
${\mathfrak Z}(\g_0)$,
we  obtain $[J,A_0]
\subseteq {\mathfrak Z}(\g_0)$.    Then  $\pi_i(J)$ is a subalgebra
of $\g_0^{[i]}$ satisfying $$[\pi_i(J),(\g_0^{[i]})^{(1)}]\subseteq \pi_i([J, A_0])\subseteq
\mathfrak Z(\g_0)\cap \pi_i(J).$$
It is easy to see that $\g_0^{[i]}$ contains a maximal
toral subalgebra $\tf^{[i]}$ such that
\begin{itemize}
\item[(1)]
$(\tf^{[i]})'\eqdef \tf^{[i]}\cap(\g_0^{[i]})^{(1)}$ is a maximal toral
subalgebra of $(\g_0^{[i]})^{(1)}$;
\item[(2)] any root space of
$(\g_0^{[i]})^{(1)}$ relative to $(\tf^{[i]})'$ is one-dimensional;
\item[(3)] no root vector of $(\g_0^{[i]})^{(1)}$ relative to $(\tf^{[i]})'$
is central in $\g_0^{[i]}$;
\item[(4)] $\g_0^{[i]}=\tf^{[i]}+(\g_0^{[i]})^{(1)}$.
\end{itemize}
Decomposing $\pi_i(J)$ into weight spaces relative to $(\tf^{[i]})'$,  one
observes readily that $\pi_i(J)\subseteq \tf^{[i]}$. From this it follows that
$\pi_i(J)$ is a solvable ideal of $\g_0^{[i]}$. Then
$\pi_i(J)\subseteq{\mathfrak Z}(\g_0)$ by our remarks earlier in the proof.
But then $J\subseteq {\mathfrak Z}(\g_0)$, as desired.  \qed

\m \smallskip

 \section{The proof of Theorem 1.63  when \,  $\boldsymbol{2 \leq q \leq r}$  \label{sec:1.14} }

\m Having established our result in the special case of depth-one
Lie algebras, we turn our attention now to graded Lie algebras $\g
=\bigoplus_{i=-q}^{r}\, \g_i$ with $q \geq 2$ satisfying the
requirements of Theorem \ref{Thm:1.48}.   In
analyzing such Lie algebras, we will make extensive
use of the depth-one graded Lie algebras $\mathcal B(\g_{-q})$
and $\mathcal B(V_{-q+1})$ described in Section \ref{sec:1.4}.

\bi {\em Suppose now that \ $\g = \bigoplus_{j=-q}^r \g_j$ \ is
a graded Lie algebra with $q \geq 2$ satisfying the assumptions
of Theorem \ref{Thm:1.48}.    Replacing $\g_{-j}$ by $\g_{-1}^j = $ \allowbreak $(\ad
\g_{-1})^{j-1}\,\g_{-1}$ for $j=2,\dots,q$ if necessary, we may
suppose that $\g_{-j} = \g_{-1}^j$.  Moreover, factoring out the
Weisfeiler radical ${\mathcal M}(\g)$ we may assume also that
${\mathcal M}(\g) = 0$ (compare (ii) of Proposition
\ref{Pro:1.10}).  Neither of these steps alters the local part
$\g_{-1}\oplus \g_0 \oplus \g_1$ (hence the assumptions of
transitivity \eqref{eq:1.3} and irreducibility \eqref{eq:1.2} are
still intact as are the constraints on $\g_0$).   In this section  we will
suppose that $q \leq r$.} \m

By Lemma \ref{Lem:1.51} (b) we know that $\g_{-q}$ is an
irreducible $\g_0$-module.  Moreover, if $F = F(\g_{-q}) =
\bigoplus_{i\ge-1}\,{F_i}$ and $\mathcal {B}={\mathcal B}(\g_{-q})
= \bigoplus_{i\ge -1} {\mathcal B}_i = F / A$, then by Lemma
\ref{Lem:1.58}, ${\mathcal B}_{1} \neq 0$.  As ${\mathcal B}_0 =
\g_0/A_0$, it follows from Lemma \ref{Lem:1.79} that all the
conditions of Theorem \ref{Thm:1.48} are satisfied by the
depth-one algebra ${\mathcal B}$. Therefore, we know from Lemma
\ref{Lem:1.77} that $(\ad b)^p - \ad b^{[p]}$ = 0 on ${\mathcal 
B}_{-1}$ for all $b \in {\mathcal B}_0^{(1)}$, or
equivalently $\left((\ad x)^{p} -\ad x^{[p]} \right)(\g_{-q}) = 0$
for all $x \in \g_0^{(1)} $.  But $(\ad x)^{p} -\ad x^{[p]} =
q\chi(x)^{p}\id$ on $\g_{-q}$ by Lemma \ref{Lem:1.71} (ii), so
either $q\equiv 0\  \mod p$ or $\chi(x) = 0$ for all $x \in
\g_0^{(1)} $.  If $q \not \equiv 0\  \mod p$, then $\chi(x) = 0$,
and we have our desired conclusion that the representation of
$\g_0^{(1)} $ on $\g_{-1}$ is restricted.  We will assume then
that $q \equiv 0 \, \mod p$. Since we are also supposing $p \geq
5$, we have $q \ge 5$ (in fact $q > 3$ will suffice for our
purposes).

\m

Now we let $V_{-q+1}$ be an irreducible submodule of the
$\g_0$-module $\g_{-q+1}$. Applying Lemma \ref{Lem:1.63} (which
assumes $3 < q \leq r$), we determine that $(\ad
V_{-q+1})^{2}\g_{q-1} \neq 0$ and hence ${\mathcal B}'_1 \neq 0$
for ${\mathcal B}' = {\mathcal B}(V_{-q+1})$.  In view of Lemma
\ref{Lem:1.77},  it follows from this that $(q-1)\chi(x) = 0$ for
all $x \in \g_0^{(1)} $.  But we are assuming $q \equiv 0 \, \mod
p$, so $\chi(x) = 0$ for all $x \in \g_0^{(1)}$.  As a
consequence, we have established

 \bi \begin{Lem} \label{Lem:1.80} \ Let $\g = \bigoplus_{i=-q}^r \g_i$ be
a graded Lie algebra satisfying the hypotheses of Theorem
\ref{Thm:1.48}, and assume $q \le r$.  Then the representation of
$\g_0^{(1)}= [\g_{0},\g_0]$ on $\g_{-1}$ is a restricted representation.
\end{Lem} \m

 \section  {The proof of Theorem 1.63 when \, $\boldsymbol{q > r}$ \label{sec:1.15}}

As in the argument above, in order to show $\g_0^{(1)}$ has a
restricted action on $\g_{-1}$, it suffices to work with the
subalgebra $\widehat{\g}$ of $\g$ generated by the local part
$\g_{-1} \oplus \g_{0} \oplus \g_{1}$ of $\g$.  We may suppose
$\widehat \g = \bigoplus_{i=-q'}^{r'} \widehat \g_{i}$ where
$\widehat \g_{-q'} \neq 0$ and $\widehat \g_{r'} \neq 0$,
$\widehat{\g}_{0} = \g_{0}$, and $\widehat{\g}_{\pm i} =
(\g_{\pm 1})^{i}$ for $i > 0$ (compare \eqref{eq:1.44}).  Then if we factor
out the Weisfeiler radical of $\widehat{\g}$ (contained in
$\sum_{j \le -2}\, \widehat{\g}_{j}$), we do not affect the local
part of $\widehat{\g}$ nor do we alter the action of $\ad \,
\widehat{\g}_{0}$ on $\widehat{\g}_{-1}$ or the transitivity
\eqref{eq:1.3}.  Under these assumptions, $\g$ is 1-transitive in
its negative part by Lemma \ref{Lem:1.15}. With these reductions,
we may suppose that the following holds.
   \m

\section {{General setup}  \label{sec:1.16}} \begin{enumerate}
\item[{\rm (i)}] $\g = \bigoplus_{i=-q}^r \g_i$ is a transitive,
irreducible graded Lie algebra generated by its local part
$\g_{-1} \oplus \g_0 \oplus \g_1$; \item[{\rm (ii)}] $\M(\g) = 0$;
\item[{\rm (iii)}] $q > r$; \item[{\rm (iv)}] $\g$ is 1-transitive
in its negative part (i.e., $[\g_{1},x] = 0$ for $x \in \g_{-}$
implies $x = 0$); \item[{\rm (v)}] the unique minimal ideal is
${\mathcal I} = \bigoplus_{i=-q}^s {\mathcal I}_i$;  \item[{\rm
(vi)}] $\g_0$ is the sum of ideals satisfying (a)--(c) of Theorem
\ref{Thm:1.48}; \ and \item[{\rm (vii)}] $[[\g_{-1},\g_1],\g_1]\ne
0$.
\end{enumerate}
  \m

 {\em Our implicit assumption in all subsequent results in 
 Section \ref{sec:1.16} 
 is that
 (i)-(vii) hold.}   \bi

\begin{Lem}\label{Lem:1.83}  We have $\Ann_{\g_{0}}(\g_{1}) = 0$  so that 
$\g$ is $1$-transitive, and $s>0$.
\end{Lem}

\pf Suppose that $Q$ $\eqdef \Ann_{\g_{0}}(\g_{1})$ $\ne 0$.  By
transitivity \eqref{eq:1.3} and irreducibility \eqref{eq:1.2},
$[\g_{-1},$ $Q]= \g_{-1}$. Therefore, $[\g_{-1},\g_{1}] =
[[\g_{-1}, Q], \g_{1}]= [[\g_{-1}, \g_{1}], Q] \subseteq Q$. But
then $[[\g_{-1},$ $\g_{1}],$ $\g_{1}]= 0$ contrary to
assumption (vii).  This gives $Q=0$, which
in conjunction with (iv) shows that  $\g$ is 1-transitive.  
Now if $s = 0$, then
$[{\mathcal I}_0, \g_1]\subseteq {\mathcal I} \cap \g_1= 0,$ hence
${\mathcal I}_0=0$. Since $0\ne [\g_{-1},\g_1]\subseteq {\mathcal 
I}_0$, this is impossible.   This demonstrates that $s > 0$ must hold. \
\qed

\m  \smallskip

\begin{Lem} \label{Lem:1.84} \ $\g_{r} \cap {\mathcal I}= 0$ if and only if
$\g_{r}$ is a trivial $[\g_{-1},\g_{1}]$-module.  In this case,
$[\g_1,\mathcal I_{r-1}] = 0$ and $\mathcal I_{r-1}$ is a nontrivial
$[\g_{-1},\g_{1}]$-module.  \end{Lem}

\pf  {F}rom Lemma \ref{Lem:1.20} it follows that the
 unique minimal graded ideal
 ${\mathcal I}$ of $\g$ satisfies ${\mathcal I}_{-i} = \g_{-i}, \ 1 \le
 i \le q$, and ${\mathcal I}_{j} = [\g_{-1},{\mathcal I}_{j+1}], \ -q
 \le j \le r - 2$.  Now if $\g_{r} \cap {\mathcal I} = 0$, then since
 $\g_{-1} \subseteq {\mathcal I}$,
\m
$$[[\g_{-1},\g_{1}],\g_{r}] \subseteq [{\mathcal I}_{0},\g_{r}] \subseteq
\g_{r} \cap {\mathcal I} = 0,$$ \m so that $\g_{r}$ is a trivial
$[\g_{-1},\g_{1}]$-module.  Conversely, if $\g_{r}$ is a
trivial $[\g_{-1},\g_{1}]$-module, then $\sum_{i\geq 1}\,(\ad \,
\g_{-1})^{i}\g_{r}$ is an ideal of $\g$ contained in ${\mathcal 
I}$, which by minimality must equal ${\mathcal I}$. Consequently,
${\mathcal I} \cap \g_r = 0$.

Now if 
$\g_{r}$ is a trivial $[\g_{-1},\g_{1}]$-module  (equivalently,
if $\g_{r} \cap {\mathcal I}= 0$),  then
$\sum_{i\geq 1}\,(\ad \,
\g_{-1})^{i}\g_{r} = \mathcal I$ as noted above,  to
force ${\mathcal I}_{r-1} = [\g_{-1},\g_{r}]$.   
As a result,   $[\g_1, \mathcal I_{r-1}] = [\g_1, [\g_{-1},\g_r]] =
[[\g_1,\g_{-1}], \g_r] + [\g_{-1}, [\g_1,\g_r]] = 0$.  
If  ${\mathcal I}_{r-1}$ happens to be a trivial
$[\g_{-1},\g_{1}]$-module,  then we would have \m
$$0 = [[\g_{-1},\g_{1}],{\mathcal I}_{r-1}] =
[[\g_{-1},\g_{1}],[\g_{-1},\g_{r}]] = [[[\g_{-1},\g_{1}],
\g_{-1}],\g_{r}] = [\g_{-1},\g_{r}]$$ \m by the transitivity
\eqref{eq:1.3}  and irreducibility \eqref{eq:1.2} of $\g$, to
contradict transitivity. \ \qed

\bi We can use the $1$-transitivity established in
Lemma \ref{Lem:1.83} to argue that $\g_{r}$ annihilates no nonzero
element in $\g_{\leq 0}$.   In fact, we can prove a  somewhat
stronger result.    Let $M$ be any graded
$\g_{0}$-submodule of $\g$ such that $[\g_{1},M]$ = $0$.  By 1-transitivity,
$M \subseteq  \g_+$.
Suppose that $[M,x] = 0$ for some nonzero $x \in \g_{-}$. Clearly,
$M$ annihilates each of the homogeneous summands $x_{-i}$ $\in
\g_{-i}$ of $x$. For any $j$ such that $x_{-j}$ $\neq$ $0$, set
$N_{-j}$ $\eqdef \left\{y \in \g_{-j} \mid [M, y] = 0\right\}$, which is a
nonzero $\g_{0}$-submodule of $\g_{-j}$.  Then by  Section \ref{sec:1.16}\,(iv),
$$(\ad \, \g_{1})^{j-1}[N_{-j},M] = [(\ad \,
\g_{1})^{j-1}N_{-j},M] = [\g_{-1},M] \ne 0,$$ contrary to the
definition of $N_{-j}$.  Consequently,
$\Ann_{\g_{-}}(M) = 0$.   In view of Lemma \ref{Lem:1.84}, we have
proven

\bi \begin{Lem}\label{Lem:1.85} \ If $M$ is any graded
$\g_{0}$-submodule of $\g$ such that $[M,\g_1] = 0$, then
$\Ann_{\g_{-}}(M) = 0$.  In particular, $\Ann_{\g_{-}}(\g_r) = 0$.
If $\g_{r}$ is a trivial $[\g_{-1},\g_{1}]$-module, then
$\Ann_{\g_{-}}( {\mathcal I}_{r-1}) = 0$.  \end{Lem}

\m We now want to show

\m

\begin{Lem}\label{Lem:1.86}  \ Suppose $M_{j}$ is a $\g_0$-submodule of
$\g_{j}$ with $0<j<q$, such that $[[\g_{-1},\g_{1}], M_j]\ne 0$
and $[\g_{1},M_{j}] = 0$. Then $[M_{j},[M_{j},\g_{-j}]] \ne 0.$
\end{Lem}

\pf  First we observe that

 $$[M_{j},(\ad \, \g_{1})^{q-j-1}\,\g_{-q}]=
 (\ad \,\g_{1})^{q-j-1}\,\left ([M_{j},\g_{-q}] \right)=\g_{-1},$$
 thanks to Section \ref{sec:1.16} (iv) and Lemma \ref{Lem:1.85}. Therefore,

 \begin{eqnarray*}  0 \ne [[\g_{-1},\g_{1}],M_{j}] &=&
 [[[M_{j},(\ad \, \g_{1})^{q-j-1}\,\g_{-q}],\g_{1}], M_{j}] \\
&=& [M_{j}, [M_{j},(\ad \, \g_{1})^{q-j}\,\g_{-q}]]\\
&\subseteq& [M_{j},[M_{j},\g_{-j}]].  \ \qquad \ \ \square
\end{eqnarray*}

\m  \begin{Lem} \label{Lem:1.87}  \ We have $[{\mathcal 
I}_{s},[{\mathcal I}_s,\g_{-s}]] \ne 0$ \end{Lem}

\pf If $\g_{r}\cap {\mathcal I}\ne 0$ (and hence $s = r$),
then by Lemma \ref{Lem:1.84}, $\g_{r}$ is a nontrivial
$[\g_{-1},\g_{1}]$-module. Since $0 < r <q$, we then apply Lemma
\ref{Lem:1.86} with $M_j=\g_r$ to conclude that
$[\g_{r},[\g_{r},\g_{-r}]] \ne 0$. On the other hand, if
$\g_{r}\cap {\mathcal I}=0$, then by Lemma \ref{Lem:1.84},
$[\g_1, {\mathcal I}_{r-1}] = 0$ and ${\mathcal I}_{r-1}$ is a
nontrivial $[\g_{-1},\g_{1}]$-module. Furthermore $q>r-1=s$ in this case,  and $s>0$ 
by Lemma \ref{Lem:1.83}, 
so Lemma \ref{Lem:1.86} applies with $M_{j}=
{\mathcal I}_{r-1} = \mathcal I_s$.    Thus, we deduce that in any event $[{\mathcal 
I}_s,[{\mathcal I}_s,\g_{-s}]]\ne 0$. \ \qed

\bi

Recall we are assuming that $\g$ is generated by its local part.
Because $[\g_{1}, {\mathcal I}_{s}] \subseteq {\mathcal I}_{s+1} = 0$,
it follows that any $\g_{0}$-submodule $R_{s}$ of ${\mathcal 
I}_{s}$ generates an ideal

$$J \eqdef \,\sum_{i \ge 0}(\ad \, \g_{-1})^{i}R_{s}$$

\noindent which is contained in ${\mathcal I}$.  By the minimality
of ${\mathcal I}$, it must be that $J = {\mathcal I}$, so that
${\mathcal I}_{s} = R_{s}$;
  i.e., ${\mathcal I}_{s}$ is
$\g_{0}$-irreducible. Thus, in view of Lemma \ref{Lem:1.87}, we
can apply Lemma \ref{Lem:1.77} to conclude that the representation
of $\mathcal B({\mathcal I}_{s})_{0}^{(1)}$ on $\mathcal B({\mathcal I}_{s})_{-1}$ is a
restricted representation. Since we are also assuming that $\g$ is
transitive, the $\g_0$-module $\g_1$ is isomorphic to a
$\g_0$-submodule of
$\text{Hom}\,(\g_{-1},\,\g_0)\cong\g_{-1}^*\otimes_{\F}\,\g_0$.
{F}rom this it follows that for every $x\in  \g_0^{(1)}$, the
endomorphism $(\ad\,x)^p-\ad\,x^{[p]}$ acts on $\g_1$ as
$-\chi(x)^p\,\id$. But then $(\ad\,x)^p-\ad\,x^{[p]}=-s\chi(x)^p\,\id$
on $\g_s=(\g_1)^s$. As in the discussion of Section \ref{sec:1.14}
leading to the proof of Lemma \ref{Lem:1.79}, either the
representation of $\g_0^{(1)}$ on $\g_{-1}$ is restricted, or
$s$ is a multiple of $p$.  In summary, we have

   \bi \begin{Lem} \label{Lem:1.88} \ The
   $\g_0$-module ${\mathcal I}_{s}$ is irreducible.  Furthermore,
   either the representation of $\g_0^{(1)}$ on $\g_{-1}$ is
   restricted, or $p\mid s$.  \end{Lem} \bigskip

{\em  Thus, in what follows, we may assume $s$ is a multiple of
$p$ so that $s > 3$}.

\m To complete the proof of Theorem \ref{Thm:1.48}, we will show
that for any nonzero irreducible $\g_{0}$-submodule $V_{s-1}$ of
${\mathcal I}_{s-1}$, the condition $[V_{s-1},[V_{s-1},\g_{-s+1}]]
\ne 0$ must hold.

\bi \begin{Lem} \label{Lem:1.89} \ $[{\mathcal 
I}_{s-1},\g_{1}]= {\mathcal I}_{s}$, and $[V_{s-1}, \g_{-s+1}]\ne
0$ for any nonzero $\g_{0}$-submodule $V_{s-1}$ of $\mathcal I_{s-1}$.
\end{Lem}

\pf   Let $V_{s-1}$ be a nonzero $\g_{0}$-submodule of
$\mathcal I_{s-1}$, and suppose initially that $[V_{s-1}, \g_{1}] = 0$.
Then, since $\g$ is assumed to be generated by its local part, it
would follow that $\sum_{j \ge 0}(\ad \, \g_{-1})^{j}V_{s-1}$ is
an ideal of $\g$ properly contained in ${\mathcal I}$;  a
contradiction to the minimality of ${\mathcal I}$.

Consequently, $0 \neq [V_{s-1},\g_1] \subseteq {\mathcal I}_{s}$.
By Lemma \ref{Lem:1.88}, ${\mathcal I}_{s}$ is
$\g_{0}$-irreducible, so that $[V_{s-1},\g_{1}] = {\mathcal 
I}_{s}$.  In particular, $[{\mathcal I}_{s-1},\g_1] = {\mathcal 
I}_{s}$.  For the second half of the lemma, we may suppose that
$[V_{s-1},{\mathcal I}_{-s+1}] = 0$. Then

$$0 =  [{\mathcal I}_{s-2},[V_{s-1},{\mathcal I}_{-s+1}]]
= [V_{s-1},[{\mathcal I}_{s-2},{\mathcal I}_{-s+1}]], $$

\noindent since we are assuming that $s > 3$, so that $[{\mathcal 
I}_{s-2},V_{s-1}]=0$. If $[{\mathcal I}_{s-2},{\mathcal I}_{-s+1}]
\neq 0$, it is the whole of $\g_{-1}$ by irreducibility
\eqref{eq:1.2}. But then $[V_{s-1},\g_{-1}] = 0$ contradicting
transitivity \eqref{eq:1.3}. This shows that $[{\mathcal 
I}_{s-2},{\mathcal I}_{-s+1}] = 0$. Therefore, by our assumption
that $[V_{s-1},{\mathcal I}_{-s+1}] = 0$,  we have

\begin{eqnarray*} 0 = [{\mathcal I}_{s-2},{\mathcal I}_{-s+1}] \supseteq
[[\g_{-1},V_{s-1}],{\mathcal I}_{-s+1}] &=& [[\g_{-1},{\mathcal 
I}_{-s+1}],V_{s-1}]\\ &=& [{\mathcal I}_{-s},V_{s-1}].
\end{eqnarray*}

\noindent Since $\g_{-1}\subset {\mathcal I}$ and
$\g_{-s}=(\g_{-1})^s$, we have ${\mathcal I}_{-s}=\g_{-s}$ forcing
$[\g_{-s},V_{s-1}]=0$. We have already established that
$[V_{s-1},\g_{1}] = {\mathcal I}_{s}$. Consequently,

\begin{eqnarray*} 0 = [V_{s-1},{\mathcal I}_{-s+1}] \supseteq
[V_{s-1},[{\mathcal I}_{-s},\g_{1}]] &=&
[V_{s-1},[\g_{-s},\g_{1}]]
\\ & =& [\g_{-s},[V_{s-1},\g_{1}]] = [\g_{-s},{\mathcal I}_{s}],
\end{eqnarray*}

\n  to contradict Lemma \ref{Lem:1.87}.  Thus, $[V_{s-1},{\mathcal 
I}_{-s+1}] \ne 0,$ as claimed.  \ \qed

\bi  Next we show how 
\begin{equation}\label{eq:1.899} [V_{s-1},[V_{s-1},\g_{-s+1}]] = 0 \end{equation}

\n  for a nonzero $\g_0$-submodule of $V_{s-1}$ of
$\mathcal I_{s-1}$ leads to a contradiction.    If  \eqref{eq:1.899} holds, then

$$W := \left\{\ad_{\g_{s-2}}\,[v,u]\, \mid \,v \, \in \,
V_{s-1},\, u \, \in \, \g_{-s+1} \right\}$$ is a weakly closed
subset of $\En\,({\g}_{s-2})$ by Lemma \ref{Lem:1.57}\,(i).

\m First we assume for all $v \, \in \, V_{s-1},\, u \, \in \,
\g_{-s+1}$ that

\begin{equation}\label{eq:1.90}
(\ad \, v)^{3}(\ad \, u)^{3}\,\g_{s-2}=0. \end{equation}

\noindent Then $0 = (\ad \, v)^{3}(\ad \, u)^{3}\,\g_{s-2} =
6(\ad\,[u, v])^{3}\,\g_{s-2}$, hence  $W$ consists of nilpotent
transformations. The associative subalgebra of $\En\,(\g_{s-2})$
generated by $W$ must act nilpotently on $\g_{s-2}$ by Lemma
\ref{Lem:1.57} (ii), so that

$$Q_{s-2}\, \eqdef\, \Ann_{\g_{s-2}}\bigl([V_{s-1},\g_{-s+1}]\bigr) \ne 0.$$

\noindent If $[Q_{s-2},\g_{-s+1}] \neq 0$, then
$[Q_{s-2},\g_{-s+1}] = \g_{-1}$ by irreducibility \eqref{eq:1.2},
and since we are assuming $s \ge 5$ (and hence $s-2+s-1 > s+1\ge
r$), we would have
$$[V_{s-1},\g_{-1}] = [V_{s-1},[Q_{s-2},\g_{-s+1}]] =
[Q_{s-2},[V_{s-1},\g_{-s+1}]] = 0$$

\noindent by the definition of $Q_{s-2}$.   Then transitivity
\eqref{eq:1.3} would force $V_{s-1}=0$,  contrary to hypothesis.
Thus, we can conclude that $[Q_{s-2},\g_{-s+1}] = 0$. As
${\mathcal I}_{2s-4}=0$ (because $s\ge 5$), we have $0 =
[{\mathcal I}_{s-2},[Q_{s-2},\g_{-s+1}]] = [Q_{s-2}, [{\mathcal 
I}_{s-2},\g_{-s+1}]]$. Now $[{\mathcal I}_{s-2},\g_{-s+1}]$ equals
0 or $\g_{-1}$, and when it is nonzero we have

$$0 = [Q_{s-2},[{\mathcal I}_{s-2},\g_{-s+1}]] =
[Q_{s-2},\g_{-1}]$$

\noindent contradicting the transitivity of $\g$. Thus, $[{\mathcal 
I}_{s-2},\g_{-s+1}] = 0$, and as a result, $$0
=[V_{s-1},[{\mathcal I}_{s-2},\g_{-s+1}]]= [{\mathcal 
I}_{s-2},[V_{s-1},\g_{-s+1}]].$$ Hence, $[V_{s-1},\g_{-s+1}]
\subseteq \Ann_{\g_{0}}({\mathcal I}_{s-2}) \cap
\Ann_{\g_{0}}(V_{s-1})$, which must be zero by Lemma
\ref{Lem:1.59}. But this contradicts Lemma \ref{Lem:1.89}.

Thus, we have proven that  if \eqref{eq:1.899} holds, then  \eqref{eq:1.90}  
cannot happen.   In other words, there exist $v\in V_{s-1}$ and $u\in
\g_{-s+1}$ such that $(\ad \, v)^{3}(\ad \, u)^{3}\,\g_{s-2}\neq
0.$ Then $U_{s-2} \eqdef (\ad \, V_{s-1})^{3}\,\g_{-2s+1}$ is a
nonzero $\g_{0}$-submodule of ${\mathcal I}_{s-2}$.  Assuming
$[\g_{1},U_{s-2}]=[\g_{2},U_{s-2}]=0$ would imply that

$$\sum_{j \ge 0}\,(\ad \,  \g_{-1})^{j}\, U_{s-2}$$

\noindent is an ideal of $\g$ properly contained in ${\mathcal 
I}$, a contradiction.  When  $[\g_{1},U_{s-2}] \ne 0$, then

$$0 \ne [\g_{1},U_{s-2}] = [\g_{1},(\ad \,
V_{s-1})^{3}\,\g_{-2s+1}] \subseteq (\ad \,
V_{s-1})^{2}\,\g_{-s+1}.$$

\noindent On the other hand, when $[\g_{2},U_{s-2}] \ne 0$, then

$$0 \ne [\g_{2},U_{s-2}] \subseteq (\ad \,
V_{s-1})^{3}\,\g_{-2s+3}.$$ In this case, transitivity
\eqref{eq:1.3} gives

$$0 \ne [\g_{-1}, (\ad \, V_{s-1})^{3}\,\g_{-2s+3}] \subseteq (\ad
\, V_{s-1})^{2}\,\g_{-s+1}.$$

\noindent  Consequently,  $(\ad \, V_{s-1})^{2}\,\g_{s+1} \neq 0$  must hold, and
that  in conjunction with \eqref{eq:1.53},
allows us to conclude the following:

\bi \begin{Lem} \label{Lem:1.91} \
$[V_{s-1},[V_{s-1},\g_{-s+1}]] \ne 0$ for any nonzero
$\g_0$-submodule $V_{s-1}$ of ${\mathcal I}_{s-1}$.  When
$V_{s-1}$ is irreducible, then for ${\mathcal B}' = {\mathcal 
B}(V_{s-1}) = \bigoplus {\mathcal B}'_i$ we have ${\mathcal 
B}'_{1} \ne 0$, and the representation of $({\mathcal 
B}'_0)^{(1)}$ on ${\mathcal B}'_{-1} = V_{s-1}$ is
restricted.
\end{Lem}

\bi As before, knowing that ${\mathcal B}_1 \neq 0$ for ${\mathcal 
B} = {\mathcal B}(\mathcal I_s)$ and ${\mathcal B}_1' \neq 0$ for
${\mathcal B}' = {\mathcal B}(V_{s-1})$ forces the representation
of $\g_0^{(1)}$ on $\g_{-1}$ to be restricted.  Thus, we may draw
the conclusion:

\bi\begin{Lem} \label{Lem:1.92} \ Let $\g = \bigoplus_{i=-q}^r
\g_i$ be a graded Lie algebra satisfying the hypotheses of Theorem
\ref{Thm:1.48}, and assume that $q > r$ and
$[[\g_{-1},\g_1],\g_1]\ne 0$. Then the representation of
$\g_0^{(1)}$ on $\g_{-1}$ is a restricted representation.
\end{Lem}
\m

\noindent \begin{pgraph} \ {\rm Combining Lemmas
\ref{Lem:1.77}, \ref{Lem:1.79}, \ref{Lem:1.80}, and \ref{Lem:1.92}
completes the proof of Theorem \ref{Thm:1.48}.} \ \qed
\end{pgraph} \m

\noindent
\section {The assumption $[[\g_{-1},\g_1],\g_1]\ne 0$ in Theorem 1.63} \label{sec:1.17}
{\rm The assumption that $[[\g_{-1},\g_1],\g_1]\ne 0$ in the
statement of Theorem \ref{Thm:1.48} cannot be dropped, as the
following example illustrates.

Let $\mathfrak{K}$ be a classical simple Lie algebra, and let $\mathfrak M$
be an irreducible nonrestricted $\mathfrak{K}$-module with basis
$x_1,\ldots,x_m$. Clearly, $\mathfrak K$ acts as homogeneous
derivations of degree zero on the graded truncated polynomial
algebra
$$\mathfrak A:=\F[x_1,\ldots,x_m]/\langle x_1^p,\ldots,x_m^p\rangle, \qquad
\deg\,x_i=1.$$ Note that $\mathfrak A=\bigoplus_{i=0}^r\, \mathfrak A_i$, where
$r=m(p-1)$ and $\mathfrak A_0=\F\,1$. Let $\mathfrak D$ denote the linear
span of the partial derivatives $\partial/\partial
x_1,\ldots,\partial/\partial x_m$ of $\mathfrak A$.

We regard $\mathfrak K$ as a subalgebra of the derivation algebra
$\text{Der}(\mathfrak A)$. Clearly, $[{\mathfrak K},{\mathfrak D}]\subseteq
\mathfrak D$ and $\mathfrak D^{(1)}=0$. It follows that
${\mathcal K}\eqdef {\mathfrak K}\oplus\mathfrak D$ is a Lie
subalgebra of $\text{Der}(\mathfrak A)$, and $\mathfrak D$ is an abelian
ideal of $\mathcal K$. The pairing ${\mathfrak D}\times
\mathfrak M\rightarrow \F$, $\,\, (d,x)\mapsto d(x),$ is nondegenerate and
$\mathfrak K$-invariant. Hence, ${\mathfrak  D}\cong \mathfrak M^*$, showing
that the $\mathfrak K$-module $\mathfrak D$ is irreducible and
nonrestricted.

Now let $\mathfrak L$ be another classical simple Lie algebra.
Then ${\mathfrak L}\ot  \mathfrak A$ carries a natural Lie algebra structure
given by $[l\ot f,l'\ot f']=[l,l']\ot ff'$ for all
$l,l'\in\mathfrak L$ and $f,f'\in \mathfrak A$. Note that
$\text{id}_{\mathfrak L}\ot \mathcal K$ is a subalgebra of the
derivation algebra of ${\mathfrak L} \ot \mathfrak A$. As a consequence, we
may regard
$$\g\eqdef({\mathfrak L}\ot \mathfrak A)\,\oplus\, (\text{id}_{\mathfrak L}\ot {\mathcal K})$$
as a semidirect product of Lie algebras. Put
\begin{eqnarray}\label{Gr}
&&\\
&&\g_{-i}:= {\mathfrak L}\ot \mathfrak A_{i},\qquad 1\le i\le r, \nonumber\\
&&\g_0\ := \left( {\mathfrak L}\ot \mathfrak A_0 \right)\,\oplus \left (\text{id}_{\mathfrak
L}\ot \mathfrak{K} \right),\nonumber \\
&&\g_1\ :=\text{id}_{\mathfrak L}\ot{\mathfrak D}\nonumber.
\end{eqnarray} Then $\g\,=\,\bigoplus_{i=-r}^1\,\g_i$ and
$[\g_i,\g_j]\subseteq \g_{i+j}$ for all $i,j$,  so that 
(\ref{Gr}) defines a grading of $\g$. It is
straightforward to see that $\g_0\cong {\mathfrak
L}\oplus\mathfrak{K}$ is classical reductive,
$\g_{-1}=\mathfrak{L}\ot \mathfrak M$ is an irreducible nonrestricted
$\g_0$-module, and $\g_1\cong \mathfrak M^*$ as $\mathfrak K$-modules.
Moreover, $\g$ is transitive, and
$$[[\g_{-1},\g_1],\g_1]\subseteq [\mathfrak{L}\ot \mathfrak A_0,\,
\text{id}_{\mathfrak L}\ot \mathfrak D]=0.$$}
\m

Transitive gradings similar to the one above were first described
by Weisfeiler who termed them {\it degenerate} (see \cite{W}).
They arise quite frequently in the classification theory of simple
Lie algebras.

%-----------------------------------------------------------------------
% Beginning of chap2.tex   9-16-05 version
%-----------------------------------------------------------------------
%
% AMS-LaTeX 1.2 sample file for a monograph, based on amsbook.cls.
% This is a data file input by chapter.tex.
%%%%%%%%%%%%%%%%%%%%%%%%%%%%%%%%%%%%%%%%%%%%%%%%%%%%%%%%%%%%%%%%%%%%%%%%%%%%%%%%%%%%%%%%%%%%%%%%%%%%%%%%%%%%%%%%%%%%%%%
%%%%%%%%%%%%%%%%%%%%%%%%%%%%%%%%%%%%%%%%%%%%%%%%%%%%%%%%%%%%%%%%%%%%%%%%

  \chapter{Simple Lie Algebras and Algebraic Groups}
\bi
  
\section {\ Introduction  \label{sec:2.1}}

\m  In finite characteristics, there are simple finite-dimensional
graded Lie algebras which resemble those in zero characteristic
and others which do not.  The former are referred to as
``classical'', as they are obtained from their characteristic zero
counterparts by reducing a Chevalley basis modulo $p$ and
factoring out the center as necessary.  The latter are the Cartan
type Lie algebras, and in characteristic 5, the Melikyan Lie
algebras.  The Cartan type and Melikyan Lie algebras possess a
grading in which the zero component $\g_0$ is isomorphic to one of
the linear Lie algebras $\mathfrak{sl}_m$ or $\mathfrak{gl}_m$ or
to one of the symplectic Lie algebras $\spm$ or $\csp$. (The
algebra $\csp$ is a central extension of $\spm$ by a
one-dimensional center.) The gradation spaces $\g_j$ are modules
for $\g_0$; indeed, in each case, $\g_{-1}$ is the natural module
(or its dual module) for $\mathfrak{sl}_m$, $\mathfrak {gl}_m$,
$\spm$, or $\csp$.  In this chapter, we gather information about
the classical, Cartan type, and Melikyan Lie algebras. Our purpose
is two-fold.  On one hand, the focus of this monograph is a
classification theorem, and these are the Lie algebras which will
appear in the theorem's conclusion.  On the other, a lot of our
discussion will consist of calculations involving these Lie
algebras, and we will include here enough facts about them to
enable the reader to follow these computations without too much
difficulty, we hope.  \m

\section {\ General  information
about the classical Lie algebras  \label{sec:2.2}}   

\m The well-established structure theory of finite-dimensional
simple Lie algebras $\dot \g$ over the field $\mathbb C$ of
complex numbers begins with a Cartan subalgebra $\dot \h$
(equivalently, a maximal toral subalgebra) and the root space
decomposition $\dot \g = \dot \h \oplus \bigoplus_{\alpha \in
\Phi} \dot \g^{\alpha}$ relative to $\dot \h$.  The root spaces $\dot
\g^{\alpha} = \{ x \in \dot \g \mid [h,x] = \alpha(h)x$ for all $h
\in \dot \h\}$ are all one-dimensional.  Because the Killing form
$\kappa(x,y) = \hbox{tr}(\ad x \, \ad y)$ is nondegenerate and
invariant, $\dot \g^{\alpha}$ and $\dot \g^{-\alpha}$ are
nondegenerately paired, and the form restricted to $\dot \h$ is
nondegenerate.  Thus, for any $\mu \in \dot \h^*$, there exists a
unique element $t_{\mu} \in \dot \h$ such that $\mu(h) =
\kappa(h,t_\mu)$ for all $h \in \dot \h$.  This allows the form to
be transferred to the dual space $\dot \h^{*}$, where $(\mu,\nu)
:= \kappa(t_{\mu},t_{\nu})$.  In particular, this form restricted
to the real span of the roots is positive definite.  For
$\alpha,\beta \in \Phi$, we set
 
\begin{equation}\label{eq:halpha}  \langle \beta,\alpha \rangle : =
{\frac {2(\beta,\alpha)} {(\alpha,\alpha)}} \qquad \hbox{and} \qquad
h_{\alpha} = {\frac {2
t_\alpha} {(\alpha,\alpha)}}\end{equation}

There is a basis of $\dot \h^*$ of {\it simple roots} $\Delta =
\{\alpha_1,\dots,\alpha_m\} \subseteq \Phi$, and every root is a
nonnegative or nonpositive integral linear combination of them.
The positive roots ($\alpha > 0$ for short) are the former, and
the negative roots ($\alpha < 0$) are the latter.  For brevity we
will write $h_i = h_{\alpha_i}$ for $i=1,\dots,m$.

\m There is a second distinguished basis for $\dot \h^*$, namely
the basis $\{\varpi_1,$ $\dots,$ $\varpi_m\}$ of {\it fundamental
weights}.  They are dual to the elements $h_i$, and so satisfy 

\begin{equation}\label{eq:2.3} \varpi_i(h_j) = \langle \varpi_i,
\alpha_j\rangle = \delta_{i,j}.
\end{equation}
The Cartan matrix $A = \Big( \langle \alpha_i,\alpha_j
\rangle\Big)$ gives the transition between these two bases:
$\alpha_i = \sum_{j=1}^m A_{i,j} \varpi_j$, \quad $A_{i,j} =
\langle \alpha_i,\alpha_j \rangle$ for $1\leq i,j \leq m$. It also
encodes all the structural information about $\dot \g$. \bi 

\begin{Thm} \label{Thm:2.4} \ Let $\dot \g$ be a
finite-dimensional complex simple Lie algebra.  Then there is a
basis $\{e_{\alpha} \mid \alpha \in \Phi\} \cup \{h_i \mid
i=1,\dots,m\}$ of $\dot \g$ such that
\begin{enumerate}
\item [{\rm(i)}] $[h_i,h_j] = 0, \qquad 1 \leq i,j \leq m$; \item
[{\rm(ii)}] $[h_i,e_{\alpha}] = \langle \alpha, \alpha_i\rangle
e_{\alpha} = \alpha(h_i)e_\alpha$ \qquad $1 \leq i \leq m, \quad \alpha \in \Phi$; \item
[{\rm(iii)}] $[e_{\alpha},e_{-\alpha}] = h_{\alpha} = \sum_{j=1}^m
\langle \varpi_j,\alpha \rangle h_j$, which is an integral linear
combination of $h_1,\dots,h_m$; \item [{\rm(iv)}] If $\beta \neq
\pm \alpha$, and $\beta-d\alpha, \dots, \beta+u \alpha$ is the
$\alpha$-string through $\beta$, then 

\begin{equation}[e_{\alpha},e_{\beta}] = \begin{cases} 0 &\qquad \hbox{\rm if} \ u
= 0 \\
\pm (d+u) e_{\alpha+\beta} &\qquad \hbox{\rm if} \ u \neq 0,
\end{cases}\nonumber\end{equation}
where the integers $d,u$ satisfy $d-u = \langle \beta,\alpha
\rangle$. 

\end{enumerate} \end{Thm}
 
\bi
 
Such a basis is referred to as a {\it Chevalley basis}.
 The $\Z$-span $\dot \g_{\Z}$ of a Chevalley basis
is a $\Z$-subalgebra of $\dot \g$.  If $\F$ is any field, then $\g = \g_\F := \F \otimes_{\Z}
\dot \g_{\Z}$ is a Lie algebra over $\F$.   It has a basis as in Theorem
\ref{Thm:2.4}, but with the structure constants reduced modulo
$p$.  If  $\F$ has characteristic $p > 3$, the Lie algebra $\g$ is simple  in all instances except when
$\Phi$ is a root system of type A$_{m}$ for $p \mid (m+1)$.  In that
particular case, $\g$ is isomorphic to $\mathfrak{sl}_{m+1}$, and $\g$
has a one-dimensional center spanned by the identity matrix $I =
h_1 + 2h_2 + \cdots + mh_{m}$   (where  $h_i $ is
identified with the difference of matrix units,  $h_i = E_{i,i}-E_{i+1,i+1}$).  Then $\g/ \F I \cong
\mathfrak{psl}_{m+1}$ is simple.  The simple Lie algebras
obtained by the process of reducing a Chevalley basis modulo $p$
(and factoring out the center if necessary) are called the {\em classical
simple Lie algebras}.      The classical simple Lie algebras  along with $\mathfrak{sl}_{m+1}$,
$\mathfrak{gl}_{m+1}$, or $\mathfrak{pgl}_{m+1}$ where $p \mid (m+1)$
and direct  sums of these algebras are known
under the general rubric of {\it classical Lie algebras}. 
A {\em classical reductive Lie algebra} 
modulo its center is  a classical Lie algebra.     Assumption (a)
of the Recognition Theorem asserts that $\g_0$ is classical reductive. 
 
\m Each classical simple Lie algebra is equipped with a
$[p]$-operation relative to which the Lie algebra is restricted:
 
\begin{equation} \label{eq:2.5} e_\alpha^{[p]} = 0 \quad \hbox{for}
\quad \alpha \in \Phi, \qquad \hbox{and} \qquad h_i^{[p]} =
h_i\quad \hbox{for} \quad i=1,\dots,m.
\end{equation} 

Mills and Seligman \cite{MS}  found a simple set of axioms which
characterize direct sums of these Lie algebras.
 
\bi
 
\begin{Thm} \label{Thm:2.6} A Lie algebra $L$ over a field  $\F$
of characteristic $p > 3$ is a direct sum of classical simple Lie
algebras if and only if  
 
\begin{enumerate}
\item [{\rm(i)}] $[L,L] = L$; \item [{\rm(ii)}] The center of $L$
is 0; \item [{\rm(iii)}] $L$ has an abelian Cartan subalgebra $H$
such that
\begin{enumerate}
\item [{\rm(a)}] $L = \bigoplus_{\alpha \in H^*} L^{\alpha}$,
where
 
$L^{\alpha} = \{x \in L \mid [h,x] = \alpha(h)x$ for all $h \in
H$\}; \item [{\rm(b)}] If $L^{\alpha} \neq 0$ and $\alpha \neq 0$,
(i.e., if $\alpha$ is a root), then $[L^{\alpha},L^{-\alpha}]$ is
one-dimensional; \item [{\rm(c)}] If $\alpha$ and $\beta$ are
roots, then not all $\beta + k \alpha$ for $k=0,1,\dots,p-1$, are
roots.
\end{enumerate}\end{enumerate}\end{Thm}
\m

\begin{Rem}\label{CSA-toral}  {\rm According to \cite[Chap.~IV, Thm.~1.2  
and  Chap.~VI, Sec.~2]{S}, 
the Cartan subalgebras of a classical simple Lie algebra $\g$ over
an algebraically closed field of characteristic $p > 3$  are all conjugate
under the automorphism group of $\g$.     Combining that result with 
Theorem \ref{Thm:2.6}, we see that any Cartan subalgebra $\h$ of $\g$
is abelian  and  that the elements of $\h$ are  ad-semisimple on $\g$.  
Hence $\h$ is a toral subalgebra of $\g$ (in fact, a maximal toral subalgebra of $\g$).
As the centralizer of a maximal toral subalgebra of a restricted Lie algebra
is a  Cartan subalgebra, and every Cartan subalgebra is such a centralizer
(for example, \cite[Thm.~4.1]{SF}),   it follows that the maximal toral subalgebras
of a classical simple Lie algebra are precisely the Cartan subalgebras.
We will use the two notions interchangeably when working with a classical simple
Lie algebra or with one of the algebras     
$\mathfrak{sl}_{m+1}$,
$\mathfrak{gl}_{m+1}$, or $\mathfrak{pgl}_{m+1}$ with $p \mid (m+1)$,
where the same result applies.  
All maximal toral subalgebras of
a classical simple Lie algebra are conjugate -- a result which follows
from the above considerations for $p > 3$, but  is known
to hold for all characteristics by \cite[Cor.~13.5]{H1}. }    \end{Rem}  

The derivations of a classical simple Lie algebra are inner
except  when $\g =\mathfrak{psl}_{m+1}$.   More specifically,
we have

\bi
\begin{Lem}\label{der}
Suppose that $\g$ is either a classical simple Lie algebra or
$\mathfrak{pgl}_{m+1}$ with $p\mid (m+1)$. Then $\Der(\g) \cong\ad
\g$ unless $\g=\mathfrak{psl}_{m+1}$ with $p\mid (m+1)$ in which
case $\Der(\g) \cong \,\mathfrak{pgl}_{m+1}$.
\end{Lem}

\pf First suppose that $\g$ is simple. We have $\g={\mathfrak
h}\oplus\bigoplus_{\alpha\in\Phi}\g^\alpha$ where $\mathfrak h$ is
a Cartan subalgebra of $\g$. Since $\ad{\mathfrak h}$
acts diagonalizably on $\g$, it also acts diagonalizably on
$\text{End}(\g)\cong\g\otimes\g^*$. It follows that $\ad\mathfrak
h$ acts diagonalizably on $\Der(\g)\subset \text{End}(\g)$. Let
$\widetilde{\mathfrak h}$ denote the centralizer of $\ad \mathfrak
h$ in $\Der(\g)$. If $D\in\Der(\g)$ is a weight vector for
$\mathfrak h$ corresponding to a nonzero weight $\mu\in{\mathfrak
h}^*$, then $\mu(h)\ne 0$ for some $h\in\mathfrak h$, forcing
$$D=\mu(h)^{-1}[\ad h, D]=-\mu(h)^{-1}\ad D(h)\in\ad\g.$$ As a consequence,
$\Der(\g)=\widetilde{\mathfrak h}+\ad\g$. If
$D\in\widetilde{\mathfrak h}$, then
$D(\g^\alpha)\subseteq\g^\alpha$ for all $\alpha\in\Phi$, as all
root spaces of $\g$ relative to $\mathfrak h$ are one-dimensional.
For each $\alpha_i\in\Delta,$ choose nonzero vectors
$e_i\in\g^{\alpha_i}$ and $f_i\in\g^{-\alpha_i}$. Since $D$
commutes with $\ad [e_i,f_i]\in\ad\mathfrak h\setminus\{0\}$, it
must be  that $D(e_i)=\lambda_ie_i$ and $D(f_i)=-\lambda_if_i$ for
some $\lambda_i\in\F,$ where $1\le i\le m$. Since in the present
case the Lie algebra $\g$ is generated by the $e_i$'s and $f_i$'s,
we may conclude that $\dim\widetilde{\mathfrak  h}\le m$.

Let $\bar{A}$ denote the matrix with entries in $\F_p$ obtained by
reducing the entries of the Cartan matrix $A$ of $\g$ modulo $p$. If
$\g\ne\mathfrak{psl}_{m+1}$ with $p\mid (m+1)$, then
$\det\bar{A}\ne 0$. Using (\ref{eq:2.3}) it is easy to observe
that in this case there exist $t_1,\ldots,t_m\in\mathfrak h$
satisfying $[t_i,e_j]=\delta_{i,j}e_j$ for all $1\le i,j\le m$.
But then $D=\sum_{i=1}^m\lambda_i\,\ad t_i\in\ad\g$, yielding
$\Der(\g)=\ad\,\g$. If $\g=\mathfrak{psl}_{m+1}$ with $p\mid
(m+1),$ then it is readily seen that there are
$t_1,\ldots,t_m\in\mathfrak{pgl}_{m+1}$ such that
$[t_i,e_j]=\delta_{i,j}e_j$ for all $1\le i,j\le m$. The above
discussion then implies that
$\Der(\mathfrak{psl}_{m+1})\cong\ad\mathfrak{pgl}_{m+1} \cong
\mathfrak{pgl}_{m+1}$ if $p\mid
(m+1)$.

Now suppose that $\g=\mathfrak{pgl}_{m+1}$ with $p\mid (m+1)$, and
let $D$ be any derivation of $\g$. Clearly, $D$ induces a
derivation of $\g^{(1)}=\mathfrak{psl}_{m+1}$. By the above, there
exists an element $x\in\g$ such that $D-\ad x$ acts trivially on
$\g^{(1)}$. We now let $y$ be any element of $\g$. Then $$0=(D-\ad
x)\big([y,\g^{(1)}]\big)=\big[(D-\ad x)(y),\g^{(1)}\big],$$
implying that $(D-\ad x)(y)$ belongs to the centralizer of
$\g^{(1)}$ in $\g$. Since $p>2$, it is straightforward to see that
the latter is trivial. Therefore, $D=\ad x$, completing the
proof.\qed

 \bi

 In Section \ref{sec:1.12},  we have used certain exponential
automorphisms. Next we establish the properties required there.

\bi

\begin{Lem} \label{Lem:2.7} \ Suppose $\g$ is
a classical simple Lie algebra, or a Lie algebra isomorphic to
$\mathfrak{sl}_{m+1}$, $\mathfrak{gl}_{m+1}$, or $\mathfrak{pgl}_{m+1}$ where
$p \mid (m+1)$.    If $x\in \g^{\al}$ is a root vector, then $\exp(\ad
x)$ is an automorphism of $\g$.
\end{Lem}

\pf  If $x\in \g^{\al}$ is a root vector, then $(\ad x)^3 = 0$
on $\g$, except when $\g$ is of type G$_2$ and $\al$ is a short root.
In this exceptional case, $\ad x$ is nilpotent
of order 4.  Thus, except for root vectors corresponding to short
roots in type G$_2$ in characteristic 5, we have $(\ad
x)^{\frac {p+1} {2}} = 0$ for all root vectors.

Now if $D$ is any derivation of $\g$ with $D^{\frac{p+1}{2}} = 0$,
then using the well-known Leibniz formula,
\begin{equation}{\frac{D^n([y,z])} {n!}} = \sum_{i = 0}^n\left[ {\frac {D^i(y)} {i! }},\,
{\frac{D^{n-i}(z)} {(n-i)! }}\right],\nonumber \end{equation}
 \n we have
\begin{eqnarray*} [\exp(D)(y),\exp(D)(z)] &=& \left [ \sum_{i =
0}^{\frac{p-1}{2}}\frac{D^i(y)}{i! },\, \sum_{j =
0}^{\frac{p-1}{2}}\frac{D^j(z)}{j! }\right ]\nonumber\\ &=& \sum_{n
= 0}^{p-1} \sum_{i = 0}^n\left [\frac{D^i(y)}{i! },
\frac{D^{n-i}(z)}
{(n-i)! }\right] = \sum_{n = 0}^{p-1}\frac{D^n([y,z])}{n! } \nonumber\\
&=& \exp(D)([y,z]), \nonumber\end{eqnarray*}
\n so that $\exp(D)$ is Lie homomorphism of $\g$.  If $D'$ is
another derivation of $\g$ with $(D')^{\frac{p+1}{2}} = 0$ and
$[D,D'] = 0$, then
\begin{eqnarray*}
\exp(D)\exp(D') &=&
\left(\sum_{i=0}^{\frac{p-1}{2}}\frac{D^{i}}{i! }\right)
\left(\sum_{j=0}^{\frac{p-1}{2}}\frac{(D')^{j}}{j! }\right) \nonumber\\
&=& \sum_{k=0}^{p-1} \frac{(D + D')^{k}}{k! } = \exp(D + D').
\nonumber\end{eqnarray*}
\n Taking $D' = -D$ shows that $\exp(D)\exp(-D) = I$, so that
$\exp(D)$ is an automorphism of $\g$. As a result, $\exp(\ad x)$
is an automorphism of $\g$ for all root vectors, except when $x$
corresponds to a short root in type G$_2$ and characteristic 5.

In the exceptional case, we suppose that $\g$ is a Lie algebra
of type G$_2$ in characteristic 5 and adopt the labelling of the
roots as in {\cite[Planche IX]{Bou1}}  so that $\al_1,\al_2$ are the
simple roots, and the positive roots are
$\{\al_1,\al_2,\al_1+\al_2,2\al_1+\al_2,3\al_1+\al_2,
3\al_1+2\al_2\}$.  Then for $x \in \g^{\al}$ ($\al$ a root) and $y
\in \g^\be$ ($\be$ a root or 0), we have $(\ad x)^{3}y= 0$ unless
$\al$ is short.  By applying a root system isomorphism, we may
assume $\al = \al_1$, $x \in \g^{\al_1}$, and $y,z$ are root
vectors in $\g$. Then
\begin{eqnarray*} [\exp(\ad x)(y), \exp(\ad x)(z)] &=& \exp(\ad
x)\big([y,z]\big) \nonumber\\
&& \quad + \left [\frac{(\ad x)^2 (y)}{2! }, \frac{(\ad x)^3
(z)}{3! }
\right] \nonumber\\
&& \quad \quad +\left [\frac{(\ad x)^3 (y)}{3! }, \frac{(\ad x)^2
(z)}{2! }\right] \nonumber\\
&& \quad \quad \quad  +\left [ \frac{(\ad x)^3 (y)}{3! },
\frac{(\ad x)^3 (z)}{3! }\right]. \nonumber\end{eqnarray*}
\n Now $(\ad x)^3(\g)$ is spanned by the root vectors
corresponding to the roots $3\al_1 + \al_2$ and $-\al_2$, and
$(\ad x)^2(\g)$ is spanned modulo $(\ad x)^3(\g)$ by root vectors
corresponding to $2 \al_1+\al_2, \al_1$, and $-\al_1-\al_2$.  From
this it is easy to see that $(\ad x)^3(\g)$ is central in $(\ad
x)^2(\g)$, so that $[\exp(\ad x)(y), \exp(\ad x)(z)]$ = $\exp(\ad
x)\big([y,z]\big)$ and $\exp(\ad x)$ is a homomorphism.  Because
\begin{eqnarray*} \exp(\ad x)\,\exp(-\ad x) &=& \left (
\sum_{i=0}^{p-1}\frac{(\ad x)^{i}}{i! } \right) \left (
\sum_{j=0}^{p-1}\frac{(-\ad x)^{j}}{j! }\right) \nonumber\\
&=& \sum_{n=0}^{p-1}\sum_{k=0}^{n}\frac{(\ad
x)^{n-k}}{(n-k)! }\frac{(-\ad x)^{k}}{k! } \nonumber\\
&& \qquad \quad + \sum_{i+j = p}^{2p-2} (-1)^j \frac{(\ad
x)^{i+j}}{i! \, j! } \nonumber\\
&=& \sum_{n=0}^{p-1}\frac{(\ad x -\ad x)^{n}}{n! } = I,
\nonumber\end{eqnarray*} \n the mapping $\exp(\ad x)$ is an
automorphism as claimed.  \ \qed

\m  Each classical Lie algebra  possesses
a finite reduced root system $\Phi$.
Corresponding to such a root system $\Phi$ is the associated
Weyl group $W$, which is generated by the reflections $s_{\alpha}
\in \En(\h^*)$, $\alpha \in \Phi$,

$$s_{\alpha}(\mu) = \mu - \langle \mu,\alpha \rangle \alpha.\nonumber$$

In fact, it can be shown that $W$ is generated by the reflections
$s_{\alpha_i}$, $1\leq i \leq m$, in hyperplanes perpendicular to
the simple roots.

\m Suppose $P = \Z \varpi_1 \oplus  \dots \oplus \Z \varpi_m$ is
the weight lattice corresponding to $\Phi$, and let $P^+ = \mathbb N
\varpi_1 \oplus \dots \oplus \mathbb N  \varpi_m$ ($\mathbb N =
\Z_{\geq 0}$)   be the set of dominant (integral) weights.  A
subset $X$ of $P$ is said to be {\it saturated} if for all
$\lambda \in X$ and $\alpha \in \Phi$, the weight $\lambda -
t\alpha \in X$ for all integers $t$ between 0 and $\langle
\lambda,\alpha \rangle,$ inclusive.  (See \cite[Chap.~VIII, Sec.~7,
No.~2]{Bou2},  or \cite[Sec.~13.4]{H3}.)

\m  The {\it minuscule weights} $\lambda \in P^+$ will play a role
in our subsequent deliberations.  They have several equivalent
characterizations:
\begin{enumerate}

\item [{\rm(i)}] The Weyl group orbit $W \lambda$ of $\lambda$ is
the smallest saturated set of weights in $P$ containing $\lambda$.

\item [{\rm(ii)}] For all roots $\alpha\in \Phi$, \ $\langle
\lambda,\alpha \rangle \in \{-1, 0, 1\}.$

\item [{\rm(iii)}] If $\mu \in P^+$ and $\mu < \lambda$,
(i.e., $\lambda-\mu = \sum_{i=1}^m k_i \alpha_i$ where $k_i \in
\mathbb N$ for all $i$), then $\lambda = \mu$.

\end{enumerate}

\noindent Relatively few weights are minuscule, and they are
displayed below:

\begin{eqnarray}\label{eq:2.8}
&& \\
\hbox{A$_m$} && \qquad \qquad \varpi_1, \dots, \varpi_m  \nonumber\\
\hbox{B$_m$} &&
\qquad \qquad \varpi_m \nonumber \\
 \hbox{C$_m$} && \qquad \qquad
\varpi_1 \nonumber \\
 \hbox{D$_m$} && \qquad \qquad  \varpi_1, \varpi_{m-1}, \varpi_m
\nonumber \\
  \hbox{E$_6$}&& \qquad \qquad \varpi_1, \varpi_6 \nonumber \\
 \hbox{E$_7$} &&
  \qquad \qquad \varpi_7. \nonumber
  \end{eqnarray}
  \m

The Borel subalgebras 

\begin{eqnarray*} \mathfrak b^+ &:= &  \mathfrak h \oplus \bigoplus_{\alpha > 0}
\g^\alpha \\
\mathfrak b^- &:=&   \mathfrak h \oplus \bigoplus_{\alpha > 0}
 \g^{-\alpha} \end{eqnarray*}  are essential in the structure and  representation theory of
$\g$.     Unlike the situation in characteristic zero where Lie's Theorem
applies, a finite-dimensional $\g$-module $M$ need not
possess a common eigenvector  for
$\mathfrak b^+$ (or for $\mathfrak b^-$).     However,  when such a common
eigenvector  $v$ exists for $\mathfrak b^+$ (resp. for $\mathfrak b^-$),   
it must be annihilated by $\mathfrak n^+ := \bigoplus_{\alpha > 0}
\g^\alpha$   (resp. $\mathfrak n^-: = \bigoplus_{\alpha > 0}
\g^{-\alpha}$)   and be a weight vector for $\mathfrak h$.   This 
motivates the following definition.  \bi

\begin{Def} \label{def:prim} Let $M$ be a finite-dimensional module for the classical algebra $\g$.
A nonzero vector $v$ in $M$ is said to be a $\mathfrak b^+$-primitive (resp. $\mathfrak b^-$-primitive)  vector of weight $\lambda$
if $x.v = 0$ for all $x \in \mathfrak n^+$ (resp. $x \in \mathfrak n^-$) and $h.v = \lambda(h)v$ for all $h \in \mathfrak h$.  \end{Def}  \m

As explained in Remark \ref{irr}, any finite-dimensional irreducible $\g$-module which
is restricted over $\g^{(1)}$  possesses  
a $\mathfrak b^+$-primitive vector  (and by symmetry, a $\mathfrak b^-$-primitive vector).  
\bi

  \section {\ Representations of algebraic groups   \label{sec:2.3} }

 \m The automorphism group of a classical simple Lie
 algebra is an algebraic group of adjoint type.
 Various proofs in the text will require us to work with
 the simply connected cover of such a group.
 In this section, we assemble the information
 we will need  about  algebraic groups.  \m

Let $G$ be a linear algebraic group, that is,  a Zariski closed
subgroup of ${\rm GL}(V)$ where $V$ is a finite-dimensional vector
space over $\F$. The  group $G$ is said to be {\em connected}  if
$G$ equals the  {\it connected component} $G^\circ$ of $G$ (the
unique irreducible component of $G$ containing the identity
element). Let $A = \F[G]$ denote the coordinate ring of $G$. Then
$G$ acts on  $A$  via left (resp. right) translation:
$(\lambda_xf)(y) = f(x^{-1}y)$ (resp. $(\rho_x f)(y) = f(yx))$.
The space $\hbox{\rm Lie}(G):= \{d \in \Der(A) \mid d \lambda_x =
\lambda_x d$ for all $x \in G\}$ of all {\em left-invariant
derivations} of the ring $A$ is a Lie algebra, since the bracket
of two derivations which commute with $\lambda_x$ obviously does
likewise. We call $\hbox{\rm Lie}(G)$ the {\em Lie algebra of G}.
Since the $p$-th power of any left-invariant derivation of $A$ is
again a left-invariant derivation, the Lie algebra $\text{Lie}(G)$
carries a canonical restricted Lie algebra structure.  \m

 We now gather a few basic  facts from Section 31 of \cite{H2}.
Let $G$ be a semisimple algebraic group.  Assume $T$ is a maximal
torus of $G$, \  $B^+ = TN^+$ is  a Borel subgroup containing $T$,
and  $B^- = TN^-$ is the opposite Borel subgroup. Suppose $\Delta$
is the base of  the root system $\Phi$ determined by $B^+$. Thus,
$\Phi$ is the set of  nonzero weights of $\hbox{\rm Ad}\, T$ in
$\hbox{\rm Lie}(G)$. The Weyl group $W$ of $G$ is defined to be
$N_{G}(T)/Z_{G}(T)$, the normalizer modulo the centralizer of the
torus $T$.  \m

Let $X(T)$ denote  the lattice of rational characters
of the torus $T$. The elements of $X(T)$ are often referred to as  {\it
weights}. Relative to
 the base $\Delta = \{\al_1,\dots,\al_m\}$, the fundamental
 weights $\{\varpi_1, \dots, \varpi_m\}\subset X(T)\ot_{\Z}\mathbb Q$
 satisfy $(\varpi_i, \alpha_j^\vee) =
 \delta_{i,j}$, where $\beta^\vee = 2(\beta,\beta)^{-1}\beta$ for any root $\beta$.
 If $\varpi_i \in X(T)$ for all $i$, then $G$ is said to be
 {\em simply connected}.  In this case, the fundamental weights form
 a basis of the $\Z$-module $X(T)$.
  (By identifying weights $\nu$ with their
 differentials  $(\mathrm{d}\nu)_e$,   we are able to use 
the same notation for the fundamental weights in $X(T)$ as for the fundamental  weights of the
 corresponding Lie algebra.)  \m

 A {\it
rational representation} of $G$ is a homomorphism of algebraic
groups $\varrho:  G \rightarrow \hbox{\rm GL}(M)$ for some
finite-dimensional vector space $M$. We say the resulting
$G$-module $M$ is rational. Every rational $G$-module $M$
admits a weight space decomposition $$M = \bigoplus_{\lambda \in X(T)} M^\lambda \nonumber$$
\noindent where $M^\lambda = \{v \in M \mid t v =
\chi_\lambda(t)v$ for all $t \in T\}$  (here $\chi_\lambda$ is the
rational character of $T$ corresponding to $\lambda \in X(T)$).
Whenever $M^\lambda \neq 0$, then $\lambda$ is said to be a {\it
weight} of $M$ with respect to $T$.    Using the fact that the
Weyl group $W=N_{G}(T)/Z_{G}(T)$, it is routine to check that $W$
permutes the weights of $M$, and $\dim M^\lambda=\dim
M^{w\lambda}$ for all $w\in W$. We let $X(M)$ denote the set of
weights of the module $M$. \m

{F}rom now on we suppose that $G$ is simply connected. Set

\begin{eqnarray*} X(T)_+ &:=& \bigoplus_{i=1}^m \mathbb N \varpi_i,
\quad \qquad \qquad (\mathbb N = \Z_{\geq 0}) \nonumber\\
 X_1(T) &:=& \left \{ \sum_{i=1}^m a_i \varpi_i \ \Bigg | \ a_i \in \mathbb N, \ 0
 \leq a_i \leq p-1 \right\} \nonumber\\
Q_+ &:=& \bigoplus_{i=1}^{m} \mathbb N \al_i.
\nonumber\end{eqnarray*} Elements of $X(T)_+$ are called the {\em
dominant weights of $T$} (with respect to the base $\Delta$).  \m
There is a partial ordering $\le$ on $X(T)$ defined by
$$\mu  \le  \lambda  \iff \lambda - \mu \in Q_+.$$
Let $X(M)$ denote the set of weights of the module $M$. Then
$\lambda \in X(M)$ is a {\it maximal} (resp. {\it minimal}) weight
of $M$ if $\mu \ge \lambda $ (resp.  $\mu \le \lambda$) for $\mu
\in X(M)$ implies $\mu = \lambda$.  If $\lambda$ is maximal weight
of $M$, then $\lambda \in X(T)_+$, and if $\lambda$ is a minimal
weight, then $\lambda \in -X(T)_+$.  An irreducible $G$-module $M$
has a unique maximal and a unique minimal weight, and the
corresponding weight spaces are one-dimensional.  Any two
finite-dimensional irreducible rational $G$-modules are isomorphic
if and only if they have the same maximal (minimal) weight. \m

Any $\lambda \in X(T)$ may be extended to all of $B^+$ by setting
$\lambda(u) = 1$ for all $u \in N^+$. Then $\lambda$ determines a
one-dimensional $B^+$-module $\F_\lambda$, where the $B^+$-action is
afforded  by $\lambda$. Let  $H^0(\lambda) = \hbox{\rm
ind}_{B^+}^G\,\, \F_\lambda \cong \{f \in \F[G] \mid f(gb) =
\lambda(b)^{-1}f(g)$ for all $g \in G$ and $b \in B^+\}$ be the
corresponding  induced module. The  socle $L(\lambda): = \hbox{\rm
soc}_G\, H^0(\lambda)$ is a simple $G$-module. Thus, for any
$\lambda \in X(T)_+$ there exists a finite-dimensional irreducible
rational $G$-module $L(\lambda)$ with maximal weight $\lambda$.
Although the dimension of $L(\lambda)$ is very difficult to
determine, in general, one knows that $\dim H^0(\lambda)$ is given
by Weyl's dimension formula:
\begin{equation}\label{eq:Weyldim} \dim
H^0(\lambda)=\prod_{\alpha\in\Phi^+}\frac{(\lambda+\rho,\alpha)}{(\rho,\alpha)},\end{equation}
where $\rho$ is half the sum of the positive roots (this is immediate from \cite[Part II, Prop.~5.10]{Ja}).
\bi

\begin{Pro} \label{Pro:2.802} {\rm (Compare \cite[Part II, Cor.~2.5]{Ja}.)}
The module dual to $L(\lambda)$ is $L(-w_0\lambda),$ where $w_0$
is the longest element of the Weyl group $W$.  \end{Pro}  \medskip

 For any
rational representation $\varrho\colon G\rightarrow \text{GL}(V)$,
the differential map
$$(\text{d}\varrho)_e\colon\,
\text{Lie}(G)\rightarrow \mathfrak{gl}(V)$$ is a homomorphism of
restricted Lie algebras. Thus, any rational $G$-module
$L(\lambda)$ can be regarded, canonically, as a restricted
$\text{Lie}(G)$-module. The following classical result is due to
Curtis \cite{C}: \bi

\begin{Pro} \label{Pro:2.803} {\rm (Compare \cite[Part II, Prop.~3.10 and Prop.~3.15]{Ja}.)}
Any irreducible restricted $\mathrm{Lie}(G)$-module is isomorphic
to exactly one module $L(\lambda)$ with $\lambda \in X_1(T)$.
Conversely, any $G$-module $L(\lambda)$ with $\lambda\in X_1(T)$
is irreducible as a $\mathrm{Lie}(G)$-module.
\end{Pro}   \medskip

The module $V(\lambda):= H^0(-w_0\lambda)^*$  which is dual to $H^0(-w_0\lambda)$
is called the {\it Weyl module} corresponding to $\lambda$.
There is an antiautomorphism $\tau$ of order 2 on $G$
which acts as the identity on $T$ and sends the root subgroup  $U_\alpha$ to $U_{-\alpha}$,
for each  $\alpha \in \Phi$
(as in \cite[Part II, Sec.~1.2]{Ja}).    Then $V(\lambda) \cong {}^\tau H^0(\lambda)$,
where the $G$-action on the right is twisted by $\tau$.   Using that fact,
one can show that the Weyl module $V(\lambda)$ is generated by
an eigenvector  $v$ for $B^+$ corresponding to $\lambda$.
Moreover, we have the following:   \bi

\begin{Pro} \label{Pro:2.804} \ {\rm (See {\cite[Part II, Sec.~2.13, 2.14(1)]{Ja}}.)}
\begin{itemize}
\item[{\rm (a)}] Any rational $G$-module $M$ which is generated by
a $B^+$-eigenvector $m$ of weight $\lambda$ is a homomorphic image
of $V(\lambda)$ via a homomorphism  sending $v$ to $m$. \item[{\rm
(b)}]  If $\hbox{\rm rad}_G V(\lambda)$ denotes the intersection
of all maximal submodules of $V(\lambda)$, then
$V(\lambda)/\hbox{\rm rad}_G V(\lambda) \cong L(\lambda)$.
\end{itemize}
\end{Pro} \medskip

The result stated below will play a crucial role in Chapter~4 of
this monograph (recall that we are assuming $G$ is simply connected):
 
\bi
\begin{Pro} \label{Pro:2.805}  {\rm (\cite[Part II, Prop.~9.24(b)]
{Ja};  see also  \cite[Part I, Sec.~9.6]{Ja}.)}  Let $M$ be a
finite-dimensional rational $G$-module such that
$$X_+(M) :=  \Big (X(M) \cap X(T)_+  \Big) \subseteq X_1(T).$$
Then any subspace of $M$ invariant under the induced action of
$\hbox{\rm Lie}(G)$ is a $G$-submodule of $M$.
\end{Pro}
\m

\section {\ Standard gradings of classical Lie algebras \label{sec:2.4}}
\m
\m In this section, we assume that $\g$ is either a classical simple Lie algebra
or $\mathfrak{pgl}_{m+1}$ with $p\mid m+1$. Our goal here is to recall
the notion of a {\it standard grading} of $\g$  and to describe
the gradings $\g=\,\bigoplus_{j\in\mathbb Z}\,\g_j$ of $\g$ such
that $\g_{-1}$ is an irreducible $\g_0$-module and $\g_-$ is
generated by $\g_{-1}$.

As explained in Section~\ref{sec:2.2}, $\g={\mathfrak
t}\oplus\bigoplus_{\alpha\in\Phi}\,\g^\alpha$, where $\mathfrak t$
is an abelian Cartan subalgebra of $\g$, and $[{\mathfrak
t},\g^\alpha]=\g^\alpha$ for all $\alpha\in\Phi$. Moreover, each
root space $\g^\alpha$ is one-dimensional, and $\exp(\ad x)\in
\text{Aut}\,\g$ for all root vectors $x\in\bigcup_{\alpha\in
\Phi}\,\g^\alpha$; see Lemma~\ref{Lem:2.7}. As in
Section~\ref{sec:2.2}, we regard $\Phi$ as a root system in an
Euclidean space. Let $\Delta=\{\alpha_1,\ldots,\alpha_m\}$ be a
basis of simple roots in $\Phi$. Recall that every root in $\Phi$
is a nonnegative or nonpositive integral linear combination of
roots in $\Delta$.

Associated to any $m$-tuple of {\it nonnegative} integers
$(a_1,\ldots,a_m)$, there is a grading of the Lie algebra $\g$ in
which a root vector corresponding to a root
$\alpha=\sum_{i=1}^mn_i\alpha_i$ is given gradation degree
$\sum_{i=1}^m a_in_i$ and all elements in $\mathfrak t$ are assigned
gradation degree $0$. If $\sum_{i=1}^ma_i=1$, i.e. if $a_k=1$ for
some $k \in \{1, \dots,  m\}$ and $a_j=0$ for $j\ne k$, then we say that  the grading is
{\it standard}.

Two gradings $L=\bigoplus_{j\in\mathbb Z}\,L_j$ and
$L=\bigoplus_{j\in\mathbb Z}\,L'_j$ of a Lie algebra $L$ are said
to be {\it conjugate} if there is an automorphism $\sigma$ of $L$
such that $\sigma(L_j)=L'_j$ for all $j\in\mathbb Z$.    In this
monograph, we identify conjugate gradings without further notice.
Our next result is known to the experts;  we provide a short proof for
the reader's convenience. \m

\begin{Pro} \label{standard}
Every grading $\g=\bigoplus_{j\in\mathbb Z}\,\g_j$ of the Lie
algebra $\g$ such that $\g_{-1}$ is an irreducible $\g_0$-module
and $\g_-$ is generated by $\g_{-1}$  is standard.
\end{Pro}

\pf Let $G\eqdef (\text{Aut}\,\g)^\circ$, the connected component of
the identity in the automorphism group of $\g$,  which is  a simple algebraic
$\F$-group of adjoint type. It is well-known that  $\text{Lie}(G)$
is isomorphic to a Lie subalgebra of $\Der(\g)$. In view of
Lemma~\ref{Lem:2.7}, we have that
$\ad\g^\alpha\subset\text{Lie}(G)$ for all $\alpha\in\Phi$. If
$\g$ is not of type ${\mathrm A}_m$ with $p\,|\,(m+1)$, then $\g$
is a simple Lie algebra and $\Der(\g)=\,\ad\,\g$. Since
$\bigcup_{\alpha\in\Phi}\ad\,\g^\alpha$ generates a nonzero ideal
of $\ad\,\g$, it follows that $\text{Lie}(G)=\ad\,\g$ in this
case. If $\g$ is of type ${\mathrm A}_m$ with $p\,|\,(m+1)$, then
either $\g=\mathfrak{psl}_{m+1}$ or $\g=\mathfrak{pgl}_{m+1}$  by
our assumption. In this case, $\Der(\g)\cong\mathfrak{pgl}_{m+1}$
implying $\text{Lie}(G)\cong\mathfrak{pgl}_{m+1}$.  Thus in any
event,  $\ad\,\g$ is an ideal of codimension $\le 1$ in
$\text{Lie}(G)$.

For ease of notation, we identify $\g$ with $\ad\,\g\subseteq
\text{Lie}(G)$. It follows from the preceding remark that $G$
contains a maximal torus $T$ such that
$\text{Lie}(T)\cap\g=\mathfrak t$. Moreover, the theory of
algebraic groups enables us to identify $\Phi$ with the root
system of $G$ relative to $T$.

Given a linear algebraic $\F$-group $H$,  we denote by $X_*(H)$ the
set of all rational homomorphisms from the one-dimensional torus
$\F^\times$ to $H$. Associated with our grading $\g=\bigoplus_{j\in\mathbb Z}\,\g_j$  is a
rational homomorphism $\lambda\in X_*(G)$ such that

\begin{equation}\label{grading}\lambda(t)\vert_{\g_j}=\,t^j\,\text{id}_{\g_j}\qquad\quad \
\big(\forall\,t\in\F^\times,\,\,\,\forall\,j\in{\mathbb
Z}\big).\end{equation} 
Clearly, $\lambda(\F^\times)$ is a
one-dimensional torus of $G$. Since any torus of $G$ is contained
in a maximal torus and all maximal tori of $G$ are conjugate, we
can assume in the rest of the proof that
$\lambda(\F^\times)\subseteq T$. In other words, we can assume
that $\lambda$ belongs to $X_*(T)$, the lattice of rational
cocharacters of $T$. Note that the Weyl group $W=N_G(T)/T$ acts on
$X_*(T)$.

 Let $X(T)$ be the lattice of characters of $T$; see Section~\ref{sec:2.3}.
 According to \cite[(8.6)]{Bo},
 the $\mathbb Z$-pairing 
 
 \begin{equation}\label{pairing}
 X(T)\times X_*(T)\rightarrow\,{\mathbb Z},\quad
 \quad (\chi,\nu)\mapsto n\eqdef\la\chi,\nu\ra\  \,\,\mbox{ if }\,\, \chi(\nu(t))=t^n,
 \end{equation}
 is $W$-invariant and nondegenerate.
 As $G$ is a group of adjoint type, the lattice $X(T)$ is spanned
over $\mathbb Z$ by the set of simple roots $\Delta$.
Consequently, the pairing (\ref{pairing}) identifies $X_*(T)$ with
$P$, the weight lattice corresponding to $\Phi$. Replacing
$\lambda$ by its $W$-conjugate if necessary, we can assume in what
follows that $\lambda$ lies in the dual Weyl chamber
corresponding to $\Delta$.  That is to say, we can assume that
$\alpha_i(\lambda(t))=t^{a_i}$ for some nonnegative integers
$a_i$, where $1\le i\le m$.

In conjunction with (\ref{grading}),  the above discussion shows
that $\g_0$ is spanned by $\mathfrak t$ and the root spaces
$\g^\alpha$ with $\alpha=\sum_{i=1}^m n_i\alpha_i$ satisfying
$n_k=0$ whenever $a_k\ne 0$. It also shows that $\g_{-1}$ is
spanned by the root spaces $\g^\gamma$ with $\gamma=\sum_{i=1}^m
r_i\alpha_i$ for which $\sum_{i=1}^m r_ia_i = -1$. All such
$\gamma$'s are negative roots in $\Phi$. Since $\g_{-1}\ne 0$, it
must be that $a_d=1$ for some $d\in \{1,\dots, m\}$.

Let $\g_{-1,d}$ denote the span of all $\g^\gamma\subset \g_{-1}$
with $\gamma=\sum_{i=1}^m r_i\alpha_i$ such that $r_d=-1$. Since
$\g^{-\alpha_d}\subseteq \g_{-1,d}$, we have $\g_{-1,d}\ne 0$.
Moreover, because  $\g_{-1}$ is an irreducible $\g_0$-module, our description
of $\g_0$ now implies that $\g_{-1}=\g_{-1,d}$. Suppose $a_j\ne 0$
for some $j\ne d$. Then $\g^{-\alpha_j}\subseteq \g_-$. On the
other hand, if $\g^{-\beta}$ lies in the Lie subalgebra generated
by $\g_{-1,d}$, then $\beta=\sum_{i=1}^m q_i\alpha_i$ with $q_d\ge
1$. It follows that $\g^{-\alpha_j}$ is not contained in the 
subalgebra generated by $\g_{-1,d}=\g_{-1}$. Since this
contradicts one of our assumptions, we derive that $a_j=0$ for
$j\ne d$. Thus, our grading is standard, as stated.\qed

\m

\section  {\  The Lie algebras of
Cartan type  \label{sec:2.5}}

\m  There are four infinite series of finite-dimensional simple,
graded
 Cartan type Lie algebras $\g = \bigoplus_{j=-q}^r \g_j$: \ the
\W, special, Hamiltonian, and contact series.  The Lie algebras in each
 series vary in ``thickness" (i.e., dimensions of corresponding
 gradation spaces $\g_j$) and ``length" (i.e., number of nonzero
 gradation spaces).  The special, Hamiltonian, and contact Lie
algebras
 are subalgebras of the \W \ Lie algebras.  The \W, 
 special, and Hamiltonian Lie algebras all have depth-one gradations
 ($q=1$) with respect to which they are transitive  \eqref{eq:1.3} and
 irreducible \eqref{eq:1.2};  while the contact Lie algebras have a
 depth-two gradation ($q=2$), which is transitive and irreducible.

\m  Each series contains both restricted and nonrestricted Lie
algebras.   Each nonrestricted algebra contains a restricted
subalgebra (its ``initial piece'').  The restricted ones, which
are exactly the initial pieces of the algebras in a given series,
are generated as Lie algebras by their local part $\g_{-1} \oplus
\g_0 \oplus \g_1$.  The local part of a nonrestricted graded
Cartan type Lie algebra coincides with the local part of its
initial piece.

 \m When $\g$ is a restricted Cartan type Lie algebra, then $x^{[p]}
= 0$
 for $x \in \g_{-1}$, since $(\ad \, x)^p$ is a derivation which
 annihilates the local part  when $p >  3$.  The restriction of the
 $[p]$-mapping of $\g$ to $\g_0$ coincides with the $[p]$-mapping on
 the classical Lie algebra $\g_0^{(1)}$, and the adjoint
representation
 of $\g_0$ on each $\g_j$ is a $p$-representation, i.e.,  is restricted.

 \m We will illustrate all these concepts by beginning with a few
simple
 examples of Cartan type Lie algebras.  Then we will introduce the
 divided power algebras and their derivations and use them to
describe
 the four infinite series of Cartan type Lie algebras.

 \m
\section {\ The Jacobson-Witt
algebras  \label{sec:2.6}}

\m The Jacobson-Witt Lie algebras were the first
finite-dimensional nonclassical simple Lie algebras to be
discovered.  Their construction starts with a polynomial algebra
$\F[x_1,\dots,x_m]$ over $\F$, a field of characteristic $p > 0$.
Let $D_i = \partial/\partial{x_i}$, $i=1,\dots,m$, be the usual
partial derivatives relative to the variables $x_i$.  The ideal
generated by the elements $x_1^p,\dots,x_m^p$ is invariant under
the derivations $D_i$, and so the $D_i$ induce derivations on the
truncated polynomial algebra $\F[x_1,\dots,x_m]/ \langle
x_1^p,\dots,x_m^p \rangle$.  The Jacobson-Witt Lie algebra is the
derivation algebra $W(m;\underline{\one}) :=
\Der(\F[x_1,\dots,x_m]/ \langle x_1^p,\dots,x_m^p \rangle)$.

\m In the special case that $m=1$, the resulting Lie algebra
$W(1;\un{1})$ is called the {\it $p$-dimensional Witt algebra} or simply the {\it Witt algebra}.  The Jacobson-Witt algebras may be viewed as ``thickenings'' of it by addition of
variables.

\m By identifying cosets with their representatives, we may assume
that the elements $x^a: = x_1^{a_1} \cdots x_m^{a_m}$ with $a =
(a_1,\dots,a_m)$ and $0 \leq a_i < p$ for all $i$ determine a
basis for $\F[x_1,\dots,x_m]/ \langle x_1^p,\dots,x_m^p \rangle$.
Then the derivations $x^a D_i$, as $a$ ranges over such $m$-tuples
and $i=1,\dots,m$, comprise a basis for $W(m;\un{1})$.  The Lie
bracket in $W(m;\un{1})$ is given by

$$[x^aD_i, x^b D_j] = b_ix^{a+b-\epsilon_i} D_j -
a_jx^{a+b-\epsilon_j} D_i,\nonumber$$

\noindent where $\epsilon_i$ is the $m$-tuple with $1$ in the
$i$th position and 0 elsewhere, and addition of tuples is
componentwise.  If any component exceeds $p-1$ or is less than 0,
the term is 0.

\m Because the Lie algebra $W(m;\un{1})$ is the derivation algebra
of an algebra, it carries a natural $[p]$-structure given by the
$p$-mapping:

$$D^{[p]} = D^{p} \qquad \hbox{\rm for all} \ D \in W(m;\un{1}). \nonumber$$

\m Rather than working with monomials as above, we could use
divided powers instead.  Thus, if we adopt the basis $x^{(a)}
=\displaystyle{ \frac{ x_1^{a_1} \cdots x_m^{a_m}}{a_1!   \cdots
a_m! }}$, where again $0 \leq a_i < p$ for all $i$, then 

\begin{equation}\label{eq:divpow} x^{(a)}x^{(b)} = {a+b \choose a} x^{(a + b)}, \qquad \hbox{\rm
where} \quad {a+b \choose a} = \prod_{j=1}^m {a_j+b_j \choose
a_j}. \end{equation}

\noindent Moreover, $W(m;\un{1})$ has a basis $\{x^{(a)}D_i \mid 0
\leq a_i < p, \ i=1,\dots,m\}$ where

\begin{eqnarray*} D_i(x^{(a)}) &=& x^{(a-\epsilon_i)} \qquad
\hbox{and} \nonumber\\
{[x^{(a)}D_i, x^{(b)}D_j]} &=& {a+b-\epsilon_i \choose
a}x^{(a+b-\e_i)}D_j - {a+b-\epsilon_j \choose b}x^{(a+b-\e_j)}D_i.
\nonumber\end{eqnarray*}

It is this form of the algebra that lends itself naturally to
further generalizations. \m

\section {\ Divided power
algebras \label{sec:2.7}}   
 
\m Let $\F$ be a field of characteristic $p > 0$, and let
$\mathcal O(m)$ denote the commutative associative algebra with 1
over $\F$ defined by generators $x_i^{(k)}$, $1 \leq i \leq m$, \ $k
\in \mathbb N = \Z_{\geq 0}$,  which satisfy the relations

\begin{eqnarray*} x_i^{(0)} &=& 1 \nonumber\\
x_i^{(k)}x_i^{(\ell)} &=& {k+\ell \choose k} x_{i}^{(k+\ell)},
\qquad 1 \leq i \leq m, \ \ k,\ell \in \mathbb N.
\nonumber\end{eqnarray*}

The ``monomials''  $x^{(a)}:= x_1^{(a_1)} \cdots x_m^{(a_m)}$, \ \
$a_i \in \mathbb N$, constitute a basis of $\mathcal O(m)$.

\m  To compute the binomial coefficients in  \eqref{eq:divpow},
it is helpful to consider base $p$-expansions.   
Assume that the base $p$ expansions of
$c$ and $d$ are given by

\begin{eqnarray} \label{eq:basep} c &=& c_k p^k+ c_{k-1}p^{k-1} + \cdots + c_1 p +
c_0\\
d &=& d_kp^k + d_{k-1}p^{k-1} + \cdots + d_1 p + d_0  \nonumber 
\nonumber\end{eqnarray}

\n where $0 \leq c_i,d_i < p$ for all $i$.  Then

\begin{equation}\label{eq:2.12} {c \choose d} = \prod_{i=0}^k {{c_i}
\choose {d_i}}.  \end{equation}

\m
\begin{Lem}\label{Lem:binom}  Suppose $c, d \in \mathbb N$ satisfy
$0\leq c,d < p^{k+1}$, but $c+d \geq  p^{k+1}$.   Then
$\displaystyle{{{c+d} \choose d} = 0}$.  \end{Lem}

\bi
\pf   This will follow from analyzing the ``carries''
in \eqref{eq:2.12} (with top number $c+d$ instead of $c$).    
Let the base $p$ expansions of $c$ and $d$ be as in
\eqref{eq:basep}.     Since $c+d \geq p^{k+1}$,   there is some smallest integer $i \geq 0$ so
that $c_i + d_i \geq p$.    Then the binomial coefficient
$\displaystyle{{c+d} \choose d}$ must equal 0, because  it contains the factor $\displaystyle{{c_i+d_i-p} \choose d_i}$,
which is $0$.   \qed

\m  \bi

 As a result of Lemma \ref{Lem:binom}, we see for any $m$-tuple $\un n =
(n_1,\dots,n_m)$ of positive integers  that

$$\mathcal
O(m;\un n) := \spa_{\F}\{ x^{(a)} \mid 0 \leq a_i <
p^{n_i}\}\nonumber$$

\noindent is a subalgebra of the divided power algebra $\mathcal
O(m)$, and we have containment $\mathcal O(m;\un n) \subseteq
\mathcal O(m;\un n')$ whenever $n_i \leq n_i'$ for all
$i=1,\dots,m$.

\m

\section {\ \ Witt Lie
algebras of Cartan type \ (the $\boldsymbol {W}$  series) \label{sec:2.8}}

\m Nonrestricted Cartan type Lie algebras can be obtained by
``lengthening" the restricted ones; that is, by allowing powers
greater than $p-1$ in the divided power variables.

\m Let $D_i$  ($1 \leq i \leq m$)  be the derivation of $\mathcal
O(m)$ defined by

\begin{equation*} D_i(x_j^{(r)}) = \delta_{i,j}x_i^{(r-1)}.\end{equation*}

The (infinite-dimensional)  {\it Witt  Lie algebra of Cartan type} is the subalgebra of
$\Der(\mathcal O(m))$ defined by

\begin{equation*} W(m):= \spa_{\F} \{ x^{(a)}D_i \mid a_i \in \mathbb N,
i=1,\dots,m\},  \end{equation*}

\noindent and having Lie bracket

\begin{equation}\label{eq:2.15}  {[x^{(a)}D_i, x^{(b)}D_j]} =
{a+b-\epsilon_i \choose a}x^{(a+b-\epsilon_i)}D_j -
{a+b-\epsilon_j \choose b}x^{(a+b-\epsilon_j)}D_i.
\end{equation}
It is a free $\mathcal O(m)$-module with basis $D_1,\dots, D_m$.

\m The fact that $\mathcal O(m;\un n)$ is a subalgebra of
$\mathcal O(m)$ for any $m$-tuple $\un n$ of positive integers
implies that

$$W(m;\un n):= \spa_{\F} \{ x^{(a)}D_i \mid 0 \leq a_i <
p^{n_i}, \ i=1,\dots,m\},\nonumber$$

\noindent is a Lie subalgebra of $W(m)$.   We refer to the Lie algebras
$W(m;\un n)$ as the {\em \W \ Lie algebras of Cartan type},  or
more succinctly,  the {\em Witt algebras}.

\m The simplest examples of Lie algebras of this kind are the
$p^n$-dimensional {\it Zassenhaus Lie algebras} $W(1;\un{n})$
obtained from the divided power algebras $\mathcal O(1;\un{n})$
for $n=1,2,\dots$.  Each Zassenhaus Lie algebra contains a copy of
the $p$-dimensional Witt algebra $W(1;\un{1})$.

\m More generally, an algebra $W(m;\un n)$ contains a copy of the
restricted Lie algebra $W(m;\un \one)$ by taking the derivations
$\sum_{i=1}^mf_i D_i$ with $f_i \in \mathcal O(m;\un {1})$.
 
\bi 
\begin{Thm}  \label{Thm:2.16} \ {\rm (See for example,
{\cite[ Prop.~2.2 and Thm.~2.4]{SF}.)}}
\begin{enumerate}
\item[{\rm (i)}] $W(m;\un {n})$ is a simple Lie algebra except
when $m = 1$ and $p =2$. \item[{\rm (ii)}] $W(m;\un{n})$ is a free
$\mathcal O(m;\un{n})$-module with basis $\{D_1,\dots,D_m\}$.
\item[{\rm (iii)}] The elements in $\{x^{(a)}D_i \mid 0 \leq a_i <
p^{n_i}, \  \ 1 \leq i \leq m\}$ determine a basis for
$W(m;\un{n})$ so that $\dim W(m;\un{n}) = mp^{n_1 + \cdots +
n_m}$.
\end{enumerate}
\end{Thm}

\bi The divided power algebra is $\Z$-graded by the subspaces $\mathcal
O(m)_{k} =$ \break  $\spa_\F \{x^{(a)}\, \big | \,   |\, a\, |
:= \sum_{j=1}^m a_j = k\}$. Thus $\mathcal O(m)$ $=
\bigoplus_{k=0}^{\infty}\mathcal O(m)_{k}$, and correspondingly,
$\mathcal O(m;\un n)$ $= \bigoplus_{k=0}^{r+1}\mathcal O(m;\un
n)_{k}$ where $r = p^{n_1}+ \cdots+ p^{n_m}- m-1$.  The associated
Lie algebras $W(m)$ and $W(m;\un n)$ inherit a grading by
assigning

\begin{equation}W(m)_k = \bigoplus_{j=1}^m \mathcal O(m)_{k+1}
D_j \ \hbox{\rm \ and} \ \ W(m;\un n)_k = \bigoplus_{j=1}^m
\mathcal O(m;\un n)_{k+1}D_j.\nonumber\end{equation}

\noindent In particular,

\begin{equation*} W(m;\un n) =\bigoplus_{j=-1}^r W(m;\un n)_j\end{equation*}

\noindent is a depth-one simple Lie algebra of height $r$.

\m Observe that the null component $W(m;\un n)_0 = W(m)_{0}$ has
the elements $x_i^{(1)}D_j$ ($1 \leq i,j \leq m$) as a basis.  Moreover,

\begin{equation}\label{eq:glm}  [x_i^{(1)}D_j, x_k^{(1)}D_\ell] = \delta_{j,k} x_i^{(1)}D_\ell - \delta_{i,\ell}
x_k^{(1)}D_j.\end{equation}

\noindent Therefore, it is evident that the null component is
isomorphic to $\mathfrak{gl}_m$ via the isomorphism that takes
$x_i^{(1)}D_j$ to the matrix unit $E_{i,j}$.

\m The
derivation
\begin{equation}\label{eq:degder}  \mathfrak D_1: = \sum_{j=1}^m x_i^{(1)}D_i  \end{equation}
(which corresponds to the identity matrix in  $\mathfrak{gl}_m$)  plays a special role in $W(m)$, as it is the ``degree derivation'':
\begin{equation*} \mathfrak D_{1}(x^{(a)}) = \sum_{j=1}^m x_j^{(1)} D_j (x^{(a)}) =
\Bigg(\sum_{j=1}^m a_j \Bigg)x^{(a)} =
\deg(x^{(a)})x^{(a)}. \end{equation*}
{F}rom the relation

\begin{equation}\label{eq:derive} [f D, gE] = f D(g) E - g E(f)D + fg [D,E], \end{equation}

\noindent which holds for all $D,E \in W(m)$ and $f,g \in \mathcal
O(m)$, we have  as a special case,

\begin{equation}\label{eq:derdeg}  \left [\mathfrak D_1, x^{(a)}D_j\right ] =
\mathfrak D_1(x^{(a)})D_j + x^{(a)}[\mathfrak D_1,D_j] = \big(\deg(x^{(a)})-1\big)x^{(a)}D_j, \end{equation}
so that $\ad \mathfrak D_1$ acts as multiplication by the  scalar
$\ell$ on $W(m)_\ell$ for each $\ell$; that is to say, $\ad \mathfrak D_1$
is the degree derivation on $W(m)$.

\m Now it follows {f}rom

\begin{equation}\label{eq:dualglm}[x_i^{(1)}D_j, D_k] = - \delta_{i,k} D_j\end{equation}

\noindent that $W(m)_{-1} = W(m;\un n)_{-1} =
\spa_\F \{D_1,\dots,D_m\}$ can be identified with the
dual module of the natural $m$-dimensional module $V$ for
$\mathfrak {gl}_m$.  Hence it is isomorphic to the space $V^*$ of
$1 \times m$ matrices over $\F$ with the $\mathfrak{gl}_m$-action
given by $y.v = - v y$ (matrix multiplication) for all $v \in V^*$
and $y \in \mathfrak {gl}_m$.

\m The mapping $\varphi: W(m)_0 \rightarrow
\mathfrak{gl}(\mathcal O(m)_1)$ given by $\varphi(D)(f) =
D(f)$ affords a representation of $W(m)_0$ (hence also
of $W(m;\un n)_0$ for any $m$-tuple $\un n$).  Moreover,
relative to the standard basis $\{x_1^{(1)},\dots,x_m^{(1)}\}$ of $\mathcal O(m)_1$, the matrix of $\varphi(x_i^{(1)}D_j)$ is exactly the matrix unit
$E_{i,j}$.

\m

\section {\ Special Lie
algebras of Cartan type \  (the $\boldsymbol {S}$ series)  \label{sec:2.9}}

 \m The Lie algebra $S(m;\un n)$ (resp. $S(m)$) is the subalgebra of
$W(m;\un n)$ (resp.  $W(m)$) of derivations whose divergence is
zero.  Thus, if $D =
 \sum_{i=1}^m f_i D_i$, then

$$\di(D) := \sum_{i=1}^m D_i(f_i) = 0.\nonumber$$

\noindent Observe that

\begin{eqnarray*}
\di ([fD_i,gD_j]) &=&
\di \bigl(f D_i(g) D_j - D_j(f) g D_i\bigr) \nonumber\\
&=& D_j(f)D_i(g) + f D_jD_i(g) - D_iD_j(f)g - D_j(f)D_i(g) \nonumber\\
&=& fD_i\bigl(D_j(g)\bigr)-gD_j\bigl(D_i(f)\bigr)\nonumber \\
&=& fD_i\bigl(\di (gD_j)\bigr) -
gD_j\bigl(\di (fD_i)\bigr), \nonumber\end{eqnarray*}

\noindent so summing over such terms will show that the relation

\begin{equation}\label{eq:divprod} \di ([D,E]) =
D\bigl(\di (E)\bigr) - E\bigl(\di (D)\bigr)
\end{equation}

\noindent holds for all $D,E \in W(m)$.  From that it is apparent
that $S(m;\un n)$ (resp.  $S(m)$) is indeed a subalgebra of
$W(m;\un n)$ (resp. $W(m)$).

\m To better understand the structure of $S(m;\un n)$, we
introduce mappings $D_{i,j}: {\mathcal O}(m;\un n) \rightarrow
W(m;\un n)$ \ defined by

\begin{equation} \label{eq:dijdef} D_{i,j}(f) = D_j(f)D_i - D_i(f) D_j, \qquad 1 \leq i,j \leq m.\end{equation}

\noindent It is easy to check that the image of $D_{i,j}$ lies in
$S(m;\un n)$.  Moreover, $D_{i,i} = 0$ and $D_{j,i} = -D_{i,j}$.

\bi
\begin{Lem} \label{Lem:2.18} \rm {(Compare {\cite[Lem.~3.2]{SF} for (a) and (b).)}}
\begin{enumerate}
\item[{\rm (a)}] For $D=\sum_{i=1}^m f_i D_i$ and $E =
\sum_{j=1}^m g_j D_j$ in $S(m;\un n)$,

$$\left [ \sum_{i=1}^m f_i D_i, \sum_{j=1}^m g_j D_j \right] =
\sum_{i,j=1}^n D_{j,i}(f_i g_j).\nonumber$$

\item[{\rm (b)}] $[D_k, D_{i,j}(f)] = D_{i,j}\bigl(D_k(f)\bigr)$
for all $1 \leq i,j,k\le m$.

\smallskip

\item[{\rm (c)}] For $1\le i,j,k,\ell\le m$ we have
\begin{eqnarray*}
[D_{i,j}(f),\, D_{k,\ell}(g)] &=&  D_{i,\ell}\Big(D_j(f)D_k(g)\Big) - D_{j,\ell}\Big(D_i(f)D_k(g)\Big) \\
&& \  \ -D_{k,j}\Big(D_i(f)D_\ell(g)\Big) +
D_{k,i}\Big(D_j(f)D_\ell(g)\Big). \end{eqnarray*}
\end{enumerate} \end{Lem}

\pf   We know  that

\begin{equation*} [D,E] = \left [ \sum_{i=1}^m f_i D_i, \sum_{j=1}^m g_j D_j \right] =
\sum_{i,j = 1}^m f_i D_i(g_j)D_j - g_j D_j(f_i)D_i \end{equation*}

\noindent from the general expression \eqref{eq:derive} for the
Lie bracket.   Then using the fact that $\di (D) = 0 =  \di (E)$,
we can see this product equals

\begin{eqnarray*} &&\sum_{i,j = 1}^m
\Big(f_iD_i(g_j)+D_i(f_i)g_j\Big)D_j -
\sum_{i,j = 1}^m \Big(g_jD_j(f_i)+D_j(g_j)f_i\Big)D_i\nonumber\\
&&\hskip 1 truein  = \sum_{i,j = 1}^m D_i(f_ig_j)D_j -
D_j(f_ig_j)D_i \nonumber\\
&&\hskip 1 truein = \sum_{i,j = 1}^m D_{j,i}\bigl(f_ig_j\bigr).
\nonumber\end{eqnarray*}

Part (b) is a simple consequence of the fact that
\begin{equation*}  [D_k,D_{i,j}(f)] =
D_k \bigl(D_j(f)\bigr) D_i - D_k\bigl(D_i(f)\bigr) D_j =
D_{i,j}\bigl(D_k(f)\bigr).
\end{equation*}

The verification of the identity in (c) is left as an exercise for
the reader.  The left side is the sum of 8 terms.  When expanded
out, the right side is the sum of 16 terms, 8 of which match
with the terms on the left, and 8 of which pairwise sum to 0.
\qed

\m

The subalgebra $S(m;\un n)$ is not simple.  However, its
commutator ideal $S(m;\un n)^{(1)}$ is, and we have the following.
 
\bi

\begin{Thm} \label{Thm:2.19} \  {\rm (Compare {\cite[Prop.~3.3,
Thm.~3.5, and Thm.~3.7]{SF}}.)} \ Suppose $m \geq 3$.
\begin{enumerate}

\item[{\rm (i)}] $S(m;\un n)^{(1)}$ is the subalgebra of $S(m;\un
n)$ generated by the elements \hfil \break
$D_{i,j}\bigl(x^{(a)}\bigr)$ \ where $0 \leq a_k < p^{n_k}$ for $1
\leq k \leq m$ and $1 \leq i,j \leq m$.

\item[{\rm (ii)}] $S(m;\un n)^{(1)}$ is simple.

\item[{\rm (iii)}] $S(m;\un n)^{(1)} = \bigoplus_{j=-1}^s (S(m;\un
n)^{(1)})_{j}$ where $s = (\sum_{j=1}^m
p^{n_j}) - m -2$ and 
$(S(m;\un n)^{(1)})_{j} = 
S(m;\un n)^{(1)} \cap W(m;\un n)_j$.

\item[{\rm (iv)}]  $S(m;\un{n}) = S(m;\un n)^{(1)}  \oplus
\bigoplus_{j=1}^m \mathbb F x^{\left(\tau -
(p^{n_j}-1)\epsilon_j\right)}D_j$,  where $\epsilon_j$ is the
$m$-tuple with $1$ in slot $j$ and 0 elsewhere,  and
$$\tau =
\tau(\un{n}) = (p^{n_1}-1, \dots, p^{n_m}-1).$$
\end{enumerate}
\end{Thm}

\m

We refer to the Lie algebras $S(m;\un n)^{(1)}$  as the simple
{\it special Lie algebras of Cartan type}. More generally, each of
the Lie algebras $S(m;\un n)$, $S(m)$, or $S(m;\un n)^{(1)}$ is
said to be a special Cartan type Lie algebra.

\m As $D_{i,j}(x_j^{(1)}) = D_i$, we have $(S(m;\un n)^{(1)})_{-1} =
W(m;\un n)_{-1}$.  The derivations $D_{j,i}\bigl(x_i^{(2)}\bigr) = x_i^{(1)}D_j$ belong to
$(S(m;\un n)^{(1)})_{0}$ for all $i \neq j$, as
do the derivations $D_{i,i+1}\bigl(x_i^{(1)}x_{i+1}^{(1)}\bigr) = x_i^{(1)}D_i -
x_{i+1}^{(1)}D_{i+1}$ for $i=1,\dots, m-1$.  Thus, $\dim (S(m;\un
n)^{(1)})_0 \geq m^2 -1$.  But the derivations $x_i^{(1)}D_i$ do not
belong to $(S(m;\un n)^{(1)})_0$, as they fail the divergence zero
test.  Thus, $\dim (S(m;\un n)^{(1)})_0 = m^2 -1$, and it is clear
that $(S(m;\un n)^{(1)})_0 \cong \mathfrak{sl}_m$.

\m
There is a Lie algebra closely aligned to $S(m;\un n)$ obtained
by adjoining the degree derivation
$\mathfrak D_1= \sum_{j=1}^m x_i^{(1)}D_i$
to the
null component.   As $[\mathfrak D_1, E] = k
E$ for all $E$ $\in W(m;\un n)_k$ (compare \eqref {eq:derdeg}), the result of adding the degree
derivation to $S(m;\un n)$ is a Lie algebra, which we denote
$CS(m;\un n)$, whose null component is isomorphic to
$\mathfrak{gl}_m$.

\m

\section{Hamiltonian
Lie algebras of Cartan type \ (the $\boldsymbol{H}$  series)   \label{sec:2.10}}

 \m In this section, we will introduce a third series of Cartan type
Lie
 algebras, but first some notation is required.  Set

\begin{eqnarray}\label{eq:2.21} \sigma(j) &=& \begin{cases} 1 & \quad
\qquad \hbox{if} \ \ 1 \leq j \leq m\\
-1 & \qquad \quad \hbox{if} \ \ m+1 \leq j \leq 2m, \end{cases} \\
j' &=& \begin{cases}j+m & \qquad \hbox{if} \ \ 1 \leq j \leq m\\
j-m& \qquad \hbox{if} \ \ m+1 \leq j \leq 2m.  \end{cases}
\end{eqnarray}

Set
\begin{eqnarray}\label{eq:2.22}  H(2m;\un n)&:=& \Biggl\{ D =
\sum_{i=1}^{2m} f_iD_{i}\in W(2m;\un n)\, \Big | \,
\sigma(j')D_{i}(f_{j'}) =
\sigma(i')D_{j}(f_{i'}), \nonumber \\
&& \hskip 2.3 truein  \quad 1 \leq i,j \leq 2m\Biggr \}.
\end{eqnarray}

\noindent  It can be argued that $H(2m;\un n) = W(2m;\un n) \cap
H(2m)$, where $H(2m)$ is the subalgebra of $W(2m)$ consisting of
all the derivations which satisfy $D(\omega_H) = 0$, where
$\omega_H$ is the differential form given by

\begin{equation*}\omega_H = \sum_{j=1}^{m} d{x_j} \wedge d{x_{j+m}}, \end{equation*}

\noindent  where $x_j = x_j^{(1)}$ for all $j$.  In fact, the defining condition in \eqref{eq:2.22} is
equivalent to the statement that $D(\omega_H) = 0$.

\m To gain further insight into the structure of $H(2m;\un n)$, we
define

\begin{equation} \label{eq:2.23} D_H: \mathcal O(2m) \rightarrow W(2m),
\qquad D_H(f) = \sum_{j=1}^{2m}\sigma(j) D_j(f)D_{j'}
\end{equation}

\noindent and  denote the image of the map $D_H$ restricted to
$\mathcal O(2m;\un n)$ by $\tilde H(2m;\un n)$.  The elements
$D_H(f)$ belong to $H(2m;\un n)$ for all $f \in \mathcal O(2m;\un
n)$, but the image of $D_H$ is proper, as the derivations
$x_j^{(p^{n_j}-1)}D_{j'}$ for $1 \leq j \leq 2m$ belong to
$H(2m;\un n)$, but fail to lie in the image.

\bi \begin{Thm} \label{Thm:2.24} \ {\rm (See {\cite[Lem.~4.1
and Thm.~4.5]{SF}}.)}
\begin{enumerate}
\item[{\rm (i)}]  Let $D = \sum_{i=1}^{2m} f_iD_i$ and $E =
\sum_{j=1}^{2m} g_j D_j$ be elements of $H(2m;\un n)$.  Then

$$[D,E] = D_H(u) \qquad \hbox{\rm where} \quad  u = \sum_{j=1}^{2m}
\sigma(j)f_jg_{j'}.\nonumber$$

\item[{\rm (ii)}] $\tilde H(2m;\un n) \supseteq  H(2m;\un
n)^{(1)}$.

\item[{\rm (iii)}] $H(2m;\un n)^{(2)}$ is a simple Lie algebra
with basis
\begin{equation*}\{ D_H(x^{(a)}) \mid a \neq (0,\dots,0) \ \ \hbox{\rm and} \ \ a \neq
(p^{n_1}-1, \dots, p^{n_{2m}}-1)\}.\end{equation*}  Thus,

$$\dim H(2m;\un n)^{(2)} = p^{n_1+ \cdots + n_{2m}}
-2.\nonumber$$

\item[{\rm (iv)}] For all $f,g \in \mathcal O(2m;\un n)$,

 \begin{eqnarray}\label{eq:2.25} {[D_H(f), D_H(g)]} &=&
 D_H\bigl(D_H(f)(g)\bigr) \\
&=& D_{H}\left ( \sum_{j=1}^{2m} \sigma(j)D_j(f)D_{j'}(g)\right).  \nonumber
\end{eqnarray}

\end{enumerate} \end{Thm}

\bi Some comments are in order.  First, it follows from assertion
(i) that $H(2m;\un n)^{(1)}$ is contained in $\tilde H(2m;\un n)$,
and hence that $\tilde H(2m;\un n)$ is an ideal of $H(2m;\un n)$.
The formula in (iv) is a consequence of (i), since

\begin{eqnarray}\label{eq:Hprod} [D_H(f), D_H(g)] &=& D_H(u) \qquad \hbox{\rm where}
\\
u &=&
\sum_{j=1}^{2m} \sigma(j)\sigma(j')D_{j'}(f)\sigma(j)D_j(g) \nonumber \\
&=& \sum_{j=1}^{2m} \sigma(j')D_{j'}(f)D_j(g) = D_H(f)(g). \nonumber
\nonumber\end{eqnarray}

\noindent As special instances of that relation we have,

\begin{eqnarray}\label{eq:2.26}&&{[D_H(x^{(a)}), D_H(x_{j'}^{(1)})]}=\sigma(j)D_H(x^{(a-\epsilon_j)})\\
&&{[D_H(x^{(a)}),D_H(x^{(b)})]}=\sum_{j=1}^{2m} \sigma(j)
{{a+b-\epsilon_j-\epsilon_{j'}}\choose{a-\epsilon_j}}
D_H(x^{(a+b-\epsilon_j-\epsilon_{j'})}). \nonumber \end{eqnarray}

    From the first equation in \eqref{eq:2.26}, we can see that
$D_H(x^{(c)}) \in H(2m;\un n)^{(1)}$ for all $c \neq (p^{n_1}-1,
\dots, p^{n_{2m}}-1)$. The coefficient of
$D_H(x^{(a+b-\epsilon_j-\epsilon_{j'})})$ in the second equation for $1 \leq j \leq m$  is

\begin{equation*} {{a+b-\epsilon_j-\epsilon_{j'}}\choose{a-\epsilon_j}}
-
{{a+b-\epsilon_j-\epsilon_{j'}}\choose{a-\epsilon_{j'}}}.\end{equation*}

\m

If
\begin{equation}\label{taudef} \tau = \tau(\un{n}): = (p^{n_1}-1, \dots, p^{n_{2m}}-1), \end{equation}
then $\displaystyle{\tau \choose c} = (-1)^{\mid c\mid }$ for all
$2m$-tuples $c$ such that $0 \leq  c \leq \tau$.  Now if $\tau =
a+b-\epsilon_j-\epsilon_{j'}$ for some $2m$-tuples $a,b$, then $0 \leq a-\epsilon_j,\, a -
\epsilon_{j'} \leq \tau$ and $\displaystyle {\tau \choose
a-\epsilon_j} - {\tau \choose a-\epsilon_{j'}} = (-1)^{\mid
a-\epsilon_j\mid} - (-1)^{\mid a-\epsilon_{j'}\mid} = 0$.  Thus,

\begin{equation*}[D_H(x^{(a)}),D_H(x^{(b)})] \in \spa_{\F}\{D_H(x^{(c)}) \mid
0 \leq  c < \tau\}. \end{equation*}

\m  As a consequence, $$\tilde H(2m;\un n)^{(1)} = \tilde H(2m;\un
n)^{(2)} = \dots = \ \spa_{\F}\{D_H(x^{(c)}) \mid 0 \leq c
< \tau\}.\nonumber$$

\m Because \ $\tilde H(2m;\un n)^{(1)} \supseteq H(2m;\un n)^{(2)}
\supseteq \tilde H(2m;\un n)^{(2)} = \tilde H(2m;\un n)^{(1)}$,  \
equality must hold, and $H(2m;\un n)^{(2)} =
\spa_{\F}\{D_H(x^{(c)}) \mid 0 \leq c < \tau\}$.

\m We may define a version of the ordinary Poisson bracket on
$\mathcal O(2m;\un n)$ by

\begin{equation}\label{eq:2.27} \{f,g\} = \sum_{i=1}^{m}
D_{i}(f)D_{i'}(g) - D_{i'}(f)D_{i}(g).
\end{equation}

Note that

\begin{equation}\label{eq:2.28} \{f,g\} = \sum_{j=1}^{2m}
\sigma(j)D_{j}(f)D_{j'}(g) = D_H(f)(g).  \end{equation}

\noindent {F}rom this and Theorem \ref{Thm:2.24} (iv) we see that
$\bar {\mathcal O}(2m;\un n) := \mathcal O(2m;\un n)/\F 1$ is a Lie
algebra under the Poisson bracket, and its derived algebra $\bar
{\mathcal O}(2m;\un n)^{(1)}$ is isomorphic to $H(2m;\un
n)^{(2)}$.

\m The Lie algebra $H(2;\un{n})^{(2)}$ is exactly
the Lie algebra $S(2;\un{n})^{(1)}$.  It is for that reason we
have excluded the $m=2$ case from Theorem \ref{Thm:2.19}.

\m The degree derivation  $\mathfrak D_1 = \sum_{j=1}^{2m} x_j^{(1)}D_{j}$  may be
adjoined to the Lie algebra $H(2m;\un n)$ or $H(2m;\un n)^{(2)}$
to produce a larger Lie algebra.  In the first case, we will write
$CH(2m;\un n)$;  while in the second we will write $H(2m;\un
n)^{(2)} \oplus \F\mathfrak D_1$ to avoid possible
confusion with $\bigl (CH(2m;\un n)\bigr )^{(2)}$, which is just the algebra
$H(2m;\un n)^{(1)}$ itself.

\m The derivations

\begin{equation} \label{eq:2.29} Q_j :=
D_H\left(x_j^{(p^{n_j})}\right) =
\sigma(j)x_j^{(p^{n_j}-1)}D_{j'},   \qquad 1 \leq j \leq 2m,
\end{equation}

\noindent belong to $H(2m;\un n)$, and expressions for their
products may be derived using (iv) of Theorem \ref{Thm:2.24}:

\begin{eqnarray} [Q_i,Q_j] &=& \delta_{i',j}\sigma(i)
D_H\left(x_i^{(p^{n_i}-1)}x_{i'}^{(p^{n_{i'}}-1)}\right),
\label{eq:2.30} \\
{[Q_i,D_H(f)]} &=& D_H(g), \quad \hbox{\rm where} \quad g =
\sigma(i)x_i^{(p^{n_i}-1)}D_{i'}(f).  \label{eq:2.31}
\end{eqnarray}

We can now state \bi 

\begin{Pro} \label{Pro:2.32} {\rm ({\cite[Chap.~1, Sec.~6,  Prop.~1]{KS}})} \
The Lie algebra $H(2m;\un n)$ is spanned over $\F$ by the
derivations $D_H(f), f \in \mathcal O(2m;\un n)$, as in
\eqref{eq:2.23}, along with the derivations $Q_j, 1 \leq j \leq
2m$, as in \eqref{eq:2.29}.  Their products are given by
\eqref{eq:2.25}, \eqref{eq:2.30}, and \eqref{eq:2.31}.
\end{Pro}
 \m

The simple Lie algebra $H(2m;\un n)^{(2)}$ of Hamiltonian type
inherits a grading from $W(2m;\un n)$:

\begin{eqnarray} H(2m;\un n)^{(2)} &=& \bigoplus_{j=-1}^r \big
(H(2m;\un n)^{(2)}\big)_j \qquad \hbox{\rm where}  \nonumber  \\
\big (H(2m;\un n)^{(2)}\big)_j &=& H(2m;\un n)^{(2)} \cap W(2m;\un
n)_j
= D_H(\mathcal O(2m;\un n)_{j+2}), \quad \hbox{\rm and} \nonumber \\
r&=& p^{n_1} + \cdots + p^{n_{2m}}- {2m}-3. \end{eqnarray}

\noindent Thus, the null component $\big (H(2m;\un
n)^{(2)}\big)_0$ is spanned by the elements

\begin{equation}\label{eq:H0basis} D_H(x_i^{(1)}x_j^{(1)}) = \sigma(j)x_i^{(1)}D_{j'} +
\sigma(i)x_j^{(1)}D_{i'}.\end{equation}

\noindent Under the representation $\varphi: W(2m)_0
\rightarrow \mathfrak{gl}(\mathcal O(2m)_1)$, $\varphi(D)(f) =
D(f)$, we see that the matrix of $\varphi\bigl(D_H(x_i^{(1)}x_j^{(1)})\bigr)$
relative to the basis $\{x_1^{(1)},\dots,x_{2m}^{(1)}\}$ is
$\sigma(j)E_{i,j'}+\sigma(i)E_{j,i'}$.  As those matrices span the
symplectic Lie algebra $\spm$,   we have  $\big
(H(2m;\un n)^{(2)}\big)_0 \cong \spm$.    Moreover,
$\big (H(2m;\un n)^{(2)}\big)_{-1}$ is the dual $V^*$ of the
natural $2m$-dimensional module $V$ for $\spm$, but
$V^* \cong V$ as $\spm$-modules.   It is well-known (and can be
seen from  \eqref{eq:H0basis} and  \eqref{eq:Hprod})  that

\begin{equation}\label{eq:spiso}   \spm \cong S^2(V)  \end{equation}
as $\spm$-modules.

\m Finally, we note that the Lie algebra $CH(2m;\un n)$ also
inherits a grading from $W(2m;\un n)$.  The null component
$CH(2m;\un n)_0$ is a one-dimensional central extension of $\big
(H(2m;\un n)^{(2)}\big)_0$ by the degree derivation; thus it is a
one-dimensional central extension of $\spm$, which
we denote by $\csp$.

 \m
  \section  {Contact
Lie algebras of Cartan type \ (the $\boldsymbol{K}$  series) \label{sec:2.11} }

 \m
 In defining the final series of Cartan type Lie algebras,
 we suppose $\sigma(j)$ and $j'$ are as in \eqref{eq:2.21} for
 $j=1,\dots,2m$.  Our approach will follow {\cite[Sec.~4.5]{SF}}
 (compare also \cite[Sec.~7.1]{KS}).  For $f \in \mathcal O(2m+1)$, let

 $$D_K(f):= \sum_{i=1}^{2m+1} f_i D_{i}, \qquad \hbox{where}\nonumber$$

\begin{eqnarray} \label{eq:2.34} f_i &=& x_i^{(1)}D_{{2m+1}}(f) +
\sigma(i')D_{i'}(f),
\ \ \ 1 \leq i \leq 2m, \\
{f_{2m+1}} &=& \Delta (f) \ := \ 2f - \sum_{j=1}^{2m} x_j^{(1)}D_j(f).
\label{eq:2.35} \end{eqnarray}
In particular,
\begin{eqnarray}
D_K(1) &=& 2 D_{2m+1} \label{eq:2.34a}\\
D_K(x_j^{(1)}) &=&  \sigma(j) D_{j'} + x_j^{(1)} D_{2m+1}  \qquad 1 \leq j \leq 2m \label{eq:2.34b}\\
D_K(x_{2m+1}^{(1)}) &=& \sum_{i=1}^{2m} x_i^{(1)}D_i + 2x_{2m+1}^{(1)}D_{2m+1}.  \label{eq:2.34c}
\end{eqnarray}

It follows (compare Proposition \ref{Pro:2.40}  below) that for all $f, g \in \mathcal O(2m+1)$,

\begin{eqnarray} \label{eq:2.36} &&{[D_K(f), D_K(g)]} =  D_K(u) \qquad
\hbox{\rm where} \\
&&\quad \quad u = \Delta(f) D_{2m+1}(g) - \Delta(g)D_{2m+1}(f)  + \{f,g\}
\ \  \hbox{and}
\nonumber \\
&&\quad \quad \{f,g\} = \sum_{j=1}^{2m} \sigma(j) D_j(f)D_{j'}(g) \quad
\hbox{\rm as in \eqref{eq:2.28}.}   \nonumber \end{eqnarray}

Thus, the elements $D_K(f), \ f \in \mathcal O(2m+1)$, form  a Lie
subalgebra $K(2m+1)$ of $W(2m+1)$.  They are precisely the
derivations $D \in W(2m+1)$ satisfying $D(\omega_K) \in \mathcal
O(2m+1)\omega_K$, where $\omega_K$ is the contact differential
form

\begin{equation*} \omega_K = d{x_{2m+1}} + \sum_{j=1}^{2m} \sigma(j) x_{j}\,d{x_{j'}}, \end{equation*}
\noindent and  $x_j = x_j^{(1)}$ for all $j$.

\begin{Rem}\label{Rem:2.37} \
Note that the form used in {\rm {\cite{Wi1}}} and in
{\rm \cite[Sec.~7.1]{KS}}  is $\omega_K' = d{x_{2m+1}} + \sum_{j=1}^{2m}
\sigma(j') x_{j}\,d{x_{j'}}$, and so the formulas in those papers
will be slightly different from the ones displayed in this
work.\end{Rem}

\m

For each $(2m+1)$-tuple $\un n$, the Lie algebra $K(2m+1;\un n)$
is by definition the intersection,  $K(2m+1;\un n) = K(2m+1) \cap
W(2m+1;\un n)$.

\m The product in \eqref{eq:2.36} can be expressed using the
modified Poisson bracket

\begin{equation}\label{eq:2.38} \langle f,g \rangle : =
\Delta(f)D_{2m+1}(g)  -  \Delta(g)D_{2m+1}(f) + \{f,g\}.
\end{equation}

\n Thus,

\begin{equation} \label{eq:2.39} [D_K(f), D_K(g)] = D_K(\langle f,g
\rangle).  \end{equation}

\m \begin{Pro} \label{Pro:2.40} \ {\rm (Compare {\cite
[Sec.~4.5, Prop.~5.3]{SF}}.)} Suppose for
$a=(a_1,\dots,a_{2m+1})$, $a_i \in \mathbb N$, that  $\parallel a
\parallel: = \mid a\mid + \,a_{2m+1} -2$.  Then
\begin{itemize}
\item[{\rm (i)}] $\langle x^{(a)}, x^{(b)} \rangle = \{x^{(a)},
x^{(b)}\}$

\hspace{-.21truein}$\displaystyle{+ \Bigg(\parallel b \parallel
{{a+b-\epsilon_{2m+1}}\choose b} -
\parallel a \parallel {{a+b-\epsilon_{2m+1}}\choose a} \Bigg)
x^{(a+b-\epsilon_{2m+1})};}$ \item[{\rm (ii)}] $\langle
1,x^{(a)}\rangle = 2 x^{(a-\epsilon_{2m+1})}$;  \item[{\rm (iii)}]
$\langle x_j^{(1)}, x^{(a)} \rangle = \sigma(j)x^{(a-\epsilon_{j'})} +
(a_{j}+1)x^{(a+\epsilon_{j}-\epsilon_{2m+1})}, \quad 1 \leq j \leq
2m$;  \item[{\rm (iv)}] $\langle x_{2m+1}^{(1)}, x^{(a)} \rangle = \,
\parallel a
\parallel x^{(a)}$;
\item[{\rm(v)}]  $\langle x_i^{(1)}x_j^{(1)}, x^{(a)}\rangle = \sigma(i)a_j x^{(a+\epsilon_j-\epsilon_{i'})}
+\sigma(j) a_i x^{(a + \epsilon_i -\epsilon_{j'})}$,

\hspace {3 truein} $1 \leq i,j \leq 2m$;
 \item[{\rm (vi)}]  $\langle x_i^{(1)}x_{i'}^{(1)}, x^{(a)} \rangle = (a_{i'}-a_i) x^{(a)}$,  \quad $1 \leq i \leq m$;
 \item[{\rm (vii)}] The spaces $K(2m+1;\un n)_j =
\spa_\F \{ D_K(x^{(a)}) \mid \, \parallel a\parallel = j\}$
give a grading of the Lie algebra $K(2m+1;\un n)$ such that
$$K(2m+1;\un n) =
\bigoplus_{j=-2}^{r} K(2m+1;\un n)_j,$$ where $r = p^{n_1} +\cdots +
p^{n_{2m}}+2p^{n_{2m+1}}-2m-3$.
\end{itemize}
\end{Pro}

\medskip
\begin{Thm} \label{Thm:2.41} \ {\rm (See {\cite[Sec.~4.5, Thm.~5.5]{SF}}.)}
\ $K(2m+1;\un n)^{(1)}$ is a simple Lie algebra, and

$$K(2m+1;\un n)^{(1)} = \begin{cases} K(2m+1;\un n) & \ \  \hbox{\rm if} \  2m+4 \not
\equiv 0 \mod p \\
\spa_{\F} \{ D_K(x^{(a)})  \mid \, a \not = \tau(\un
{n}) \} & \ \  \hbox{\rm if} \ 2m+4 \equiv 0 \mod p, \end{cases}\nonumber$$

\n where $\tau(\un{n}) = (p^{n_1}-1, \dots, p^{n_{2m+1}}-1)$. Thus, $\dim K(2m+1;\un n)^{(1)}$ has dimension
$p^{n_1+\cdots+n_{2m+1}}$ in the first case and dimension
$p^{n_1+\cdots+n_{2m+1}}-1$ in the second.  Moreover, $K(2m+1;\un
n)^{(1)}$ is isomorphic to $\mathcal O(2m+1;\un n)^{(1)}$, where
$\mathcal O(2m+1;\un n)$ is viewed as a Lie algebra under the
product $f \times g \mapsto \langle f, g \rangle$.
\end{Thm}

\medskip

{F}rom parts (ii)-(vi) of Proposition \ref{Pro:2.40}, we see that
$K(2m+1;\un n)_{-2}$ is spanned by $d_{2m+1}: = D_K(1) =
2D_{2m+1}$; \ and $K(2m+1;\un n)_{-1}$ is spanned by the elements
$d_j: = \sigma(j')D_K(x_{j'}^{(1)}) = D_{j}+
\sigma(j')x_{j'}^{(1)}D_{2m+1}$ for $j=1,\dots, 2m$.  These
elements satisfy the rule

\begin{equation}\label{eq:2.42} [d_i,d_j] = \delta_{i,j'}\sigma(i)d_{2m+1} \qquad
1 \leq i,j \leq  2m. \end{equation}

The elements $D_K(x^{(a)})$ with $\parallel a\parallel = 0$ span
the subalgebra  $K(2m+1;\un n)_0$, which is isomorphic to the Lie
algebra $\csp = \spm  \oplus \F I$.
The space $K(2m+1;\un n)_{-1}$ is its natural $2m$-dimensional
module. \bigskip

\begin{Rem} \label{rem:natural}   {\rm For any of the simple Lie algebras $\g = \bigoplus_i \g_i$ of Cartan type discussed in Sections \ref{sec:2.8}-\ref{sec:2.11}, the sum
$\mathfrak m_0: = \bigoplus_{i \geq 0}  \g_i$  is a maximal subalgebra of $\g$
which is  invariant under
the automorphisms of $\g$.      This uniqueness property of $\mathfrak m_0$ is justification for calling the grading we have described in those
sections  the {\em natural grading} of $\g$.    Moreover,  the spaces 
$m_\ell: = \bigoplus_{i \geq \ell} \g_i$ for $\ell \geq -1$   afford a filtration
$\g \supseteq \mathfrak m_{-1} \supset  \mathfrak m_0  \supset  \mathfrak m_1 \supset  \cdots$
of $\g$, which is often referred to as the {\it natural filtration} (see for example, \cite[Defn.~4.2.8]{St}.) }
 \end{Rem}

\medskip

In closing this subsection, we comment  that very few of the simple
Cartan type Lie algebras described above are restricted.  In fact,
the restricted ones may be characterized as having their defining
$m$-tuple $\un n = \un 1$, the tuple of all ones.

\bi

\begin {Pro} \label{Pro:2.43} \  {\rm{(See \cite[Thm.~2]{K3} and \cite[(7.2)]{St}.})}
\ Let $\g$ be a restricted
simple Lie algebra of Cartan type.  Then $\g$ isomorphic to one of
the following: \ \ $W(m;\un {1})$ $ (m \geq 1)$, $S(m;\un
{1})^{(1)}$ $ (m \geq 3)$, $H(2m;\un {1})^{(2)}$ $ (m \geq 1)$, or
$K(2m+1;\un {1})^{(1)}$ $ (m \geq 1)$.  \end{Pro}

\m
\section  {The Recognition Theorem with stronger hypotheses  \label{sec:2.12}}

By imposing strong assumptions on the non-positive homogeneous components, we can deduce
the following version of the Recognition Theorem.     The next two chapters will be devoted
to showing that these additional hypotheses must hold when the conditions of the Main
Theorem are fulfilled.     A similar result,  phrased in the language of
filtered Lie algebras,  appears in \cite[Cor.~5.5.3 (Weak Recognition Theorem)]{St}
for {\em simple} Lie algebras.  

\bi
 \begin{Thm} \label{Thm:2.81}  \ Let
$\g = \bigoplus_{j=-2}^r  \g_j$  be a finite-dimensional graded
Lie algebra over an algebraically closed field ${\mathbb F}$ of
characteristic $p > 3$.   Assume that: \smallskip

\begin{itemize}
\item[{\rm (a)}]  $\g_0$ is isomorphic to $\mathfrak{gl}_m,\,
\mathfrak{sl}_m, \,\spm$, or $\csp =
\spm \oplus {\mathbb F} I$. \smallskip \item[{\rm
(b)}]  $\g_{-1}$ is a standard $\g_0$-module (of dimension $m$ or
$2m$,  depending on the simple classical Lie algebra $\g_0^{(1)} =
[\g_0,\g_0]).$ \smallskip \item[{\rm (c)}] If $\g_{-2} \neq 0,$
then it is one-dimensional and equals $[\g_{-1}, \g_{-1}],$ and
$\g_0$ is isomorphic to $\csp.$
\smallskip \item[{\rm (d)}]\   If  $x \in \bigoplus_{j \geq 0} \g_j$ and $[x,\g_{-1}] = 0$,  then
$x = 0$ (transitivity);

\item[{\rm (e)}]\  If  $x \in \bigoplus_{j \geq 0} \g_{-j}$  and $[x,\g_{1}] = 0$,  then
$x = 0$ (1-transitivity).
  \end{itemize}
 \smallskip \noindent Then either $\g$ is a
Cartan type Lie algebra with the natural grading:

\begin{equation*} X(\mathfrak m;{\un n})^{(2)} \subseteq \g \subseteq X(\mathfrak m;{\un n}) \end{equation*}

\noindent where
\begin{equation*} X = \begin{cases} W, S, \ \hbox{\rm or} \ CS & \ \hbox{\rm and } \ \ \mathfrak m = m, \\
H  \ \hbox{\rm or} \ CH
& \ \hbox{\rm and } \ \ \mathfrak m  = 2m, \\
K & \ \hbox{\rm and } \ \ \mathfrak m= 2m+1; \end{cases}\end{equation*}
or $\g$ is a classical simple Lie algebra:
\begin{itemize}
\item[{\rm a)}]  $ \g \cong \mathfrak{sl}_{m+1}  \hookrightarrow
W(m;\un n),$ or
\item[{\rm b)}] $ \g \cong \mathfrak{sp}_{2(m+1)}  \hookrightarrow
K(2m+1;{\un n}).$
\end{itemize}

\end{Thm}

\pf \  By \cite[Prop.~2.7.3]{St},  the Lie algebra $\g$ is isomorphic to
a subalgebra of a \W \ Lie algebra of Cartan type.  Then the theorem
follows from \cite[Lem.~5.2.3]{St}.  \qed

\m

\section {$\boldsymbol{\g_\ell}$ as a $\boldsymbol{\g_0}$-module for 
Lie algebras $\boldsymbol{\g}$ of Cartan type  \label{sec:2.13}} \m

Next we will  investigate the structure of $\g_\ell$ ($\ell < p-1$)   as a
$\g_0$-module when $\g$ is a depth-one  Lie algebra of Cartan
type. Here we will  follow Section 10 of Chapter I of \cite{KS}.
See also {\cite[Sec.~5.2]{St}} for related results. For this
purpose, it is convenient to introduce the following derivations
in $W(m)$:

\begin{eqnarray}\label{eqn:Df} \mathfrak D_f &:=& f \mathfrak D_1 \qquad  \hbox{\rm for} \ \ f
\in \mathcal O(m),  \end{eqnarray}

\noindent  where $\mathfrak D_1 = \sum_{j=1}^m x_j^{(1)}D_j$ is
the degree derivation.  Then  $\deg(\mathfrak D_f)
= \deg(f)$ for all homogeneous $f \in \mathcal O(m)$,
and these derivations  multiply according to the rule,
\begin{eqnarray}\label{eq:degdermult}&& \\
 {[\mathfrak D_f, \mathfrak D_g]} &=&
[f \mathfrak D_1, g \mathfrak D_1]\nonumber \\ &=&
\Big(f \mathfrak D_1(g) - g \mathfrak D_1(f)\Big)\mathfrak D_1\nonumber  \\
&=& \Big(\deg(g)-\deg(f)\Big) fg \mathfrak D_1 = \Big(\deg(g)-\deg(f)\Big) \mathfrak D_{fg}. \nonumber
\end{eqnarray}
In deriving this expression,  we have applied the general relation
$$[f D, gE] = f D(g) E - g E(f)D + fg [D,E]$$ in \eqref{eq:derive}.
Recall {f}rom \eqref{eq:derdeg} that  $\ad\mathfrak D_1$ acts as multiplication by the  scalar
$\ell$ on $W(m)_\ell$ for each $\ell$.  In particular, when $\g= W(m)$ and $m = 1$,  then $\g_0 =
\F \mathfrak D_1$,  and for each $\ell \geq -1$, the space
$\g_\ell = \F x_1^{(\ell+1)}D_1$ is an irreducible $\g_0$-module on which $\mathfrak D_1$
acts as multiplication by the scalar $\ell$.
The $\g_0$-module structure of $\g_\ell$  for  $m \geq 2$  is the topic of
the next result.

\bi
\begin{Thm} \label{Thm:2.82} \  Assume $\g = W(m)$ for $m \geq 2$,
and let $\g_\ell = W(m)_\ell$ for  all
$\ell \geq -1$.  Let  $\mathfrak h$ be the Cartan subalgebra
of $\g_0$ with basis $x_j^{(1)}D_j$ ($1 \leq j \leq m$), and
let $\varepsilon_i$ ($1 \leq i \leq m$) be the dual basis in $\mathfrak h^*$ so
that $\varepsilon_i(x_j^{(1)}D_j) = \delta_{i,j}$.   Let  $\mathfrak b^+ = \mathfrak h \oplus
\mathfrak n^+$, where $\mathfrak n^+$ is  the $\F$-span of all $x_i^{(1)}D_j$
with $1\le i<j\le m$,  and set  

\begin{eqnarray} \g_{\ell}^\dagger &:=&\{ D \in \g_\ell \mid \di (D) = 0\}  \label{eq:2.82a}\\
\g_{\ell}^\sharp &:=&  \spa_{\mathbb F}\Big\{ \mathfrak D_f = f
\sum_{j=1}^m x_j^{(1)}D_j = f \mathfrak D_1 \, \Big | \, \deg(f) = \ell \Big\}. \label{eq:2.82b}
\end{eqnarray}  Then for $\ell \leq  p-2$ we have:

\begin{itemize}
\item[{\rm(i)}]  $ \g_\ell^\dagger$ is an irreducible
$\g_0$-submodule of $\g_\ell$  with a $\mathfrak b^+$-primitive vector of 
weight 
$(\ell+1)\varepsilon_1-\varepsilon_m$ when $m+\ell \not \equiv 0
\mod p$;  \smallskip

\item[{\rm(ii)}]  $\g_\ell^\sharp$ is an irreducible
$\g_0$-submodule of $\g_\ell$ with a $\mathfrak b^+$-primitive vector of weight $\ell
\varepsilon_1$, and $\g_\ell^\sharp$ is isomorphic as a
$\g_0$-module to $S^\ell(V)$, where $V = \mathcal O(m)_1$ is the
natural $m$-dimensional module of $\g_0$; \smallskip

\item[{\rm(iii)}] $\g_\ell = \g_\ell^\sharp \oplus
\g_\ell^\dagger$ when   $m+\ell \not \equiv 0  \mod p$; \smallskip

\item[{\rm(iv)}] When $m \geq 3$ and $m+\ell \equiv 0 \mod p$, then $\g_\ell \supset
\g_\ell^\dagger \supset \g_\ell^\sharp \supset 0$ is the unique
composition series of the $\g_0$-module $\g_\ell$; \smallskip

\item[{\rm(v)}] When $m=2$ and $m+\ell \equiv 0 \mod p$ (i.e. $\ell = p-2$), then 
$\g_{p-2}^\dagger$   is spanned modulo $ \g_{p-2}^\sharp$ by $x_1^{(p-1)}D_2$
and $x_2^{(p-1)}D_1$.  Moreover,  $\g_{p-2}^\dagger /\g_{p-2}^\sharp$  $\cong L(0) \oplus L(0)$ and
$\g_{p-2}/\g_{p-2}^\dagger \cong \g_{p-2}^\sharp \cong  L((p-2)\varpi_1) \cong S^{p-2}(V)$
as modules  for $\g_0^{(1)} \cong \mathfrak{sl}_2$.  
\end{itemize}   \end{Thm}

\pf   By \eqref{eq:divprod}  we have
\begin{eqnarray}\label{eq:divhom}
\di\left(\left[x_i^{(1)}D_j, D\right]
\right)&=& x_i^{(1)}D_j\big( \di(D)\big) - D\bigl(\di(x_i^{(1)}D_j)\bigr) \\
&=& x_i^{(1)}D_j\big( \di(D)\big) -  D(\delta_{i,j}1)\nonumber  \\
&=& x_i^{(1)}D_j\big( \di(D)\big), \nonumber \end{eqnarray}
which implies that $\g_\ell^\dagger$ is a $\g_0$-submodule of $\g_\ell$ for all $\ell$.
This is also evident from the fact that $\g_0 = \g_0^{(1)} \oplus \F \mathfrak D_1$,
where $\g_0^{(1)} = S(m)_0$,  and $\g_\ell^\dagger = S(m)_\ell$.

\m Applying \eqref{eq:derive}, we obtain \vskip -.2 truein
\begin{eqnarray}\label{eq:sharp} && \\ \left [x_i^{(1)}D_j, \,\mathfrak D_f \right] &=
& x_i^{(1)}[D_j, \mathfrak D_f] - \mathfrak D_f\left(x_i^{(1)}\right) D_j   \nonumber  \\
&=& x_i^{(1)}[D_j, f \mathfrak D_1]  - f\mathfrak D_1\left(x_i^{(1)}\right) D_j  \nonumber  \\
&=& x_i^{(1)}D_j(f) \mathfrak D_1 +  x_i^{(1)}f D_j  - x_i^{(1)}f D_j  =
x_i^{(1)}D_j(f) \mathfrak D_1.  \nonumber  \end{eqnarray}

\noindent  Thus, $\g_\ell^\sharp$ is also a $\g_0$-submodule of
$\g_\ell$.  Moreover, if $V = \mathcal O(m)_1$, then the
calculation in \eqref{eq:sharp}  shows that $\g_\ell^\sharp \cong
\mathcal O(m)_\ell = S^\ell(V)$ as $\g_0$-modules via the
identification  $f \mapsto \mathfrak D_f$.   The vector $\mathfrak
{D}_{x_1^{(\ell)}} = x_1^{(\ell)} \mathfrak D_1$ has zero product
with all $x_i^{(1)}D_j$  such that $i< j$, and $[x_i^{(1)}D_i,
\mathfrak {D}_{x_1^{(\ell)}}] = \ell \delta_{i,1}\mathfrak
{D}_{x_1^{(\ell)}}$. Therefore,  $\mathfrak {D}_{x_1^{(\ell)}}$ is
a $\mathfrak b^+$-primitive vector of $\g_\ell^\sharp$ and its weight is
$\ell \varepsilon_1$ relative to the Cartan subalgebra $\mathfrak
h$. Relative to $\g_0^{(1)} \cong \mathfrak {sl}_m$ and its Cartan
subalgebra $\mathfrak h \cap \g_0^{(1)}$,  the module $S^\ell(V)$ has a 
$\mathfrak b^+$-primitive vector $x_1^{(\ell)}$ of weight $\ell
\varpi_1$.    Since $\ell \leq  p-2$,  the irreducible
$\mathfrak{sl}_m$-module $L(\ell \varpi_1)$  is isomorphic to the Weyl module  $V(\ell
\varpi_1)$, which has dimension given by Weyl's dimension formula (see \eqref{eq:Weyldim}) so 
that
$$\dim V(\ell \varpi_1) = {{m+\ell-1} \choose \ell} = \dim S^\ell(V).$$
 It follows that  $\g_\ell^\sharp$ is an irreducible $\g_0^{(1)}$-module (hence an
 irreducible $\g_0$-module)  isomorphic to $S^\ell(V)$
for all $\ell \leq p-2$,  as asserted in (ii) .

Consider the sequence of $\g_0$-module maps,

\begin{eqnarray*}  \g_\ell &\stackrel{\sim}{\rightarrow}& S^{\ell+1}(V)
\otimes V^* \twoheadrightarrow \,S^\ell(V) \\
\quad  f D_k  &\mapsto& f \otimes D_k \ \mapsto \ \di(f D_k) =
D_k(f). \end{eqnarray*}

\noindent The kernel is $\g_\ell^\dagger$.    Since
\begin{equation}\label{eq:2.825} \di \bigl(\mathfrak D_u\bigr) = (m+\ell)u\end{equation}
 for all $u \in \mathcal O(m)_\ell$,  we see that $\g_\ell^\sharp \subset \g_\ell^\dagger$
 if $m+\ell \equiv 0$, and $\g_\ell^\sharp \cap \g_\ell^\dagger = 0$
 if $m+\ell \not \equiv 0$  mod $p$.

\m  Suppose $x^{(a)} D_k \in \g_\ell$ and set

\begin{equation*}  E_j: = (a_j+1)D_{k,j}\bigl(x^{(a+\epsilon_j)}\bigr)
= (a_j+1)x^{(a)}D_k - x^{(a-\epsilon_k)}x_j^{(1)}D_j - \delta_{j,k}x^{(a)}D_j \in \g_\ell^\dagger.
\end{equation*}   Then

\begin{eqnarray*}(m+\ell)x^{(a)}D_{k} &=&
\Big(\textstyle{\sum_{j=1}^m \bigl (a_j+1\bigr)}\Big)x^{(a)}D_{k} - x^{(a)}D_k \\
&=& \sum_{j=1}^m E_j + x^{(a-\epsilon_k)} \sum_{j=1}^m x_j^{(1)}D_j  \\
&=&  \sum_{j=1}^{m}E_j  + x^{(a-\epsilon_k)}  \mathfrak D_1.  \end{eqnarray*}

\noindent  Consequently,  when $m+\ell \not \equiv 0 \mod p$, \ $\g_\ell \subseteq
 \g_\ell^\dagger+ \g_\ell^\sharp$ so that $\g_\ell = \g_\ell^\dagger \oplus \g_\ell^\sharp$
 as asserted in (iii).    We also see that the above map
$x^{(a)}D_k \mapsto D_k(x^{(a)}) = x^{(a-\epsilon_k)}$
 is essentially (up to a factor of $(m+\ell)^{-1}$) the projection of $\g_\ell$ onto the
 $\g_0$-submodule $\g_\ell^\sharp \cong S^\ell(V)$.

Recall that  ${\mathfrak n}^+$ is the $\F$-span of all $x_i^{(1)}D_j$
with $1\le i<j\le m$, a maximal nilpotent subalgebra of
$\g_0^{(1)}$, and  ${\mathfrak b}^+:=\h\oplus{\mathfrak n}^+$.    We
claim that any ${\mathfrak b}^+$-primitive vector of $\g_\ell$
 is a scalar multiple of $x_1^{(\ell+1)}D_m$ or of
${\mathfrak D}_{x_1^{(\ell)}}$.    The claim certainly holds when
$\ell=0$ (see \eqref{eq:sharp}).     Assume that it holds for all $\ell<s$ where $1\le s\le
p-2$, and let $u$ be a ${\mathfrak b}^+$-primitive vector of $\g_s$. Since
$[\g_{-1},u]\ne 0$, by transitivity, $[D_k,u]$ is a ${\mathfrak b}^+$-primitive
vector of $\g_{s-1}$   for some $k\le m$. In view of our assumption, 
this means that either $[D_k,u]= \zeta \, x_1^{(s)}D_m$ or
$[D_k,u]=\zeta \,{\mathfrak D}_{x_1^{(s-1)}}$ where $\zeta \in
\F^\times$. No generality will be lost by assuming that
$\zeta=1$.

Suppose $[D_k,u]={\mathfrak D}_{x_1^{(s-1)}}.$ If $k>1$,  then
$$u\,=\,x_{1}^{(s-1)}x_k^{(1)}\Big(\sum_{i\ne k}
x_i^{(1)}D_i+\frac{1}{2}x_k^{(1)}D_k\Big)+u_0,\ \ \quad
u_0\in\g_{s}\cap\ker\ad D_k.$$      {F}rom this it is immediate that
$$[x_{k-1}^{(1)}D_k,u]\,\equiv\, x_{1}^{(s-1)}x_{k-1}^{(1)}\sum_{i\ne
k} x_i^{(1)}D_i\mod {\mathcal O}(m)D_k.$$    Since $[{\mathfrak
n}^+,u]=0$,  we have reached a contradiction,  showing that $k=1$. Then
$$u\,=\,sx_{1}^{(s+1)}D_1+
x_{1}^{(s)}\sum_{i=2}^m x_i^{(1)}D_i+u'_0,\ \ \quad
u'_0\in\g_{s}\cap\ker\ad D_1,$$ implying
\begin{eqnarray*}
[x_{1}^{(1)}D_2,u]&=& -sx_{1}^{(s+1)}D_2+x_1^{(s)}x_1^{(1)}D_2+
x_1^{(1)}[D_2,u_0']-u_0'(x_1^{(1)})D_2\\
&=&x_{1}^{(s+1)}D_2+x_1^{(1)}\sum_{i=1}^mf_iD_i+gD_2
\end{eqnarray*}
where $g$ and the $f_i$ belong to the subalgebra of ${\mathcal
O}(m)$ generated by $x_j^{(r)}$ with $j\ge 2$ and $r\ge 0$.    Since
$s\ge 1$ it follows that $[x_1^{(1)}D_2,u]\ne 0$, a contradiction.
We conclude that $[D_k,u]=x_1^{(s)}D_m.$

If $k>1$ then  $u=x_1^{(s)}x_k^{(1)}D_m+u_1$ where
$u_1\in\g_{s}\cap\ker\ad D_k$. Then, for $i<k$,
$$0=\,[x_{i}^{(1)}D_k,u]\,=\, x_{1}^{(s)}x_{i}^{(1)}D_m-u_1(x_{i}^{(1)})D_k.$$
This shows that $k=m$ and
$u_1=aD_m+\sum_{i=1}^{m-1}\,x_1^{(s)}x_i^{(1)}D_i$ for some $a\in
{\mathcal O}(m)_{s+1}$ with $D_m(a)=0$. But then
$u=aD_m+{\mathfrak D}_{x_1^{(s)}}$ implying $[{\mathfrak
n}^+,aD_m]=0$. As a consequence, $(x^{(1)}_iD_j)(a)=0$ whenever
$1\le i<j\le m-1$. From this it is easy to deduce that $a=\mu
x_1^{(s+1)}$ for some $\mu\in\F$.

If $k=1$ then  $u=x_1^{(s+1)}D_m+u'_1$ where
$u'_1\in\g_{s}\cap\ker\ad D_1$. Then, for $i>1$,
$$0=\,[x_{1}^{(1)}D_i,u]\,=\, x_{1}^{(1)}[D_i,u'_1]-u'_1(x_{1}^{(1)})D_i.$$
Since $u'_1(x_{1}^{(1)})\in{\mathcal O}(m)$ does not involve
$x_1^{(1)}$, it follows that $[D_i,u'_1]=0$ for all $i$. But then
$u_1'=0$ by transitivity. This proves our claim.

Now assume $m \geq 3$.   Set $G_0\eqdef \hbox{\rm SL}(V)$, a simply connected algebraic
$\F$-group, and identify $\text{Lie}(G_0)$ with $\g_0^{(1)}$. Let
$T$ be the maximal torus of $G_0$ such that $\text{Lie}(T)=\h$,
and let $B^+$ be the Borel subgroup of $G_0$ with
$\text{Lie}(B^+)={\mathfrak b}^+$.       By abuse of notation, we will
identify the weights $\mu$ of $T$ with their differentials
$({\mathrm d}\mu)_e\in \h^*$.       In particular, the fundamental
weights in $X(T)_+$ will be denoted by $\varpi_1,\varpi_2,\ldots, 
\varpi_{m-1}$ (this will cause no confusion).     Note that
$\g_\ell$ is isomorphic to $S^{\ell+1}(V)\otimes V^*$, which is a rational $G_0$-module.
Since $\ell\le p-2$,   it is easy to see that
$X_+\big(S^{\ell+1}(V)\otimes V^*\big)\subset X_1(T)$. So the
$G_0$-module $\g_\ell$ satisfies all the conditions of Proposition
\ref{Pro:2.805}. It follows that every $\g_0^{(1)}$-submodule of
$\g_{\ell}$ is $G_0$-stable.

We claim that the $\g_0^{(1)}$-submodule of $\g_\ell$ generated by
$x_1^{(\ell+1)}D_m$ coincides with $\g_\ell^\dagger$. Let us
denote this submodule by $\g_\ell'$ and suppose for a contradiction
that $\g_\ell'\ne \g_\ell^\dagger$. Since we are assuming $m\ge 3$, Theorem
\ref{Thm:2.19}\,(i) applies yielding $D_{i,j}(f)\not\in\g_\ell'$ for
some $f\in {\mathcal O}(m)_{\ell+2}$. Since $\g_\ell'$ is stable
under the normalizer $N_{G_0}(T)$ of $T$,  which acts doubly
transitively on the set of lines $\{\F x_i^{(1)} \mid 1\le i\le
m\}$, we may assume that $(i,j)=(m,1)$ and $f$ is a weight vector
for $T$. Then $f=\xi \,x_1^{(a_1)}x_2^{(a_2)}\cdots x_m^{(a_m)}$
for some $a_i\in{\mathbb N}$ with $\sum_{i=1}^m\,a_i=\ell+2$ and
$\xi\in\F^\times$. Because $\ell+1<p$, the
$\mathfrak{sl}_{m-1}$-module generated by $x_1^{(\ell+1)}D_m$
coincides with ${\mathcal O}(m-1)_{\ell+1}D_m$. Consequently,
$a_m\ge 1$. From (\ref{eq:dijdef}) it follows that $[x_i^{(1)}D_j,
D_{m,1}(g)]=D_{m,1}\big((x_i^{(1)}D_j)(g)\big)$ whenever $2\le
i,j\le m-1.$ As $\ell+2-a_1-a_m<p$, the
$\mathfrak{sl}_{m-2}$-module ${\mathcal O}(m-2)_{\ell+2-a_1-a_m}$
is irreducible. Thus,  no generality will be lost in assuming that
$a_i=0$ for $i\not\in\{1,2,m\}$. As $[x_2^{(1)}D_1,
D_{m,1}(g)]=D_{m,1}\big((x_2^{(1)}D_1)(g)\big)$ and $a_1+a_2<p$,
we may assume further that $a_2=0$. Then
$f=x_1^{(a_1)}x_m^{(a_m)}$ where $a_1+a_m=\ell+2$.

Let $t$ be the minimal integer with the property that
$D_{m,1}\big(x_1^{(t)}x_m^{(\ell+2-t)}\big)\not\in\g_\ell'$. If
$t>0$,  then $\ell+2-t\not\equiv 0\mod p$. In this case Lemma
\ref{Lem:2.18}(c) implies
$$\big[D_{1,m}\big(x_m^{(2)}),D_{m,1}(x_1^{(t+1)}x_m^{(\ell+1-t)}\big)\big]\,=\,
(\ell+2-t)D_{m,1}\big(x_1^{(t)}x_m^{(\ell+2-t)}\big)\in
\g_\ell'\setminus\{0\},$$ a contradiction.      Hence $t=0$;  that is,
$x_m^{(\ell+1)}D_1=-D_{m,1}\big(x_m^{(\ell+2)}\big)\not\in\g_\ell'$.
On the other hand, $\g_\ell'$ is $G_0$-stable and there is an
element in $N_{G_0}(T)$ which permutes the lines spanned by
$x_1^{(\ell)}D_m$ and $x_m^{(\ell)}D_1$. But then
$x_1^{(\ell+1)}D_m\not\in\g_\ell'$,   which is false. By
contradiction, the claim follows.    In conjunction with our
description of the $\mathfrak b^+$-primitive vectors in $\g_\ell$,  this yields that for
$\ell+m\not\equiv \,0 \mod p$ the $\g_0^{(1)}$-module
$\g_\ell^\dagger$ is irreducible, proving (ii).

Now suppose $\ell+m\equiv \,0 \mod p$.    Then our earlier remarks
show that $\g_\ell^\sharp$ coincides with the $\g_0^{(1)}$-socle
of $\g_\ell^\dagger$. We denote by $M$ the maximal
$\g_0^{(1)}$-submodule of $\g_\ell^\dagger$. Since all
$\g_0^{(1)}$-submodules of $\g_\ell$ are $G_0$-stable, so is $M$.
Let $\mu$ be a maximal $T$-weight in $X(M)$ and let $v$ be a
weight vector of weight $\mu$ in $M$. Then $[{\mathfrak
n}^+,v]=0$.      Our description of the $\mathfrak b^+$-primitive vectors in $\g_\ell$
shows that a nonzero multiple of $v$ lies in
$\{x_1^{(\ell+1)}D_m,\,{\mathfrak D}_{x_1^{(\ell)}}\}$.  Since
$x_1^{(\ell+1)}D_m$ generates the $\g_0^{(1)}$-module
$\g_\ell^\dagger$, it must be that $\mu=\ell\varpi_1$, the weight
of  $\mathfrak D_{x_1^{(\ell)}}$.   Hence,
$\nu\le \ell\varpi_1$ for all $\nu\in X(M)$. It also follows that
$-\ell\varpi_{m-1}$ is the only minimal weight in $X(M)$ and $\dim
M^{\ell\varpi_1}=1$.

Next we look at the dual $G_0$-module $M^*$. Since $X(M^*)=-X(M)$,
the preceding remark yields that $\ell\varpi_{m-1}$ is the only
maximal weight of $M^*$ and $\dim (M^*)^{\ell\varpi_{m-1}}=1$. Let
$N\eqdef\{\varphi\in M^*\mid \varphi(u)=0\ \text{for all}\
u\in\g_\ell^\sharp\}$. Clearly, $N$ is a $G_0$-submodule of $M^*$
and $M^*/N\cong (\g_\ell^\sharp)^*$. Recall that
$\g_\ell^\sharp\cong V(\ell\varpi_1)\cong L(\ell\varpi_1)$.    Then
$M^*/N\cong L(\ell\varpi_{m-1})\cong V(\ell\varpi_{m-1})$; see Proposition
\ref{Pro:2.802}. Since $\dim (M^*)^{\ell\varpi_{m-1}}=1$,   it must be that
$\ell\varpi_{m-1}\not\in X(N)$. Let $\psi\in
(M^*)^{\ell\varpi_{m-1}}\setminus\{0\}$ and denote by $N'$ the
$G_0$-submodule generated by $\psi$.      By Proposition
\ref{Pro:2.804}(a), the $G_0$-module $N'$ is a homomorphic image
of the Weyl module $V(\ell\varpi_{m-1})$.     As the latter is
irreducible, so is $N'$.           On the other hand, $N'\not\subset N$
(because $\psi\not\in N$). As a consequence, $M^*=N\oplus N'$ and
$M\cong N^*\oplus N'^*$. Since the subspace of  ${\mathfrak b}^+$-primitive vectors
of $M$ is one-dimensional, the $G_0$-module $M$ is indecomposable.
As $N'\ne 0$, this forces $N=0$. But then $M=\g_\ell^\sharp$, and
our proof is complete in the $m\geq 3$ case.

When $m=2$ and $\ell = p-2$, set

\begin{equation*} v_k: = x_1^{(p-1-k)}x_2^{(k)}D_2 - x_1^{(p-k)}x_2^{(k-1)}D_1\in \g_{p-2}^\dagger \end{equation*}
for $k=0,1,\dots,p$,  where by convention $x_j^{(-1)} = 0$ for $j=1,2$.   The elements  $e= x_1^{(1)}D_2$,
$f = x_2^{(1)}D_1$, and $h = x_1^{(1)}D_1 - x_2^{(1)}D_2$ form a canonical
basis of $\g_0^{(1)} \cong \mathfrak {sl}_2$, and we have

\begin{eqnarray}\label{eq:sl2action}  {[f, v_k]}  &=& (k+1) v_{k+1}, \nonumber  \\
{[e,v_k]} & = & (p+1-k)v_{k-1},  \\
{[h, v_k]} &=& (p-2k) v_k, \nonumber \end{eqnarray}
where $v_{p+1} = 0 = v_{-1}$.    The vectors
$v_k = k^{-1} x_1^{(p-1-k)} x_2^{(k-1)} \mathfrak D_1$ for $k=1,\dots,p-1$,
determine a basis of $\g_{p-2}^\sharp$, and we can see from the expressions in
\eqref{eq:sl2action} that $\g_{p-2}^\sharp \cong L((p-2)\varpi_1)$ as a module
for $\g_0^{(1)}$.     As $[f,v_0] = v_1$ and $[e,v_p] = v_{p-1}$, 
and the other basis elements of $\g_0^{(1)}$ act trivially on $v_0$ and $v_p$, 
we have that  $v_0$ and $v_p$ span trivial $\g_0^{(1)}$-modules of $\g_{p-2}^\dagger$
modulo $\g_{p-2}^\sharp$.     Moreover, because  $[e,x_1^{(p-1)}D_1] = -v_0$
and $[h,x_1^{(p-1)}D_1] = (p-2)x_1^{(p-1)}D_1$,  the 
$\g_0^{(1)}$-submodule of
$\g_{p-2}/\g_{p-2}^\dagger$  generated by $x_1^{(p-1)}D_1$ is isomorphic
to  $L((p-2)\varpi_1)$.   But since $\dim \g_{p-2} = 2p$,  it follows that
the $v_k$ ($0 \leq k \leq p$) comprise a basis of $\g_{p-2}^\dagger$
and that $\g_{p-2}/\g_{p-2}^\dagger \cong L((p-2)\varpi_1)$.    Consequently, all  the assertions
in (v) hold.   \qed
\bi

The submodules in \eqref{eq:2.82a} and \eqref{eq:2.82b}  are useful
in describing the subalgebras of the restricted Lie algebras of Cartan type
which contain $\g_{-1} \oplus \g_0$.
\m

\begin{Pro} {\rm (Compare  \label{Pro:2.845} \ {\rm {\cite[Chap.~1,
Sec.~10, Prop.]{KS}} and \cite[Lem.~5.2.3]{St}}).} \  Let $\g =
\bigoplus_{j=-1}^r \g_j$ be one of the simple restricted Lie
algebras  $W(m;\un 1),$ $S(m;\un 1)^{(1)}$, or $H(2m;\un
1)^{(2)}$,  and let $L$ be a (not necessarily graded)  subalgebra
of $\g$ containing $\g_{-1} \oplus \g_0$.   Then only the
following possibilities can occur:
\begin{itemize}
\item[{\rm (a)}] $L= \g$,

\item[{\rm (b)}] $L =  \g_{-1} \oplus \g_0$,

\item[{\rm (c)}] $\g = W(m;\un 1)$ and $L = \g_{-1} \oplus  \g_0 \oplus \g_1^\sharp$,   where
$\g_1^\sharp$ is  as in \eqref{eq:2.82a}, and $L \cong
\mathfrak{sl}_{m+1}$ when $m+1 \not \equiv 0$ mod $p$.

\item[{\rm (d)}] $\g = W(m;\un 1)$ and $L =
\g_{-1} \oplus \g_0 \oplus \g_1^\dagger \oplus L_2 \oplus \cdots
\oplus L_{r-1}$,   where  $L_\ell = \g_\ell^\dagger$ for $\ell
\neq m(p-1)-p$, and $\bigl(S(m;\un 1)^{(1)}\bigr)_\ell  \subseteq  L_\ell \subseteq
\g_\ell^\dagger$ for $\ell = m(p-1)-p$;
 thus  $S(m;\un 1)^{(1)} \oplus \F \mathfrak D_1   \subseteq
 L \subseteq S(m;\un{1}) \oplus \F \mathfrak D_1$
 where $\mathfrak D_1 = \sum_{j=1}^m x_j^{(1)}D_j$;

\item[{\rm (e)}] $\g = S(m;\un 1)^{(1)}$,  where $m + 1 \equiv 0
\mod p,$ and $L = \g_{-1} \oplus \g_0 \oplus \g_1^\sharp \cong
\mathfrak{psl}_{m+1}.$
\end{itemize}
\end{Pro}

\m

\begin{Pro} \label{Pro:2.848}
 \ {\rm (See  \cite[Chap.~1, Sec.~10]{KS} and {\cite[Lem.~5.2.2]{St}}.)}
\begin{itemize}
\item[{\rm (a)}] The space  $H(2m)_\ell$  is irreducible as a
module for $H(2m)_0$  whenever $1 \leq \ell \leq  p - 3.$
Moreover, $H(2m)_1 \cong \mathcal O(2m)_3 \cong S^{3}(V)$
where $V$ is the natural $2m$-dimensional module for $H(2m)_0
\cong \mathfrak{sp}_{2m}$. \item[{\rm (b)}] The space
$H(2m;\underline{1})_\ell$ is irreducible as a module for
$H(2m)_0$ for all $\ell$. 
\end{itemize}
\end{Pro} \m

The next lemma is an immediate consequence of Theorem
\ref{Thm:2.82}. \m

\begin{Lem} \label{Lem:2.85} \
\begin{itemize}\item [\rm{(a)}] Suppose that $\g = \bigoplus_{j\geq -1} \g_j$
is a Lie algebra of Cartan type $S(m)$ or $CS(m)$ with the natural grading.
Relative to the Cartan subalgebra $\mathfrak h:= \spa_\F\{D_{\ell,\ell+1} (x_\ell^{(1)}x_{\ell+1}^{(1)})
\mid \ell=1,\dots, m-1\}$ of $\g_0^{(1)}$ and corresponding
Borel subalgebras $\mathfrak b^+$ and $\mathfrak b^-$ of $\g_0^{(1)}$,  the
element $D_1$ (resp. $D_m$)  is a $\mathfrak b^-$-primitive (resp. $\mathfrak b^+$-primitive)  vector of $\g_{-1}$ of weight
$-\varpi_{1}$ (resp. $\varpi_{m-1}$).
For $k=2,3$,  the $\g_0^{(1)}$-module $\g_k$ is generated by the $\mathfrak b^-$-primitive vector $D_{1,m}(x_m^{(k+2)})$
of weight
$-\varpi_1 - (k+1)\varpi_{m-1}$. \smallskip

\item [\rm{(b)}] Suppose that $\g = \bigoplus_{j\geq -1} \g_j$
is a Lie algebra of Cartan type  $H(2m)$ or
$CH(2m)$ with the natural grading.
Relative to the Cartan subalgebra $\mathfrak h:= \spa_\F\{h_\ell =
-D_H(x_\ell^{(1)}x_{\ell+m}^{(1)}), \ 1 \leq \ell \leq  m\}$ of $\g_0^{(1)}$
and corresponding
Borel subalgebras $\mathfrak b^+$ and $\mathfrak b^-$ of $\g_0^{(1)}$,
the vector $D_{1}= D_H(x_{1+m}^{(1)})$ is
a $\mathfrak b^+$-primitive vector of ${\g}_{-1}$ of weight $\varpi_1$.
If  $p > 5$ and  $m \geq 2$,  the irreducible  $\g_0^{(1)}$-module $\g_k$ is generated by
the $\mathfrak b^-$-primitive vector
$e^{-(k+2)\varpi_1} = D_H(x_{1+m}^{(k+2)})$,  $k = 2, 3$.    If $p = 5$,  then the
${\g}_0^{(1)} $-module ${\g}_2$ is generated by a
$\mathfrak b^-$-primitive vector $e^{-4\varpi_1},$ and the
${\g}_0^{(1)} $-module ${\g}_3$ is generated by a
$\mathfrak b^-$-primitive vector $e^{-3\varpi_1- \varpi_2}$ (resp. $e^{-3\varpi_1}$)  when $m \geq 2$
(resp. $m=1$).
  \end{itemize}
\end{Lem}

\pf  From Theorem \ref{Thm:2.19}, we know that
$S(m)_k$  is spanned by derivations
of the form $D_{i,j}(f),$ $f \in \mathcal O(m)_{k+2}$ for $k=-1,2,3$.
Since by Lemma \ref{Lem:2.18}(c) and \eqref{eq:dijdef},

\begin{eqnarray*}
[ D_{1,2}\bigl(x_1^{(1)} x_2^{(1)}\bigr),  D_{1,m}\bigl(x_m^{(j)}\bigr)]
 &=&  - D_{1,2}\bigl(x_2^{(1)}x_m^{(j-1)}\bigr) \\
&=& - x_m^{(j-1)}D_{1}  =  -   D_{1,m}\bigl(x_m^{(j)}\bigr)\end{eqnarray*}

\noindent and

\begin{eqnarray*}
[ D_{m-1,m}\bigl(x_{m-1}^{(1)} x_m^{(1)} \bigr),   D_{1,m}\bigl(x_m^{(j)}\bigr)]
&=&  - j  D_{1,m}\bigl(x_m^{(j)}\bigr) +   D_{1,m-1}\bigl(x_{m-1}^{(1)}x_m^{(j-1)}\bigr) \\
&=& - j  x_m^{(j-1)}D_1 + x_m^{(j-1)}D_{1}  \\
&=&  -(j-1)  D_{1,m}\bigl(x_m^{(j)}\bigr),  \end{eqnarray*}

\noindent  and since $[D_{\ell,\ell+1}(x_\ell^{(1)}x_{\ell+1}^{(1)}),  D_{1,m}\bigl(x_m^{(j)}\bigr)]
= 0$ for all $\ell \neq 1,m-1$, we see that
$D_{1,m}\bigl(x_m^{(j)}\bigr)$ is a weight vector relative
to $\mathfrak h$  of weight $-(\varpi_1 + (j-1)\varpi_{m-1})$.

Now it follows from  Theorem  \ref{Thm:2.82}  that for $k=2,3$ the
$S(m)_0$-module $\g_k = S(m)_k = \spa_{\mathbb F}\{ D \in W(m)_k
\mid \di (D) = 0\}$ is generated by a $\mathfrak b^+$-primitive vector of
weight $(k+1)\varpi_1 + \varpi_{m-1}$ relative to $\h$. The vector
$D_{1,m}\bigl(x_m^{(k+2)}\bigr)$ is  a $\mathfrak b^-$-primitive vector of
${\g}_k$ of weight $-(\varpi_1 + (k+1)\varpi_{m-1})$, as required.

    We turn our attention now to $\g$ of Cartan type  $H(2m)$ or $CH(2m)$.
In this case, $\g_k$ is seen to be irreducible for $k \leq  p-3$
by \cite[Chap.~1, Sec.~10]{KS}  (compare also
\cite[Lem.~5.2.2]{St}).
We have the Cartan subalgebra $\mathfrak h$ of $\g_0^{(1)}$
spanned by vectors
$$h_\ell= -D_H(x_\ell^{(1)}x_{\ell+m}^{(1)}) =
x_\ell^{(1)}D_{\ell} - x_{\ell+m}^{(1)}D_{\ell+m} \qquad  (1 \leq \ell \leq  m)$$
 (see
\eqref{eq:H0basis}),  and corresponding to the simple roots relative to $\h$,  we have the root vectors
\begin{eqnarray*} e_\ell &=&
-D_H(x_\ell^{(1)}x_{\ell+1+m}^{(1)}) = x_\ell^{(1)}D_{\ell+1} - x_{\ell+1+m}^{(1)}D_{\ell+m}
\qquad \ \  (1 \leq \ell < m), \\
f_\ell &=& -D_H(x_{\ell+1}^{(1)}x_{\ell+m}^{(1)})=
x_{\ell+1}^{(1)}D_{\ell} - x_{\ell+m}^{(1)}D_{\ell+1+m} \qquad \qquad  (1\leq \ell < m), \\
e_m &=& D_H(x_m^{(2)}) = x_m^{(1)}D_{2m} \\
f_m &=& -D_H(x_{2m}^{(2)})
= x_{2m}^{(1)}D_m. \end{eqnarray*}      The
$\mathfrak b^-$-primitive vector of ${\g}_k $ is realized by
$D_H(x_{1+m}^{(k+2)})$, which has weight $-(k+2)\varpi_{1}$ for $k \leq  p-3$.   It remains only to show
that when $p = 5$, the ${\g}_0^{(1)} $-module ${\g}_3$
is generated by a $\mathfrak b^-$-primitive vector $e^{-(3\varpi_1 + \varpi_2)}$ of weight $-3\varpi_1-\varpi_2$.
But that follows from Lemma \ref{Lem:2.87} below. \qed

\m
\begin{Lem} \label{Lem:2.87} \
Let  $V$ be the natural $2m$-dimensional module for  $\mathfrak{sp}_{2m}$.
If $p > 5$, then $S^5(V)$ is an irreducible $\mathfrak{sp}_{2m}$-module.
If $p = 5$, then $S^5(V)$ has a trivial submodule $Y := \spa_\F\{\big(x_i^{(1)}\big)^5 \mid i=1,\dots,2m\}$.
 The quotient module
$S^5(V)/Y$   is an irreducible module for  $\mathfrak{sp}_{2m}$
with a ${\mathfrak b}^+$-primitive vector of weight  $3\varpi_1 + \varpi_2  \  (3\varpi_1$ if $m=1)$.
\end{Lem}

\pf   We identify $V$ with $\mathcal O(2m;\un 1)_1$
and  the  action of
$\mathfrak{sp}_{2m}$  with that of $H(2m)_0$ on $\mathcal O(2m;\un 1)_1$.      
Then when $p > 5$,   we can identify  $S^5(V)$
with $\mathcal O(2m;\un 1)_5$.   
By \cite[Lem.~5.2.2]{St},  the $H(2m)_0$-module
$\mathcal O(2m;\un 1)_k$ is irreducible for all $k$.     Thus, the
result holds for $p > 5$.       Now when $p = 5$,  the space $Y:= \spa_\F\{\big(x_i^{(1)}\big)^5 \, \big | \, i=1,\dots,2m\}$ is a trivial  $\mathfrak{sp}_{2m}$-submodule of $S^5(V)$.
The image of $\big(x_1^{(1)}\big)^4x_2^{(1)}$ in $S^5(V)/Y$ is a 
$\mathfrak b^+$-primitive vector
  of weight $3 \varpi_1 + \varpi_2$
if $m \geq 2$  ($3\varpi_1$ if $m = 1$) relative to the Cartan
subalgebra $\h$  in the proof of Lemma \ref{Lem:2.85}(b). 
Since $\dim S^5(V)/Y = \dim H(2m;\un 1)_3$,   and
$H(2m;\un 1)_3$ is an irreducible module for $\mathfrak{sp}_{2m} \cong H(2m)_0$
by Proposition \ref{Pro:2.848}(b), it follows from Lemma \ref{Lem:2.85}(b)
that $S^5(V)/Y$ is a irreducible $\mathfrak{sp}_{2m}$-module
with a ${\mathfrak b}^+$-primitive vector
  of weight $3 \varpi_1 + \varpi_2$
if $m \geq 2$  ($3\varpi_1$ if $m = 1$).    \qed
\m

\begin{Lem} \label{Lem:2.88} \  Assume  $p = 5$. Then

\begin{equation*}H(2;\un 1)^{(1)}=  \tilde H(2;\un 1) = \left
\{ D_H(f) \mid f \in \mathcal O(2;\un 1)\right\} = H(2;\un 1)^{(2)} \, \oplus\, \F
D_H(x_1^{(4)}x_2^{(4)}), \end{equation*}  where $D_H(x_1^{(4)}x_2^{(4)})  \in$
$\left(H(2;\un 1)^{(1)} \right)_6$.   If $\g =
\bigoplus_{i=-1}^{r}\g_i$ is a restricted Lie algebra of Cartan
type $H(2;\un 1)$ or $CH(2;\un 1)$, and $\widehat \g$ is
the subalgebra  of $\g$
generated by  the local part $\g_{-1} \oplus \g_0 \oplus \g_1$,
then
\begin{itemize}
\item[{\rm (a)}] $\widehat \g_3 = \g_3 \,=\,\spa_\F \left \{D_H(f)\ \big | \  f =
x_1^{(4)}x_2^{(1)}, x_1^{(3)}x_2^{(2)}, x_1^{(2)}x_2^{(3)},
x_1^{(1)}x_2^{(4)}\right\}$  is an irreducible
${\g}_0^{(1)}$-module generated by the $\mathfrak b^-$-primitive vector  \break
$D_H(x_1^{(1)}x_2^{(4)})$ of weight  $-3 \varpi_1$;

\item[{\rm(b)}]  $(\widehat \g)^{(1)}\,=\,H(2;1)^{(2)}$ is a
simple Lie algebra;

\item[{\rm (c)}] $\widehat \g_i \subset (\widehat \g)^{(1)}$ for
all $i\ne 0$;

\item[{\rm (d)}] $\widehat\g_6 = 0$.
\end{itemize}
\end{Lem}

\pf  This result is apparent from  Theorem \ref {Thm:2.24} (iii)
and equations \eqref{eq:2.30} and \eqref{eq:2.31}.     \qed

\m
 \section {Melikyan Lie
algebras  \label{sec:2.14}}

 \m In characteristic 5, there exist finite-dimensional simple Lie
algebras
 which are neither classical nor of Cartan type.  Here we quote from
 {\cite[Sec.~4.3]{St}}, which in turn follows the presentation of \cite{Ku2}.

 \m Let $\widetilde {W}(2;\un n)$ denote a second copy of the vector
 space $W(2; \un n)$ and set

\begin{equation}\label{eq:2.45} M(2;\un n):= \mathcal O(2;\un n)
\oplus W(2;\un n) \oplus \widetilde {W}(2;\un n).  \end{equation}

\m \noindent The space $M(2;\un n)$ can be given the structure of
a Lie algebra by defining an anticommutative bracket multiplication
as follows:
\begin{gather} \label{eq:2.46}\\ {[D, \widetilde E]} = \widetilde{[D,E]} + 2
\di (D) \widetilde E \nonumber \\
{[D, f]} = D(f) - 2 \di (D)f \nonumber   \\
{[f_1 \widetilde{D_1} + f_2 \widetilde{D_2}, g_1 \widetilde{D_1} +
g_2
\widetilde{D_2}]} = f_1g_2 - f_2g_1 \cr [f, \widetilde{E}] = fE \nonumber  \\
{ [f,g]}
= 2(f \widetilde{D}_g - g \widetilde{D}_f ) \qquad \hbox{\rm where} \nonumber \\
\widetilde {D}_f = D_1(f)\widetilde {D_2} - D_2(f)\widetilde {D_1}
.\nonumber  \end{gather}

\noindent for all $D \in W(2;\un n), \widetilde{E} \in
\widetilde{W}(2;\un n)$, and $f, g, f_i,g_i \ (i=1,2)$ in
$\mathcal O(2;\un n)$. The product $[D,E]$ of two elements of
$W(2;\un n)\subset M(2;\un{n})$ is the same as in the \W \  Lie
algebra of Cartan type $W(2;\un n)$.

Although the above makes sense for any prime $p$, it is only for
$p=5$ that the multiplication defined by (\ref{eq:2.46}) satisfies
the Jacobi identity; the proof of this fact can be found in
\cite[Lem.~4.3.1]{St}. The Lie algebras $M(2;\un n)$ are called
the {\it Melikyan algebras}. From the construction we see that
\begin{equation}\label{eq:2.47} \dim M(2;\un n) = 5^{n_1+n_2}+ 2 \cdot
5^{n_1+n_2}+2 \cdot 5^{n_1+n_2} = 5^{n_1+n_2+1}.  \end{equation}
The subspace $W(2;\un{n})$ of $M(2;\un{n})$ is a Lie subalgebra of
$M(2;\un{n})$, and both $\mathcal O(2;\un n)$ and $\widetilde
{W}(2;\un n)$ are $W(2;\un{n})$-modules.

 We assign gradation degrees $\deg_M(\cdot)$ to
elements of $M(2;\un n)$, by using the natural gradation degrees
$\deg_W$ in the Lie algebra $W(2;\un n)$ and the degrees
$\deg_{\mathcal O}$ in the associative algebra $\mathcal O(2;\un
n)$:

\begin{eqnarray*} \deg_M(D) &=& 3 \deg_W(D) \\
\deg_M(\widetilde{E}) &=& 3 \deg_W(E) + 2 \\
\deg_M(f) &=& 3 \deg_{\mathcal O}(f) - 2,
\end{eqnarray*}

\noindent where $D, \widetilde E,$ and $f$ are as above.  Thus,
$M(2;\un n) = \bigoplus_{j=-3}^r M(2;\un n)_j$, where $r=
3(5^{n_1}+ 5^{n_2}) -7$.    We refer to this as
the {\em natural grading} of $M(2;\un n)$, a phrase 
motivated  by the fact that $\bigoplus_{j\geq 0}^r M(2;\un n)_j$
is the unique maximal subalgebra of codimension 5 and depth $\geq 3$
(see \cite[Thm.~4.3.3]{St},  which is based on \cite{Ku2}).  
The natural grading  determines a $\Z/3\Z$-grading
$M(2;\un n) = M_{\ov {-2}} \oplus M_{\bar 0} \oplus M_{\bar 2}$
where $M_{\ov{-2}} = \mathcal O(2;\un n)$, $M_{\bar 0} = W(2;\un
n)$, and $M_{\bar 2} ={\widetilde W}(2;\un n)$.
  
\m The natural gradation on $M(2;\un n)$ is inspired  by a certain
gradation of a classical simple Lie algebra $\mathcal G$ of type G$_2$.
Suppose $\mathcal H$ is a Cartan subalgebra of $\mathcal G$ and $\{\alpha,\beta\}$
is a base of simple roots with $\alpha$ short and $\beta$ long. We
assign $\alpha$ degree 1, $-\alpha$ degree $-1$, and $\pm \beta$
degree 0.  Then $\mathcal G = \bigoplus_{j=-3}^3 \mathcal G_{j}$, where

\begin{gather*} \mathcal G_{-3} = \mathcal G^{-3\alpha-2\beta} \oplus
\mathcal G^{-3\alpha-\beta}, \quad \mathcal G_{-2} = \mathcal G^{-2\alpha - \beta}, \quad
\mathcal G_{-1} =
\mathcal G^{-\alpha -\beta} \oplus \mathcal G^{-\alpha} \\
\mathcal G_0 = \mathcal G^{-\beta} \oplus \mathcal H \oplus \mathcal G^{\beta} \\
\mathcal G_{1} = \mathcal G^{\alpha+\beta} \oplus \mathcal G^{\alpha}, \ \ \qquad \mathcal G_{2} =
\mathcal G^{2\alpha +\beta},\ \ \qquad \mathcal G_{3} = \mathcal G^{3\alpha+2\beta} \oplus
\mathcal G^{3\alpha+\beta}.
\end{gather*}

Indeed for the Melikyan algebra $M = M(2;\un n)$, we have

\begin{gather*} M_{-3} = \F D_1 \oplus \F D_2, \qquad M_{-2} = \F 1,
\qquad M_{-1} = \F \widetilde {D_1} \oplus \F \widetilde {D_2} \\
M_0 = \spa_\F \{x_i^{(1)} D_j \mid i,j = 1,2 \}   \\
M_{1} = \F x_1^{(1)}  \oplus \F x_2^{(1)},
\end{gather*}
so that $\bigoplus_{j \leq 0} M_j \cong \bigoplus_{j\leq 0}
\mathcal G_j$ as graded Lie algebras.

\m The Melikyan algebras can be characterized by their gradation.

\m
\begin{Pro} \label{Pro:2.48} {\rm (See \cite{Ku2} and \cite[Thm.~5.4.1]{St}.)}
Let $\g = \bigoplus_{j=-3}^{r} \g_j$  be a finite-dimensional,
transitive and $1$-transitive graded Lie algebra of height $r$, and suppose that
$\bigoplus_{j=-3}^0\,\g_j \cong\, \bigoplus_{j=-3}^0\,M(2;\un
1)_j$ as graded Lie algebras.   If $r\le 3$, then $\g$ is isomorphic
to a classical Lie algebra of type ${\mathrm G}_2$. If $r>3$, then
$\g$ is isomorphic as a graded Lie algebra to a Melikyan algebra
$M(2;\un n)$ with its natural grading.
\end{Pro}

\pf By our assumption, $\g_0\cong\mathfrak{gl}_2$, both $\g_{-1}$
and $\g_{-3}$ are two-dimensional, irreducible faithful
$\g_0$-modules, and $\g_{-}$ is generated by $\g_{-1}$. Since $\g$
is $1$-transitive, Proposition~\ref{Pro:1.17} shows that $\g$
contains a unique minimal ideal $\mathcal I$ which is graded and
contains $\g_{-1}$.

Since $\dim \g_{-2}=1$, the Lie product on $\g$ induces a
$\g_{0}^{(1)}$-equivariant pairing $\g_{-3}\times \g_1\rightarrow
\F$. The $1$-transitivity of $\g$ along with the above remark implies
that this pairing is nondegenerate. As a result,
$\g_1\cong\g_{-3}\cong\g_{-1}$ as $\g_0^{(1)}$-modules.

The Lie product on $\g$ induces a $\g_0$-module isomorphism
$\g_{-3}\cong\g_{-2}\otimes \g_{-1}$. Choose $u_1\in\g_1$ and
$u_{-3},v_{-3}\in \g_{-3}\setminus\{0\}$ such that
$[u_1,u_{-3}]\ne 0$ and $[u_1,v_{-3}]=0$. Choose
$u_{-1},v_{-1}\in\g_{-1}$ such that $[u_{-1},[u_1,u_{-3}]]=u_{-3}$
and $[v_{-1},[u_{1},u_{-3}]]=v_{-3}$. Then $[[u_1,u_{-1}],
u_{-3}]=[[u_1,u_{-3}],u_{-1}]=-u_{-3}$ and $[[u_1,u_{-1}],
v_{-3}]=0$.  Hence, $[\g_{-1},\g_1]$ is a noncentral ideal of
$\g_0$, and $[\g_{-1},\g_1]\not\subseteq \g_0^{(1)}$. But then
$\g_0=[\g_{-1},\g_1]$ implying ${\mathcal I}_0=\g_0$ and
${\mathcal I}_1=\g_1$.

By Lemma \ref{Lem:1.20}, we have ${\mathcal
I}=\bigoplus_{j=-3}^s\,{\mathcal I}_j$ where $s\in\{r-1,r\}$.
Clearly, $[\g_i,{\mathcal I}_s]=0$ for all $i>0$. Since any
$\g_0$-submodule $\mathcal K_s$ of ${\mathcal I}_s$ generates an ideal
$\bigoplus_{i\ge 0}\big(\ad \g_{-1}\big)^i(\mathcal K_s)$ of $\g$ contained
in $\mathcal I$, the minimality of $\mathcal I$ shows that
${\mathcal I}_s$ is an irreducible $\g_0$-module. If
$[\g_0,{\mathcal I}_s]=0$, then $\bigoplus_{i\ge 1}\big(\ad
\g_{-1}\big)^i({\mathcal I_s})$ is an ideal of $\g$ contained in
$\bigoplus_{j=-3}^{s-1}\,{\mathcal I}_j$, a contradiction.
Therefore, ${\mathcal I}_s=[\g_0,{\mathcal I}_s]$.

Let ${\mathcal J}$ be an arbitrary nonzero ideal of the Lie
algebra ${\mathcal I}$. Then $\g_{-3}\subset {\mathcal J}$ in view
of transitivity and the fact that $\g_{-1}\subset{\mathcal I}$.
Since $\g_1\subset{\mathcal I}$ and $\g$ is $1$-transitive, it
follows that $\g_{-1}\subset{\mathcal J}$. But then
$\g_0=[\g_{-1},\g_1]\subset\mathcal J$ yielding ${\mathcal
I}_s=[\g_0,{\mathcal I}_s]\subset {\mathcal J}$. This implies
${\mathcal I}=\bigoplus_{i\ge 0}\big(\ad \g_{-1}\big)^i({\mathcal
I}_s)\subseteq \mathcal J$. We conclude that ${\mathcal I}$ is a
simple Lie algebra.

For $k\ge -3$, set ${\mathcal I}_{(k)}\eqdef \bigoplus_{j\ge
k}\,{\mathcal I}_j$. Clearly, $[{\mathcal I}_{(i)},{\mathcal
I}_{(j)}]\subseteq {\mathcal I}_{(i+j)}$ for all $i,j\in\mathbb
Z$, that is ${\mathcal I}\,=\,{\mathcal
I}_{(-3)}\supset\cdots\supset{\mathcal I}_{(s)}\supset 0$ is a
filtration of the Lie subalgebra $\mathcal I$. By construction,
the corresponding graded Lie algebra $\bigoplus_{j\ge
3}\big({\mathcal I}_{(j)}/{\mathcal I}_{(j+1)}\big)$ is isomorphic
to the graded Lie algebra $\mathcal I$. It is straightforward to
check that ${\mathcal I}_{(0)}$ is a maximal subalgebra of
$\mathcal I$. Since $\mathcal I$ is a simple Lie algebra we now can
apply \cite[Thm.~5.4.1]{St} to deduce that either $s=3$ and
${\mathcal I}$ is isomorphic as a graded algebra to a simple Lie
algebra of type ${\mathrm G}_2$,  or $s>3$ and $\mathcal I$ is
isomorphic as a graded algebra to a Melikyan algebra $M(2,\un{n})$
with the natural grading.

The adjoint action of $\g$ on $\mathcal I$ gives rise to a Lie
algebra homomorphism $\pi\colon\,\g\rightarrow\Der({\mathcal I})$.
If $\ker\,\pi\ne 0$, then $\ker\,\pi\supseteq\mathcal I$, by the
minimality of $\mathcal I$. Since $[{\mathcal I},{\mathcal I}]\ne
0$, this is impossible. So $\pi$ is injective. If ${\mathcal I}$
is a Lie algebra of type ${\mathrm G}_2$, then all derivations of
${\mathcal I}$ are inner, forcing $\g=\mathcal I$. Now suppose
${\mathcal I}\cong M(2;\un{n})$, and identify $\g$ with a
subalgebra of $\Der\big(M(2;\un{n})\big)$. Note that the  natural  
gradation of $M(2;\un{n})$ gives rise to a grading of the
derivation algebra $\Der\big(M(2;\un{n})\big)$ relative to which
$\pi\colon\g\rightarrow\Der\big(M(2;\un{n})\big)$ becomes a homomorphism
of graded Lie algebras. According to \cite[Thm.~7.1.4]{St},
the Lie algebra $\Der\big(M(2;\un{n})\big)$ is spanned by
$\ad M(2;\un{n})$ and the $p$th powers of $D_1,D_2\in
W(2;\un{n})=M(2;\un{n})_{\bar{0}}$. Since
$\deg\,(D_i)^{p^j}=-3p^j$ for all $j\ge 0$, this description
implies that $$\pi(\g)\,\subseteq \,{\bigoplus_{j\ge
-3}}\big(\Der\big(M(2;\un{n})\big)\big)_j\,=\,\ad M(2;\un{n}).$$  Thus,
$\g={\mathcal I}\cong M(2;\un{n})$, completing the proof.\qed

%-----------------------------------------------------------------------
% Beginning of chap3.tex   9-15-05 version
%-----------------------------------------------------------------------
%
% AMS-LaTeX 1.2 sample file for a monograph, based on amsbook.cls.
% This is a data file input by chapter.tex.
%%%%%%%%%%%%%%%%%%%%%%%%%%%%%%%%%%%%%%%%%%%%%%%%%%%%%%%%%%%%%%%%%%%%%%%%%%%%%%%%%%%%%%%%%%%%%%%%%%%%%%%%%%%%%%%%%%%%%%%
%%%%%%%%%%%%%%%%%%%%%%%%%%%%%%%%%%%%%%%%%%%%%%%%%%%%%%%%%%%%%%%%%%%%%%%%
 
  \chapter{The Contragredient Case}  
\bi 

\section {\ Introduction  \label{sec:3.1}}  

 \m This chapter is devoted to a proof of the Main Theorem for graded Lie algebras
$\g = \bigoplus_{j=-q}^r \g_j$  in which $\g_{-1}$ and $\g_1$ are
dual modules for $\g_0$ (the so-called contragredient case). We
begin by focusing on three-dimensional subalgebras of $\g$  and
their representations.    We investigate conditions under which
certain pairs of such subalgebras, which share a common
one-dimensional Cartan subalgebra, generate an
infinite-dimensional algebra.   Next we focus on extreme vectors in
$\g_1$ and $\g_{-1}$ relative
to Borel subalgebras ${\mathfrak b}^+$ and ${\mathfrak b}^-$ of
$\g_0$  and show that a ${\mathfrak b}^+$-primitive vector of
$\g_1$ and a ${\mathfrak b}^-$-primitive vector  of  $\g_{-1}$  generate a
three-dimensional simple Lie algebra. We then proceed to show that
$\g$ contains a subalgebra with a ``balanced'' gradation and use
that fact to conclude that either  $\g$ is  a classical Lie algebra or 
the characteristic of $\F$ is 5 and $\g$ is isomorphic to a Melikyan Lie algebra.

\m
\section {\ Results on
modules for three-dimensional Lie algebras\label{sec:3.2}}

In this section, we create computational tools for dealing with
modules for three-dimensional simple Lie algebras and Heisenberg
Lie algebras and apply them to show that  certain pairs of
such Lie algebras sharing a common Cartan subalgebra  generate
an infinite-dimensional Lie algebra.

The first result, which can be proved easily by induction, concerns modules
for a Lie algebra $\mathfrak s$ spanned by elements $e,f,h$
satisfying the commutation relations 

\begin{equation}\label{eq:3.1}[e,f] = h, \quad \quad [h,e] = \xi e, \quad \quad  [h,f]
= -\xi f. \end{equation}

\noindent When $\xi$ is specialized to be 2 (or equivalently, any nonzero scalar as $p >2$), the result gives
information about $\mathfrak{sl}_2$-modules; when $\xi = 0$,  the result applies to modules for a Heisenberg Lie
algebra.
 
 \bi
\begin{Pro}\label{Pro:3.4}
\begin{itemize}
\item[{\rm (a)}]   Let $\mathcal U$ be a module for ${\mathfrak s}
= \spa_{\mathbb F}\{e,f,h\}$, where $e,f,h$ satisfy the
commutation relations in \eqref{eq:3.1}, and suppose $\mathcal U$
contains a vector $u_0 \neq 0$ such that

\begin{equation*}h.u_0 = \mu u_0 \quad \quad f . u_0 = 0. \end{equation*}

\noindent Set $u_{-1} = 0$ and $u_j = e.u_{j-1}$ for $j = 1,2, \dots$ \ .  Then

\begin{equation*}f.u_j = -\Big(j\mu + \frac {j(j-1)} {2}\xi \Big) u_{j-1}\end{equation*}

\noindent for all $j = 0,1, \dots$.   If $f.u_j = 0$ for some $j$, then $j \equiv 0 \mod p$
or  $2\mu =  -(j-1)\xi$. \item[{\rm (b)}]   Let $\mathcal V$ be a
module for ${\mathfrak s} =\spa_{\mathbb F}\{e,f,h\}$, and suppose
$\mathcal V$ contains a vector $v_0 \neq 0$ such that

\begin{equation*}h.v_0 = \lambda v_0 \quad \quad e.v_0 = 0. \end{equation*}

\noindent Set $v_{-1}=0$ and $v_j = f .v_{j-1}$ for $j = 1,2, \dots$\ .  Then

\begin{equation*}e. v_j = \Big (j\lambda - \frac {j(j-1)} {2}\xi \Big) v_{j-1}\end{equation*}

\noindent for all $j = 0,1, \dots$.   If $e.v_j = 0$ for some $j$, then $j \equiv 0 \mod p$
or $2\lambda = (j-1)\xi$.
\end{itemize}

\end{Pro}

\bi In the next result and subsequent ones, we will adopt the
convention of using the corresponding capital letter for the
adjoint mapping.  Thus, $E = \ad e$, etc.

\bi
\begin{Pro} \label{Pro:3.5} \ Suppose ${\mathfrak s} = \spa_{\mathbb F}\{e,f,h\}$  is a
subalgebra inside a Lie algebra $\g$,  and the elements $e,f,h$
satisfy the commutation relations in \eqref{eq:3.1}.  Assume
$u_0,v_0$ are elements of $\g$ such that

\begin{equation}\begin{array}{ccc} \label{eq:3.4} Fu_0 = 0&&\hspace*{-.5cm}Hu_0 = -\lambda u_0, \\
Ev_0 = 0&&\hspace*{-.55cm}Hv_0= \lambda v_0,  \\
&[u_0,v_0] = h.&\end{array}
 \end{equation}

\noindent Then for  $\ell  \in \{2, \dots, p-1\}$, $[E^\ell u_0,
F^\ell v_0]$ is a nonzero multiple of $h$ if and only if $2\lambda
\neq (k -1)\xi$ for $k= 2, \dots, \ell$ and $\lambda \neq \ell
\xi$.
\end{Pro}

\pf  Observe that

\begin{equation*} \hbox{\rm ad}\bigl(E^\ell u_0\bigr) = [\underbrace{E, \dots [E,[E}_\ell,U_0]]\dots] =
\sum_{k = 0}^\ell (-1)^\ell {\ell \choose k}E^{\ell-k}U_0 E^k.
\end{equation*}

\noindent Thus

\begin{equation}[E^\ell u_0, F^\ell v_0] = \sum_{k = 0}^\ell (-1)^k
{\ell \choose k}E^{\ell-k}U_0 E^k F^\ell v_0.\nonumber
\end{equation}

\noindent Now by Proposition \ref{Pro:3.4}\,(b), $E^k$ acting on
$F^\ell v_0$ is a multiple of $F^{\ell-k}v_0$,  and
$U_0 F^{\ell-k}v_0=F^{\ell-k}U_0v_0$ $= F^{\ell-k}h$,  which
is 0 if $\ell- k$ $\geq 2$. Thus,

\begin{eqnarray}   \label{eq:3.5}&& \\
{[E^\ell u_0, F^\ell v_0]}  &=& (-1)^{\ell-1}\ell
EU_0E^{\ell-1}F^\ell v_0 +
(-1)^\ell U_0E^\ell F^\ell v_0 \nonumber  \\
 &=& (-1)^{\ell-1}\ell  \prod_{k= 2}^\ell  \Big(k\lambda -\frac {k(k-1)}{2}\xi\Big)
EU_0F v_0 \nonumber \\
&& \quad \quad +
 (-1)^{\ell}\prod_{k = 1}^\ell  \Big(k\lambda -\frac {k(k-1)}{2}\xi\Big) U_0 v_0
\nonumber \\
 &=&  (-1)^{\ell-1}\prod_{k= 2}^\ell \Big(k\lambda -\frac {k(k-1)}{2}\xi\Big)
(\ell\xi-\lambda)h. \nonumber
\end{eqnarray}

\noindent We obtain a nonzero multiple of $h$ precisely when
$2\lambda \neq (k-1)\xi$  for any $k= 2, \dots, \ell$ and $\lambda
\neq \ell  \xi$.   \qed
 
\m
\begin{Pro}\label{Pro:3.8} Suppose $\g$ is a Lie algebra
containing nonzero elements $e_i,f_i, i = 1,2$, and  $h$
such that

$$\begin{gathered}
\left[h,e_1\right] = 2e_1, \quad \quad \quad [h,f_1] = -2f_1, \nonumber \\
[h,e_2] = ae_2, \quad \quad \quad [h,f_2] = -af_2, \nonumber \\
[e_i,f_j]
= \delta_{i,j}h, \nonumber \\
 \end{gathered}\nonumber $$

\noindent for some $a\ne 0$.  If $a\in\{1,2,\ldots,p-1\} = \F_p\setminus \{0\}$, set

\begin{eqnarray*}
e_1' = F_1^{p-a}f_2,&& \quad \quad f_1' = E_1^{p-a} e_2,   \\
e_2' = F_2^2e_1' = F_2^2F_1^{p-a}f_2, && \quad \quad f_2' = E_2^2 f_1' =
E_2^2E_1^{p-a}e_2.
 \end{eqnarray*}

If $a\not\in\mathbb F_p$, set

\begin{eqnarray*}
e_1' = E_1^{p-2}e_2,&& \quad \quad f_1' = F_1^{p-2} f_2,   \\
e_2' = E_2^{p-2}E_1^2e_2, && \quad \quad f_2' = F_2^{p-2}F_1^2f_2.
\end{eqnarray*}

\noindent Then we have

\begin{eqnarray*}
\left[h,e_1'\right] = be_1',&&\qquad   [h,f_1'] = -bf_1,   \\
\left[h,e_2'\right] = -be_2',&&\qquad   [h,f_2'] = bf_2',  \\
\left[e_i,f_j\right]
&=& \delta_{i,j}\zeta_i h,
\end{eqnarray*}

\noindent where $\zeta_1$ and $\zeta_2$ are nonzero scalars,  and  
$b=a$ for $a\in\mathbb F_p$ and $b=a-4$ for $a\not\in\mathbb F_p$.
\end{Pro}

\pf   First we assume that $a\in\{1,2,\ldots,p-1\}$. Suppose $a
\neq p-1$. Then $[e_2,e_1'] = 0 = [f_2,f_1']$, from which it
follows that

\begin{eqnarray*}
\ [e_1^{\prime},f_2^{\prime}] &=& [e_1^{\prime},E_2^2f_1^{\prime}] \,\in\, \F \,E_2^2h = 0 \\
\ [e_2^{\prime},f_1^{\prime}] &=& [F_2^2e_1^{\prime},f_1^{\prime}]
\,\in\, \F \, F_2^2 h = 0.
\end{eqnarray*}

\noindent Applying Proposition \ref{Pro:3.5} with ${\mathfrak s} =
\spa_{\mathbb F}\{e_1,f_1,h\}$, $u_0 = e_2$,  $v_0 = f _2$,
$\lambda = p-a$, $\xi=2$, and $\ell = p-a$ shows that $[e_1',f_1']
= -[E_1^{p-a}e_2,\,F_1^{p-a}f_2]$ is a nonzero multiple of $h$,
say $\mu h$.

Next take ${\mathfrak s} = \spa_{\mathbb F}\{e_2,f_2,h\}$, $u_0 =
f_1'$, and $v_0 = -\mu^{-1} e_1'$. Notice that $\lambda = \xi = a$
in this case. Applying Proposition \ref{Pro:3.5} with $\ell=2$, we
see that $[E_2^2 u_0, F_2^2 v_0] = -\mu^{-1}[e_2',f_2']$ is a
nonzero multiple of $h$.   This finishes the proof of the
proposition for $a\in \{1,2,\dots,p-2\}$.

Suppose $a = p-1$. Then $e_1' = F_1f_2$, $f_1' = E_1e_2$, $e_2' =
F_2^2 e_1' = -F_2^3 f_1$ and $f_2' = E_2^2 f_1' = -E_2^3e_1$.
Direct calculation shows that

$$\begin{gathered}
\left[e_1',f_2'\right] = 0 = [e_2',f_1']\nonumber \\
[e_1',f_1'] = -h.
 \end{gathered}\nonumber $$

\noindent Using Proposition \ref{Pro:3.5} once again with
${\mathfrak s} = \spa_{\mathbb F}\{e_2,f_2,h\}$, $u_0 = e_1$, $v_0
= f_1$, $\lambda = p-2$, $\xi = p-1$, and $\ell = 3$ allows us to
conclude that $[e_2',f_2']$ is a nonzero multiple of $h$ also.

Now assume that $a\not\in\mathbb F_p$ and take ${\mathfrak s} =
\spa_{\mathbb F}\{e_1,f_1,h\}$, $u_0 = e_2$, and $v_0 = f_2$. Then
$\lambda = -a$ and $\xi=2$. Applying Proposition \ref{Pro:3.5}
with $\ell=p-2$ shows that $[e_1',f_1']=[E_1^{p-2} u_0, F_1^{p-2}
v_0]$ is a nonzero multiple of $h$, say $\nu h$.

 Finally,  take ${\mathfrak s} = \spa_{\mathbb F}\{e_2,f_2,h\}$, $u_0 =
E_1^2e_2$, $v_0 = \nu^{-1} F_1^2 f_2$, and $\ell = p-2$ in Proposition
\ref{Pro:3.5},  where $\nu = 2(a+1)(a+4)$.    Since $\lambda=-(a+4)$ and $\xi=a$ in this case,  we
see that $[E_2^{p-2} u_0, F_2^{p-2} v_0] = \nu^{-1}[e_2',f_2']$ is
a nonzero multiple of $h$.  Moreover, because $[e_1',f_2]=0=[f_1',e_2]$ and
$p-2\ge 3$, we have
\begin{eqnarray*}
[e_{1}^{\prime},f_{2}^{\prime}] &=&
[e_{1}^{\prime},F_{2}^{p-2}F_{1}^{2}f_{2}]\,=\,
F_{2}^{p-2}[E_{1}^{p-2}e_{2},F_{1}^{2} f_{2}] = 0 \\
\ [e_{2}^{\prime},f_{1}^{\prime}] &=&
[E_2^{p-2}E_{1}^{2}e_{2},f_{1}^{\prime}] \,=\,
E_{2}^{p-2}[E_{1}^{2}e_{2},F_{1}^{p-2}f_{2}] = 0.
\end{eqnarray*}
This completes the proof. \qed

\bi

\begin{Thm} \label{Thm:3.9}  Suppose $\g = \bigoplus_{t
\in \Z} \g_t$ is a $\Z$-graded Lie algebra containing homogeneous
elements $e_i,f_i, i = 1,2,$ and $h$ satisfying the assumptions of
Proposition \ref{Pro:3.8}. Suppose $e_i \in \g_{\ell_i}$, $f_i \in
\g_{-\ell_i}$, where $\ell_1\ell_2 \geq 0$ and not both $\ell_1$ and
$\ell_2$ are 0. Then $\g$ is infinite-dimensional.
\end{Thm}

\pf We may suppose that $e_i', f_i'$ for $i = 1,2$ are as
in the statement of Proposition \ref{Pro:3.8}, and that
$[e_i',f_i'] = \zeta_i h$ where $\zeta_i \neq 0$.  Define

\begin{eqnarray*}  \widetilde e_1 = 2b^{-1}e_1',&&\qquad
\widetilde f_1 = \zeta_1^{-1}f_1', \\
\widetilde e_2 = 2b^{-1}e_2',&&\,\qquad \widetilde f_2 =
\zeta_2^{-1}f_2' ,\\ \widetilde h = [\widetilde e_1,\widetilde
f_1] &=&2b^{-1}h\ \,  = [\widetilde e_2, \widetilde f_2].
\end{eqnarray*}

\noindent Then

\begin{eqnarray*} {[\widetilde h,\widetilde e_1]} = 2 \widetilde e_1,&&\qquad \quad
{[\widetilde h,\widetilde f_1]} = -2 \widetilde f_1, \\
{[\widetilde h, \widetilde e_2] = -2 \widetilde e_2}, &&\qquad \quad
[\widetilde h, \widetilde f_2] = 2 \widetilde f_2.
\end{eqnarray*}

\noindent Moreover, $\widetilde e_1 \in \g_{k_1}$, $\widetilde e_2
\in \g_{k_2}$, $\widetilde f_1 \in \g_{-k_1}$, and
 $\widetilde f_2
\in \g_{-k_2}$, where $k_1 = -(p-b)\ell_1 - \ell_2$, $k_2 =
-(p-b)\ell_1 - 3\ell_2$ if $b\in\mathbb F_p$ and $k_1 =
(p-2)\ell_1 + \ell_2$, $k_2 = 2\ell_1 +(p-1)\ell_2$ if
$b\not\in\mathbb F_p$. Since $k_1k_2\ge 0$, we may replace the
initial set of elements with this new set and continue the
process.  Since the elements constructed lie in spaces $\g_t$ with
$|\,t\,|$ increasing at each stage, $\g$ must be
infinite-dimensional. \qed

\m
\begin{Thm} \label{Thm:3.10} {\rm (Compare  \cite[Lem.~20]{K1})}.
Suppose $\g = \bigoplus_{t \in \Z} \g_t$ is a $\Z$-graded Lie
algebra containing nonzero homogeneous elements $e_i,f_i, h_i, \,
i = 1,2,$  which satisfy $e_i \in \g_{\ell_i}$, $f_i \in
\g_{-\ell_i}$, where $\ell_1\ell_2 \geq 0$ and not both $\ell_1$
and $\ell_2$ are 0. Assume

\begin{equation*}[e_i,f_j] = \delta_{i,j}h_i, \quad [h_i,e_j] = a_{i,j}e_j,
\quad [h_i,f_j] = -a_{i,j}f_j,  \end{equation*}

\noindent where $a_{i,j}$ is the $(i,j)$-entry of the matrix

\begin{equation}A = \left(\begin{matrix} 2 & 0  \\ c & b  \\ \end{matrix} \right)
\nonumber \end{equation}

\noindent and $c \neq 0$. Then  $\g$ is infinite-dimensional.
\end{Thm}

\pf We have $[f_2,[e_1,e_2]] = -[e_1,h_2] = ce_1 \neq
0$, so that $[e_1,e_2] \neq 0$.    Similarly, $[f_1,f_2] \neq 0$.
Set $\overline{e}_1 = e_1, \overline{e}_2 = [e_1,e_2],
\overline{f}_1 = f_1$, $\overline{f}_2 = c^{-1}[f_1,f_2]$ and
$\overline{h} = h_1$. Then we claim these elements satisfy the
hypotheses of Proposition \ref{Pro:3.8} (with $a = 2$).  All this
is apparent except perhaps for the following calculations:

\begin{equation*}
\left[\overline{h},\overline{f}_2\right]  = c^{-1}[[h_1, f_1],f_2] =
 -2c^{-1}[f_1,f_2] = -2 \overline{f}_2,  \end{equation*}

\noindent and
\begin{eqnarray*} \left[\overline{e}_2,\overline{f}_2\right]
&=& c^{-1}[[e_1,e_2],[f_1,f_2]]\\ &=& c^{-1}
[[[e_1,e_2],f_1],f_2] + c^{-1}[f_1,[[e_1,e_2],f_2]]  \\
&=& c^{-1}[[h_1,e_2],f_2] + c^{-1}[f_1,[e_1,h_2]] = -[f_1,e_1] =
h_1 = \overline{h}.
\end{eqnarray*}

\noindent Since $[e_1,e_2] \in \g_{\ell_1+\ell_2}$ and
$\ell_1(\ell_1+\ell_2) \geq 0$, we may apply Theorem \ref{Thm:3.9}
to obtain the desired conclusion.  \qed

\m
\begin{Cor} \label{Cor:3.11} \  Suppose that
there are  nonzero elements $e_i,f_i\in\g$, $i = 1,2$,  and $h$ such that the relations in
Proposition \ref{Pro:3.8} hold with $a = 0$. Assume
further that $e_i \in \g_{\ell_i}$, $f_i \in \g_{-\ell_i}$, where
$\ell_1\ell_2 \geq 0$ and not both $\ell_1$ and $\ell_2$ are 0.
Then $\g$ is infinite-dimensional.
\end{Cor}

\pf This is a direct consequence of Theorem
\ref{Thm:3.10} with $h_1 = h =  h_2$, $c = 2$, and  $b = 0$. \qed

\m
\begin{Pro} \label{Pro:3.12}  \ Suppose $e,f,h$ span a Heisenberg
Lie subalgebra of a Lie algebra $\g$.  Assume $u_0,v_0$ are
elements of $\g$ such that
\begin{equation}\label{eq:3.11}\begin{array}{cccc}
{[f,u_0]} = 0, &\quad &{[h,u_0]} = -\lambda u_0,&\\
{[e,v_0]} = 0, & \quad &{[h,v_0]}= \lambda v_0,& \\
{[u_0,v_0]} = z, & \quad &{[z,v_0]} = 2v_0,&\quad [z,u_0] = -2u_0, \\
{[z,f]} = -\eta f,&\quad&{[z,e]} = \eta e& \ \ \, (\hbox{\rm some}
\ \eta \in \mathbb F).  \end{array}\end{equation}
 \noindent Then
\begin{eqnarray}
&& \left [E^\ell u_0, \, F^\ell v_0\right ]  =  (-1)^\ell  \ell!\,
\lambda^{\ell-1} (\lambda z - \ell \eta h) \label{eq:3.12}\\
&&\left [ \bigl[E^\ell u_0, F^\ell v_0 \bigr],
E^\ell u_0\right]  = 2(-1)^\ell \ell! \, \lambda^{\ell }(\ell \eta -
1)E^\ell u_0\label{eq:3.13}\\
&&\left [ \bigl[E^\ell u_0, F^\ell v_0 \bigr],
F^\ell v_0\right]  = -2(-1)^\ell \ell! \, \lambda^{\ell }(\ell \eta -
1)F^\ell v_0\label{eq:3.14}
\end{eqnarray} for all $\ell \in \{1, \dots, p-1\}$.
\end{Pro}

\pf  Following the argument in \eqref{eq:3.5},  we have

\begin{eqnarray*}
\left [E^\ell u_0,\,F^\ell v_0\right]  &=& \sum_{k = 0}^\ell (-1)^k {\ell \choose k} E^{\ell-k}U_0E^kF^\ell v_0
 \\
  &=& (-1)^{\ell-1}\ell\,\ell!\, \lambda^{\ell-1} EU_0F v_0 +
 (-1)^{\ell}\ell!\,\lambda^\ell U_0 v_0 \\
 &=& (-1)^\ell\, \ell!\, \lambda^{\ell-1} (\lambda z - \ell \eta h).
\end{eqnarray*}

\noindent Since $[\lambda z - \ell \eta h, E^\ell u_0] =\bigl
((\ell\eta -2)\lambda + \ell\eta \lambda\bigr) E^\ell u_0 = 2
\lambda (\ell\eta - 1)E^\ell u_0$, equation \eqref{eq:3.13} follows,
and \eqref{eq:3.14} is completely analogous.
 \qed

\m
\section  {\ Primitive vectors
in $\boldsymbol{\g_1}$ and $\boldsymbol{\g_{-1}}$\label{sec:3.3}} \m

    In this section, we show that if $e$ is a ${\mathfrak b}^+$-primitive
    vector in $\g_1$, and $f$ is a ${\mathfrak b}^-$-primitive vector
    in $\g_{-1}$, then $[e, \, [e, \, f]] \neq 0$,
    and $e$ and $f$ generate a three-dimensional simple Lie
    algebra. \bigskip

{\it Henceforth in this section  we assume the following:  \m
\begin{itemize}
\item[{\rm (a)}] $\g = \bigoplus_{j \in \Z} \g_j$ is a
finite-dimensional $\Z$-graded Lie algebra over an algebraically
closed field of characteristic $p>3$.
\item[{\rm(b)}] $\g_0$ is
classical reductive,  and $\mathfrak t$ is a maximal toral subalgebra of $\g_0$.
\item[{\rm(c)}] $\{\alpha_1, \dots, \alpha_m\}$ is a base of
simple roots relative to $\mathfrak t$, \item[{\rm(d)}]  $e_i \in
\g_0^{\alpha_i},$ $f_i \in \g_0^{-\alpha_i}$  denote nonzero root
vectors corresponding to these simple roots. \item[{\rm(e)}]
$\mathfrak b^+ = \mathfrak t
\oplus \mathfrak n^+$  where $\mathfrak n^+ = \bigoplus_{\alpha > 0} \g_0^{\alpha}$,  and  
$\mathfrak b^- = \mathfrak t \oplus \mathfrak n^-$ where $\mathfrak n^- = \bigoplus_{\alpha > 0} \g_0^{-\alpha}$.
\item[{\rm(f)}] the representation of $\g_0^{(1)}:=[\g_0,\g_0]$ on
$\g_{-1}$ is irreducible and restricted.
\end{itemize}}
\m

Assumption (f) implies that $\g_{-1}$ possesses $\mathfrak b^-$-primitive vector $f \neq 0$.   Thus, there exists $\delta \in \mathfrak t^*$
such that $[t,f] = \delta(t)f$ for all $t \in \mathfrak t$ and $[y,f] = 0$
for all $y \in \mathfrak n^-$.  If $\g$ is transitive \eqref{eq:1.3}, then $\g_k$ is a
restricted module for $\g_0^{(1)}$ for all $k \geq 0$.  In particular,
$\g_1$ will possess a $\mathfrak b^+$-primitive vector $e \neq 0$.  Thus,
there will be $\gamma \in \mathfrak t^*$ such that $[t,e] = \gamma(t)e$ for
all $t \in  \mathfrak t$ and $[x,e] = 0$ for all $x \in \mathfrak n^+$.

\bi
\begin{Lem} \label{Lem:3.18} \ Assume $\g$ is transitive \eqref{eq:1.3}
and $e$ and $f$ are as in Section \ref{sec:3.3}. Then $[e,f] \neq 0$.  \end{Lem}

\pf  Suppose the contrary.  Then since $\g_{-1}$ is
irreducible, it is spanned by vectors of the form $[e_{i_k}
\cdots,[e_{i_2},[e_{i_1},f]] \dots ]$, $k  \geq 0$,
 and

\begin{equation}\left[e,[e_{i_k} \cdots,[e_{i_2},[e_{i_1},f]] \dots ]\right] =
\left [e_{i_k} \cdots,[e_{i_2},[e_{i_1},[e,f]]]\dots\right] =
0.\nonumber
\end{equation}

\noindent By transitivity,  $e = 0$, a
contradiction.  \qed
 
\bi
\begin{Lem} \label{Lem:3.19}  If $f$ and $e$ are as in Section \ref{sec:3.3}
and $\delta = -\gamma$, then $[[e,f],f]$ = $0$ if and only if
$[[e,f],e] = 0$.
\end{Lem}

\pf Clearly, $[[e,f],f] = 0$ if and only if
$\gamma([e,f]) = 0$ if and only if $[[e,f],e] = 0$.  \qed \bi

\begin{Thm}\label{Thm:3.20}   Let $\g$
be a graded Lie algebra satisfying  (a)-(f) of 
Section \ref{sec:3.3},  and assume that $\g$ is  transitive \eqref{eq:1.3}.   Let  $e$
and $f$ be as in
Section \ref{sec:3.3}, 
and suppose that $\delta = -\gamma$.  Then $[[e,f],e] \neq 0$.
\end{Thm}

\pf  We begin by assuming the conclusion is false.   
Observe that then both $[[e,f],f]$ and $[[e,f],e]$ are 0 by Lemma
\ref{Lem:3.19},  and $h = [e,f]$ is nonzero by Lemma \ref{Lem:3.18}, so that $e,f,h$ span a Heisenberg Lie algebra.  If $\alpha_j(h) = 0$ for all simple roots
$\alpha_j$, then $h \in \mathfrak t$ belongs to the center of
$\g_0$.   But then $h$
acts as a nonzero scalar $\zeta$ on $\g_{-1}$ by Schur's Lemma and
transitivity. That would imply
 $[[e,f],f] = \zeta f \neq 0.$     Consequently,
there must exist some simple root $\alpha_k$ such that
$\alpha_k(h) \neq 0$.   We suppose that $e_k$ and $f_k$ are root
vectors corresponding to $\alpha_k$ and  $-\alpha_k$ respectively,  chosen so that
$e_k, f_k, h_k = [e_k,f_k]$ form a canonical basis for
$\mathfrak{sl}_2$.

Let us set $u_0 = -f_k$, $v_0 = e_k$,  and $z = h_k$.  Then the relations
in  \eqref{eq:3.11} hold,

\begin{equation*}
{\begin{array}{cccc}
{[f,u_0]} = 0,&\quad& [h,u_0] = -\lambda u_0, & \\
{[e,v_0]} = 0,&\quad&[h,v_0] = \lambda v_0,& \\
{ [u_0,v_0] = z},&\quad&[z,v_0] = 2v_0, & \qquad[z,u_0] = -2u_0, \\
{[z,f]} = -\eta f, &\quad&[z,e] = \eta e, &
\end{array}}\end{equation*}

\noindent where $\lambda = \alpha_k(h) \neq 0$ and $\eta =
\gamma(h_k)$.  Note that by replacing  $f$ by  $\lambda^{-1}f$ and $h$
by $\lambda^{-1} h$,   we may suppose that $\lambda = 1$.  We claim $\eta
\neq 0$. If $\eta = 0$, then for

\begin{equation*}\begin{array}{ccc}
 \widetilde e_1 = -u_0 = f_k, & \quad \widetilde f_1 = v_0 = e_k, &
\quad \widetilde h_1 = -h_k = [f_k,e_k] = -z \\
 \widetilde e_2 = e, &\quad \widetilde f_2 = f, & \quad \widetilde
h_2 = h = [e,f],   \end{array}\end{equation*}

\noindent  we have the relations

\begin{equation}\label{eq:3.19}[\widetilde e_i,\widetilde f_j] = \delta_{i,j}
\widetilde h_i, \qquad [\widetilde h_i,\widetilde e_j] = a_{i,j}\widetilde e_j,
\qquad [\widetilde h_i,\widetilde f_j] = -a_{i,j}\widetilde f_j,
\end{equation}

\noindent  where $a_{i,j}$ is the $(i,j)$-entry of the matrix

\begin{equation}A = \left(\begin{matrix} 2 & -\eta  \\  -1 & 0  \\ \end{matrix}
\right)\nonumber \end{equation} 

\noindent When $\eta = 0$, then by Theorem  \ref{Thm:3.10}  we would
have that $\g$ is infinite-dimen-sional.  Therefore, $\eta \neq 0$.
Note also that since $\g_{-1}$ is a restricted $\g_0^{(1)}$-module,
$\eta \in \F_p$.

Because the relations in \eqref{eq:3.11} hold, we can now apply Proposition \ref{Pro:3.12}  to get
\eqref{eq:3.12}-\eqref{eq:3.14}:

\begin{eqnarray*} &&\\
&&{\left[E^\ell u_0, \,F^\ell v_0\right]}=
 (-1)^\ell  \ell! \,(z - \ell \eta h),  \nonumber \\
&&{ \left[\bigl[E^\ell u_0, \, F^\ell v_0\bigr], E^\ell
u_0\right]} = 2(-1)^\ell  \ell!\,(\ell\eta -
1)E^\ell u_0,  \quad  \hbox{\rm and} \nonumber\\
&&\left [ \bigl[E^\ell u_0, F^\ell v_0 \bigr],
F^\ell v_0\right] = -2(-1)^\ell \ell!\,(\ell \eta -
1)F^\ell v_0\nonumber \end{eqnarray*}

\noindent for all $\ell \in \{1, \dots, p-1\}$.

    We suppose first that $\eta = 1$.   Here we exploit the fact that now

\begin{equation*}[[E u_0,\,Fv_0], Eu_0]  = -2(\eta - 1)E u_0 = 0, \end{equation*}

\noindent so that $Eu_0$  and  $Fv_0$ span a Heisenberg
algebra.    We set

\begin{equation*}\begin{array}{ccc}
e^\prime= Eu_0,&\quad f^\prime= Fv_0,&\quad
h^\prime= h - z,  \\
u_0^\prime  = E^\ell u_0, &\quad v_0^\prime= cF^\ell v_0, &\quad
z^\prime = (-1)^\ell\,\ell ! c(z - \ell h), \end{array} \end{equation*}

\noindent  for an arbitrary but fixed nonzero $c \in \mathbb F$.
Then  in view of the relations above, we have for $3 \leq \ell \leq p-1$ that \allowbreak

\begin{eqnarray*}
{[f^{\prime},u_0^{\prime}]} &=& \left [ [f,v_0], E^\ell u_0\right]  = 0 \\
{[e^{\prime}, v_0^{\prime}]} &=& 0\\
{[h^{\prime},u_0^{\prime}]}&=& [h -  z, \, E^\ell u_0] = (1 - \ell)u_0^\prime \neq 0\\
{[h^{\prime}, v_0^{\prime}]} &=&(\ell - 1)v_0^\prime \hfil \\
{[ u_0^{\prime}, v_0^{\prime} ]} &=&z^\prime \\
{[ z^{\prime}, u_0^{\prime}]} &=& 2(-1)^\ell\,\ell! \,c(\ell - 1)u_0^{\prime} \\
{[ z^{\prime}, v_0^{\prime}]}&=& -2(-1)^\ell\,\ell! \,c(\ell - 1) v_0^\prime  \\
{[ z^{\prime},f^{\prime}]} &=&(-1)^\ell\,\ell! \,c \left[z - \ell h,\, [f, v_0]\right] =
(-1)^\ell\,\ell!\,c(1-\ell)f^{\prime} \\
{[ z^{\prime},e^{\prime}]} &=& (-1)^\ell\,\ell!\,c(\ell-1)e^\prime .
\end{eqnarray*}

Taking  $c$ so $c^{-1} = (-1)^{\ell-1}\ell!\,(\ell-1)$,   we get
that the primed letters satisfy all the relations in \eqref{eq:3.11} with
$\lambda' = \ell- 1$ and $\eta' = -1 \equiv p - 1 \mod p.$   
Consequently, by replacing  $e,f,h,u_0,v_0,z$
by the corresponding primed letters, we may  
suppose that  $\eta \neq 0,1$.

Then setting  $\ell = \eta^{-1}$, we see that the right side of
  \eqref{eq:3.13}  becomes 0.   The expression in \eqref{eq:3.12}  is
nonzero since $z$ and $h$ are linearly independent due to the fact
$[z,e] = \eta e \neq 0$, while $[h,e] = 0$.   Since $\ell > 1$,
the vector $\widetilde e = E^\ell u_0$ (resp. $\widetilde f =
F^\ell v_0$) is a $\mathfrak b^+$-primitive (resp. $\mathfrak b^-$-primitive) vector. Moreover
the properties that pertain to $e$ and $f$, namely $[\widetilde
e,\widetilde f] \neq 0$ and $[[\widetilde e,\widetilde
f],\widetilde f] = 0 = [[\widetilde e, \widetilde f],\widetilde
e]$ hold for them. Note also that $\widetilde h: =[\widetilde e, \widetilde f]$ is not central.    Indeed, if $\widetilde h \in
\mathfrak Z(\g_0)$, then it acts as a scalar, say $\kappa \neq 0$,
on $\g_{-1}$, and as $\ell \kappa$ on $(\g_{-1})^\ell$.  But then
$[\widetilde h, \widetilde f] = \ell\kappa \widetilde f$, a
contradiction, since $\kappa \neq 0$ and $\ell = \eta^{-1} \neq
0,$ but $[[\widetilde e, \, \widetilde f], \,  \widetilde f] = 0.$
We may continue the argument with this new triple $\widetilde e,
\widetilde f, \widetilde h$ of elements. Since when $f \in
\g_{-i}$, then $\widetilde f \in \g_{-i}^n$ for some $n \in \{2,
\dots, p-1\}$, the elements produced belong to the spaces $\g_t$
with $\mid t\mid$ increasing.
 Thus, $\g$ must be infinite-dimensional.  We have reached a contradiction.  Thus,
$[[e,f],e] \neq 0.$ \qed

\m
\section {\ Subalgebras with a balanced grading \label{sec:3.4}} \m

We  consider  graded Lie algebras $\g = \bigoplus_{j=-q}^r \g_j$   with  $\g_1$  
contragredient (as a $\g_0$-module) to
$\g_{-1}$.   When $\g$
is generated  by  its  local part 
$\g_{-1} \oplus \g_0 \oplus \g_1$, we prove  that $\g$ must have a balanced
gradation (i.e. $q = r$) and then argue that $\g$ must be classical.  For
example, the Lie subalgebra 
of the Melikyan algebra $M(2;\un n)$ generated by the local part 
is a simple classical Lie algebra of type G$_2.$

 \bi  {\it Our assumptions for the remainder of Chapter 3 are the following: \m
 \begin{itemize}
\item[{\rm (1)}] $\g = \bigoplus_{j = -q}^r \g_j$  is a graded Lie
algebra over an algebraically closed field ${\mathbb F}$ of
characteristic $p  > 3$;
\item[{\rm(2)}] $\g_0$ is
classical reductive;
\item[{\rm(3)}] $\g$ is transitive \eqref{eq:1.3} and 1-transitive \eqref{eq:1.4}; 
\item[{\rm(4)}] $\g_{-1}$ is an
irreducible restricted module for $\g_0^{(1)}= [\g_0,\g_0]$; 
\item[{\rm(5)}] $\g_1$ is a  $\g_0$-module isomorphic to the dual module $\g_{-1}^*$.
\item[{\rm(6)}] $\mathfrak t$ is a maximal toral subalgebra of $\g_0$;
$\mathfrak b^+ = \mathfrak t
\oplus \mathfrak n^+$  where $\mathfrak n^+ = \bigoplus_{\alpha > 0} \g_0^{\alpha}$;   and  
$\mathfrak b^- = \mathfrak t \oplus \mathfrak n^-$ where $\mathfrak n^- = \bigoplus_{\alpha > 0} \g_0^{-\alpha}$.
\end{itemize}}
\m

  Let  $f^{\Lambda} \in \g_{-1}$ be a $\mathfrak b^+$-primitive vector of
  weight $\Lambda$  in the irreducible restricted
$\g_0^{(1)}$-module $\g_{-1}$.   We assume $e^{-\Lambda} \in \g_1$ is
a $\mathfrak b^-$-primitive vector of weight $-\Lambda$ in the irreducible $\g_0^{(1)}$-module $\g_1.$
Since the Lie algebra $\g_0$ is classical reductive, it has a
presentation by generators $e_i,f_i, i \in {\tt I}$, $h_j, j\in {\tt J}$, and Serre
relations, (where
${\tt I =J}$ if no ideal summands of the form $\mathfrak Z(\g_0)$, or
$\mathfrak{gl}_n$, $\mathfrak {sl}_n$ or $\mathfrak{pgl}_n$ with $p \mid n$ occur).
There is an
automorphism $\sigma$ on $\g_0$ such
that $\sigma(f_i) = e_i$, $\sigma(e_i) = f_i$, $i \in {\tt I}$, and
$\sigma(h_j) = -h_j$, $j \in {\tt J}$.
\m

Since the $\g_0$-modules $\g_1$ and
$\g_{-1}$ are dual to one another,  there is a
$\g_0$-homomorphism $\phi: \g_{-1} \otimes \g_1 \rightarrow \g_0$
given by $\phi(u \otimes v) = [u,v]$.   Let ${\mathfrak F}^\pm$ be
the free Lie algebra generated by $\g_{\pm 1}$,  and set ${\mathfrak
F} = {\mathfrak F}^- \oplus \mathfrak F_0 \oplus {\mathfrak F}^+$,
where ${\mathfrak F}_0 = \g_0.$ Then ${\mathfrak F}$ is a Lie
algebra (compare \cite{BKM}),  and it is generated by the elements
$e_i,f_i, i \in {\tt I}$, $h_j, j \in {\tt J}$, $f^\Lambda$, $e^{-\Lambda}.$
We may assign a grading to ${\mathfrak F}$ by setting 0 =
\hbox{deg}($e_i$) = \hbox{deg}($f_i$) = \hbox{deg}($h_j$) for all
$i,j$ and deg$(f^\Lambda) = -1 = -$deg$(e^{-\Lambda})$ so that
${\mathfrak F}^\pm = \sum_{t = 1}^\infty {\mathfrak F}_{\pm t}$
where ${\mathfrak F}_{\pm t}$ is the span of the elements of
degree $\pm t$. We claim we can extend the automorphism $\sigma$
on $\g_0 = {\mathfrak F}_0$ to one on ${\mathfrak F}$. Indeed,
first define a new $\g_0$-module structure on
 $\g_{-1}$  by $x \ast v = [\sigma(x),v]$.  Then $f_i \ast f^\Lambda
= [e_i,f^\Lambda] = 0$ and $h_j \ast f^\Lambda = -[h_j,f^\Lambda]
= -\Lambda(h_j)f^\Lambda$.  It follows from these calculations
that $\g_{-1}$ with this new $\g_0$-action is an irreducible
$\g_0$-module with a $\mathfrak b^-$-primitive vector of weight $-\Lambda$.  Since such a module
is unique up to isomorphism, we can define a $\g_0$-module
isomorphism  $\sigma: \g_{-1} \rightarrow \g_1$ so that
$\sigma(f^\Lambda) = e^{-\Lambda}$.    The same argument shows that
$\g_1$ under the $\g_0$-action  $x \ast u = [\sigma(x), u]$
is a $\g_0$-module with a $\mathfrak b^+$-primitive vector
of weight $\Lambda$.  Thus we can define a $\g_0$-module
isomorphism  $\sigma: \g_1 \rightarrow \g_{-1}$ so that
$\sigma(e^{-\Lambda}) = f^\Lambda$.  Note then that

\begin{equation}\label{eq:3.24}\sigma([x,w]) = \sigma(\sigma(x) \ast w) =
[\sigma(x),\sigma(w)] \end{equation}

\noindent for all $x \in \g_0$ and $w \in \g_{-1}$ or $w \in
\g_1$.  We may now extend $\sigma$ to be an automorphism on the
free Lie algebras $\mathfrak F^+$ and $\mathfrak F^{-}$ so that
\eqref{eq:3.24} holds for all $x \in {\mathfrak F}_0 = \g_0$ and
$w \in \mathfrak F^+$ or $w \in \mathfrak F^-$.

Let $\widehat \g$ denote the Lie subalgebra of $\g$ generated by the local
part $\g_{-1} \oplus \g_0 \oplus \g_1$. There is a natural
epimorphism $\phi: {\mathfrak F} \rightarrow \widehat \g$ taking
the generators to the corresponding elements in $\widehat \g$. The
kernel is a homogeneous ideal,  ${\mathfrak K} = \bigoplus_{t \in \Z}
{\mathfrak K}_t$  of $\mathfrak F$, which trivially intersects the
local part ${\mathfrak F}_{-1} \oplus {\mathfrak F}_0 \oplus
{\mathfrak F}_1$.

Suppose $m \geq 2$ is minimal with $\mathfrak K_m \neq 0$.  Let
$y \in \mathfrak K_m$. Then $[y,\mathfrak F_{-1}] = 0$ by the
minimality of $m$.   Hence $[\sigma(y),\mathfrak F_1] = 0$, and
applying $\phi$ we have $[\phi(\sigma(y)),\g_1] = 0$.  By the
1-transitivity  of $\g$, we must have
$\phi(\sigma(y)) = 0$;  that is, $\sigma(y) \in \mathfrak K_{-m}$. Thus,
$\sigma: \mathfrak K_m \rightarrow \mathfrak K_{-m}$. Arguing with
the minimal $n$ such that $\mathfrak K_{-n} \neq 0$ shows that
$\sigma: \mathfrak K_{-n} \rightarrow \mathfrak K_n$.  Thus, these
two minimal values $m,n$ must be equal, and $\sigma$ must
interchange the subspaces $\mathfrak K_m$ and $\mathfrak K_{-m}$.
Suppose we have shown that $\sigma$ interchanges $\mathfrak K_s$
and $\mathfrak K_{-s}$ for all $2 \leq m \leq s < t$. Then
$[\mathfrak K_t, \mathfrak F_{-1}] \subseteq \mathfrak K_{t-1}$,
and $[\sigma(\mathfrak K_t), \mathfrak F_1] \subseteq \mathfrak
K_{-(t-1)}$. Applying $\phi$ gives $[\phi(\sigma(\mathfrak K_t)),
\g_1] \subseteq \phi(\mathfrak K_{-(t-1)}) = 0$.  Once again by
1-transitivity, it must be that $\sigma(\mathfrak
K_t) \subseteq \mathfrak K_{-t}$. Replacing $t$ with $-t$ and $-1$
with $1$ in the argument gives $\sigma(\mathfrak K_{-t}) \subseteq
\mathfrak K_t$. Since $\sigma$ is an automorphism, equality must
hold.  Thus, $\mathfrak K$ is $\sigma$-invariant, so that there is
an induced automorphism $\sigma$ on ${\mathfrak F}/{\mathfrak K}
\cong \widehat \g$. This says that the gradation of $\widehat \g$
is balanced. We now summarize what we have just shown.

\bi
\begin{Thm} \label{Thm:3.24} Assume $\g = \bigoplus_{j = -q}^r
\g_j$  is a graded Lie algebra satisfying (1)-(5) of Section \ref{sec:3.4}. 
Then the gradation of the subalgebra  $\widehat \g$
of $\g$ generated by the local part $\g_{-1}\oplus \g_0 \oplus
\g_1$ is balanced.
\end{Thm}

\m We now come to the first main result of this section.

\bi

\begin{Thm} \label{Thm:3.25} \ Assume $\g = \bigoplus_{j = -q}^r
\g_j$  is a graded Lie algebra satisfying (1)-(5) of Section \ref{sec:3.4}. Then the graded
Lie subalgebra $\widehat{\mathfrak{g}}$ generated by the local
part $\g_{-1} \oplus \g_0 \oplus \g_1$ of $\g$ is either classical
simple or isomorphic to $\mathfrak{pgl}_{kp}$ for some
$k\in\mathbb N$.
\end{Thm}

\pf  The Lie algebra $\widehat{\mathfrak{g}}$ generated
by $\g_{-1} \oplus \g_0 \oplus \g_1$ satisfies all the hypotheses
of the theorem, so we will assume from now that $\g=\widehat{\g}$.
Let $S\eqdef\bigoplus_{i\ge 0}\,(\ad \g_1)^i\g_{-q}$. Since $\g$
is generated by its local part, Lemma \ref{Lem:1.20} shows that
$S$ is the unique minimal ideal of $\g$. This, in turn, implies
that $\text{Ann}_{\g}\,S=0$. In particular, the center ${\mathfrak Z}(S)=0$
and $\g$ can be regarded as a Lie subalgebra of $\Der(S)$. Let $G$
denote the connected component of the automorphism group of the Lie algebra $S$, and let
$\mathcal L$ denote the Lie algebra of the algebraic group $G$.
Then $\mathcal L$ is a Lie subalgebra of $\Der(S)$. \m

\noindent (a) Suppose we have established that $\ad S\subseteq
\mathcal L$. Then $\ad S$ is an ideal of $\mathcal L\subseteq
\Der(S)$. Let $R$ denote the solvable radical of the algebraic
group $G$ and ${\mathcal R}= \text{Lie}(R)$, a solvable subalgebra of
$\Der(S)$. Since $\mathcal R$ is an ideal of $\mathcal L$, the
subspace $[{\mathcal R},\ad S]$ is a solvable ideal of $\ad S$. Let
$\mathfrak r$ denote the inverse image of $[{\mathcal R},\ad S]$
under the isomorphism
$\ad\colon\,S\stackrel{\sim}{\longrightarrow} \ad S$. \m

Suppose ${\mathcal R}\ne 0$. Then $[{\mathcal R},\ad S]\ne 0$ and
hence ${\mathfrak r}\cap\g_{-1}\ne 0$ by the transitivity of $\g$.
The irreducibility of $\g_{-1}$ shows that the $\g_0$-module
$\{x\in \g_{-1}\mid [x,\g_1]=0\}$ is zero. This, in turn, shows
that $J\eqdef [{\mathfrak r}\cap \g_{-1}, \g_1]$ is a nonzero
solvable ideal of $[\g_{-1},\g_1]$. Since $\g_0$ acts irreducibly
and faithfully on $\g_{-1}$, Schur's lemma implies that ${\mathfrak
Z}(\g_0)$ acts on $\g_{-1}$ as scalar operators, so that
$\dim\,{\mathfrak Z}({\g}_0)\le 1$ (see the discussion in Section \ref {sec:1.10}).   Lemma \ref{Lem:1.79} now says
that $J\subseteq {\mathfrak Z}(\g_0)$ (one should keep in mind
that $[\g_{-1},\g_1]$ is an ideal of $\g_0$). Since $J\ne 0$,
there is $z\in\mathcal R$ such that $z(x)=x$ for all
$x\in\g_{-1}$. Since $\g_1$ is isomorphic to a $\g_0$-submodule of
$\text{Hom}(\g_{-1},\g_0)$, we have $z(x)=-x$ for all $x\in\g_1$.
But then $\ad \g_{\pm1}\subset\mathcal R$, contradicting the
solvability of $\mathcal R$. We thus deduce that ${\mathcal R}=0$. Consequently, $G$ is a semisimple (connected) algebraic
group.

\m

Because the center  $Z(G)$ of the group $G$ acts trivially on ${\mathcal L}$,
it also acts trivially on $\text{ad}\,S\subseteq \mathcal L$.
Since ${\mathfrak Z}(S)=0$, it follows
that the scheme $Z(G)$ is trivial, i.e., $G$ is
a semisimple group of adjoint type. It follows that $G$ is
isomorphic (as an algebraic group) to the direct product of its
simple components $G_1,\ldots,G_d$. As a consequence, ${\mathcal
L}={\mathcal L}_1\oplus\dots \oplus {\mathcal L}_d$ where
${\mathcal L}_i= \text{Lie}(G_i)$ (a direct sum of Lie algebras).
Since $[{\mathcal L}_i, \ad\,S]\subseteq \ad\,S$ is a  {\it
nonzero} ideal of $\ad\, S$ and ${\mathcal L}_i\cap {\mathcal
L}_j=0$ for $i\ne j$, the simplicity of $S$ forces $d=1$. In other
words, $G$ is a simple algebraic group of adjoint type. But then
either $\mathcal L$ is a classical simple Lie algebra or
${\mathcal L}\cong \mathfrak{pgl}_{kp}$ for some $k\in\mathbb N$.
Since $\ad\, S$ is a simple ideal of $\mathcal L$, we derive that
$S$ is classical simple. If $S\not\cong\mathfrak{psl}_{kp}$, then
$\ad S=\Der(S)$,  hence $\g\hookrightarrow\Der(S)$ is classical
simple. If $S\cong\mathfrak{psl}_{kp}$ for some $k\in\mathbb N$,
then $\ad S$ has codimension $1$ in $\Der(S)\cong
\mathfrak{pgl}_{kp}$. Therefore, in any event either $\g$ is
classical simple  or $\g\cong\mathfrak{pgl}_{kp}$, as wanted.

\bi

\noindent (b) By Theorem \ref{Thm:3.24}, $r=q$. The argument used
in the proof of Lemma \ref{Lem:2.7} shows that $\exp\, (t \ad
x)\in G$ for all $x\in S$ such that $(\ad x)^{(p+1)/2}=0$ and all
$t\in\mathbb{F}$. From this it follows that for all such $x$ we
have that $\ad\,x\in \mathcal L$. Since $p>3$ and $(\ad\,\g_{\pm
q})^3(S)=0$, this implies that $\ad\,\g_{\pm q}\subset \mathcal
L$.

\m

If $q=1$, then
$$\ad S=\ad \g_{-1}\oplus [\ad \g_{-1},\ad \g_1]\oplus
[[\ad \g_{-1},\ad \g_1],\ad \g_1] \subseteq {\mathcal L},$$ and
the statement follows from part (a) of this proof. Thus in what
follows we may assume that $q\ge 2$.

Suppose $\ad \g_{q-1}\subset \mathcal L$. Then $\ad \g_{-1} =[\ad
\g_{-q}, \g_{q-1}]\subset \mathcal L$,  and hence $\ad S \subseteq
\mathcal L$ thanks to Lemmas \ref{Lem:1.20} and \ref{Lem:1.55}.

\bi

\noindent (c) Suppose $q\ge 3$. If $p>5$, then
$$(\ad \g_{q-1})^{(p+1)/2}(\g)\subseteq (\ad \g_{q-1})^{4}(S)\subseteq
\g_{3q-4}=0,$$ which shows that $\ad \g_{q-1}\subset {\mathcal
L}$. Then our remarks in part~(b) imply that $\g$ is either
classical simple or is isomorphic to $\mathfrak{pgl}_{kp}$. If
$p=5$ and $q\ge 4$, then $(\ad \g_{q-1})^{(p+1)/2}(S)\subseteq
\g_{2q-3}=0$. Again the result follows by part~(b).

Suppose $p=5$ and $q=3$.  Given $x\in \g_{q-1} = \g_2$ set $X\eqdef \ad x$.
Clearly, $X^4=0$. As $X^2(S)\subset \g_1\,\oplus\, \g_2\,\oplus\,
\g_3$ and $X^3(S)\subset \g_3$, it must be that
$[X^3(S),X^2(S)]=0.$ By (the proof of) Lemma \ref{Lem:2.7}, $\exp
tX \in G$ for all $t\in\mathbb F$. But then again $\ad
\g_{q-1}\subset {\mathcal L}$ and we are done by part~(b).

Suppose $p>5$ and $q=2$. For $x\in \g_1$ we now have
$X^4(S)\subseteq \g_2$ and $X^3(S)\subset \g_1\oplus\g_2$. Also,
$X^5=0$. Hence $\exp tX \in G$ for all $t\in\F$ provided that
$p>7$. If $p=7$, then $[X^i(S), X^{p-i}(S)]\subseteq [\g_1\oplus
\g_2,\g_2]=0$, and again $\exp tX \in G$ for all $t\in\F$. Thus
this case can be argued as in part~(b).

\bi

\noindent (d) Finally, suppose $q = 2$ and $p=5$. This case is
{\it much} more complicated. We now have
\begin{equation}\g
= \g_{-2} \oplus \g_{-1} \oplus \g_0 \oplus \g_1 \oplus
\g_2.\nonumber
\end{equation}
Recall that $\g$ is irreducible and transitive. Therefore,
$\g_{-2}$ is an irreducible $\g_0$-module and
$[[\g_{-2},\g_2],\g_{-2}]=\g_{-2}$; see Lemmas \ref{Lem:1.51}\,(b) and
\ref{Lem:1.58}.

Let $A_0 \eqdef \Ann_{\g_0}(\g_{-2})$ and $A_2\eqdef \{x\in
\g_2\mid [x,\g_{-2}]\subseteq A_0\}.$ By Lemma \ref{Lem:1.57}\,(i), the
set $\{\ad_{\g_{-1}\,}[u,v]\,\vert \, u\in\g_{-2},\,v\in A_2\}$ is
weakly closed in $\text{End}(\g_{-1})$. Since
$$0=(\ad\,u)^2(\ad v)^2 \g_{-1}=2(\ad [u,v])^2\g_{-1}$$ for all $u\in \g_{-2},\, v\in A_2$,
Lemma \ref{Lem:1.57}\,(ii)  shows that $[\g_{-2},A_2]$ annihilates
$\g_{-1}$. Since $\g$ is transitive, it must be that
$[\g_{-2},A_2]=0$. But then $A_2\subseteq \text{Ann}_{\g}\,S=0.$
Since $[A_0,\g_2]\subseteq A_2$, it follows that $A_0$ is an ideal
of the Lie algebra $\g_{-2}\oplus\g_0\oplus \g_2$.

Let ${\mathcal B}\eqdef \bar{\g}_{-1}\oplus \bar{\g}_0\oplus
\bar{\g}_1$,  where $\bar{\g}_i=(\g_{2i}+A_0)/A_0$. The graded Lie
algebra $\mathcal B$ is irreducible and transitive, but  a priori it
is not clear whether it is $1$-transitive. Nevertheless,  we will
proceed as in part (a) of this proof.

Let $I=\bar{\g}_{-1}\oplus [\bar{\g}_{-1},\bar{\g}_1]\oplus
\bar{\g}_1.$ Since any nonzero ideal of $\mathcal B$ contains
$\bar{\g}_{-1}$, we have ${\mathfrak Z}(I)=0$. Let $\bar{G}$ denote
the connected component of the automorphism group of $I$ and
$\bar{R}$ the solvable radical of the algebraic group $\bar{G}$.
Let $\bar{\mathcal R}:= \text{Lie}(\bar R)$, a solvable subalgebra of
$\Der(I)$. Since $(\ad\bar{\g}_{\pm 1})^3(I)=0$, we have
$\ad I\subseteq \text{Lie}(\bar{G}).$ Hence $\ad I$ is an ideal of
$\text{Lie}(\bar{G})$. Since $\bar{\mathcal R}$ is an ideal of
$\text{Lie}(\bar{G})$, we deduce that $[\bar{{\mathcal R}},\ad I]$ is
a solvable ideal of $\ad I$. Let $\bar{\mathfrak r}$ denote the
inverse image of $[{\bar{\mathcal R}},I]$ under the isomorphism
$\ad\colon\,I\stackrel{\sim}{\longrightarrow} \ad I$. \m

Suppose $\bar{{\mathcal R}}\ne 0$. Then $[\bar{{\mathcal R}},\ad
I]\ne 0$ and hence $\bar{\mathfrak r}\cap\bar{\g}_{-1}\ne 0$ by the
transitivity of $\mathcal B$. The irreducibility of
$\bar{\g}_{-1}$ shows that
$\text{Ann}_{\bar{\g}_{-1}}\,\bar{\g}_1=0$. This, in turn, implies
that $J\eqdef [\bar{\mathfrak r}\cap \bar{\g}_{-1}, \bar{\g}_1]$ is a nonzero solvable ideal of
$[\bar{\g}_{-1},\bar{\g}_1]$.  Lemma \ref{Lem:1.79} says that
$\bar{\g}_0$ is classical reductive and $J\subseteq {\mathfrak
Z}(\bar{\g}_0)$. Arguing as in part (a) of this proof we  reach a
contradiction. Thus $\bar{G}$ is a semisimple algebraic group.

Since the center  $Z(\bar{G})$ acts trivially on
$\text{Lie}(\bar{G})$, it also acts trivially on
$\text{ad}\,I\subseteq \text{Lie}(\bar{G})$. Since ${\mathfrak
Z}(I)=0$, this yields that $\bar{G}$ is a group of adjoint type.
The structure theory of  algebraic $\mathbb F$-groups now shows
that $\text{Lie}(\bar{G)}$ is classical reductive with ${\mathfrak
Z}(\text{Lie}(\bar{G}))=0$. Since any nonzero ideal of $I$
contains $\bar{\g}_{-1}$ we can argue as in part (a) to deduce
that $\bar{G}$ is a simple algebraic group. But then $I\cong \ad
I$ is classical reductive, too. Moreover, its derived subalgebra
$I^{(1)}$ is classical simple. Since
$I\hookrightarrow \Der(I^{(1)})$, the inequality $I\ne
I^{(1)}$ would imply that $I^{(1)}\cong\mathfrak{psl}_{kp}$ for
some $k\in\mathbb N$ and $I\cong\mathfrak{pgl}_{kp}\cong\Der
\big(\mathfrak{psl}_{kp}\big)$. But then all {\it nilpotent} derivations
of $I^{(1)}$ would be inner, forcing $\bar{\g}_{1}\subset
I^{(1)}$. This shows that $I$ is classical simple. Since $\mathcal
B$ is irreducible and transitive, it embeds into $\Der(I)$.
Therefore, either ${\mathcal B}=I$ or ${\mathcal B}\cong
\mathfrak{pgl}_{kp}$ and $I={\mathcal B}^{(1)}$. In any event, the
simple algebraic group $\bar{G}$ has the property that
$\text{Lie}(\bar{G})=\Der(I)$.

\bi \noindent (e) We now turn our attention to the Lie algebra
$\widetilde{\mathcal B}\eqdef \g_{-2}\oplus\g_0\oplus\g_2$ and its
ideal $\widetilde{I}\eqdef \g_{-2}\oplus[\g_{-2},\g_2]\oplus\g_2$. Let
$\varphi\colon\,\widetilde{\mathcal B}\twoheadrightarrow \mathcal B=
\widetilde{\mathcal B}/A_0$ denote the canonical homomorphism. Since $A_0$
annihilates  $\g_2$ and $\g_{-2}$, it commutes with $\widetilde{I}$.
So $\widetilde{I}\cap A_0$ is an abelian ideal of $\g_0$. Since $\g_0$
is classical reductive, $\widetilde{I}\cap A_0\subseteq {\mathfrak
Z}(\g_0)$. Since ${\mathfrak Z}(\g_0)$ acts faithfully on
$\g_{-2}=[\g_{-1},\g_{-1}]$, we get $\widetilde{I}\cap A_0=0$.
Therefore, $\widetilde{I}\cong \varphi(\widetilde{I})=I$, and we can
identify $\Der(\widetilde{I)}$ with $\text{Lie}(\bar{G})$. Since
${\mathfrak Z}(\Der(\widetilde{I}))=0$, the natural $p$th power
map on $\Der(\widetilde{I})$ then gets identified with the
canonical $(\text{Ad}\,\bar{G})$-equivariant $p$th power map on
$\text{Lie}(\bar{G})$.

\m

The adjoint action of $\widetilde{\mathcal B}$ on $\widetilde{I}$ induces
a natural Lie algebra homomorphism $\psi\colon \widetilde{\mathcal
B}\rightarrow \Der(\widetilde{I})=\text{Lie}(\bar{G})$. Since
$\g_0$ is classical reductive, it contains a diagonalizable Cartan
subalgebra, $\h$ say. We denote by $\mathcal H$ the image of $\h$
under $\psi$. Since $\g_{-1}$ is an irreducible $\g_0$-module,
$\h$ acts diagonalizably on $\g_{-1}$; see Lemma \ref{Lem:1.73}.
Therefore, $\h$ acts diagonalizably on
$\g_{-1}\hookrightarrow\text{Hom}(\g_{-1},\g_{0})$. Since $\g$ is
generated by its local part, $\h$ acts diagonalizably on $\g$. As
a consequence, $\mathcal H$ is a diagonalizable subalgebra of
$\text{Lie}(\bar{G})$, i.e., consists of pairwise commuting
$[p]$-semisimple elements of $\text{Lie}(\bar{G})$. Let $\mathcal C$ 
(respectively, ${\mathcal C}'$) denote the centralizer of
$\mathcal H$ in $\Der(\widetilde{I})$ (respectively, in
$\ad\widetilde{I}$). Clearly, ${\mathcal C}'$ is an ideal of
codimension $\le 1$ in $\mathcal C$. It is easy to see that
$\mathcal C'$ is isomorphic to the centralizer of $\h$ in
$\widetilde{I}$. The group $\bar{G}$ acts on its Lie algebra
$\text{Lie}(\bar{G})=\Der(\widetilde{I})$ via the adjoint
representation $\text{Ad}$, and we have $(\text{Ad}\,g)(D)=g\circ
D\circ g^{-1}$ for all $g\in\bar{G}$ and
$D\in\Der(\widetilde{I})$. Let $$Z_{\bar{G}}(\mathcal H)\eqdef
\{x\in\bar{G}\mid (\text{Ad}\,g)(x)=x \  \, \mbox{for all}\ \,
x\in\mathcal H\},$$  the centralizer of $\mathcal H$ in $\bar{G}$.
Thanks to \cite[Sec.~9.1]{Bo} we have ${\mathcal C}=\text{Lie}(Z_{\bar{G}}(\mathcal H))$, while 
\cite[Sec.~13.19]{Bo} says that the connected component of
$Z_{\bar{G}}(\mathcal H)$ is reductive. As a consequence,
$\mathcal C$ is isomorphic to a Lie algebra of reductive group.
Let
$$\mathfrak c'\eqdef \g_{-2}^\h\oplus \big(\h\cap [\g_{-2},\g_2]\big)\oplus
\g_{2}^\h$$ where ${\g}_{\pm 2}^\h=\{x\in\g_{\pm
2} \mid [x,\h]=0\}$. Then $\mathfrak c'\cong\mathcal C'$ by our earlier
remarks. The set
$$\big\{(\ad x)\vert_{\mathfrak c'}\mid x\in\g_{-2}^\h\cup\big(\h\cap[\g_{-2},\g_2]\big)\cup
\g_2^\h\big\}$$ is weakly closed and consists of nilpotent
endomorphisms of $\mathfrak c'$.  Applying \cite[p. 33]{J} (compare Lemma \ref{Lem:1.57}\,(ii)),  we now deduce
that the Lie algebra ${\mathcal C}'$ is nilpotent. Since
$\mathcal C'$ is an ideal of codimension $\le 1$ in $\mathcal 
C$, the Lie algebra $\mathcal C=\text{Lie}(Z_{\bar{G}}(\mathcal
H))$ is solvable.  Since the group $Z_{\bar{G}}(\mathcal H)$ is
reductive, the connected component $Z_{\bar{G}}(\mathcal H)^\circ$ must be a torus. 
 It follows that $\mathcal C$ is
abelian and consists of $[p]$-semisimple elements of
$\text{Lie}(\bar{G})$. This, in turn, implies that $\mathfrak c'$ is a
diagonalizable subalgebra of $\widetilde{I}$, yielding $\g_{\pm
2}^\h=0$. 
\m

Given $t\in{\mathbb F}^\times$ define
$\lambda(t)\in\text{Aut}\,\widetilde{I}$ by setting
$\lambda(t)(x)=t^ix$ for all $x\in \g_i\cap \widetilde{I}$. Since
$\lambda({\mathbb F}^\times)$ is a one-dimensional torus,
$\lambda(\mathbb F^\times)$ is contained in $\bar{G}$. Since
$\mathcal H$ preserves the graded components of $\widetilde{I}$, we
have $\lambda(\mathbb F^\times)\subseteq Z_{\bar{G}}(\mathcal
H)^\circ$. Also,
$$\text{Lie}(\bar{G})=\text{Lie}(\bar{G})_{-2}\oplus\text{Lie}(\bar{G})_0\oplus
\text{Lie}({\bar{G}})_2$$ where
$\text{Lie}(\bar{G})_i=\{x\in\text{Lie}({\bar{G}}) \mid 
\big(\text{Ad}\,\lambda(t)\big)\,x=t^i x\ \ \mbox{for all}\ \
t\in\mathbb F^\times\}$. Let $Z_{\bar{G}}(\lambda)$ denote the
centralizer of $\lambda(\mathbb F^\times)$ in $\bar{G}$. By
\cite[Sec.~8.4]{Sp}, there exists a parabolic subgroup $P$ of
$\bar{G}$ with Levy decomposition $P=Z_{\bar{G}}(\lambda)\,R_u(P)$  (where
$R_u(P)$ is the unipotent radical of $P$) 
such that $\text{Lie}(R_u(P))=\text{Lie}(\bar{G})_2$. Let $T$ be a
maximal torus of $\bar{G}$ containing $Z_{\bar{G}}({\mathcal
H})^\circ$. Then $T\subset Z_{\bar{G}}(\lambda)$ and $${\mathcal
H}\subseteq{\mathcal C}= \text{Lie}(Z_{\bar{G}}(\mathcal
H))^\circ\subseteq \text{Lie}(T).$$    Let $B$ be a Borel subgroup of
$P$ containing $T$. Then $B$ is a Borel subgroup of $\bar{G}$ and
$R_u(P)$ is a normal subgroup of $B$. Therefore,
$\text{Lie}(R_u(P))=\text{Lie}({\bar{G}})_2$ is
$(\text{Ad}\,B)$-stable.    But then $\ad\g_2=(\ad\widetilde{I})\cap
\text{Lie}({\bar{G}})_2$ is $(\text{Ad}\,B)$-stable as well. It
follows that $\g_2$ contains a nonzero element fixed by the
derived subgroup $(B,B)$  of $B$, say $e$.

\m

Since $\widetilde{I}$ is an irreducible $\bar{G}$-module, the subspace
of $(B,B)$-invariants in $\widetilde{I}$ is one-dimensional, hence
equals $\F e$. Since $B=T\,(B,B)$, the line ${\mathbb F}\, \ad e$
is a root space for $\text{Ad}\,T$. Let $\Psi$ denote the root
system of $\bar{G}$ relative to $T$ and $\Psi_+(B)$  the positive
system of $\Psi$ corresponding to $B$. Note that all root spaces of
$\text{Lie}(\bar{G})$ lie in $\ad \widetilde{I}$ (for
$\text{Lie}({\bar{G}})/\ad \widetilde{I}$ is a trivial
$\bar{G}$-module). Let $\gamma$ denote the root of $\ad\,e$. It is
immediate from the above discussion that $\gamma$ is the highest
root in $\Psi_+(B)$. As $\gamma\big(\lambda(t)\big)=t^2$ for all
$t\in\F^\times$, there is a nonzero $f\in \g_{-2}$ such that $ \ad
f$ is a root vector corresponding to $-\gamma$.  Since ${\mathcal
H}\subseteq\text{Lie}(T)$, there is a linear function $\eta$ on
$\h$ such that  $$[h,e]=\eta(h)e,\qquad [h,f]=-\eta(h)f
\qquad(\forall\,h\in\h)$$ (to be more precise,
$\eta(h)=((\text{d}\gamma)_e\circ \psi)(h)$ where
$(\text{d}\gamma)_e$ is the differential of the morphism
$\gamma\colon\,T\rightarrow \F^\times$ at the identity element of
$T$). Since $\ad e$ and $\ad f$ belong to the Lie algebra of a
semisimple subgroup of rank one in $\bar{G}$, they generate an
$\mathfrak{sl}_2$-triple in $\text{Lie}(\bar{G})$. Let $h=[e,f]$.
No generality is be lost by assuming that $\{e,h,f\}$ is a
standard basis of ${\mathfrak s}\eqdef \F f\oplus\F h\oplus\F
e\cong \mathfrak{sl}_2$. Since  all root spaces of $\text{Ad}\,T$
are one-dimensional and $2\gamma+\alpha > \gamma$ for all
$\alpha\in\Psi\setminus\{-\gamma\}$, we have the equalities $(\ad
e)^2(\widetilde{I})=\F e$ and $(\ad f)^2(\widetilde{I})=\F f$.

\m

\noindent (f) Let  $E$, $H$, and $F$ denote the adjoint
mappings on \ $S= \bigoplus_{i\ge 0}\,(\ad \g_1)^i\g_{-q}$ \   $(q=2)$ corresponding to $e$, $h$, and $f$,
respectively.   Then $S$ is an $\mathfrak{s}$-module with
$E^3=F^3=0$ so that  $$E^{p-1}(H + 1) = 0 = (H + 1)F^{p-1}.$$   Therefore, 
\cite[Sec.~5.2]{S} implies that the $\mathfrak s$-module $S$ is
completely reducible with composition factors $L(0)$, $L(1)$,
$L(2)$, where $L(k) = L(k\varpi_1)$.   In particular, all eigenvalues of $H$ 
are in  $\F_p$. Let
$V_k$ denote the eigenspace of $H$ in $S$  corresponding to $k\in\F_p$.
Since $E^2(S)=(\ad e)^2(\widetilde{I})=\F e$ by part (e),  and $\dim
V_2=\dim\,E^2(S) $ equals the multiplicity of $L(2)$ in $S$, we
have $V_2=\F e$. Similarly, $V_{-2}=\F f$.     As $p=5$,   we have
$E(V_1) \subseteq V_{-2}$  and $F(V_{-1})\subseteq V_2$.   The
equalities $E^3=F^3=0$ now yield $[V_1,V_2]=[V_{-1}, V_{-2}]=0$.
But then the decomposition
 \begin{equation} S=V_{-2} \oplus V_{-1} \oplus V_0 \oplus V_1
\oplus V_2\nonumber\end{equation} is a $\mathbb Z$-grading of $S$.
Since $h\in \h$, all graded components $V_i$ are $\h$-stable.

As explained in part (e), the subalgebra $\h$ acts diagonalizably
on $S$.  Thus,  any $\h$-stable subspace $M$ of $S$ decomposes as
$M=\bigoplus_{\alpha\in\h^*}\, M^\alpha$ where $M^\alpha=\{x\in
M\mid [h,x]=\alpha(h)x\ \,\mbox{for all}\ \,h\in \h\}$. By our
remarks in part (e), $V_2=V_2^\eta$ and $V_{-2}=V_{-2}^{-\eta}$.

\m

\noindent (g) The graded component $V_1$ is spanned by the set
${\mathcal X}\eqdef \bigcup_{i\in{\mathbb Z},\,\,\alpha(h)=1}\,
\g_i^{\alpha}.$ We claim that $(\ad x)^4=0$ for all $x\in\mathcal
X$. If $x\in  \g_i^{\alpha}$ and $\alpha\ne \eta/2$, then
$$(\ad x)^4(S)=(\ad x)^4(V_{-2}^{-\eta})\subseteq
V_2^{-\eta-\alpha}=0,$$ so our claim is valid in this case. Assume
for a contradiction that there is a $y\in \g_s^{\eta/2}$,  for some
$s\in\Z$, such that $(\ad y)^4\ne 0$. Notice that
$\g^{-\eta}=\g_{-2}^{-\eta}= \F f$ and $\g^\eta=\g_2^{\eta}=\F e$.
Since $y\in V_1$, we have $(\ad y)^4(S)= (\ad y)^4(\F f)\subseteq
\g_2$. As $f\in\g_{-2}$, it must be that $s=1$. We now consider
$$\g(\eta):=\g_{-2}^{-\eta}\oplus \g_{-1}^{-\eta/2}\oplus \h\oplus
\g_{1}^{\eta/2}\oplus \g_{2}^\eta,$$ a graded subalgebra of $\g$
(one should take into account that
$[\g_{-1}^{\eta/2},\,\g_1^{-\eta/2}]\subseteq \h$ as $\h$ is
self-centralizing in $\g_0$). Let $\h_0$ denote the kernel of
$\eta$, a subspace in ${\mathfrak Z}\big(\g(\eta)\big)$. Consider
the pairing $$\xi\colon\, \g_{-1}^{-\eta/2}\times
\g_1^{\eta/2}\longrightarrow \,\,\h/\h_0,\quad \ \
(x,y)\longmapsto \,\,[x,y]\ \,(\mbox{mod}\ \h_0).$$ Let $W_{i}$
denote the  orthogonal subspace to $\xi$ in $\g_{i}^{i\eta/2}$, $\,i=-1,1$.
By our assumption, $y\not\in W_1$.    If $W_1$ had codimension $\ge
2$ in $\g_{1}^{\eta/2}$, we would be able to find
$a_1,a_2\in\g_1^{\eta/2}$ and $b_1,b_2\in \g_{-1}^{-\eta/2}$ such
that $\xi(a_i,b_j)=\delta_{ij}$ for all $i,j\in \{1,2\}$. But then
the images of $a_1, a_2, \,b_1, b_2$ and $h$ would generate an
infinite-dimensional  subalgebra of $\g(\eta)/\h_0$; see Theorem
\ref{Thm:3.9}. Therefore, $W_i$ has codimension one in $\g_i,\,$
$i=-1,1$. Since $E(\g_1)=F(\g_{-1})=0$ and $\eta(h)=2$, the  map
$E$ induces a linear isomorphism between $\g_{-1}^{-\eta/2}$ and
$\g_{1}^{\eta/2}$.   On the other hand, our assumptions on $y$ and
$s$ imply that $e\in[\g_1^{\eta/2},\g_1^{\eta/2}]$. Hence
$$E(W_{-1})\subseteq [\g_1^{\eta/2},[\g_1^{\eta/2},W_{-1}]]\subseteq
[\g_1^{\eta/2},\h_0]=0,$$ yielding $W_{-1}=0$. As a result,
$\dim\,\g_1^{\eta/2}=1$, that is $\g_1^{\eta/2}=\F y$. But then
$(\ad y)^4(f)=0$, a contradiction. The claim follows.

\bi

\noindent (h) We now follow \cite[Sec.~3]{P1} to show that for all
$x \in \mathcal X$ and  $t \in \mathbb{F}$ the linear map
$\exp\,(t \ad x)$ belongs to $G$, the connected component of the automorphism group
$\text{Aut}(S)$ of $S$.
Let $X=\ad x$. Since $X^4=0$ by part (g), this amounts to proving
that
 $$[X^2(u), \, X^3(v)] = 0\qquad
(\forall\, u,  v \in S).$$ Let $u\in V_i$ and $v\in V_j$ where
$-2\le i,j\le 2$ (we can view $i$ and $j$ as integers here). Since
$X^4\big([X(u),  v]\big) = 0$ and $p=5$, it follows that
$$[X^3(u),  X^2(v)] = [X^2(u), \, X^3(v)]=
\frac{1}{2}\,X\big([X^2(u), X^2(v)]\big).$$

Note that $X^3$ annihilates $V_0\oplus V_1\oplus V_2.$ So if
$[X^2(u),  X^3(v)] \neq 0,$ then $u,  v \in V_{-1} \cup V_{-2}$.
Since $\dim\, V_{-2}=1$, we also have $\{u,v\} \not \subset
V_{-2},$ while the equality $[V_1,V_2]=0$ implies $\{u,v\} \not
\subset V_{-1}$. In view of what is displayed above, we may thus
assume that $u=f$ and $v \in V_{-1}.$  But then $X^3(v)=\mu e$ for
some $\mu \in \mathbb F,$ and
$$[X^2(u),  X^3(v)] = \mu [[x,  [x,  f]],  e] = -\mu [x, [x, h]] = 0.$$

    Thus, we have shown that $\exp tX \in G$ for all $x \in
\mathcal X$ and $t\in \F$. Since $\mathcal X$ spans $V_1$, it
follows that $\ad V_1\subset \mathcal L = \text{Lie}(G)$. Since $E,F\in
\mathcal L$  too, and $\ad V_{-1}=[F,\ad V_1]$, we derive that
$\mathcal L$ contains a nonzero ideal of $\g$. But then $\ad
S\subseteq \mathcal L$.     By part (a) of this proof the result
follows. \qed
\medskip

\section {\ Algebras with an unbalanced grading}\label{sec:3.5}

\m

The previous section dealt with balanced gradings,  so what remains to be considered are 
the Lie algebras $\g$  with $\g_1  \cong
({\g}_{-1})^\ast$ as ${\g}_0$-modules
but with an unbalanced grading. 
We will address such algebras in this section
and complete the proof of  our  final (main) result on
the contragredient case:

\bi
\begin{Thm} \label{Thm:3.34} \ Assume $\mathfrak{g} = \bigoplus_{j =
-q}^r \mathfrak{g}_j$  is a graded Lie algebra
 over an algebraically closed field ${\mathbb F}$ of characteristic
$p > 3$.  Suppose $\mathfrak{g}_0$ is classical reductive with
$\mathfrak{g}_0^{(1)} \neq 0$, $\mathfrak{g}$ is transitive
\eqref{eq:1.3} and $1$-transitive \eqref{eq:1.4}, and
$\mathfrak{g}_{-1}$ is an irreducible restricted
$\mathfrak{g}_0^{(1)}$-module.  Assume further that
$\mathfrak{g}_1$ is an irreducible $\g_0$-module generated by a
$\mathfrak b^-$-primitive vector of weight $-\Lambda$. Then either
$\mathfrak{g}$ is a classical Lie algebra with one of its standard
gradings (see Section \ref{sec:2.4}), or else $p = 5$ and $\mathfrak{g}$ is a
Melikyan algebra.  In the latter case, the grading of
$\mathfrak{g}$ is the natural grading if $r \geq q$ and is
opposite to the natural grading if $r < q$.
\end{Thm}

\pf  We assume that

\begin{itemize}
\item[{\rm (i)}] $\mathfrak{g}_0 = \mathfrak{g}_0^{[1]} \oplus
\cdots \oplus \mathfrak{g}_0^{[\ell]}$ is a decomposition of
$\mathfrak{g}_0$ into ideals $\mathfrak{g}_{0}^{[i]}$ which are
classical simple, $\mathfrak{sl}_n$, $\mathfrak{gl}_n$,
$\mathfrak{pgl}_n$ where $p \mid n$, or one-dimensional;

 \item[{\rm (ii)}] $\mathfrak t$ is a fixed maximal toral
subalgebra in $\mathfrak{g}_0$ and ${\mathfrak t}^{[i]} =
{\mathfrak t} \cap \mathfrak{g}_{0}^{[i]}$;

\item[{\rm
(iii)}] $\Phi = \Phi^{[1]} \cup \cdots \cup \Phi^{[\ell]}$, where
$\Phi^{[i]}$ is the root system of $\mathfrak{g}_{0}^{[i]}$ with
respect $\mathfrak t^{[i]}$; \; $\Delta = \{\alpha_1, \dots,
\alpha_m\}$ is a system of simple roots in $\Phi$; and  $\Phi^+$
and $\Phi^-$ are the positive and negative roots respectively
relative to $\Delta$.

 \item[{\rm (iv)}] for $\alpha \in \Phi$, the root vector $e_{\alpha}$
spans the root space  $\mathfrak{g}_0^{\alpha}$, and the triple of
vectors $(e_{\alpha},e_{-{\alpha}},h_{\alpha})$, where $h_\alpha
\in \mathfrak{t}$, form a standard basis for a copy of
$\mathfrak{sl}_2$:

\begin{equation} [e_{\alpha},e_{-{\alpha}}] = h_{{\alpha}},\quad [h_{{\alpha}},e_{{\alpha}}] = 2e_{{\alpha}},
\quad [h_{{\alpha}},e_{-{\alpha}}] = -2e_{-{\alpha}}.
\nonumber\end{equation}

\smallskip \item[{\rm (v)}] $\mathfrak{n}^+ = \bigoplus_{\alpha \in \Phi^+}
\mathfrak{g}_0^{\alpha}$,  \  $\mathfrak{n}^- = \bigoplus_{\alpha
\in \Phi^-} \mathfrak{g}_0^{\alpha}$, \  $\mathfrak b^+ = \mathfrak t \oplus \mathfrak n^+$,
and $\mathfrak b^- = \mathfrak t \oplus \mathfrak n^-$.  
\end{itemize}

 \m

For $j\ge 0$, the $\g_0$-module $\g_{j+1}$ embeds into
$\text{Hom}(\g_{-1},\g_j)$ by transitivity \eqref{eq:1.3}. For
$j\le 0$, the $\g_0$-module $\g_{j-1}$ embeds into
$\text{Hom}(\g_{1},\g_j)$ by $1$-transitivity \eqref{eq:1.4}. Easy
induction on $\mid j\mid$ now shows that all homogeneous pieces
$\mathfrak{g}_j$ are restricted $\mathfrak{g}_0^{(1)}$-modules.

\m

Let $e^{-\Lambda}$ be a nonzero $\mathfrak b^-$-primitive vector of weight $-\Lambda$ in 
$\mathfrak{g}_1$,  and let $\widetilde{\mathfrak{g}}$ be the graded
subalgebra of $\mathfrak{g}$ generated by $\mathfrak{g}_0$,
$f^\Lambda$, and $e^{-\Lambda}$.  Set $\widetilde{\mathfrak{g}}_i
= \widetilde{\mathfrak{g}} \cap \mathfrak{g}_i$. As
$\widetilde{\mathfrak{g}}_{\pm 1} = \mathfrak{g}_{\pm 1}$  by our
initial assumption, the graded Lie algebra
$\widetilde{\mathfrak{g}}$ is transitive and $1$-transitive.
 It follows that the Weisfeiler radical $\mathcal
M(\widetilde{\mathfrak{g}})$ of $\widetilde{\g}$ is zero. Since
our assumption on $\g_1$ implies that $\g_1\cong \g_{-1}^*$ as
$\g_0^{(1)}$-modules, $\widetilde{\mathfrak{g}}$ satisfies the
assumptions of Theorem \ref{Thm:3.25}. By that result,
$\widetilde{\mathfrak{g}}$ is a classical Lie algebra with a
standard grading (where $e^{-\Lambda}$ has degree 1, $f^\Lambda$
degree $-1$, and elements of $\g_0$ have degree 0).

Suppose now that $\widetilde{\mathfrak{g}}_{t}\neq {\mathfrak{g}}_t$ for some $t \geq 2$, but
$\widetilde{\mathfrak{g}}_{i} = \mathfrak{g}_i$ for all
$1 \leq i < t$.    Let $x_{\beta} \in \g_t$ be such that
$x_\beta+\widetilde{\g}_t$ is  a $\mathfrak b^-$-primitive vector of weight $\beta$ for the
$\mathfrak{g}_0$-module $\g_t/\widetilde{\g}_t$.  Then
$[f^{\Lambda}, x_{\beta}] \in \widetilde{\g}_{t-1}$   and $[f_j,
x_{\beta}] \in \widetilde{\g}_t$ for all $f_j := e_{-\alpha_j}$ \
($j=1,\dots,m$). Set $f_0: = f^{\Lambda}$ and $\alpha_0 :=
-\Lambda.$     The elements $f_0,
f_1, \dots, f_m$ generate a maximal nilpotent subalgebra
$\widetilde{\mathfrak{n}}^- \subset \widetilde{\g}$. Since
$\widetilde{\mathfrak{n}}^- \supset \mathfrak{g}_{-1}$ and
$\mathfrak{g}$ is transitive, there exists an $l  \in \{0,1
\dots, m\}$ for which $[f_l, x_{\beta}] \neq 0$.   Thus,
$\beta-\alpha_l  \in \Phi(\widetilde{\g})$,   the set of roots of
$\widetilde{\g}$.    \m

Let us suppose first  that $\beta$ is a root of $\widetilde{\g}$
and that $[f_l, e_{\beta}] \neq 0$.     Write $\beta =
\sum_{i=0}^m b_i \alpha_i$. Then  $b_0 = \deg(e_\beta) =
\deg(x_{\beta}) = t \geq 2$. Since twice a root is never a root in
a classical Lie algebra (\cite[Lem.~ II.2.5]{S}), it follows that
$\beta$ has height $\hbox{\rm ht}(\beta) = \sum_{i=0}^m b_i \geq
3$.

\m   Since root spaces of a classical Lie algebra are
one-dimensional (Lemma II.2.6 of \cite{S}),  there is some scalar
$a \in {\mathbb F}^\times$ for which $[f_l, x_{\beta} -
ae_{\beta}] = 0$.   If $\beta-\alpha_k \not\in \Phi(\widetilde
\g)$ for all $k \neq l$,  then  $[f_k, x_{\beta} - ae_{\beta}] =
0$ for all such $k$.      But then
$[\widetilde{\mathfrak{n}}^-,x_{\beta} - ae_{\beta}] = 0,$
contradicting the transitivity of $\mathfrak{g}$. Thus, there must
exist a  $k \neq l$ such that $\beta-\alpha_k$ belongs to
$\Phi(\widetilde \g)$.

\m Suppose $\Phi(\widetilde{\g})$ has rank two. Then
$\widetilde{\Delta}\eqdef \{\alpha_k,\alpha_l\}$ is a basis of
$\Phi(\widetilde{\g})$ and $\beta$ is a root of height $\ge 3$ in
$\Phi(\widetilde{\g})$ with the property that $\beta-\alpha\in
\Phi(\widetilde{\g})$ for all $\alpha\in\widetilde{\Delta}$. Since
$\Z\gamma\cap \Phi(\widetilde{\g})=\pm\gamma$ for all $\gamma\in
\Phi(\widetilde{\g})$, this is impossible. Therefore, the rank of
the root system $\Phi(\widetilde{\g})$ is $\geq 3$.

    Set

    $$\widetilde{\Delta}(\beta) \eqdef \{\alpha_i \mid \beta - \alpha_i \in
    \Phi(\widetilde{\g})\}.$$

If all of the roots in $\Phi(\widetilde{\mathfrak{g}})$ have the
same length (i.e.,  if the root system is {\it simply laced}),  and if $\alpha_i$ and $\alpha_j$ are distinct simple
roots in $\widetilde{\Delta}(\beta)$,  then $(\alpha_i,
\alpha_j) = 0$.    To see this note that
when all roots have equal length, then root strings
through linearly independent roots have length
at most 2.    Thus, since $\alpha_i \in \widetilde \Delta(\beta)$,
we have

$$ \frac{2(\beta,\alpha_i)}{(\alpha_i,\alpha_i)} = 1,$$
and similarly for $\alpha_j$.   If $(\alpha_i,\alpha_j) \neq 0$,
then  $\alpha_i+\alpha_j$ is a root, and

\begin{equation} \frac{2(\beta,\,\alpha_i + \alpha_j)}
{(\alpha_i+\alpha_j,\,\alpha_i+\alpha_j)} = \frac{2(\beta, \,
\alpha_i)}{(\alpha_i,\alpha_i)} + \frac{2(\beta,\,
\alpha_j)}{(\alpha_j, \alpha_j)} = 2.\nonumber \end{equation}

\noindent  This forces
$\beta$ to equal $\alpha_i + \alpha_j$,     which  contradicts  the fact that
 ht$\beta \geq 3$.   Therefore, when $\Phi(\widetilde{\mathfrak{g}})$ is simply laced,  we must
have $(\alpha_i,  \alpha_j) = 0$ for all $\alpha_i,\alpha_j \in \widetilde{\Delta}(\beta)$.
 \m

 We have argued earlier that there exists a $k \neq l$ so that
 $[f_k, x_\beta-a e_\beta] \neq 0$.    Since $\alpha_k$ and $\alpha_l$ both belong
to $\widetilde \Delta(\beta)$, we must have $(\alpha_k,\alpha_l) = 0$
in the simply laced case,  and $[f_k,f_l] = 0$.  But then

\begin{equation*}(\beta - \alpha_k,\alpha_l) = (\beta, \alpha_l) > 0,
\end{equation*} \noindent so that  $[f_l,[f_k,x_\beta-a e_\beta]]
\neq 0$.    On the other hand,

\begin{equation} \label{xbeta} [f_l,[f_k,x_\beta-a e_\beta]]   =
[[f_l,f_k],x_\beta-a e_\beta]+ [f_k,[f_l,x_\beta-a e_\beta]] = 0.
\end{equation}

\noindent
This contradiction shows that all roots cannot have equal length,
so it must be that $\Phi(\widetilde{\mathfrak{g}})$ is of type B$_{m+1}$,
C$_{m+1}$,  or F$_4$ with  $m \geq 2$,  since the rank is at least 3.  \m

Suppose first that $\Phi(\widetilde{\mathfrak{g}})$ is of type B$_{m+1}$.
Let $\dot \alpha_0, \dot \alpha_1, \dots, \dot \alpha_{m}$
denote the labeling of the simple roots of B$_{m+1}$ compatible with that in  \cite[Planche II]{Bou1}
so that  the  highest
root  is

\begin{equation}\dot \alpha_0 + 2\dot \alpha_1 + \cdots  +2\dot \alpha_{m-1} + 2\dot \alpha_m\end{equation}

\noindent  Since the coefficient of $\alpha_0$ in the expression for $\beta$ is
$b_0 \geq 2$, it must be that  $\alpha_0 = \dot \alpha_s$ for some
$1 \leq s \leq m$.   Thus, $\beta$ must be
 of the form

\begin{equation} \sum_{j  \leq i < l}  \dot \alpha_i + 2 \sum_{l  \leq n \leq m}
\dot \alpha_n,  \end{equation}

\noindent   for $0 \leq j < l \leq m-1$,  and $\widetilde {\Delta}(\beta) = \{\dot \alpha_{j},
\dot \alpha_{l}\}.$ \m

Likewise, if $\Phi(\widetilde{\mathfrak{g}})$ is of type C$_{m+1}$,  then  the highest
root is

\begin{equation}2\dot \alpha_0 + 2\dot \alpha_1 + \cdots  +2\dot \alpha_{m-1} + \dot \alpha_m.  \end{equation}

\noindent  The root
$\beta$ must have  the form

\begin{equation}  \sum_{j \leq i < l} \dot \alpha_i + 2\left( \sum_{l \leq n < m}
\dot \alpha_n\right) + \dot \alpha_m  \quad \ \hbox{\rm or} \quad \
2\left(\sum_{j \leq l < m}  \dot \alpha_l \right) + \dot \alpha_m. \end{equation}

\noindent The second of these possibilities cannot happen, since $\widetilde \Delta(\beta) =
\{\dot \alpha_{j}\}$ in that case.    For the first,    $\widetilde \Delta(\beta) =
\{\dot \alpha_{j},\dot \alpha_{l}\}$ just as in the B$_{m+1}$-case.    \m

In either  the B$_{m+1}$-  or  C$_{m+1}$-case,   if  $l-j   \, >
1$,   then $(\dot \alpha_{j}, \dot \alpha_{l}) = 0$,
and we can use an argument as in (\ref{xbeta})  to arrive at a contradiction.  If
instead $l-j = 1$,   then   $\widetilde \Delta(\beta) = \{\dot \alpha_{j},
\dot \alpha_{j+1}\}$.    Since $j \leq m-2$,
the Cartan matrix entry $A_{j,j+1} = -1$, so  that the
length of the $\dot \alpha_{j+1}$-string through $\dot \alpha_j$ is
1.    If $\dot \alpha_j= \alpha_k$, then set $\dot f_j = f_k$, and adopt a similar
convention for $\dot \alpha_{j+1}$.

  Consequently\,  $[\dot f_{j+1}[\dot f_{j+1}, \dot f_{j}] ] = 0$, \,  and both 
$\beta - \dot \alpha_{j}- \dot \alpha_{j+1}$  and $\beta -$
$\dot \alpha_j -2\dot \alpha_{j+1}$ are roots of
$\widetilde{\mathfrak{g}}$.    As root spaces of a classical Lie
algebra are one-dimensional, we can choose  an appropriate $c \in
\F$, so that   $[\dot f_{j+1}, x_{\beta}-c e_\beta] = 0$. Since
$\mathfrak{g}_{-1}$ is generated as a $\mathfrak{g}_0$-module by
$\{f_0,f_1, \dots, f_m\}$,   it follows from transitivity  that
$[\dot f_j,  x_{\beta}-ce_\beta] = \zeta e_{\beta-\dot\alpha_j}
\neq 0$.  Then

$$[ \dot f_{j+1}, [
\dot f_{j+1},[\dot f_{j}, x_\beta - c e_\beta]]]
= \zeta[ \dot f_{j+1},[\dot f_{j+1},e_{\beta - \dot \alpha_i}] ] \neq
0,$$  since $\beta - \dot \alpha_j -$ $2\dot \alpha_{j+1}$ is a root of
$\widetilde{\mathfrak{g}}_0$.   However,

\begin{equation}\label{xbetai}[ \dot f_{j+1}, [\dot f_{j+1},[\dot f_{j}, x_\beta - c e_\beta]]]
= [[\dot f_{j+1}, [\dot f_{j+1}, \dot f_j]], x_\beta -c e_\beta]  = 0.
\end{equation}

\noindent  This contradiction
shows that $\Phi(\widetilde{\mathfrak{g}})$ is neither B$_{m+1}$ nor
C$_{m+1}$.   \m

    Suppose finally  that $\Phi(\widetilde{\mathfrak{g}}) \cong$ F$_4.$
We will adopt the labeling of roots
compatible with  \cite[Planche VIII] {Bou1}.
Because the coefficient \,$b_0$\,  of  \,$\alpha_0$\,  in $\beta$ is
$\geq 2$,  we need only
consider  roots having some coefficient larger than 1.  Below we list
the quadruple of coefficients of such roots  first relative to the simple roots
$\dot \alpha_j$, then relative to
the fundamental weights $\varpi_j$. Translation
between the two can be accomplished using the fact that
$\dot \alpha_i = \sum_j A_{i,j} \varpi_j$.

\begin{equation*}\begin{array} {cccccccc}
\hbox{\rm root}&\hbox{\rm0,1,2,0} & \hbox{\rm1,1,2,0} & \hbox{\rm 0,1,2,1} & \hbox{\rm1,2,2,0} &
\hbox{\rm1,1,2,1} &\hbox{\rm0,1,2,2} & \hbox{\rm 1,2,2,1} \\
\hbox{\rm weight}&\hbox{\rm -1,0,2,-2}&  \hbox{\rm 1,-1,2,-2} &\hbox{\rm -1,0,1,0}&
\hbox{\rm 0,1,0,-2}&\hbox{\rm 1,-1,1,0}&\hbox{\rm -1,0,0,2}& \hbox{\rm 0,1,-1,0}\\
&&&&&&&\\
\hbox{\rm root}& \hbox{\rm1,1,2,2} & \hbox{\rm1,2,3,1} & \hbox{\rm1,2,2,2} & \hbox{\rm1,2,3,2} &
\hbox{\rm1,2,4,2} & \hbox{\rm1,3,4,2} & \hbox{\rm2,3,4,2} \\
\hbox{\rm weight}&\hbox{\rm 1,-1,0,2}& \hbox{\rm 0,0,1,-1} & \hbox{\rm 0,1,-2,2}
& \hbox{\rm0,0,0,1} & \hbox{\rm 0,-1,2,0} & \hbox{\rm  -1,1,0,0} & \hbox{\rm1,0,0,0}\end{array}
\end{equation*}

\noindent  Now
$\beta-\dot\alpha_j$ is a root if and only if either the $j$th coordinate
in the weight expression of $\beta$ is positive or $\beta+\alpha_j$ is a
root and the $j$th coordinate in the weight expression of $\beta$ is zero.    Thus, the
only roots having $| \widetilde \Delta(\beta)| > 1$ are

$$\begin{array} {cccccc}
\hbox{\rm root}& \hbox{\rm 1,2,3,2}& \hbox{\rm 1,2,2,2}& \hbox{\rm1,1,2,2} &  \hbox{\rm 1,2,2,1}
 &
 \hbox{\rm1,1,2,1} \\
\widetilde \Delta(\beta) & \{\dot \alpha_3,\dot \alpha_4\}&\{\dot \alpha_2,\dot
\alpha_4\}&\{\dot \alpha_1,\dot \alpha_4\}&\{\dot \alpha_2,\dot \alpha_4\}
&\{\dot \alpha_1, \dot \alpha_3, \dot \alpha_4\}\\
\\ \hbox{\rm root}&\hbox{\rm 0,1,2,1} &  \hbox{\rm 1,1,2,0} &&& \\
\widetilde \Delta(\beta) & \{\dot \alpha_3,\dot \alpha_4\}&\{\dot \alpha_1,\dot \alpha_3\}.&&&
\end{array}$$ \m

Since there is a coefficient equal to 0 in each  of the last two roots, we may view the calculations
for them  as occurring in a root system
of type C$_3$ or B$_3$.  So by the same reasoning as in those cases, those possibilities
for $\beta$ can be eliminated.  For the roots corresponding to the tuples (1,2,2,2), (1,1,2,2),
(1,2,2,1), (1,1,2,1),  there are two orthogonal roots in $\widetilde \Delta(\beta)$.  Thus, an argument
as in (\ref{xbeta}) using the $f_k,f_l$ corresponding to those
two roots can be employed in those cases.  \m

Finally,  if the root $\beta$ has coefficient tuple $(1,2,3,2)$,
$\widetilde \Delta(\beta) = \{\dot \alpha_3,
\dot \alpha_4 \}$,  and both  $\beta - \dot \alpha_3 - \dot \alpha_4$ and $\beta -
\dot \alpha_3 - 2\dot \alpha_4$ are roots of $\widetilde{\mathfrak{g}}$.
However,  $\dot \alpha_3 + 2\dot \alpha_4$ is not a root, so that
$[\dot f_4,[\dot f_4,\dot f_3]] = 0$.    Then a calculation such as that in
(\ref{xbetai})  can be done to rule out this case.
Thus,  $\Phi(\widetilde{\mathfrak{g}})$ cannot be F$_4$ either.  \m

We have proven that if $[f_l, x_{\beta}] \neq 0$,  then either
$\beta \not \in \Phi(\widetilde{\mathfrak{g}})$  or $[f_l,
e_{\beta}] = 0$.    Because  $[f_l, x_{\beta}] = e_{\gamma}$, where
 $\gamma= \beta-\alpha_l$,   we have $[[f_l, x_{\beta}], e_{-\gamma}] =
h_{\gamma}$, and it follows that

\begin{equation}h_{\gamma} = [[f_l, e_{-\gamma}], x_{\beta}]
+ [f_l, [x_{\beta}, e_{-\gamma}]].\nonumber \end{equation}

If $\beta \not \in \Phi(\widetilde{\mathfrak{g}}_0)$ (the root system of $\widetilde \g_0$), then since  $-\alpha_l - \gamma = -\beta$,
it must be that $[f_l, e_{-\gamma}] = 0,$ and the first summand in the
expression for $h_\gamma$ vanishes.
If instead  $\beta \in \Phi(\widetilde{\mathfrak{g}}_0)$,   then  $[f_l,e_\beta] = 0$, and
we can replace $x_{\beta}$ by a suitable linear combination
$x_{\beta} + ce_{\beta}, \, c \in {\mathbb F}$ to  get $[[f_l,
e_{-\gamma}], x_{\beta}] = 0.$  (Such a replacement does not
render the equality $[[f_l, x_{\beta}], e_{-\gamma}] = h_{\gamma}$
invalid, because we have shown that $[f_l, e_{\beta}] = 0.)$  \m

    Thus, we may assume that $h_{\gamma} = [f_l,  [x_{\beta}, e_{-\gamma}]].$
As $[f_l, [x_{\beta}, e_{-\gamma}]] \in  \F h_{\alpha_l}$, it must
be that $h_{\gamma} \in \F h_{\alpha_l}.$ But then
$\varpi_i(h_{\gamma}) = 0$ for all $i \neq l.$  Since $p>3$,
Bourbaki's tables show that either $\gamma = \alpha_l$ or $\gamma
= -\alpha_l$; see \cite{Bou1}. Since $\gamma \in
\Phi(\widetilde{\mathfrak{g}})^+$, we conclude that $\gamma =
\alpha_l$.  As $\deg(e_{\gamma}) \geq t-1 \geq 1$, it follows that
$l = 0;$ consequently, $\beta = n\alpha_0$ for some integer $n$.
However, since $[f_0, x_{\beta}]$ is a root vector of a classical
Lie algebra, $n$ must be two; i.e.,
\begin{equation}\beta = 2\alpha_0\quad \hbox{\rm (and $t = 2$).} \nonumber \end{equation}

\noindent This implies that $x_{\beta}$ is a $\mathfrak b^-$-primitive vector of $\g_t$
(because $2\alpha_0 - \alpha_j \not \in \Phi(\widetilde{\mathfrak{g}}_0),
\, 1 \leq j \leq m).$  Now take any $\nu \in \{1, \dots,m \}$ for
which $(\alpha_{\nu},\alpha_0) \neq 0.$  Then, since $[f_0,
x_{\beta}] = e_0$ and $[f_{\nu}, x_{\beta}] = 0,$ we have

\begin{eqnarray}
\big[[f_0, [f_0, f_{\nu}]], x_{\beta}\big] &=&\Big( \ad[f_0, [f_0, f_{\nu}]]\Big)  (x_{\beta})\nonumber \\
&=& [F_0, [F_0, F_{\nu}]] (x_{\beta})  \qquad (\hbox{\rm where} \ F_0 = \ad f_0, \ F_\nu = \ad f_\nu) \nonumber \\
&=& \sum_{i=0}^2 (-1)^i{2 \choose i} \,F_0^{2 -i} F_{\nu} F_0^i\,(x_\beta)  \nonumber\\
&=& F_{\nu} F_0^2(x_{\beta}) = F_{\nu}F_0 (e_0)\nonumber \\ &=& [h_0,f_{\nu}] =  - \alpha_{\nu}(h_{\alpha_0})f_{\nu}\nonumber \\
&=& - \frac{2(\alpha_{\nu},\alpha_0)}{(\alpha_0,\alpha_0)}f_{\nu}\,\, \neq \,0,  \nonumber  \end{eqnarray}

\noindent whence $2\alpha_0 + \alpha_{\nu} \in
\Phi(\widetilde{\mathfrak{g}}_0).$ Since $\alpha_0$ and $\alpha_{\nu}$
are simple roots, the classical Lie algebra
$\widetilde{\mathfrak{g}}$ must have a root system of type B$_{m+1}$,
C$_{m+1}, \, m \geq 2$,   F$_4$, or G$_2$ .   Thus, we must consider, respectively,

\begin{equation}\label{eq:3.2805}
{\beginpicture \setcoordinatesystem
units <0.45cm,0.3cm> % sets scale
 \setplotarea x from -4 to 14, y from -1 to 1    % sets plot size up
 \linethickness=0.03pt                          % sets line thickness
  \put{$\circ$} at 0 0
 \put{$\circ$} at 2 0
 \put{$\cdots$} at 5 0
 \put{$\circ$} at  8 0
 \put{$\bullet$} at 10  0   \plot .15 .1 1.85 .1 /
 \plot 2.15 .1 3.85 .1 /
 \plot 6.15 .1 7.85 .1 /
 \plot 8.15 -.2 9.85 -.2 /
  \plot 8.15 .2 9.85  .2 /
  \put {$>$} at 9 0
 \put{$\alpha_1$} at 0 -1.5
  \put{$\alpha_2$} at 2 -1.5
    \put{$\alpha_{m}$} at 8 -1.5
        \put{$\alpha_0$} at 10 -1.5
        \put {$(m \geq 2)$} at 15.5 0  \endpicture}\end{equation}

\begin{equation}\label{eq:3.2810}
{\beginpicture \setcoordinatesystem
units <0.45cm,0.3cm> % sets scale
 \setplotarea x from -4 to 12, y from -1 to 1    % sets plot size up
% \linethickness=0.03pt                          % sets line thickness
  \put{$\circ$} at 0 0
 \put{$\circ$} at 2 0
 \put{$\cdots$} at 5 0
 \put{$\circ$} at  8 0
 \put{$\bullet$} at 10  0 \put{$\circ$} at  12 0
 \plot .15 .1 1.85 .1 /
 \plot 2.15 .1 3.85 .1 /
  \plot 6.15 .1 7.85 .1 /
  \plot 10.15 -.2 11.85 -.2 /
  \plot 10.15 .2 11.85  .2 /
  \put {$<$} at 11 0   \plot 8.15 .1  9.85 .1 /
  \put{$\alpha_1$} at 0 -1.5
  \put{$\alpha_2$} at 2 -1.5
   \put{$\alpha_{m-1}$} at 8 -1.5
      \put{$\alpha_{0}$} at 10 -1.5
        \put{$\alpha_m$} at 12 -1.5
        \put {$(m \geq 3)$} at 16 0  \endpicture}\end{equation}

\begin{equation}\label{eq:3.2815}
{\beginpicture \setcoordinatesystem
units <0.45cm,0.3cm> % sets scale
 \setplotarea x from -6 to 12, y from -1 to 1    % sets plot size up
 \linethickness=0.15pt                          % sets line thickness
  \put{$\circ$} at 0 0
 \put{$\circ$} at 2 0
  \put{$\bullet$} at 4 0
\put{$\circ$} at 6 0
 \plot .15 .1 1.85 .1 /
 \plot 2.15 -.2 3.85 -.2 /
  \plot 2.15 .2 3.85 .2 /
  \plot 4.15 .1 5.85 .1 /
  \put {$>$} at 3 0
  \put{$\alpha_1$} at 0 -1.5
  \put{$\alpha_2$} at 2 -1.5
  \put{$\alpha_0$} at 4 -1.5
  \put{$\alpha_3$} at 6 -1.5
  \endpicture}
\end{equation}

\begin{equation}\label{eq:3.2820}
{\beginpicture \setcoordinatesystem
units <0.45cm,0.3cm> % sets scale
 \setplotarea x from -6 to 8, y from -1 to 1    % sets plot size up
\put{$\circ$} at 0 0
 \put{$\bullet$} at 2 0
 \plot .15 -.3 1.85 -.3 /
  \plot .15 0 1.85  0 /
  \plot .15 .3 1.85 .3 /
  \put {$>$} at 1 0
   \put{$\alpha_1$} at 0 -1.5
  \put{$\alpha_0$} at 2 -1.5 \endpicture}
\end{equation}

\noindent Since we showed that $(\alpha_{\nu}, \alpha_0) \neq
0$ implies that $2\alpha_0 + \alpha_{\nu} \in
\Phi(\widetilde{\mathfrak{g}}_0)$,  neither  \eqref{eq:3.2810} nor
 \eqref{eq:3.2815} can occur, since in   \eqref{eq:3.2810},
we can take $\nu = m- 1,$ and in \eqref{eq:3.2815}, we can
take $\nu = 3$.   \m

Let us now consider  \eqref{eq:3.2805}, where the root system is of type B$_{m+1}$.
 Recall that in the Cartan matrix for B$_{m+1}$  (with roots numbered $\dot \alpha_0,\dot \alpha_1,
 \dots, \dot \alpha_m$),
we have $A_{m-1,m} = -2$,   so that $\alpha_m(h_{\alpha_0}) =
\dot \alpha_{m-1}(h_{\dot \alpha_{m}}) = A_{m-1,m} = -2$.     Then since
$\beta = 2\alpha_0$ and neither $2\alpha_0 - \alpha_m$ nor
$\alpha_0 - \alpha_m$ is a root, we have

\begin{eqnarray}
\big[[f_0,[f_0,[f_0,f_m]], x_{\beta}\big] &=&
\sum_{i=0}^3 (-1)^{i}{3 \choose i}  F_0^{3-i} F_m F_0^{i}\,(x_{\beta})\nonumber \\
&=& \Big(3 F_0 F_m F_0^{2}- F_m F_0^{3}\Big)(x_{\beta})\nonumber \\
&=& \Big(3 F_0 F_mF_0 - F_m F_0^{2}\Big)(e_{\alpha_0})\nonumber \\
&=& \Big( F_m F_0  - 3 F_0F_m\Big)(h_{\alpha_0})\nonumber \\
&=& \big(2-3(-2)\big)[f_m,f_0]\nonumber \\
&=& 8[f_0,f_m].  \nonumber\end{eqnarray}

\noindent  However, because $-A_{m-1,m} = 2$,   the length of the $\dot \alpha_m =
\alpha_0$ root string through $ \dot \alpha_{m-1} = \alpha_m$ is 2,    that so $[f_0,f_m] \neq 0$, but
$[f_0,[f_0,[f_0,f_m]]] = 0$.    This contradiction shows that  the configuration in
\eqref{eq:3.2805} is impossible.  \m

 It remains to consider   \eqref{eq:3.2820}, where the root system is of type G$_2.$
Then  $\alpha_1(h_{\alpha_0}) = \dot \alpha_{2}(h_{\dot \alpha_{1}})= A_{2,1} = -3$.    We calculate

\begin{eqnarray*}
\big[[f_0,[f_0,[f_0,[f_0,f_1]]], x_{\beta}]\big]&\hspace{-.1truein}=&\sum_{i=0}^4 (-1)^{i}{4 \choose i}
F_0^{4-i}F_1 F_0^{i}(x_{\beta})\nonumber \\
&\hspace{-.1truein}=&\Big(6 F_0^2F_1F_0^{2} - 4F_0F_1 F_0^{3}
+ F_1 F_0^{4}\Big)(x_{\beta})\nonumber \\
&\hspace{-.1truein}=& \Big(-6 F_0^2 F_1+ 4F_0F_1F_0 -
F_1 F_0)^{2}\Big)(h_{\alpha_0})\nonumber \\
&\hspace{-.1truein}=& - 6\alpha_1(h_{\alpha_0})[f_0,[f_0,f_1]] + 8[f_0,[f_1,f_0]]\nonumber \\
&\hspace{-.1truein}=& (18 -8)[f_0,[f_0,f_1]] = 10 [f_0[f_0,f_1]].  \end{eqnarray*}

\noindent Since the length of the $\ \dot \alpha_1 = \alpha_0$ root string
through $\dot \alpha_2 = \alpha_1$ is  $-A_{2,1} =3$,    $[f_0,[f_0,f_1]]\neq 0$,   but
$[f_0,[f_0,[f_0,[f_0,f_1]]]] = 0$.    It follows that
\eqref{eq:3.2820} is possible only if $p = 5$.   \m

    Let $\mathfrak{g}^{\geq -1}$ denote the
subalgebra of $\mathfrak{g}$ generated by the graded components
$\mathfrak{g}_i, \, i \geq -1$ and let ${\mathcal G}=\bigoplus_{j=-3}^3\,{\mathcal G}_j$ be the
graded Lie algebra of type ${\mathrm G}_2$ from Section
\ref{sec:2.14}.   The above discussion shows that
$\bigoplus_{j \le
0}\,\g^{\ge -1}_j \cong \, \bigoplus_{j\le 0}\,{\mathcal G}_j$  as graded Lie algebras. Clearly,
$\mathfrak{g}^{\geq -1}$ is both transitive and $1$-transitive.  Recall
that $x_\beta\not\in\widetilde{\g}_2$.
Therefore,  $\dim \g^{\ge -1}_2>\dim {\mathcal G}_2=1$, implying
$\g^{\ge -1}\not\cong\mathcal G$. Now Proposition \ref{Pro:2.48}
enables us to deduce that $\mathfrak{g}^{\geq -1}$ is isomorphic
to a Melikyan algebra $M: = M(2;\un{n})$ with its
natural grading (as in Section \ref{sec:2.14}). \m

Suppose $\mathfrak{g} \neq \mathfrak{g}^{\geq -1}.$  Let $\mathfrak{g}^{\leq 1}$
denote the subalgebra of $\mathfrak{g}$ generated by the graded
components $\mathfrak{g}_i, \, i \leq 1$.      Set
$\mathfrak{g}_i^{\leq 1} \eqdef \mathfrak{g}_{-i} \cap
\mathfrak{g}^{\leq 1}.$  Then, since $\mathfrak{g}$ is
1-transitive,  the graded Lie algebra $\mathfrak{g}^{\leq 1} =
\bigoplus_{i \in {\mathbb Z}}\mathfrak{g}_i^{\leq 1}$ satisfies
conditions (i) through (v).  Furthermore, because  $\mathfrak{g} \neq
\mathfrak{g}^{\geq -1},$ it follows that $\mathfrak{g}^{\leq
1} \neq \widetilde{\mathfrak{g}}_0$.   Reasoning as above, we now
conclude that the graded Lie algebra $\mathfrak{g}^{\leq 1} =
\bigoplus_{i \in {\mathbb Z}}\mathfrak{g}_i^{\leq 1}$ is
isomorphic to a Melikyan algebra $M' = M(2;\un{n}')$ with its
natural  grading.

\m Set
\begin{equation*}M_{\overline{\,0}}
\eqdef \bigoplus_{i \equiv 0 \mod 3}\mathfrak{g}_i^{\geq
-1} \quad \hbox{\rm and} \quad M_{\overline{\pm 2}} \eqdef \bigoplus_{i \equiv \pm 2 \mod
3}\mathfrak{g}_i^{\geq -1}.\end{equation*}

\noindent Then, as in Section \ref{sec:2.14},  $
M_{\overline{\,0}}$ is isomorphic to the \W \  Lie algebra
$W(2;\un{n})$,   and the subspaces $M_{\overline{-2}}$ and
$M_{\overline{\,2}}$ are irreducible $M_{\overline{\,0}}$-modules.
More precisely, $M_{\overline{\,2}} \,\cong\,
\widetilde{W}(2;\un{n})= \{\widetilde{E} \mid E \in W(2;\un{n})\}$
and $M_{\overline{-2}} \,\cong\,  \mathcal O(2; \un{n})$, and the
multiplication is as in \eqref{eq:2.46}:
\begin{gather*} {[D, \widetilde E]} = \widetilde{[D,E]} + 2
\di(D) \widetilde E \\
{[D, f]} = D(f) - 2 \di(D)f   \\
{[f_1 \widetilde{D_1} + f_2 \widetilde{D_2}, g_1 \widetilde{D_1} +
g_2
\widetilde{D_2}]} = f_1g_2 - f_2g_1 \cr [f, \widetilde{E}] = fE  \\
{ [f,g]}
= 2(f \widetilde{D}_g - g \widetilde{D}_f ) \qquad \hbox{\rm where}  \\
\widetilde D_f = D_1(f)\widetilde {D_2} - D_2(f)\widetilde {D_1} .
\end{gather*}

As in  Section \ref{sec:2.14},  the algebra $M$ may be assigned
a $\mathbb Z$-gradation according to:

\begin{gather*} M_{-3} = \F D_1 \oplus \F D_2, \qquad M_{-2} = \F 1,
\qquad M_{-1} = \F \widetilde {D_1} \oplus \F \widetilde {D_2} \\
M_0 = \spa_\F \{x_i D_j \mid i,j = 1,2 \}   \\
M_{1} = \F x_1^{(1)} \oplus \F x_2^{(1)},
\end{gather*}
Using these results, it is easy to check that the following hold:

\begin{eqnarray*}
&&\left[ x_1^{(1)} D_{1} - x_2^{(1)} D_{2}, \,x_1^{(1)} D_{2}\right] = 2x_1^{(1)} D_{2},  \\
&&\left[ 3x_2^{(1)} D_{2}, x_2^{(1)} \right] = 3x_2^{(1)}  - 2\cdot 3 \cdot x_2^{(1)}
\equiv 2x_2^{(1)}  \mod 5, \\
&&\left[ x_1^{(1)} D_{1} - x_2^{(1)} D_{2},\, x_2^{(1)} D_{1}\right] = -2x_1^{(1)} D_{2},  \\
&&\left[ 3x_2^{(1)} D_{2}, 3\widetilde{D_{2}}\right] = -9\widetilde{D_{2}} + 2\cdot 3 \cdot
3\widetilde{D_{2}} \equiv -\widetilde{D_{2}} \mod 5,  \\
&&\left[ x_1^{(1)} D_{2}, x_2^{(1)} D_{1}\right] = x_1^{(1)} D_{1} - x_2^{(1)} D_{2},  \\
&&\left[ x_2^{(1)} , 3\widetilde{D_2}\right] = 3x_2^{(1)} D_{2}, \hbox{ and } \\
&&\left[ x_1^{(1)} D_{2}, \widetilde{D_{2}}\right] = [x_2^{(1)} D_{1},x_2^{(1)} ] = 0.
\end{eqnarray*}
\smallskip

\noindent Thus, we can put
\begin{eqnarray}
e_1&\eqdef& x_1^{(1)} D_{2},\label{melg2}  \\
f_1&\eqdef& x_2^{(1)} D_{1},\nonumber \\
h_1&\eqdef& x_1^{(1)} D_{1} - x_2^{(1)} D_{2}\nonumber \\
e_0&\eqdef& x_2^{(1)} ,\nonumber \\
f_0&\eqdef& 3\widetilde{D_{2}},\nonumber \\
h_0&\eqdef& 3x_2^{(1)} D_{2}. \nonumber\end{eqnarray}

\noindent Set $x_{\beta} \eqdef x_2^{(1)} \widetilde{D_1}.$  Then
\begin{eqnarray}
[f_1, x_{\beta}] &=& 0,\label{eq:3111}  \\
{[h_1, x_{\beta}]} &=& - 2x_{\beta}, \nonumber\\ {
[f_0,x_{\beta}]} &=& [3\widetilde{D_{2}},\,x_2^{(1)} \widetilde{D_{1}}] = -
3x_2^{(1)}   \equiv 2e_0
\mod 5, \nonumber \\ { [h_0, x_{\beta}]} &=& 9x_2^{(1)} \widetilde {D_{1}} \equiv 4 x_\beta \ \mod 5,\nonumber .  \end{eqnarray}

    Since $\widetilde{\mathfrak{g}}_0 \neq \mathfrak{g}^{\leq1}$,  there exists
by symmetry an $x_{-\beta} \in \mathfrak{g}_{-2}$ satisfying the
following conditions

\begin{eqnarray}\label{eq:3222} [e_1, x_{-\beta}] &=& 0,  \\
{[h_1, x_{-\beta}]} &=& 2x_{-\beta},\nonumber \\
{[e_0,x_{-\beta}]} &=& -2f_0,\nonumber \\  {[h_0, x_{-\beta}]} &=&
- 4x_{-\beta}. \nonumber  \end{eqnarray}

\noindent (Recall that $\mathfrak{g}_{-2} = \mathfrak{g}_2^{\leq
1}$  and that the role of the $f_i$'s is played in this case by
the $e_i$'s, $i = 1,2.)$  \m

    We claim that $[x_{\beta},$ $x_{-\beta}]$ = $0.$ Indeed, let
$\tf'\eqdef  \F x_1^{(1)} D_{1} \oplus \F x_2^{(1)} D_{2}$. Then
$h\eqdef [x_{\beta},x_{-\beta}]$ is in $\tf',\,$ $[x_1^{(1)}
\widetilde{D_{1}},$ $x_{\beta}] = 0,\,$ $[x_1^{(1)}
\widetilde{D_1}, \, x_{-\beta}]\in \mathfrak{g}_0,$ and
$$[x_1^{(1)} D_{1}- x_2^{(1)} D_{2},\,[x_1^{(1)}
\widetilde{D_{1}},\,x_{-\beta}]]= 2[x_1^{(1)} \widetilde{D_{1}},\,
x_{-\beta}]$$ (see \eqref{eq:3222} above). By comparing
eigenvalues relative to $x_1^{(1)} D_{1}-x_2^{(1)} D_{2}$,   we
see that $[x_1^{(1)} \widetilde{D_{1}}, \,x_{-\beta}]$ $=
ax_1^{(1)} D_{2},$ where $a$ $\in {\mathbb F}.$ Therefore,

\begin{eqnarray}
[x_1^{(1)} \widetilde{D_{1}},\,[x_{\beta}, x_{-\beta}]] &=& [x_{\beta}, [x_1^{(1)}
\widetilde{D_{1}},\,x_{-\beta}]]\nonumber \\
&=& [x_2^{(1)} \widetilde{D_{1}}, \,ax_1^{(1)} D_{2}]\nonumber \\
&=& a\big(x_2^{(1)} \widetilde{D_{2}} - x_1^{(1)} \widetilde{D_{1}}\big).
\nonumber\end{eqnarray}

\noindent On the other hand, $\big[x_1^{(1)}
\widetilde{D_{1}},\,[x_{\beta}, x_{-\beta}]\big]$ = $-[h,x_1^{(1)}
\widetilde{D_{1}}]\,\in \,{\mathbb F}x_1^{(1)} \widetilde{D_{1}}.$
This implies that $a = 0$,  so that $[x_1^{(1)}
\widetilde{D_{1}}, x_{-\beta}] = 0 = [h,\,x_1^{(1)}
\widetilde{D_{1}}].$ Since $h \in \tf'$ and $[x_1^{(1)}
D_{1},\,x_1^{(1)} \widetilde{D_{1}}] = [x_2^{(1)} D_{2},\,x_1^{(1)}
\widetilde{D_{1}}] = 2x_1^{(1)} \widetilde{D_{1}}$,  we obtain
that $h = b(x_1^{(1)} D_{1} - x_2^{(1)} D_{2})$ for some $b \in
{\mathbb F}.$  \m

    Suppose that $b \neq 0,$ and set $e_1= x_1^{(1)} D_{2}$,
     $f_1 = x_2^{(1)} D_{1}$ as above, 
$e_2 \eqdef b^{-1}x_{-\beta},$ and $f_2 \eqdef x_{\beta}.$  Then
$[e_i, f_j]$ $= \delta_{i,j}h_{\alpha_1}, \, 1 \leq i,j \leq 2,$
and $[h_{\alpha_1}, e_i] = 2e_i, \, i = 1,2,$ and $[h_{\alpha_1},
f_i] = - 2f_i, \, i = 1,2$.   Then by Theorem \ref{Thm:3.9}, the
Lie algebra  generated by the elements $e_i, \, f_i$ \, $i=1,2$,  would be
infinite-dimensional, so that we must have that $b = 0$ and $h =
[x_{\beta}, x_{-\beta}] = 0$.    \m

Set $e_1'\eqdef [x_{\beta},e_0]$ and $f_1'\eqdef
[x_{-\beta},f_0]$.   By \eqref{eq:3111} we have
$[[x_{\beta},e_0],f_0]\in \F x_{\beta}.$ This yields
\begin{eqnarray*}
[e'_1,f'_1]&=&[[x_\beta,e_0],[x_{-\beta}, f_0]]\,=\,
[[[x_\beta ,e_0],x_{-\beta}],f_0]\\
&=& [[x_{\beta},[e_0,x_{-\beta}]],f_0]\,=\, -2[[x_{\beta},f_0],f_0]\\
&=&4[e_0,f_0]\, =\, 4h_{\alpha_0}.
\end{eqnarray*}
 Using the results of Section \ref{sec:2.14}, one
sees  that the $\mathfrak{g}_0^{(1)}$-modules
$\mathfrak{g}_4^{\geq -1}$ and $\mathfrak{g}_4^{\leq 1}$ map
into the graded component $\mathcal{O}(2)_2$, and so, as in the
proof of Lemma \ref{Lem:2.85}\,(b),  are isomorphic to the
3-dimensional $\mathfrak{g}_0^{(1)}$-module $L(2)$. Inasmuch as
$[e_1',\, e_0]$ $\in \mathfrak{g}_4^{\geq -1}$ and $[f_1',\,
f_0] \in \mathfrak{g}_4^{\leq 1},$ and 
$[h_{\alpha_1},\,e_0] = - e_0$ (see \eqref{eq:3111}),  we have
\begin{equation}[h_{\alpha_1},[e_1', e_0]]
= [h_{\alpha_1},[[x_{\beta}, e_0], e_0]] = -4[e_1', e_0]] \equiv
[e_1', e_0]] \mod 5,\nonumber \end{equation} and, by symmetry,
\begin{equation*}[h_{\alpha_1},[f_1', f_0]] = - [f_1', f_0]\mod 5.  \end{equation*}

\smallskip

\noindent Comparing eigenvalues, we see that $[e_1', e_0] = [f_1',
f_0] = 0.$ \m

  Now set $e_2' \eqdef e_0,$ $f_2' \eqdef f_0,$ $h' \eqdef h_{\alpha_0},$ $e_3' \eqdef f_1',$
and $f_3' \eqdef e_1'.$  Then $[e_i', f_j'] = \delta_{i,j}h', \, j
= 2,3,$ and $[h', e_2'] = 2e_2',$ $[h', f_2'] = -2f_2'$, $[h', e_3']$ $= [h', f_1']$ $= [h_{\alpha_0},
[x_{-\beta}, f_0]]$ $= -4e_3'  = e_3' \mod 5,$ and, by symmetry,
$[h', f_3'] = -f_3'.$  Therefore, Theorem \ref{Thm:3.9} applies,
showing that the Lie algebra generated by  $e_i',\, f_i',\,  i = 2,3,$ is
infinite-dimensional.    \m

    Thus, we have proved that if $\mathfrak{g}_t \neq \widetilde{\mathfrak{g}}_t$ for some
    $t \geq 2$,  then $p = 5,$
and $\mathfrak{g}$ is isomorphic to a Melikyan algebra with the
natural grading.   If $\mathfrak{g}_t = \widetilde{\mathfrak{g}}_t$
for all $t \geq 2$ and $\mathfrak{g}_{-s} \neq
\widetilde{\mathfrak{g}}_{-s}$ for some $s
> 0$,   we deal with the opposite grading

\begin{equation*}\mathfrak{g} = \bigoplus_{i = -r}^q \mathfrak{g}_i',  \end{equation*}

\noindent where $\mathfrak{g}_i' = \mathfrak{g}_{-i}$.    Since
$\mathfrak{g}_1$ is an irreducible restricted $\g_0^{(1)}$-module,
this graded Lie algebra satisfies the conditions of Theorem
\ref{Thm:3.34}. By the same reasoning as above, one can prove now
that again $p = 5, \ \mathfrak{g}$ is a Melikyan algebra, and the
grading of $\mathfrak{g}$ is the opposite of the natural grading.  \m

    This completes the proof of Theorem \ref{Thm:3.34} and with it
the contragredient case of the Main Theorem. \qed

%-----------------------------------------------------------------------
% Beginning of chap4.tex 9-15-05 version
%-----------------------------------------------------------------------
%
% AMS-LaTeX 1.2 sample file for a monograph, based on amsbook.cls.
% This is a data file input by chapter.tex.
%%%%%%%%%%%%%%%%%%%%%%%%%%%%%%%%%%%%%%%%%%%%%%%%%%%%%%%%%%%%%%%%%%%%%%%%%%%%%%%%%%%%%%%%%%%%%%%%%%%%%%%%%%%%%%%%%%%%%%%
%%%%%%%%%%%%%%%%%%%%%%%%%%%%%%%%%%%%%%%%%%%%%%%%%%%%%%%%%%%%%%%%%%%%%%%%

 \chapter{The Noncontragredient Case}
\bi

\section  {\ General assumptions and
notation  \label{sec:4.1}}
\m

Throughout this chapter, in which we prove the noncontragredient
case of the Main Theorem,  {\em our blanket assumptions are}  \m
\begin{itemize} \item[{\rm (1)}] $\mathfrak{g} = \bigoplus_{j = -q}^r \mathfrak{g}_j$ is a
finite-dimensional  graded Lie algebra over an algebraically
closed field ${\mathbb F}$ of characteristic $p>3$ with
$\mathfrak{g}_1 \neq 0$;

\smallskip \item[{\rm (2)}] $\mathfrak{g}_0$ is a classical reductive Lie algebra, and
$\mathfrak{g}_0 = \mathfrak{g}_0^{[1]} \oplus \cdots \oplus
\mathfrak{g}_0^{[\ell]}$ is a decomposition of $\mathfrak{g}_0$ into
ideals $\mathfrak{g}_{0}^{[i]}$ which are
\begin{itemize}
\item[{(i)}] classical simple,
\item[{(ii)}] $\mathfrak{sl}_n$, $\mathfrak{gl}_n$, or  $\mathfrak{pgl}_n$ where $p
\mid n$,  or
\item[{(iii)}] one-dimensional. \end{itemize}
\end{itemize}
\m

Also throughout {\em our conventions are}
\begin{itemize}
\item[{\rm (a)}] $\mathfrak{t}$ is a fixed maximal
toral subalgebra in $\mathfrak{g}_0$ and $\mathfrak{t}^{[i]}
= \mathfrak{t} \cap \mathfrak{g}_{0}^{[i]}$;

\smallskip \item[{\rm (b)}] $\Phi$ is the root system of $\mathfrak{g}_0$ with
respect to $\mathfrak{t}$; \; $\Delta = \{\alpha_1, \dots,
\alpha_m\}$ is a system of simple roots in $\Phi$; and  $\Phi^+$
and $\Phi^-$ are the positive and negative roots respectively
relative to $\Delta$;

\smallskip \item[{\rm (c)}] for $\alpha \in \Phi$, the root vector $e_{\alpha}$
spans the root space  $\mathfrak{g}_0^{\alpha}$, and the triple of
vectors $(e_{\alpha},e_{-{\alpha}},h_{\alpha})$, where $h_\alpha
\in \mathfrak{t}$, form a standard basis for a copy of
$\mathfrak{sl}_2$:

\begin{equation} [e_{\alpha},e_{-{\alpha}}] = h_{{\alpha}},\quad
[h_{{\alpha}},e_{{\alpha}}] = 2e_{{\alpha}},
\quad [h_{{\alpha}},e_{-{\alpha}}] = -2e_{-{\alpha}}.
\nonumber\end{equation}

\noindent (If $(e_{\alpha},e_{-{\alpha}},h_{\alpha})$ is the triple
corresponding to $\alpha \in \Phi^+$, then we assume that the
triple $(e_{-{\alpha}},e_{\alpha},h_{-{\alpha}}= -h_{\alpha})$ is
the one corresponding to $-\alpha \in \Phi^-$. One of the vectors
$e_{\alpha}$ or $e_{-{\alpha}}$ may be chosen to be an arbitrary
nonzero root vector, and then the remaining vectors in the triple
are chosen accordingly.)

\item[{\rm (d)}] $\mathfrak{n}^+ = \bigoplus_{\alpha \in
\Phi^+} \mathfrak{g}_0^{\alpha}$, \   $\mathfrak{n}^- =
\bigoplus_{\alpha \in \Phi^-} \mathfrak{g}_0^{\alpha}$, \
$\mathfrak b^+ = \mathfrak t \oplus \mathfrak n^+$,
and $\mathfrak b^- =\mathfrak t \oplus \mathfrak n^-$.

\item[{\rm (e)}]  $(\,,\,)$ is a symmetric bilinear form on the dual space to
${\mathfrak t}\cap\g_0^{(1)}$ which is invariant under the Weyl group $W$
of $\g_0^{(1)}$ and whose restriction to each
subspace $\big({\mathfrak t}^{[i]}\cap \g_0^{(1)}\big)^*$ of $\big({\mathfrak
t}\cap\g_0^{(1)}\big)^*$ is nonzero.
  (Such a form can be
obtained by mod $p$ reduction of the $W$-invariant form used
in the tables at the end of  \cite{Bou1}.)
\end{itemize}

\bi   Weight
vectors in $\g$ will be relative to
the toral subalgebra $\mathfrak t$.    Let
$\mathfrak{g}_i^{\lambda} = \{x \in \mathfrak{g}_i \mid [t,x] =
{\lambda}(t)x$ for all $t \in \mathfrak{t}\}$ denote  the
weight space in $\mathfrak{g}_i$ corresponding to $\lambda \in \mathfrak t^\ast$.
Let $\varpi_1, \dots, \varpi_m$  be
the fundamental weights relative to our system of simple roots
$\Delta$ so that $\varpi_i(h_{\alpha_j}) = \langle \varpi_i, \alpha_j\rangle = \delta_{i,j}$. Then any
weight of $\mathfrak{t}$ is a linear combination of the
$\varpi_i$'s.

\m  The classical Lie algebras of types A$_m$ and C$_m$ play a
special role in what follows. Often we identify an algebra of type
A$_m$ with $\mathfrak{sl}(V)$ (or $\mathfrak{sl}(V)$ modulo its
center) and an algebra of type C$_m$ with $\mathfrak{sp}(V) \subseteq
\mathfrak{sl}(V)$.   In
the A-case, the natural module $V$ has dimension $m+1$ and weights $ \e_i \ (1
\leq i \leq m+1)$ relative to the toral subalgebra of
diagonal matrices, where $\e_i$ denotes the projection of
a matrix onto its $(i,i)$-entry.    Thus, $\e_1 + \cdots + \e_{m+1} = 0$
on $\mathfrak{sl}(V)$, and the weights of $V$ are 
$\varpi_1, \varpi_2-\varpi_1, \dots, \varpi_m-\varpi_{m-1}, -\varpi_m$
in the language of the fundamental weights.  
In the C-case,  $V$ has dimension $2m$ and weights $ \pm \e_i \ (1
\leq i \leq m)$, or in terms of the fundamental weights, 
$\pm \varpi_1, \pm(\varpi_2-\varpi_1), \dots, \pm(\varpi_m-\varpi_{m-1})$. 
We will adopt the notation in \cite{Bou1} giving the
expressions for the roots  in terms of the weights $\e_i$
whenever we consider the root systems of type A$_m$ and C$_m$.
Thus for type A$_m$,  we have
\begin{eqnarray}\label{eq:4.2}&&\\
&& \Phi = \{\e_i - \e_j \mid 1 \leq i,j \leq m+1, \, i \neq j \}, \nonumber \\
&& \Delta = \{\alpha_i = \e_i -\e_{i+1} \mid 1 \leq i \leq m\},
\nonumber
\end{eqnarray}

\noindent and for type C$_m$, $m \geq 2$,
\begin{eqnarray}\label{eq:4.3}&&\\
& \Phi = \{\pm \e_i \pm \e_j \mid 1 \leq i,j \leq m \}, \nonumber \\
& \Delta = \{\alpha_i = \e_i -\e_{i+1} \mid 1 \leq i \leq m-1\} \cup \{
\alpha_m =2\e_m \}.
\nonumber
\end{eqnarray}
\noindent (Compare the discussion in Section \ref{sec:2.1}.)
The module $V$ has  a unique highest weight, namely $\varpi_1$, 
(relative to the usual partial order in which $\nu \leq \mu$ if $\mu-\nu
= \sum_{i=1}^m k_i \alpha_i$, where  $k_i \in \mathbb N$ for all $i$),
and a unique lowest weight, which is $-\varpi_{m}$ in the A$_m$-case
and is $-\varpi_1$ in the C$_m$-case.     We use 
 the term {\it standard representation} to refer to either $V$ or its dual module
$V^*$.    In the C$_m$-case, $V$ is isomorphic to  $V^*$, but this
is not true in the A$_m$-case,  since $V^*$ has highest weight
$\varpi_m$ and lowest weight $-\varpi_1$.   

\m
\section {\ Brackets of weight vectors
in opposite gradation spaces  \label{sec:4.2}} 
 
 \m
    Here we obtain information about the weights of certain pairs of weight
 vectors lying in homogeneous spaces of opposite
 degree.  This information will later help us to analyze the root
 structure of $\g_0$ and the
 structure of the $\g_0$-module $\g_{-1}$.
 \bi

\begin{Lem} \label{Lem:4.5}  \ Suppose for $i \neq 0$ that there are
weight vectors $x_{\lambda} \in \mathfrak{g}_{-i}^{\lambda}$ and
$x_\gamma \in \mathfrak{g}_i^\gamma$ corresponding to weights
${\lambda},\gamma \in \mathfrak{t}^\ast,$ which satisfy the
following conditions:

\begin{equation}[x_{\lambda},x_\gamma] = e_{-\delta}, \quad [x_{\lambda},e_\delta] = 0 =
[x_\gamma,e_{-\delta}] \;\; \text{for some}\; \delta \in
\Phi.\nonumber \end{equation}

\noindent Then either ${\lambda}(h_{\delta}) = 1$ or
${\lambda}(h_{\delta}) = 0$.
\end{Lem}

\pf  Assume that
${\lambda}(h_{\delta})({\lambda}(h_{\delta}) - 1) \neq 0$ and set

\begin{eqnarray} e_1  &=& [x_\gamma,e_\delta], \quad e_2 = [[[x_\gamma,e_\delta],e_\delta],e_\delta],
\quad f_1 = x_{\lambda} \nonumber \\
f_2  &=& \Big ({\lambda}(h_{\delta})({\lambda}(h_{\delta})
-1)\Big)^{-1}[[x_{\lambda},e_{-\delta}],e_{-\delta}]]\nonumber\end{eqnarray}

\noindent If  $\xi \eqdef \Big
({\lambda}(h_{\delta})({\lambda}(h_{\delta}) -1)\Big)^{-1}$, then
calculation shows that

\begin{eqnarray}
\left[e_1,f_1\right]  &=& [[x_\gamma,x_{\lambda}],e_\delta] = h_\delta \nonumber \\
\left[e_1,f_2\right] & \in& {\mathbb F}\bigl[\bigl[x_\gamma, [[x_{\lambda},e_{-\delta}],
e_{-\delta}]\bigl]\,,e_\delta\bigr]
+ {\mathbb F}\bigl[x_\gamma,\bigl[e_\delta,\,[[x_{\lambda},e_{-\delta}],e_{-\delta}]
\bigr ]\bigr]\nonumber \\
  && \subseteq {\mathbb F} [x_\gamma,[x_{\lambda},e_{-\delta}]] = 0 \nonumber \\
\left[e_2,f_1\right]  &=& [[[[x_\gamma,x_{\lambda}],e_\delta],e_\delta],e_\delta] = 0\nonumber \\
\left[e_2,f_2\right]  &=& \xi \Big (X_\gamma E_\delta^3 - 3
E_\delta X_\gamma E_\delta^2 + 3E_\delta^2 X_\gamma E_\delta  -
E_\delta^3X_\gamma \Big )
\big([[x_{\lambda},e_{-\delta}],e_{-\delta}]\big) \nonumber \\
&=& -3\xi E_\delta X_\gamma E_\delta^2\big([[x_{\lambda},e_{-\delta}],e_{-\delta}]\big) \nonumber \\
&=& -6 E_\delta X_\gamma(x_{\lambda}) = 6 h_\delta,
\nonumber\end{eqnarray}

\noindent where $E_\delta = \ad e_\delta$ and $X_\gamma = \ad
x_\gamma$.  Furthermore, since $e_1 = [x_{\gamma}, \, e_{\delta}]$
and $\lambda = -(\gamma + \delta),$

\begin{equation*}\begin{array}{cc}
\left[ h_\delta,e_1\right]  =  -{\lambda}(h_\delta)e_1, &\qquad \left[ h_\delta,e_2\right] =
(-{\lambda}(h_\delta) + 4)e_2 \\
\left[ h_\delta,f_1\right]  = {\lambda}(h_\delta)f_1, &\qquad \left[ h_\delta,f_2\right] =
({\lambda}(h_\delta) - 4)f_2 .
\end{array}\end{equation*}

\noindent Since ${\lambda}(h_\delta) \neq 0$, this contradicts
Theorem \ref{Thm:3.9}, which we apply by replacing $h_{\delta}$ by
$-2(\lambda(h_{\delta}))^{-1}h_{\delta}$ and $f_2$ by
$\frac{1}{6}f_2$.     \qed

\m \smallskip \smallskip
\section {\ Determining $\boldsymbol{\g_0}$
and its representation on $\boldsymbol{\g_{-1}}$ \label{sec:4.3}}  

 \m
Under a few additional assumptions,  we are able to prove here that  the commutator
$\g_0^{(1)}$  must be of type A$_m$ or C$_m$, and that the representation of
$\g_0^{(1)}$ on $\g_{-1}$  is a {standard representation}.

\bi \begin{Thm} \label{Thm:4.7}  \ Under hypotheses (1)-(2) of
Section  \ref{sec:4.1},   suppose that there exist a
$\mathfrak b^+$-primitive vector $f^\Lambda \in \mathfrak{g}_{-1}$ of weight
$\Lambda$  and a
$\mathfrak b^-$-primitive vector $e^{\Gamma} \in \mathfrak{g}_1$ of weight $\Gamma$ such that
\begin{itemize}
\item[{\rm (i)}] $[f^\Lambda,e^\Gamma] = e_{-{\alpha}}$\, for some
$\alpha \in \Phi$;   \item[{\rm (ii)}] $\Lambda(\mathfrak{t}^{[i]})
\neq 0$ for all $i = 1, \dots, \ell$\,  if $\alpha \in \Phi^+$, and
$\Gamma(\mathfrak{t}^{[i]}) \neq 0$ for all $i = 1, \dots, \ell$ if
$\alpha \in \Phi^-$.
\end{itemize}
\noindent Then the commutator ideal  $\g_0^{(1)}$ of
$\mathfrak{g}_0$ consists of a single summand, which is of type
$\hbox{\rm A}_m$ or $\hbox{\rm C}_m$.   Moreover,  if $\alpha \in \Phi^+$, then ${\alpha}$
is the highest root of $\Phi$, and $\Lambda$ is the highest weight
of a standard representation of $\g_0^{(1)}$.   If $\alpha \in \Phi^-$, then
${\alpha}$ is the lowest root of $\Phi$ and $\Gamma$ is the lowest
weight of a standard representation of $\g_0^{(1)}$.
\end{Thm}

\m
\begin{Rem}\label{Rem:4.8}  \ {\rm The condition $\Lambda(\mathfrak{t}^{[i]}) \neq 0$
for any $i = 1, \dots, \ell$  holds automatically when
$\mathfrak{g}_{-1}$ is an irreducible $\g_0$-module and
$\mathfrak{g}$ is transitive. Indeed, if
$\Lambda(\mathfrak{t}^{[i]}) = 0$ for some $i$,  then
$[\g_0^{[i]},f^\Lambda] = 0$, because in this case $\g_{-1}$ is a
completely reducible $\g_0^{[i]}$-module with all composition
factors being simple modules generated by $\mathfrak b^+$-primitive
vectors of weight zero. Then
$\g_{-1} = \sum_{j \neq i} \mathfrak{U}\big(\g_0^{[j]}\big)
f^\Lambda$ by irreducibility. But this implies $\g_{-1}$ is
annihilated by $\g_0^{[i]}$, which contradicts transitivity.}
\end{Rem}
\m

\pf  It suffices to prove the result when $\alpha \in
\Phi^+$, for the argument in the case $\alpha \in \Phi^-$ is completely
symmetrical and uses the hypothesis $\Gamma(\mathfrak t^{[i]}) \neq 0$ for any
$i$.

First  we assume that ${\alpha}$ is a highest root of $\Phi$,
hence a highest root of some summand $\g_0^{[k]}$,  and we show
that then $\g_0^{(1)}$ is of type  A$_m$ or C$_m$ and $\Lambda$ is
the highest weight of a standard representation of $\g_0^{(1)}$.
We begin by proving that $\Lambda(h_{\alpha}) \neq 0$. Indeed,
suppose $\Lambda(h_{\alpha}) = 0$.   If $\g_0^{[k]}$ is isomorphic
to $\mathfrak {sl}_2$, then $\Lambda(\mathfrak t^{[k]}) = 0$,
contrary to assumption. Therefore, we may assume that the rank of
$\g_0^{[k]}$ is $\geq 2$  and  may choose $\beta \in
\Phi^+\setminus \{{\alpha}\}$ such that $(\alpha,\beta)>0$. Then
$[e_{\alpha}, \, e_{-\beta}] \neq 0$ in $\g_0$. Consider the
weight vectors

\begin{equation}x_{\lambda} = [f^\Lambda,e_{-\beta}], \quad \quad x_\gamma = [e^\Gamma,e_\beta],
\quad \quad e_{-\delta} = (\Gamma +
{\alpha})(h_\beta)e_{-{\alpha}} = -\Lambda
(h_\beta)e_{-{\alpha}}.\nonumber \end{equation}

\noindent  Since $f^\Lambda$ is a $\mathfrak b^+$-primitive vector of weight $\Lambda$, and $e^\Gamma$ is a $\mathfrak b^-$-primitive vector of weight $\Gamma$, and $\Lambda + \Gamma = -\alpha,$ a lowest root, we have

\begin{eqnarray*}
[x_{{\lambda}},x_{\gamma}]
&=& \bigl[[f^{\Lambda}, \, e_{-\beta}], \, [e^{\Gamma}, \, e_{\beta}]\bigr]\\
&=& [[[f^{\Lambda}, \, e_{-\beta}], \, e^{\Gamma}], \, e_{\beta}] + [e^{\Gamma}, \,
[[f^{\Lambda}, \, e_{-\beta}], e_\beta]] \\
&=& - [e^{\Gamma}, \, [f^{\Lambda}, \, [e_{\beta}, \, e_{-\beta}]]\\
&=& \Lambda(h_{\beta})[e^{\Gamma}, \, f^{\Lambda}] = -\Lambda (h_{\beta})e_{-{\alpha}} = e_{-\delta}.
\end{eqnarray*}

\noindent In addition, since $\delta = \alpha,$ the highest root
of $\g_0^{(1)}$,

\begin{equation}[x_{\lambda},e_\delta] = 0 = [x_\gamma,e_{-\delta}].\nonumber \end{equation}

\noindent Therefore, Lemma \ref{Lem:4.5}  applies and gives
$\lambda(h_{\delta}) =(\Lambda-\beta)(h_\alpha) = 0$ or 1. Since
we are assuming that $\Lambda(h_{\alpha}) = 0$,   this forces
$\beta(h_{\alpha}) = 0$ or $\beta(h_{\alpha}) = -1$.   But
$\beta(h_{\alpha}) \in \{0,1,2,3\}$ as ${\alpha}$ is a highest
root of $\Phi$. Therefore, $\beta(h_{\alpha}) = 0$, which implies
$(\beta,\alpha) = 0$,  a contradiction to our choice of $\beta$.
We conclude that $\Lambda(h_{\alpha}) \neq 0$.   But then applying
Lemma \ref{Lem:4.5}  to the weight vectors $f^\Lambda$, $e^\Gamma$
and $e_{-{\alpha}}$ gives $\Lambda(h_{\alpha}) = 1$.

Our next goal is to prove that $\Lambda(h_\beta) = 0$ for any
$\beta \in \Phi^+$ such that $\alpha - \beta \not \in \Phi$.
Indeed, if on the contrary, $\Lambda(h_\beta) \neq 0$ for some
$\be \in \Phi^+$ with $\alpha - \be \not \in \Phi$,  we may
consider the elements:

\begin{equation*}
\begin{array}{ccc} e_1  = [e^\Gamma,e_{\alpha}], &\qquad  e_2 = [[e^\Gamma,e_{\alpha}],e_\beta],
&\qquad  \\
f_1  = f^\Lambda, &\qquad  f_2 =
-\Lambda(h_\beta)^{-1}[f^\Lambda,e_{-\beta}],  &\qquad  h =
h_{\alpha}. \end{array} \end{equation*}

\noindent Then because ${\alpha}-\beta \not \in \Phi$ (so that
$[e_{-\beta}, \, e_{\alpha}] = 0)$ and ${\alpha}$ is a highest
root, we have $\alpha(h_\be) = 0$, and hence $\be(h_{\alpha}) = 0$.   Moreover, since $\Gamma + \alpha =
-\Lambda$ and $f^\Lambda$ (resp. $e^\Gamma$) is a $\mathfrak b^+$-primitive vector of weight $\Lambda$ (resp. $\mathfrak b^-$-primitive vector of weight $\Gamma$),  we have  the following relations:

\begin{eqnarray}
\left[e_1,f_1\right]  &=& h_\alpha = h, \nonumber \\
\left[e_1,f_2\right]  &=& -\Lambda(h_\beta)^{-1}[[e^\Gamma,e_{\alpha}], \,
[f^\Lambda,e_{-\beta}]]\nonumber \\
&=& -\Lambda(h_\beta)^{-1}\left( [[[e^\Gamma,e_{\alpha}], f^\Lambda],\,
e_{-\beta}] +[ f^\Lambda, \, [[e^\Gamma, e_{\alpha}],\,e_{-\beta}] ]\right) = 0 \nonumber \\
\left[e_2,f_1\right]  &=&
[[[e^\Gamma,f^\Lambda],e_{\alpha}],e_\beta] =
\beta(h_{\alpha})e_\beta =
0, \nonumber \\
\left[e_2,f_2\right]  &=&
\Lambda(h_\beta)^{-1}[[f^\Lambda,[e_{-\beta},[[e^\Gamma,e_{\alpha}],e_\beta]] \nonumber \\
 &=&
\Lambda(h_\beta)^{-1}(\Gamma+{\alpha})(h_\beta)[f^\Lambda,[e^\Gamma,e_{\alpha}]]
= h_\alpha = h.  \nonumber\end{eqnarray}

\noindent Furthermore, since $\Lambda(h_{\alpha}) = 1, \, \Gamma +
\alpha = -\Lambda,$ and $\beta(h_{\alpha}) = 0,$ we have

\begin{equation}[h,e_i]  = -e_i \quad \quad [h,f_i] = f_i, \quad i = 1,2.\nonumber \end{equation}

\noindent  These calculations show that the elements $h' = -2h$,
$e_i' = e_i$, $f_i' = -2f_i$, for $i = 1,2$, satisfy all the
conditions of Theorem \ref{Thm:3.9} and so generate an infinite-dimensional algebra.
Since $\g$  is
finite-dimensional, this is impossible.  Hence $\Lambda(h_\beta) =
0$ for any $\beta \in \Phi^+$ such that $\alpha - \beta \not \in
\Phi$.

 Take a
simple root $\alpha_n$  in $\Phi^{[i]}$ for some $i \neq k$. Since
$\alpha -\alpha_n \not \in \Phi$, we have $\Lambda(h_{\alpha_n}) =
0$ by what we have just shown.    But
$\Lambda(\mathfrak{t}^{[i]}) \neq 0$ for all $i = 1, \dots, \ell$.
Consequently, $\Phi$ must be an irreducible root system.

If  $\be$ is a simple root for which $\alpha - \beta
\in \Phi$ and $\alpha - 2 \beta \not \in \Phi$,  then necessarily $\beta \in \Phi^{[k]}$
and $\alpha-\beta \in (\Phi^{[k]})^+$.   In this case, set

$$x_{\lambda} = f^\Lambda, \quad \quad x_\gamma = [e^\Gamma,e_{{\alpha}-\beta}].\nonumber$$

\noindent Then  $[x_{\lambda},e_\beta]$ $= 0$ $= [x_\gamma,$
$e_{-\beta}]$. Moreover, 
$$[x_{\lambda},x_\gamma] =
[[f^\Lambda, e^\Gamma], e_{{\alpha}-\beta}] =
[e_{-{\alpha}}, e_{{\alpha}-\beta}] = \zeta e_{-\beta}$$ for
some $\zeta$ $\neq 0$.
 Thus, $\Lambda(h_\beta)$ $= 0$ or $\Lambda(h_\beta)$
$= 1$ by Lemma \ref{Lem:4.5}  for such a root $\beta$.

Let $\alpha = \sum_{i = 1}^m  k_i \alpha_i$.  Examination of the
root systems (see \ref{sec:4.5}) shows that when $\Phi \neq
\hbox{\rm A}_m$ or $\hbox{\rm C}_m$, there exists $j \in \{1,
\dots, m\}$ such that
\smallskip

\begin{itemize} \item[{\rm (i)}] $\alpha - \alpha_j \in \Phi$, \smallskip \item[{\rm (ii)}]
 $\alpha - 2 \alpha_j \not \in \Phi$,   \item[{\rm (iii)}]
$({\alpha},{\alpha}) = (\alpha_j, \alpha_j),$
\item[{\rm (iv)}] ~$k_j = 2$
 \item[{\rm (v)}] $\alpha - \alpha_i \not \in \Phi$ for all $i \in
\{1, \dots, m\}$ such that  $i \neq j$.
\end{itemize}

\m \noindent Then by (v) and what we have shown, $\Lambda(h_{\alpha_i}) = 0$ for $i \neq j$. Since
$\mathfrak{g}_{-1}$ is not a trivial $\g_0^{(1)}$-module, we obtain {from} the preceding
paragraph that $\Lambda(h_{\alpha_j}) =1$.    But
then in view of Theorem \ref{Thm:2.4} we have that

$$\Lambda(h_{\alpha}) = \Lambda\Big(\sum_{i=1}^n\langle \varpi_i,\alpha \rangle h_{\alpha_i}\Big) =
\langle \varpi_j,\alpha \rangle\Lambda(h_{\alpha_j}) = k_j\Lambda(h_{\alpha_j}) = k_j = 2,\nonumber$$

\noindent contradicting the fact that $\Lambda(h_{\alpha}) = 1$.
It must be then that  $\Phi = \hbox{\rm A}_m$ or $\Phi = \hbox{\rm C}_m$.

If $\Phi = $ C$_m$, then $\alpha = 2\e_1$ and  $\alpha - \alpha_i
\not \in \Phi$ for $i > 1$, so that $\Lambda(h_{\alpha_i}) = 0$
for $i \neq 1$. (We are assuming the standard indexing of the
simple roots as in 4.2.)  Thus, $\Lambda = a \varpi_1$ and $\alpha_j = \alpha_1$.   But since
$({\alpha},{\alpha}) = 2(\alpha_1,\alpha_1),$   $k_1 = 2$,
and $\Lambda(h_{\alpha}) = 1$, we have using \eqref{eq:halpha} that  
$$1 = \Lambda(h_\alpha) = a \langle \varpi_1, \alpha\rangle =  a \frac{2(\varpi_1, \alpha)}{(\alpha,\alpha)} =
   a \frac{4(\varpi_1,\alpha_1)}{2 (\alpha_1,\alpha_1)} = a \langle \varpi_1, \alpha_1\rangle = a,$$
so we can conclude that $\Lambda =
\varpi_1.$

For $\Phi =$ A$_m$, we have $\alpha =
\varepsilon_1-\varepsilon_{m+1} = \alpha_1 + \cdots + \alpha_m$; \
$\alpha - \alpha_1, \alpha - \alpha_m \in \Phi$; \ $\alpha - 2
\alpha_1 \not \in \Phi$, $\alpha - 2 \alpha_m \not \in \Phi$;  and
$\alpha - \alpha_i \not \in \Phi$ for $i \in \{2, \dots, m-1\}$.
This forces $\Lambda = a_1 \varpi_1 + a_m \varpi_m$.    However,
since $1 = \Lambda(h_{\alpha}) = a_1 + a_m$, we see that $\Lambda = \varpi_1$
or $\Lambda = \varpi_m$.

The proof of Theorem \ref{Thm:4.7} in the case that $\alpha \in
\Phi^+$ will be complete once we show that ${\alpha}$ is a highest
root of $\Phi$.  Let $\widetilde{\mathfrak{g}}$ denote the graded
subalgebra of $\mathfrak{g}$ generated by $f^\Lambda$,
$\mathfrak{g}_0$, and $e^\Gamma$.  Set $\widetilde{\mathfrak{g}}_i
= \mathfrak{g}_i \cap \widetilde{\mathfrak{g}}$.  Let $J$ be the
sum of all graded ideals of $\widetilde{\mathfrak{g}}$ having zero
intersection with the local part $\widetilde{\mathfrak{g}}_{-1}
\oplus \mathfrak{g}_0 \oplus \widetilde{\mathfrak{g}}_1$, and set
$\overline{\mathfrak{g}} = \widetilde{\mathfrak{g}}/ J$.   Clearly
$\overline{\mathfrak{g}} = \bigoplus_{i = -{q'}}^{r'}
\overline{\mathfrak{g}}_i$ is a graded Lie algebra satisfying all
the hypotheses of Theorem \ref{Thm:4.7}.   Indeed, we can identify
$\overline \g_0$ with $\g_0$ and $\overline{\mathfrak{g}}_{\pm 1}$
with $\g_{\pm 1}$, as $J$ intersects the local part trivially. The
ideals $\mathcal A^-$ and $\mathcal A^+$ (see Proposition
\ref{Pro:1.45}) must be contained in $J$. As a result, the
following conditions hold: \m

\begin{itemize}
\item[{\rm (a)}] If $[\overline{x},\overline{\mathfrak{g}}_{-1}] =
0$ for some $\overline{x} \in \overline{\g}_i$, $i \geq 0$, then
$\overline{x} = 0$.
\smallskip \item[{\rm (b)}] If $[\overline{x},\overline{\mathfrak{g}}_{1}] = 0$
for some $\overline{x} \in \overline{\g}_i$, $i \leq 0$, then
$\overline{x} = 0$.
\end{itemize}

 \m Thus, in what follows we will assume that $\mathfrak{g} = \overline{\mathfrak{g}}$,
and hence that $\g$ is generated by $\g_0$,  $f^\Lambda$,
$e^\Gamma$,  that $\g$ is transitive and 1-transitive (but not
necessarily irreducible), and that there exists $\be \in \Delta$
for which $[e_{\alpha},e_\beta] \neq 0$ (as otherwise $\alpha$ is
a highest root, and we are done). Again we suppose that
$\Phi^{[k]}$ is the irreducible component of the root system
$\Phi$ containing ${\alpha}$.  As before, applying Lemma
\ref{Lem:4.5}  to the weight vectors $f^\Lambda$, $e^\Gamma$ and
$e_{-{\alpha}}$ we get $\Lambda(h_{\alpha}) = 0$ or
$\Lambda(h_{\alpha}) = 1$. In either case,  $\Gamma(h_{\alpha}) =
-(\Lambda + {\alpha})(h_{\alpha}) \neq 0$ since $p>3$.

Put $V_{-1}={\mathfrak U}(\mathfrak{n}^-) f^\Lambda$ and
$V_1={\mathfrak U}(\mathfrak{n}^+)\, e^\Gamma$. For $i\ge 2$,
define $V_{\pm i}$ inductively by setting $V_{\pm i}=[V_{\pm
1},\,V_{\pm (i-1)}]$. Each $V_i$ is a $\g_0$-submodule of $\g_i$,
and the graded subspace $I\eqdef \g_0\,\oplus\,\sum_{i\ne 0}
\,V_i$ is invariant under the endomorphisms $\ad f^\Lambda$ and
$\ad e^\Gamma$. Therefore, $I$ is an ideal of $\g$ containing
$f^\Lambda,\,$ $e^\Gamma$, and $\g_0$. But then $I=\g$, yielding
$\g_{-1}={\mathfrak U}(\mathfrak{n}^-) f^\Lambda$.

Set $f^{\Lambda'} = [[e_{-\be},f^\Lambda],f^{\Lambda}]$ and
$e^{\Gamma'} = [[e_{\alpha},e^\Gamma],e^\Gamma]$, and let
$F^\Lambda = \ad f^\Lambda$ and $E_{-\be} = \ad e_{-\be}$ denote
the adjoint mappings.    Then since $\beta + \alpha$ is a root,
and  twice a root is never a root, we must have $\beta \neq
\alpha$ and

\begin{eqnarray} [f^{\Lambda'},e^{\Gamma'}]  &=&
\Big((F^\Lambda)^2 E_{-\be} -2 F^\Lambda E_{-\be} F^\Lambda +
E_{-\be}(F^\Lambda)^2 \Big)
\Big( [[ e_{\alpha},e^\Gamma],e^\Gamma]\Big) \nonumber \\
&=& E_{-\be}(F^\Lambda)^2 \Big( [[e_{\alpha},e^\Gamma],e^\Gamma]\Big) \nonumber \\
&=& 2 E_{-\be}F^\Lambda \Big(
[[e_{\alpha},e_{-{\alpha}}],e^\Gamma]\Big)
  = 2\Gamma(h_{\alpha})[e_{-\be},e_{-{\alpha}}]. \nonumber \\
\nonumber\end{eqnarray}

\noindent As $\be \in \Delta$, it can be easily checked that
$f^{\Lambda'} \in \mathfrak{g}_{-2}$ is a $\mathfrak b^+$-primitive vector
of weight $\Lambda'$.  
We claim that $e^{\Gamma'} \in
\mathfrak{g}_2$ is a $\mathfrak b^-$-primitive vector of weight $\Gamma'$.    Indeed, if this were not true, then there would be a
sequence $F_{i_1}, \dots, F_{i_n}$ of length $n \geq 1$,  where $F_i = \ad
e_{-\alpha_i}$,   such that $e' \eqdef
F_{i_1} \cdots F_{i_n}e^{\Gamma'} \neq 0$ and $[\mathfrak{n}^-,
e'] = 0$.  But
$$e' = [[F_{i_1} \cdots F_{i_n}e_{\alpha}\,,\,e^\Gamma],e^\Gamma]\nonumber$$
\noindent implies $0\neq e' = \xi [[e_\sigma,e^\Gamma],e^\Gamma]$
for some $\xi \in {\mathbb F}^\times$ and $\sigma \in \Phi^+ $.
Since $n\geq 1$, we must have
$$[f^\Lambda,e'] = 2 \xi [[e_\sigma,e_{-{\alpha}}],e^\Gamma] = 0.\nonumber$$
\noindent Since $\g_{-1} = \mathfrak{U}(\mathfrak n^-)f^\Lambda$,
this yields $[e',\mathfrak{g}_{-1}] = 0$, contradicting
transitivity and proving the claim.

  Next we demonstrate that $(2 \Lambda -\be)(\mathfrak{t}^{[k]}) \neq 0$.
Assume that $\theta$ is any positive root for which
$[e_{-{\alpha}},e_{-\theta}] \neq 0$.  Set

$$x_\lambda= [e_{-\theta},f^\Lambda], \quad \quad x_\gamma = e^\Gamma,
\quad \quad  e_{-\delta} = [e_{-{\alpha}},e_{-\theta}].\nonumber$$

\noindent An easy calculation shows that these weight vectors
satisfy all the conditions of Lemma \ref{Lem:4.5}. Therefore,
$(\Lambda -\theta)(h_{\alpha + \theta}) \in \{0,1\}$. Setting
$\theta = \be$, we get $(\Lambda -\be)(h_{\alpha + \be}) \in
\{0,1\}$. If $(2\Lambda - \be)(\mathfrak{t}^{[k]}) = 0$, then in
particular $\frac {1}{2}\be(h_{\alpha}) = \Lambda(h_{\alpha}) \in
\{0,1\}$ and $1 = \frac {1}{2}\be(h_\beta) = \Lambda(h_\beta)$.
If $(\alpha,\beta) = 0$, then

\begin{eqnarray} (\Lambda-\beta)(h_{\alpha+\beta}) &=& -\textstyle{\frac {1}{2}}
\beta(h_{\alpha + \beta})\nonumber \\
&=&-\frac{(\beta,\alpha + \beta)}{(\alpha+\beta,\alpha+\beta)}
=  -\frac{(\beta,\beta)}{(\alpha+\beta,\alpha+\beta)} \in \{0,1\}.  \nonumber \end{eqnarray}

\noindent  This forces $(\alpha+\beta,\alpha+\beta) = - (\beta,\beta)$, contradicting
the assumption $p > 3$.   Hence it must be that $(\alpha,\beta) \neq 0$
and $\Lambda(h_\alpha)  = \frac {1}{2}\be(h_\alpha) = 1$.
This implies that $(\beta,\alpha) = (\alpha,\alpha)$ so that
$(\beta,\alpha+\beta) = (\alpha,\alpha) + (\beta,\beta) \neq 0$.
Consequently, $ (\Lambda-\beta)(h_{\alpha+\beta})  =
\displaystyle{-\frac{(\beta,\alpha + \beta)}{(\alpha+\beta,\alpha+\beta)}   = 1}$.
This means that $2(\alpha,\alpha)+(\beta,\beta) \equiv 0$  $\mod p$.
If $\alpha$ is a long root, then $\beta(h_\alpha) \in \{-1,0,1\}$ (see \cite{Bou1}).  However,
we know that $\beta(h_\alpha) = 2$.   As $p > 3$, we conclude $\alpha$ is a short root.
Then from the relation $2(\alpha,\alpha)+(\beta,\beta) \equiv 0$ it follows that only
the case that  $p = 5$ and
$\Phi^{[k]}$ is of type G$_2$ is possible.     We assume
that the elements of the base $\{\alpha_1,\alpha_2\}$ of G$_2$ are numbered so  $\alpha_1$ is  short
and $\alpha_2$ is long.   Since $\beta$ is simple and long, it must
be that $\beta = \alpha_2$.    Then from the fact that $\alpha$ is a short positive root,
and $\alpha + \beta \in \Phi^{[k]}$, we see that $\alpha = \alpha_1$.
Set $\theta = \alpha_1 + \alpha_2$.    Then $[e_{-\alpha},e_{-\theta}] \neq 0$, \  $\Lambda(h_{\alpha+\theta})
= \frac{1}{2}\alpha_2(h_{2\alpha_1+\alpha_2}) = 0$, and $\theta(h_{\alpha+\theta})
= (\alpha_1+\alpha_2)(h_{2\alpha_1+\alpha_2}) = 1$.   Therefore,
$(\Lambda-\theta)(h_{\alpha+\theta}) = -1  \not \in \{0,1\}$.

Thus, we have proved that $(2 \Lambda - \be)(\mathfrak{t}^{[k]})
\neq 0$, that $f^{\Lambda'}$ is a $\mathfrak b^+$-primitive vector, and that
$e^{\Gamma'}$ is a $\mathfrak b^-$-primitive vector.    Hence,  after scaling $e^{\Gamma'}$ by
a suitable scalar factor if necessary so that $[f^{\Lambda'},e^{\Gamma'}] =
e_{-\alpha-\beta}$, we have that  the triple
$(f^{\Lambda'},\mathfrak{g}_0,e^{\Gamma'})$ satisfies all the
conditions of Theorem \ref{Thm:4.7}.  Repeating the argument with
this new triple, we will eventually arrive at the highest root
of $\Phi^{[k]}$.  Thus, we may assume that $\alpha + \be$ is the
highest root of $\Phi^{[k]}$, and that $[f^{\Lambda'},e^{\Gamma'}] =
e_{-\alpha-\beta}$.     Applying the previous part of the
proof to the triple $(f^{\Lambda'},\mathfrak{g}_0,e^{\Gamma'})$,
we obtain that $\Phi = \Phi^{[k]}$, and either $\Phi =$ A$_m$,
$\Lambda' \in \{\varpi_1,\varpi_m\}$, or else $\Phi =$ C$_m$,
$\Lambda' = \varpi_1$.

We adopt the notation in \eqref{eq:4.2}
and \eqref{eq:4.3}.   Let $\Phi =$ C$_m$, $m \geq 2$. Then $\alpha
+ \be = 2\e_1$, and since $\be \in \Delta$, $\be = \e_1 - \e_2$ and
$\alpha = \e_1 + \e_2$. Because $\Lambda' = 2 \Lambda - \be = \varpi_1 =
\e_1$ in this case, we obtain $2 \Lambda(h_{\alpha}) = 1$. But
$\Lambda(h_{\alpha}) \in \{0,1\}$ by our previous considerations.
This shows that the case $\Phi =$ C$_m$ cannot occur; i.e., there
cannot be any such $\beta$, so that  $\alpha$ must be the highest
root.

Now suppose that $\Phi =$ A$_m$ and $\Lambda' =
2\Lambda-\be = \varpi_1 = \e_1$.   Since we are assuming $\alpha \in \Phi^+$
is not the highest root, we have  $m\geq 2$.
 As ${\alpha}+\beta = \alpha_1 +
\cdots + \alpha_m$, \;$\beta \in \Delta$, and $\alpha \in \Phi$, it must be
$\beta = \alpha_1$ or $\beta = \alpha_m$.  Then  $\Lambda(h_{\alpha}) \in \{0,1\}$ forces  $\beta
= \alpha_m$, \ $\Lambda(h_{{\alpha}+\beta}) = 1$, and
$\Lambda(h_{\alpha}) = 0$. Set

$$\begin{gathered} x_{{\lambda}'} =
[[[e_{-\be},f^\Lambda],f^\Lambda],f^\Lambda], \quad
\quad x_{\gamma'} = [[[[e^{\Gamma},e_{\alpha}],e_{\alpha}],e^\Gamma],e^\Gamma],  \nonumber \\
e_{-\delta'} = -8 [e_{-\be},e_{-{\alpha}}]. \nonumber \\
\end{gathered}\nonumber$$

\noindent Then
\begin{eqnarray*}
[x_{{\lambda}'},x_{\gamma'}]  &=&
 \Big( \sum_{j=0}^3 (-1)^{j} {3 \choose j}(F^\Lambda)^{j}E_{-\beta}(F^\Lambda)^{3-j}\Big)
\big([[[[e^\Gamma,e_{\alpha}],e_{\alpha}],e^\Gamma],e^\Gamma]\big) \\
&=& E_{-\be}(F^\Lambda)^3\Big([[[[e^\Gamma,e_{\alpha}],e_{\alpha}],e^\Gamma],e^\Gamma]\Big) \\
&=&
E_{-\be}(F^\Lambda)^2\Big([[[[e_{-{\alpha}},e_{\alpha}],e_{\alpha}],e^\Gamma],e^\Gamma]
+ 2[[[[e^\Gamma,e_{\alpha}],e_{\alpha}],e_{-{\alpha}}],e^\Gamma]\Big)  \nonumber \\
&=& (-4\Gamma-3{\alpha})(h_{\alpha})E_{-\be}(F^\Lambda)^2([[e_{{\alpha}},e^\Gamma],e^\Gamma])\\
&=& (2\Gamma)(h_{\alpha})(4\Lambda+{\alpha})(h_{\alpha})[e_{-\be},e_{-{\alpha}}] =
-8[e_{-\be},e_{-{\alpha}}], \end{eqnarray*}

\noindent because $\Lambda + \Gamma = -{\alpha}$ and
$\Lambda(h_{\alpha}) = 0$. The relations
$[x_{{\lambda}'},e_{\delta'}] = [x_{\gamma'},e_{-\delta'}] = 0$
are easy to verify.  Applying Lemma \ref{Lem:4.5}, we get that
$(3\Lambda - \be)(h_{{\alpha}+\beta}) \in \{0,1\}$. But it follows
from the above calculation that $(3\Lambda - \beta)(h_{{\alpha}+\beta}) = 2$.
This shows that the case we are considering ($\Phi =$ A$_m$,
$\Lambda' = \varpi_1$) cannot occur.  The case $\Phi =$ A$_m$, $m
\geq 2$, $\Lambda' = \varpi_m$ can be handled similarly.
This completes the proof of Theorem \ref{Thm:4.7}.
\qed  \m

\section {\ Additional assumptions  \label{sec:4.4}}  
\m 

In the remainder of the chapter,  we will often impose
some or all of the following  conditions in addition to our blanket
assumptions (1)-(2).   \m

\begin{itemize}
\item[{\rm (3)}] $\mathfrak{g}_{-1}$ is an irreducible
$\g_0$-module;
\item[{\rm(4)}]  $\g$ is transitive \eqref{eq:1.3};
\item[{\rm (5)}]  $\mathfrak{g}$ is 1-transitive \eqref{eq:1.4};
\item[{\rm(6)}]  $\g_1$ is a faithful $\g_0$-module.

\end{itemize} \m

\begin{Rem}\label{Rem:4.10} {\rm
When $\g$ is a graded Lie algebra  satisfying constraints (1)-(4) along with (6),  then
$[\g_1,\g_{-1}] \neq 0$ by
transitivity.  Since $\g_0$ is assumed to act faithfully on  $\g_1$, we have   $[[\g_{1},\g_{-1}],\g_1] \neq 0$.
Thus all the hypotheses of Theorem \ref{Thm:1.48} hold,  so $\g_{-1}$ is a restricted $\g_0^{(1)}$-module.
Condition (5) (1-transitivity)  implies (6), so whenever constraints (1)-(5)  are assumed
(as in the statement of the Main Theorem,  for example),
$\g_{-1}$ must be a restricted $\g_0^{(1)}$-module.} \end{Rem}

\m
\section {\ Computing weights of $\mathfrak b^-$-primitive vectors in $\boldsymbol{\g_1}$ 
\label{sec:4.5}} 
\m

    In Remark \ref{Rem:4.8},  we have argued
    that when $\g$ is transitive and $\g_{-1}$ is irreducible,
  one of the hypotheses of Theorem \ref{Thm:4.7} holds automatically--namely,  
   $\Lambda(\mathfrak t^{[i]}) \neq 0$ for all $i$  whenever $\Lambda$ is the weight of a $\mathfrak b^+$-primitive vector of $\g_{-1}$.    Here we
    show that under a few additional assumptions  on the $\g_0$-module $\g_{-1}$, the
    other assumption in Theorem \ref{Thm:4.7} holds as well;
$\Gamma(\mathfrak t^{[i]}) \neq 0$ for all $i$ whenever
$\Gamma$ is the weight of a $\mathfrak b^-$-primitive vector of
$\g_1$.     \bi

\begin{Thm} \label{Thm:4.12}  \  Let  $\g$  be a graded Lie algebra
satisfying conditions (1)-(4), and assume $\mathfrak{g}_{-1}$ is a
restricted $\g_0^{(1)}$-module with  a $\mathfrak b^+$-primitive vector of
weight $\Lambda$.
Suppose also that $\mathfrak{g}_1$ contains an irreducible
$\g_0$-submodule generated by a $\mathfrak b^-$-primitive vector $e^{-\Lambda}$ corresponding to the weight $-\Lambda \in \mathfrak{t}^*$.
Then the following hold: \begin{itemize} \item[{\rm (i)}]
$\Gamma(\mathfrak{t}^{[i]}) \neq 0$, $i = 1, \dots, \ell$ whenever 
$\g_1$ has a $\mathfrak b^-$-primitive vector of 
weight $\Gamma$; \item[{\rm (ii)}]   If all $\mathfrak b^-$-primitive vectors
of $\g_1$ have weight
$-\Lambda$, then $\g_1$ is an
irreducible $\g_0^{(1)}$-module.
\end{itemize}
\end{Thm}

\pf Given a linear function $\psi$ on $\tf$ we denote by
$\psi^{[i]}$ the restriction of $\psi$ to
$\tf^{[i]}\cap\g_0^{(1)}$. Assume $e^\Gamma \in \g_1$ is a 
$\mathfrak b^-$-primitive vector of
weight $\Gamma$  with $\Gamma(\mathfrak t^{[k]}) = 0$ for some $k$.
 As $\mathfrak{g}$ is transitive, we can regard $\mathfrak{g}_1$ as
a submodule of $\Hom(\mathfrak{g}_{-1},\mathfrak{g}_0)$ $\cong
\mathfrak{g}_{-1}^* \ot \mathfrak{g}_0$.   In particular,
$\mathfrak{g}_1$ is a restricted $\g_0^{(1)}$-module. Let
$f^\Lambda$ be a $\mathfrak b^+$-primitive vector of
$\mathfrak{g}_{-1}$ corresponding to the weight $\Lambda$. By
transitivity and irreducibility, we have from  Remark
\ref{Rem:4.8} that $\Lambda(\mathfrak t^{[i]}) \neq 0$ for any
$i$,  so $\Gamma \neq -\Lambda$. Then by transitivity,
$[f^\Lambda, e^\Gamma] = \zeta  e_{-\alpha}$ for some $\zeta  \in
\F^\times$ and some $\alpha \in \Phi$.    It follows that
$\Lambda^{[k]}$ is a root of $\big(\g_0^{[k]}\big)^{(1)}$.

According to Theorem \ref{Thm:3.25}, the subalgebra
$\widetilde{\mathfrak{g}}$ generated by $f^\Lambda$,
$\mathfrak{g}_0$, and $e^{-\Lambda}$ modulo its Weisfeiler radical
$\mathcal M(\widetilde{\mathfrak{g}})$  is isomorphic  to a
classical Lie algebra with a standard grading. Let $\widetilde
\Delta = \Delta\cup\{-\Lambda\}$ be the canonical base for the
root system of the classical Lie algebra
$\widetilde{\mathfrak{g}}/\mathcal M(\widetilde{\mathfrak{g}})$.
Using the Dynkin diagrams in \cite{Bou1}, it is straightforward to
see that ${\Lambda}^{[i]}$ is a minuscule weight (see Section
\ref{sec:2.2}) of the root system $\Phi^{[i]}$ except in the
following cases:

\bi

(I) $\Phi = \Phi^{[1]}$ $\cup$ $\Phi^{[2]}$ where $\Phi^{[1]} =$
A$_{m-1}$, $\Lambda^{[1]} = \varpi_{m-1}^{[1]}$, $\Phi^{[2]} =$
A$_1$, $\Lambda^{[2]} = 2\varpi_1^{[2]}$, and the resulting root
system is of type B$_{m+1}$,

\begin{equation}\label{eq:4.13}
{\beginpicture \setcoordinatesystem
units <0.45cm,0.3cm> % sets scale
 \setplotarea x from -4 to 16, y from -1 to 1    % sets plot size up
 \linethickness=0.15pt                          % sets line thickness
  \put{$\circ$} at 0 0
 \put{$\circ$} at 2 0
 \put{$\cdots$} at 5 0
 \put{$\circ$} at  8 0
 \put{$\bullet$} at 10  0 \put{$\circ$} at  12 0
 \plot .15 .1 1.85 .1 /
 \plot 2.15 .1 3.85 .1 /
 \plot 6.15 .1 7.85 .1 /
 \plot 10.15 -.2 11.85 -.2 /
  \plot 10.15 .2 11.85  .2 /
  \put {$>$} at 11 0
 \plot 8.15 .1  9.85 .1 /
 \put{$\alpha_1$} at 0 -1.5
  \put{$\alpha_2$} at 2 -1.5
    \put{$\alpha_{m-1}$} at 8 -1.5
      \put{$-\Lambda$} at 10 -1.5
        \put{$\alpha_m$} at 12 -1.5
        \put {$(m \geq 1)$} at 16 0  \endpicture}\end{equation}

\smallskip (II) $\Phi =\Phi^{[1]}=$ A$_m$, $\Lambda^{[1]} = 2\varpi_m$, and
the resulting root system is of type C$_{m+1}$,

\begin{equation}\label{eq:4.14}
{\beginpicture \setcoordinatesystem
units <0.45cm,0.3cm> % sets scale
 \setplotarea x from -4 to 12, y from -1 to 1    % sets plot size up
 \linethickness=0.15pt                          % sets line thickness
  \put{$\circ$} at 0 0
 \put{$\circ$} at 2 0
 \put{$\cdots$} at 5 0
 \put{$\circ$} at  8 0
 \put{$\circ$} at 10  0 \put{$\bullet$} at  12 0
 \plot .15 .1 1.85 .1 /
 \plot 2.15 .1 3.85 .1 /
  \plot 6.15 .1 7.85 .1 /
  \plot 10.15 -.2 11.85 -.2 /
  \plot 10.15 .2 11.85  .2 /
  \put {$<$} at 11 0   \plot 8.15 .1  9.85 .1 /
  \put{$\alpha_1$} at 0 -1.5
  \put{$\alpha_2$} at 2 -1.5
      \put{$\alpha_{m-1}$} at 8 -1.5
      \put{$\alpha_m$} at 10 -1.5
        \put{$-\Lambda$} at 12 -1.5
        \put {$(m \geq 2)$} at 16 0  \endpicture}
\end{equation}

\smallskip (III) $\Phi = \Phi^{[1]}$ $\cup$ $\Phi^{[2]}$ where $\Phi^{[1]} =$A$_1$,
$\Lambda^{[1]} = \varpi_{1}^{[1]}$, $\Phi^{[2]} =$ A$_2$,
$\Lambda^{[2]} = 2\varpi_1^{[2]}$, and the resulting root system
is of type F$_4$,

\begin{equation}\label{eq:4.15}
{\beginpicture \setcoordinatesystem
units <0.45cm,0.3cm> % sets scale
 \setplotarea x from -6 to 12, y from -1 to 1    % sets plot size up
 \linethickness=0.03pt                          % sets line thickness
  \put{$\circ$} at 0 0
 \put{$\bullet$} at 2 0
  \put{$\circ$} at 4 0
\put{$\circ$} at 6 0
 \plot .15 .1 1.85 .1 /
 \plot 2.15 -.2 3.85 -.2 /
  \plot 2.15 .2 3.85 .2 /
  \plot 4.15 .1 5.85 .1 /
  \put {$>$} at 3 0
  \put{$\alpha_1$} at 0 -1.5
  \put{$-\Lambda$} at 2 -1.5
  \put{$\alpha_2$} at 4 -1.5
  \put{$\alpha_3$} at 6 -1.5
  \endpicture}\end{equation}

 \smallskip (IV) $\Phi=\Phi^{[1]} =$ A$_1$, $\Lambda^{[1]} = 3\varpi_1$, and the resulting
root system is of type G$_2$,

\begin{equation}\label{eq:4.16}
{\beginpicture \setcoordinatesystem
units <0.45cm,0.3cm> % sets scale
 \setplotarea x from -6 to 8, y from -1 to 1    % sets plot size up
 \linethickness=0.03pt                          % sets line thickness
  \put{$\circ$} at 0 0
 \put{$\bullet$} at 2 0
 \plot .15 -.3 1.85 -.3 /
  \plot .15 0 1.85  0 /
  \plot .15 .3 1.85 .3 /
  \put {$<$} at 1 0
   \put{$\alpha_1$} at 0 -1.5
  \put{$-\Lambda$} at 2 -1.5 \endpicture}\end{equation}

\smallskip (V)  $\Phi=\Phi^{[1]} =$ C$_3$, $\Lambda^{[1]} = \varpi_3$, and the
resulting root system is of type F$_4,$

\begin{equation}\label{eq:4.18}
{\beginpicture \setcoordinatesystem
units <0.45cm,0.3cm> % sets scale
 \setplotarea x from -6 to 12, y from -1 to 1    % sets plot size up
 \linethickness=0.15pt                          % sets line thickness
  \put{$\circ$} at 0 0
 \put{$\circ$} at 2 0
  \put{$\circ$} at 4 0
\put{$\bullet$} at 6 0
 \plot .15 .1 1.85 .1 /
 \plot 2.15 -.2 3.85 -.2 /
  \plot 2.15 .2 3.85 .2 /
  \plot 4.15 .1 5.85 .1 /
  \put {$<$} at 3 0
  \put{$\alpha_1$} at 0 -1.5
  \put{$\alpha_2$} at 2 -1.5
  \put{$\alpha_3$} at 4 -1.5
  \put{$-\Lambda$} at 6 -1.5
  \endpicture}\end{equation}

\bi If $\Lambda^{[k]}$ is a minuscule weight, then
$\Lambda^{[k]}(h_\gamma) \in \{-1,0,1\}$ for all $\gamma\in
\Phi^{[k]}$. Since $\gamma(h_\gamma) = 2$, this implies that
$\Lambda^{[k]}$ it is not a root of $\big(\g_0^{[k]}\big)^{(1)}$,
contrary to our earlier remarks. Thus, $\Lambda^{[k]}$ is not
minuscule, and hence it remains to handle Cases (I)--\,(V).

We can rule out Case (II),  since here  $\Phi = \Phi^{[1]}$ =
A$_m$ and $\Lambda^{[1]} = 2 \varpi_m$ is not a root of
$\big(\g^{[1]}\big)^{(1)}$ for $m \geq 2$. In Case (III), $\Phi =
\Phi^{[1]} \cup \Phi^{[2]}$, where $\Phi^{[1]}$ = A$_1$ and
$\Phi^{[2]}$ = A$_2$, and $\Lambda^{[1]} = \varpi_1^{[1]}$,
$\Lambda^{[2]} = 2 \varpi_1^{[2]}$.    Hence, $\Lambda^{[i]}$ is
not a root of $\big(\g^{[i]}\big)^{(1)}$ for $i=1,2$. Likewise in
Case (V), $\Phi = \Phi^{[1]}$ = C$_3$, and $\varpi_3$ is not a
root of $\big(\g^{[1]}\big)^{(1)}$.

In Case (I), if $m \geq 2$, then $\Phi = \Phi^{[1]} \cup
\Phi^{[2]}$, where $\Phi^{[1]}$ = A$_{m-1}$, $\Phi^{[2]}$ = A$_1$,
\ $\Lambda^{[1]} = \varpi_{m-1}^{[1]}$, and $\Lambda^{[2]} = 2
\varpi_1^{[2]}$.     If the root $\alpha$ in the first paragraph of the
proof belongs to
$\Phi^+$, then since $\Lambda^{[i]} \neq 0$ for $i=1,2$, Theorem
\ref{Thm:4.7} applies and gives that $\Phi$ is irreducible.  Since
this is false, we must have $\alpha \in \Phi^-$. If
$\Gamma(\mathfrak t^{[1]}) = 0$, then  $\Lambda^{[1]} $ is a root
of  $\big(\g^{[1]}\big)^{(1)}$, contradicting the fact that
$\Lambda^{[1]}$ is a minuscule weight for $\Phi^{[1]}$. Hence, we
may assume $\Gamma(\mathfrak t^{[2]}) = 0$,  so that $\alpha =
-\alpha_1^{[2]}$.    (It is easy to check that   
$\alpha =
-\alpha_1^{[2]}$ also holds when $m=1$.)   But then

$$[e^\Gamma, [f^\Lambda,[f^\Lambda, e_{\alpha}]]] = -2 \zeta  [f^\Lambda, [e_{-\alpha}, e_{\alpha}]]
= 4 \zeta f^\Lambda \neq 0,$$

\noindent where $\zeta$ is
as in the first paragraph of the proof,  which implies $[f^\Lambda,[f^\Lambda, e_{\alpha}]] \neq 0$.  On the other hand,
$f^\Lambda$ corresponds to a long root of $\widetilde \g_0$, and
  $(\ad f^\Lambda)^2 (\widetilde \g_0) = \F  f^\Lambda \subseteq \g_{-1}$.
As $[[f^\Lambda,[f^\Lambda, e_{\alpha}]] \in \g_{-2}$, this shows that Case (I) cannot occur.

In Case (IV),  we have  $\Phi=\Phi^{[1]} =$ A$_1$, $\Lambda^{[1]}
= 3\varpi_1$, and the resulting root system is of type G$_2$. Then
$\alpha = \pm \alpha_1$.
 If $\alpha = \alpha_1$, then Theorem \ref{Thm:4.7}
gives $\Lambda = \varpi_1$ (the highest weight of the natural
module), which is false.    Hence $\alpha = -\alpha_1$.  But in
this case, $\Gamma^{[1]} = \alpha_1 - \Lambda = -\varpi_1$, and
$\Gamma(\mathfrak{t}) \neq 0$. This completes the proof of part
(i).

\smallskip

{F}rom now on we assume that all $\mathfrak b^-$-primitive
vectors of $\g_1$ have weight  $-\Lambda$. 
Set $L=\widetilde{\g}/{\mathcal M}(\widetilde{\g})$.
It follows from Theorem \ref{Thm:3.25} that
$\widetilde{G}\eqdef\text{Aut}(L)^\circ$ is a simple algebraic
group (of adjoint type), and $\ad L\cong L$ is a
$\widetilde{G}$-stable ideal of codimension $\le 1$ in
$\text{Lie}(\widetilde{G})=\Der(L)$ (the equality
$\text{Lie}(\widetilde{G})=L$ holds if and only if
$L\not\cong\mathfrak{psl}_{kp}$). There exists a homomorphism of
algebraic groups $\phi \colon\,\F^\times \rightarrow
\widetilde{G}$ such that for all $i\in\Z$ we have
$$L_i=\big(\widetilde{\g}/{\mathcal M}(\widetilde{\g})\big)_i=\{x\in
L \mid \phi(t)(x)=t^i\, x\ \,\text{for all}\ \,\,
t\in\F^\times\}$$ (see the proof of Theorem \ref{Thm:3.25} for
more detail). Let $\widetilde{G}_0$ be the centralizer  in
$\widetilde{G}$ of the one-dimensional torus $\phi(\F^\times)$.
This is a connected, reductive subgroup of $\widetilde{G}_0$, and
it is immediate from \cite[Sec. 9.4]{Bo} that
$\text{Lie}(\widetilde{G}_0)$ consists of all derivations of $L$
commuting with $\phi(\F^\times)$. Since $\phi(\F^\times)$ acts
trivially on $\text{Lie}(\widetilde{G})/\ad L$, the Lie algebra
$\g_0\cong \big(\widetilde{\g}/{\mathcal
M}(\widetilde{\g})\big)_0\cong \text{ad}_L(L_0)$ can be identified
with a $\widetilde{G}_0$-stable ideal of codimension $\le 1$ in
$\text{Lie}(\widetilde{G}_0)$. Moreover, $
(\text{Lie}(\widetilde{G}_0))^{(1)}\cong \g_0^{(1)}$ as Lie
algebras. Identifying the $\g_0$-module $L_{-1}$ with $\g_{-1}$,
we see that the adjoint action of $\g_0$ on $\g_0$ and $\g_{-1}$
is induced by the differential of the (rational) action of
$\widetilde{G}_0$ on $L$. There exist maximal unipotent subgroups
$N^\pm$ and a maximal torus $\widetilde{T}$ in $\widetilde{G}_0$
such that ${\rm Lie}(N^\pm) = \mathfrak{n}^\pm$ and ${\rm
Lie}(\widetilde{T})\supseteq {\mathfrak t}$.

Let $G_0$ be the derived subgroup of $\widetilde{G}_0$, and set
$T=\widetilde{T}\cap G_0$. It follows from the general theory of
algebraic groups that $G_0$ is semisimple (but not necessarily
simply connected), that $T$ is a maximal torus of $G_0$, and that
$\widetilde{G}_0=G_0\cdot Z(\widetilde{G}_0)^\circ$; see
\cite{Bo}, for example. Note that $N^\pm\subset G_0$, and
$\g_0^{(1)}=(\text{Lie}(\widetilde{G}_0))^{(1)}=(\text{Lie}(G_0))^{(1)}.$
Let $\widehat{G}_0$ be the simply connected cover of $G_0$. It is
well-known that there exists a surjective homomorphism of
algebraic groups $\iota\colon\,\widehat{G}_0\rightarrow G_0$ whose
kernel is finite (and central in $\widehat{G}_0$) and whose
restriction to any maximal unipotent subgroup $\widehat{U}$ of
$\widehat{G}_0$ induces an isomorphism of algebraic groups
$\widehat{U}\stackrel{\sim}{\rightarrow}\iota(\widehat{U})$. It is
straightforward to see that $\ker
(\mathrm{d}\iota)_e\subseteq{\mathfrak
Z}\big(\text{Lie}(G_0)\big)$. The inverse image $\widehat{T}\eqdef
\iota^{-1}(T)$ is a maximal torus of $\widehat{G}_0$. Since
$\widehat{G}_0$ is simply connected, the Lie algebra
$\text{Lie}(\widehat{G}_0)$ is generated by root vectors relative
to $\widehat{T}$ ( this is due to the fact that $\widehat{T}$ has
the smallest possible group of rational cocharacters). As
$\mathfrak{sl}_2$ is simple for $p>2$, we then have
$(\text{Lie}(\widehat{G}_0))^{(1)}=\text{Lie}(\widehat{G}_0)$. In
conjunction with our earlier remarks, this shows that
$(\mathrm{d}\iota)_e$ maps $\text{Lie}(\widehat{G})$ onto
$\g_0^{(1)}$. As a consequence, the adjoint action of $\g_0^{(1)}$
on $\g_{-1}$ and $\g_0$ is induced by the differential of a
rational action on $\g_{-1}\oplus \g_0$ of the semisimple, simply
connected group $\widehat{G}_0$.
 We retain our notation associated with the
lattice of weights of $\widehat{T}$.

To prove that $\g_1$ is irreducible, we identify it with a
$\g_0^{(1)}$-submodule of $\text{\rm Hom}(\g_{-1},$ $\g_0)$. By
the preceding remark, the action of $\g_0$ on $\Hom(\g_{-1},$
$\g_0)$ $\cong \g_{-1}^*$ $\ot$ $\g_0$ is induced by the
differential of the natural rational action of $\widehat{G}_0$ on
$\Hom(\g_{-1},\g_0).$ Since $\iota(\widehat{G}_0)=G_0$ and
$\widetilde{G}_0=G_0\cdot Z(\widetilde{G})^\circ$, the
$\widehat{G}_0$-module $\g_{-1}$ is irreducible.  We denote by
$\lambda$ the maximal weight of the representation of
$\widehat{G}_0$ on $\g_{-1}$. By the above discussion,
$(\mathrm{d}\iota)_e$ sends $\text{Lie}(\widehat{T})$ onto
$\tf\cap\g_0^{(1)}$. To simplify notation we will identify
$\tf\cap\g_0^{(1)}$ with
$\text{Lie}(\widehat{T})/\ker(\mathrm{d}\iota)_e$. As
$\widehat{G}_0$ is simply connected, the differential map
$\mathrm{d}\colon\,X(\widehat{T})\rightarrow
\text{Lie}(\widehat{T})^*,\ \,\nu\mapsto (\mathrm{d}\nu)_e,\,$ has
kernel $p\,X(\widehat{T})$ and induces an isomorphism
$\big(X(\widehat{T})/p\,X(\widehat{T})\big)\otimes_{\Z}\F\,
\stackrel{\sim}{\longrightarrow} \text{Lie}(\widehat{T})^*$. As
$\mathrm{d}\lambda$ vanishes on $\ker(\mathrm{d}\iota)_e$, we may
regard it as a linear function on $\tf\cap\g_0^{(1)}$.
 Note that $\lambda\in X_1(\widehat{T})$, because $\g_{-1}$ is
 an irreducible $\g_0^{(1)}$-module (see Remark \ref{irr} and
 Proposition \ref{Pro:2.803}). Since $\mathrm{d}\lambda$ is the weight of 
 a $\mathfrak b^+$-primitive vector in the
 $\g_0^{(1)}$-module $\g_{-1}$, it must coincide with the
 restriction of $\Lambda$ to $\tf\cap\g_0^{(1)}$.

Suppose $\Lambda^{[k]}$ is a minuscule weight for all
$k\in\{1,\dots,\ell\}$. Then $\lambda$ is a minuscule weight of
$X(\widehat{T})$, which means that
$X\big(\Hom(\g_{-1},\g_0)\big)\,= \,X'\cup X''$ where $X'=
\{-w\lambda \mid w\in W\}$ and $X''= \{-w \lambda+ \gamma \mid 
w\in W,\,\gamma\in \Phi\}$ (here $W=
N_{\widehat{G}_0}(T)/Z_{\widehat{G}_0}(T)$ is the Weyl group of
$\widehat{G}_0$). As $\Phi$ does not have any components of type
$\mathrm{G}_2$, this implies that
$\displaystyle{\frac{2(\nu,\gamma)}{(\gamma,\gamma)}}\in \{0, \pm
1,\pm 2,\pm 3\}$ for any $\nu\in X'\cup X''$ and $\gamma\in \Phi$.
As $p$ $> 3$, it follows that all dominant weights of the
$\widehat{G}_0$-module $\Hom(\g_{-1},\g_0)$ belong to
$X_1(\widehat{T})$.  Hence Proposition \ref{Pro:2.804} applies and
we obtain that $\g_1$ is $\widehat{G}_0$-stable.

Let $\mu$ be a minimal weight of the $\widehat{G}_0$-module
$\g_1$. Then $w_0\mu$ is a dominant weight of $\g_1$, hence it
belongs to $X_1(\widehat{T})$.     It follows that $\mu\in
-X_1(\widehat{T})$. As the weight component $\g_1^\mu$ consists of
$\mathfrak b^-$-primitive vectors, the weight $\mathrm{d}\mu$
coincides with the restriction of $-\Lambda$ to
$\tf\cap\g_0^{(1)}$. As $-\lambda\in -X_1(\widehat{T})$ too,  it
must be that $\mu=-\lambda$.

Let $\nu \in X(\g_1)$.  Then there exists $w \in W$ such that $w
\nu \in -X(T)_+$.  Because $X(\g_1)$ is $W$-stable, and $-\lambda$
is the only minimal weight of $\g_1$, we have that $w \nu \geq
-\lambda$. But then $\lambda \geq -w \nu$ with $-w \nu \in
X(T)_+$.  As $\lambda$ is a minuscule weight, this implies $w\nu =
-\lambda$. Hence all weights of $\g_1$ are conjugate under $W$. If
$e^\mu$ and $e^{-\Lambda}$ are linearly independent vectors of
$\g_1^\mu$, then $[f^\Lambda,[te^{-\Lambda}-e^\mu, f^\Lambda]]
\neq 0$ for any $t \in {\mathbb F}$, as  $t e^{-\Lambda}-e^\mu$ is
a $\mathfrak b^-$-primitive vector, $[t e^{-\Lambda}-e^\mu,f^\Lambda] \neq 0$, 
and the triple $\{f^\Lambda, \g_0, t e^{-\Lambda} - e^\mu\}$
satisfies all the conditions of Theorem \ref{Thm:3.20}).  But
$[\g_1^\mu, f^\Lambda] \subseteq \mathfrak t$ and $[\mathfrak
t,f^\Lambda] = {\mathbb F} f^\Lambda$.  This implies that
$[f^\Lambda, [t_0 e^{-\Lambda}-e^\mu,f^\Lambda]] = 0$ for a
suitable $t_0 \in F$. Consequently, $\dim \g_1^\mu = 1$.

Summarizing, we have that $\g_1$ is an irreducible
$\widehat{G}_0$-module with maximal weight $-w_0\lambda$.  As
$-w_0\lambda \in X_1(T)$, Proposition \ref{Pro:2.803} shows that
$\g_1$ remains irreducible under the action of
$(\mathrm{d}\iota)_e\big(\hbox{\rm
Lie}(\widehat{G}_0)\big)=\g_0^{(1)}$. Therefore, if
$\Lambda^{[k]}$ is minuscule for any $k \in \{1, \dots, \ell\}$,
then we have the desired conclusion that $\g_1$ is an irreducible
$\g_0^{(1)}$-module.

It remains to consider the case where the classical algebra
$L=\widetilde \g/\mathcal M(\widetilde \g)$ corresponds to one of
the diagrams (I--V) pictured above. In Case (I), we have that
$\Phi = \Phi^{[1]} \cup \Phi^{[2]}$ where $\Phi^{[1]}=
\mathrm{A}_{m-1},\,$ $\Phi^{[2]} =\mathrm{A}_1,\,$ $\Lambda^{[1]}
= \varpi_{m-1}^{[1]},$ and $\Lambda^{[2]} = 2 \varpi_1^{[2]}$.
(When $m= 1$, we assume that $\Phi^{[1]} = \emptyset$ and
$\Lambda^{[1]} = 0$.) Since $\lambda\in X_1(\widehat{T})$, this
yields $\lambda=\varpi_{m-1}^{[1]}+2\varpi_1^{[2]}$, where
$\varpi_{m-1}^{[1]}$ and $\varpi_1^{[2]}$ are now viewed as
fundamental weights in $X(\widehat{T})$. In our situation,
$\g_{-1}\cong L(\lambda) \cong L\big(\varpi_{m-1}^{[1]}\big)
\otimes L\big(2 \varpi_m^{[2]}\big)$, and $\Phi$ has the following
Dynkin diagram:

\begin{equation}\label{eq:4.2205}
{\beginpicture \setcoordinatesystem
units <0.45cm,0.3cm> % sets scale
 \setplotarea x from 0 to 12, y from -2 to 1    % sets plot size up
 \linethickness=0.15pt                          % sets line thickness
  \put{$\circ$} at 0 0
 \put{$\circ$} at 2 0
 \put{$\cdots$} at 5 0
 \put{$\circ$} at  8 0
 \put{$\circ$} at 10  0
 \put{$\circ$} at  12 0
 \plot .15 .1 1.85 .1 /
 \plot 2.15 .1 3.85 .1 /
 \plot 6.15 .1 7.85 .1 /
 \plot 8.15  .1 9.85 .1 /
 %\plot 10.15 .1  11.85 .1 /
 \put{$\alpha_1$} at 0 -1.5
  \put{$\alpha_2$} at 2 -1.5
    \put{$\alpha_{m-2}\,\,$} at 8 -1.5
      \put{$\,\,\,\,\alpha_{m-1}$} at 10 -1.5
        \put{$\,\,\,\,\alpha_m$} at 12 -1.5
          \endpicture}\end{equation}

\bi\noindent Obviously, $\dim L(\lambda) = 3(m+1)$, and
$X(\g_{-1})=W\lambda\cup W\varpi_{m-1}^{[1]}$. As $p>3$, it is
straightforward to see that all dominant weights of $\Hom(\g_{-1},
\g_0)$ belong to $X_1(\widehat{T})$. Applying Proposition
\ref{Pro:2.804} we now obtain that the subspace $\g_1$ is
$\widehat{G}_0$-stable. The above reasoning then shows that
$-\lambda$ is the only minimal weight of the
$\widehat{G}_0$-module $\g_1$, and $\dim \g_1^{-\lambda} = 1$.
This, in turn, yields that $\lambda$ is the only maximal weight of
the dual $\widehat{G}_0$-module $\g_1^*$, and $\dim
(\g_1^*)^{\lambda} = 1$. Let $M$ be the $\widehat{G}_0$-submodule
of $\g_1^*$ generated by a nonzero element $\xi \in
(\g_1^*)^\lambda$. If $M \neq \g_1^*$, then
$$M^\perp = \{ x \in \g_1 \mid  \psi(x) = 0 \ \, \text{ for all
}\ \, \psi \in M\}$$ is a nonzero $\widehat{G}_0$-submodule of
$\g_1$. As $M^\perp$ is $\text{Lie}(\widehat{G}_0)$-stable, it
contains a $\mathfrak b^-$-primitive vector.   By our
assumption, this forces $M^\perp \cap \g_1^{-\Lambda} \neq 0$. As
$\g_1^{-\Lambda} = \g_1^{-\lambda}$, we must have that
$\g_1^{-\lambda} \subset M^\perp$.  It follows that
$\xi(\g_1^{-\lambda}) = 0$.  But $\xi(\g_1^\nu) = 0$ for all $\nu
\neq -\lambda$ by the definition of $(\g_1^*)^\lambda$. Therefore,
$\xi = 0$ contrary to our assumption.

As a result, $\g_1^*$ is generated by a $\mathfrak b^+$-primitive vector of
weight $\lambda$.  By Proposition \ref{Pro:2.805}, there exists a
surjective $\widehat{G}_0$-module homomorphism $\psi: V(\lambda)
\rightarrow \g_1^*$, where $V(\lambda)$ is the Weyl module with
maximal weight $\lambda$.  In our situation, $V(\lambda) \cong
V\big(\varpi_{m-1}^{[1]}\big) \otimes V\big(2
\varpi_1^{[2]}\big)\cong L(\lambda)$. This implies that in Case
(I), the $\g_0^{(1)}$-module $\g_1$ is irreducible. A similar
argument also works in Cases (II) and (III). (In seeing this, one
should keep in mind that for groups of type $\mathrm{A}_m$, the
Weyl modules $V(2 \varpi_1)$ and $V(2 \varpi_m)$ are irreducible
provided $p > 2$.) In Case (V), we have
$\Phi=\Phi^{[1]}=\mathrm{C}_3$ and $\Lambda^{[1]}=\varpi_3^{[1]}$.
It is well-known that in this situation 
$V\big(\varpi_3^{[1]}\big)\cong L\big(\varpi_3^{[1]}\big)$ for
$p>2$, and
$X(\g_{-1})=X\big(L(\varpi_3^{[1]})\big)=W\varpi_3^{[1]}\cup
W\varpi_{1}^{[1]}$. Direct verification shows that for $p>3$ all
dominant weights of $\Hom(\g_{-1},\g_0)$ are in
$X_1(\widehat{T})$.  Reasoning as above we derive again that
$\g_1$ is an irreducible $\g_0^{(1)}$-module.

It remains to consider Case (IV).  In this case, $\g_0\cong
\mathfrak{gl}_2$, $\Lambda^{[1]}=3\varpi_1,\,$ $\mathfrak{n}^- =
{\mathbb F} f$, and $\mathfrak{n}^+ = {\mathbb F} e$. Hence,
$\g_{-1}\cong L(3)$ as $\g_0^{(1)}$-modules (here we adopt
the abbreviated notation $L(k)$ for $L(k\varpi_1)$, $0 \leq k \leq p-1$).  Since $\g_1$ is a
restricted $\g_0^{(1)}$-module, we have $(\ad f)^p(\g_1) = 0,$
while the above reasoning yields $\dim\big(\g_1 \cap \ker \, \ad
f\big) = 1.$ As a consequence, the $\mathfrak{sl}_2$-module $\g_1$
is indecomposable of dimension $\leq p.$ Note that $L(3)\cong
L(3)^*$ and $\g_0\cong L(2)\oplus\F$ as $\g_0^{(1)}$-modules.
Therefore, $$\g_1\,\hookrightarrow\,\Hom(\g_{-1},\g_0)\cong L(3)
\oplus \big(L(2) \ot L(3)\big).$$ If $p>5$, then the
$\mathfrak{sl}_2$-module $L(2)\ot L(3)$ is completely reducible;
see \cite{BO} for example. But then $\g_1$ is irreducible, and we
are done. So we may assume that $p=5$. Recall that the graded
component $\widetilde{\g}_1$ of $\widetilde{\g}$ is an irreducible
$\g_0$-module generated by a $\mathfrak b^-$-primitive vector of weight
$-\Lambda$. In our case, this says that $\widetilde{\g}_1\cong
L(3)$ as $\mathfrak{sl}_2$-modules. If $\dim \g_1<5$, then
$\g_1=\widetilde{\g}_1$, because $\dim L(3)=4$.

 Suppose then that $\dim \g_1 = p=5$. As $\widetilde{\g}_1 \cong L(3)$,
 the quotient module
$\g_1/\widetilde{\g}_1$ is trivial.  As the endomorphism $(\ad
h)\vert_{\g_1}$ is semisimple, $\g_1$ contains a nonzero weight
vector of weight zero, $u_0$ say. By transitivity, we have
$[u_0,\g_{-1}]\ne 0$. Thus, $\g_0$ contains a nonzero weight
vector with weight equal to a weight of $\g_{-1}.$ But we know
that $\g_{-1} \cong L(3)$ and $\g_0 \cong L(0) \oplus L(2).$ The
$h$-eigenvalues of $L(3)$ are $3, 1, -1$, and $-3$, while the
$h$-eigenvalues of $L(0) \oplus L(2)$ are $2$, $0$, and $-2$.
Since $p=5$, it must be that either the bracket of a 
$\mathfrak b^+$-primitive vector $v_3$ of $\g_{-1}$ with $u_0$ is a nonzero scalar
multiple of $f$,  or the bracket of a $\mathfrak b^-$-primitive vector $v_{-3}$
of $\g_{-1}$ with $u_0$ is a nonzero scalar multiple of $e$. For
definiteness, we will assume that we are in the former situation.
The proof in the latter case is similar. We can assume that
$$[v_3,u_0] = f.$$
Since $v_3$ is a $\mathfrak b^+$-primitive vector of $\g_{-1}$, we have
$$h =  [e, f] =  [e, [v_3, u_0]]  =  [v_3, [e,
u_0]].$$

\m

\noindent Set $x_{-3}=[e,u_0]$. Then $x_{-3}$ is an
$h$-eigenvector in $\g_1$ with eigenvalue $2\equiv -3\mod 5$. It follows
that $x_{-3}\in\widetilde{\g}_1$. Since $\widetilde{\g}_1\cong
L(3)$ as $\g_0^{(1)}$-modules, it must be that $[f,x_{-3}]=0.$ We
thus have
$$h = [v_3, x_{-3}],\quad\ \,  [e,v_3]=[f,x_{-3}]=0.$$
But then Theorem \ref{Thm:3.9} shows that the subalgebra generated
by
\begin{eqnarray}
e_1 &\eqdef& e,\quad \quad\ \,\, f_1 \eqdef f,\nonumber\\
e_2&\eqdef&x_{-3},\quad\ \,f_2 \eqdef v_{3},\nonumber
\end{eqnarray}

\smallskip

\noindent is infinite dimensional. This shows that $\g_1 =
\widetilde{\g}_1$, completing the proof. \qed

\m   \smallskip    

\begin{Pro} \label{Pro:4.20}   \  Assume $\g$ is a graded Lie algebra
satisfying assumptions (1)-(4) and (6).   Let $e^\Gamma \in
\mathfrak{g}_1$ be a $\mathfrak b^-$-primitive  vector corresponding to
a weight $\Gamma \neq -\Lambda$. Then $\Gamma(\mathfrak t^{[i]})
\neq 0$ for all $i = 1, \dots, \ell$.
\end{Pro}    

\pf    Observe that  $\g_{-1}$ is a restricted $\g_0$-module by
Remark \ref{Rem:4.10}. The irreducibility and transitivity of $\g$
imply that ${\mathfrak Z}(\g_0)$ acts on $\g_{-1}$ as scalar
operators and this action is faithful. The transitivity of $\g$
allows us to identify the $\g_0$-module $\g_1$ with a submodule of
$\text{Hom}(\g_{-1},\g_0)\,\cong\, \g_{-1}^*\otimes\g_0$. As a
consequence, ${\mathfrak Z}(\g_0)$ acts faithfully on $\g_1$.

 Now suppose $\Gamma(\mathfrak t^{[i]}) = 0$ for some $i$. Without
 loss of generality we may assume $i=1$.
 Since ${\mathfrak Z}(\g_0)$ acts on $\g_1$ as scalar operators,
 we have $\Gamma(z)\ne 0$ for any nonzero $z\in{\mathfrak
 Z}(\g_0)$. Then ${\mathfrak Z}(\g_0)\cap \g_0^{[1]}=0$,
 showing that $\g_0^{[1]}$ is centerless (in particular, nonabelian).
 We denote by
 $\g_{0,1}$ the derived ideal of $\g_0^{[1]}$. The above remark shows that $\g_{0,1}$
 is classical simple.
Let  $$\mathfrak{g}_1 = Q_N \supset Q_{N-1} \supset \cdots \supset
Q_1 \supset Q_{0} = 0$$ be a composition series for the
$\mathfrak{g}_0$-module $\mathfrak{g}_1$, and let $\Gamma_k$
denote  the weight of a $\mathfrak b^-$-primitive vector of $Q_k/Q_{k-1}$. Since $\mathfrak{g}_1$
is a submodule of the restricted $\g_0^{(1)}$-module
$\mathfrak{g}_{-1}^* \ot \mathfrak{g}_0$, all the composition
factors $Q_k/Q_{k-1}$ are irreducible, restricted
$\g_0^{(1)}$-modules. It follows that $Q_k/Q_{k-1}$ is a
completely reducible $\g_{0,1}$-module generated by $\mathfrak b^-$-primitive
vectors of weight
${\Gamma_k}\mid_{\,\mathfrak{t}\,\cap\,\g_{0,1}}$. If
$\Gamma_k(\mathfrak{t}^{[1]}) = 0$ for all $k$, then $Q_k/Q_{k-1}$
is a trivial $\mathfrak{g}_{0,1}$-module for all $k$.  However, if
$Q_{k-1}$ and $Q_k/Q_{k-1}$ are trivial
$\mathfrak{g}_{0,1}$-modules, then so is $Q_k$.  Indeed,
 for any  $v \in Q_k$,  we have   $x. v \in Q_{k-1}$ for all $x
\in \mathfrak{g}_{0,1}$.  But then $[x,y]v = x.(y . v) - y.(x.v) =
0$ for all $x,y \in \mathfrak{g}_{0,1}$  so that $Q_k$ is a
trivial $\mathfrak{g}_{0,1}$-module too, because $\g_{0,1}$ is a
simple Lie algebra. This would lead to $\mathfrak{g}_1$ being a
trivial $\mathfrak{g}_{0,1}$-module, which cannot happen since
$\mathfrak{g}_1$ is assumed to be a faithful
$\mathfrak{g}_0$-module.   Consequently,
 $\Gamma_k(\mathfrak{t}^{[1]}) \neq 0$ for some
$k$.    We may suppose that  $k$ is chosen so
$\Gamma_l(\mathfrak{t}^{[1]}) = 0$ for all $l < k$ but  $\Gamma_k(\mathfrak{t}^{[1]}) \neq 0$.

 Let  $e^{\Gamma_k}$ be a weight vector of $Q_k$ whose image in
$Q_k/Q_{k-1}$ is a nonzero $\mathfrak b^-$-primitive vector.  We claim that
$e^{\Gamma_k}$ itself is a $\mathfrak b^-$-primitive vector.  First note that
$Q_{k-1}$ is a trivial $\g_0^{[1]}$-module, so that
$\mu(\mathfrak{t}^{[1]}) = 0$ for any weight $\mu$ of the
$\mathfrak{g}_0$-module $Q_{k-1}$ (an immediate consequence of our
choice of $k$).  This implies that $[e_\theta, e^{\Gamma_k}] = 0$
for any $\theta \in (\Phi^{[i]})^-$ with $i > 1$.   Suppose that
$u \eqdef [e_{-\alpha_j},e^{\Gamma_k}] \neq 0$ for some simple
root  $\alpha_j$ of $\Phi^{[1]}$.  As $u \in Q_{k-1}$ is a weight
vector corresponding to the weight $\Gamma_k -\alpha_j$, we have
$(\Gamma_k-\alpha_j)(\mathfrak{t}^{[1]}) = 0$.  Further,
$[e_\theta,u] = [e_{-\alpha_j},[e_\theta,e^{\Gamma_k}]] = 0$ for
any $\theta \in (\Phi^{[i]})^-$ with $i > 1$.     Thus,
$[(\mathfrak{n}^{[i]})^-,u] = 0$ for all $i =1, \dots,\ell$, which is to say that $u$ is a
${\mathfrak b}^{-}$-primitive vector.   Then it follows that 
$[f^\Lambda,u] \neq 0$ by the transitivity  of $\mathfrak{g}$.  As
$\Lambda(\mathfrak t^{[1]}) \neq 0$ by transitivity (compare
Remark \ref{Rem:4.8}), we have $[f^\Lambda,u] = \zeta e_\beta$ for
some $\zeta \in {\mathbb F}^\times$ and $\be \in \Phi^{[1]}$.
Since  $u \in \g_{0}^{[1]}$, which is centerless,
and $[e_{{\alpha}},u] = 0$ for any $\alpha \in (\Phi^{[1]})^+$,
we see that $\be$ is the highest root of $\Phi^{[1]}$.   Observe
that $ \Lambda + \Gamma_k =  \be + \alpha_j$ as functions on
$\mathfrak t$. Since $p>3$, the function
$\be+\alpha_j\in\mathfrak{t}^*$ is not a root of $\g_0^{[1]}$.
Therefore, $[f^\Lambda,e^{\Gamma_k}] = 0$.

Now if $[e_{-\alpha_l},e^{\Gamma_k}]$ were nonzero for some $l
\neq j$, the above argument would yield
 $\Lambda+\Gamma_k  =  \be + \alpha_l$, which would force $\alpha_l =
\alpha_j$.     Hence, $[e_{-\alpha_l},e^{\Gamma_k}] = 0$ for all $l \neq j$. Using
that property and the fact that
$[e_{-\alpha_l},[e_{-\alpha_j}, e^{\Gamma_k}]] = 0$ for all $l$,  it is easy to verify that
$[e_{-\nu},e^{\Gamma_k}] =
0$ for all $\nu \in \Phi^+$, $\nu  \neq \alpha_j$.

Suppose $\be \neq \alpha_j$.    Using the relations obtained, we
get

\begin{eqnarray*}
[[[[f^\Lambda,e_{-\alpha_j}],e_{-\be}],e_{-\be}], e^{\Gamma_k}]
&
=& [[[f^\Lambda,u],e_{-\be}],e_{-\be}] \\
&=& \zeta[[e_\be,e_{-\be}],e_{-\be}] = -2 \zeta e_{-\be} \neq 0.
\end{eqnarray*}

\noindent  This yields that
$[[[f^\Lambda,e_{-\alpha_j}],e_{-\be}],e_{-\be}] \neq 0$.  Since
$\mathfrak{g}_{-1}$ is an irreducible $\mathfrak{g}_0$-module and
$f^\Lambda$ is a $\mathfrak b^+$-primitive vector of $\mathfrak{g}_{-1}$, the
subspace $P:= \mathfrak{U}\bigl((\mathfrak{n}^{[1]})^- \bigr).f^\Lambda
\subseteq \mathfrak{g}_{-1}$ generated by $f^\Lambda$ is an
irreducible $\mathfrak{g}_0^{[1]}$-module with a $(\mathfrak b^{[1]})^+$-primitive vector
of weight
$\Lambda^{[1]} = \Lambda \mid_{\,\mathfrak{t}^{[1]}}$.    It follows
from $(\Gamma_k - \alpha_j)(\mathfrak{t}^{[1]})= 0$  that
$\Lambda^{[1]} = \be$. Since the Lie algebra $\g_{0,1}$ is simple,
the $\mathfrak{g}_0^{[1]}$-modules $\mathfrak{g}_{0,1}$ and $P$
are isomorphic.   Moreover, because $\beta$ is a long root of $\Phi^{[1]}$, we
must have

$$\begin{gathered} {[[P,e_{-\be}],e_{-\be}]} \subseteq
P^{-\be} \quad \text {(the } -\be \text {\; weight space
of\;} P).\end{gathered}\nonumber$$

\noindent However, this contradicts the fact that $0$ $\neq
[[[f^\Lambda,$ $e_{-\alpha_j}]$ $,e_{-\be}],$ $e_{-\be}]$ $\in
P^{-\be -\alpha_j}.$  Consequently,  the case $\beta \neq
\alpha_j$ cannot occur.

So let us suppose that $\beta = \alpha_j$.  Then

\begin{eqnarray*}
[[[[f^\Lambda,e_{-\alpha_j}],e_{-\alpha_j}],e_{-\alpha_j}], e^{\Gamma_k}]
&=& 3[[[f^\Lambda,u],e_{-\alpha_j}], e_{-\alpha_j}] \\
&=& 3 \zeta [[e_{\alpha_j},e_{-\alpha_j}],
e_{-\alpha_j}]\\
&=& -6\zeta e_{-\alpha_j} \neq 0,\\
\end{eqnarray*}
\noindent whence  $0 \neq
[[[f^\Lambda,e_{-\alpha_j}],e_{-\alpha_j}],e_{-\alpha_j}] \in
P^{-2\be}$. But since  $P$ is isomorphic to $\mathfrak{g}_{0,1}$
as $\g_0^{[1]}$-modules, it is impossible for $-2\be$ to be a
weight of $P$.  This contradiction establishes the claim that
$e^{\Gamma_k}$ is a $\mathfrak b^-$-primitive vector.  As a consequence,
$[f^\Lambda,e^{\Gamma_k}] \neq 0$ by the transitivity of
$\mathfrak{g}$.

Suppose $\Gamma_k = -\Lambda$. Let $\widetilde{\g}_1$ denote the
$\g_0$-module generated by $e^{-\Lambda}$. If $\widetilde{\g}_1$
is irreducible then $\Gamma(\mathfrak{t}^{[1]}) \neq 0$ by Theorem
\ref{Thm:4.12}\,(i), contradicting our initial assumption.
Consequently, $\widetilde{\g}_1\cap Q_{k-1}\ne 0$. We may thus
assume without loss of generality that $e^{\Gamma_1}\in
\widetilde{\g}_1$. For all positive roots $\alpha\in \Phi^{[1]}$, 
choose $e_\alpha\in\g_0^\alpha$, $e_{-\alpha}\in\g_0^{-\alpha}$,
and $h_\alpha\in\mathfrak{t}$ such that $(e_\alpha,h_\alpha,
e_{-\alpha})$ form an $\mathfrak{sl}_2$-triple in $\g_{0,1}$.
Recall that $[e_{-\alpha},e^{\Gamma_1}]=0$ for all
$\alpha\in\Phi^{[1]}$. Since $[f^{\Lambda}, e^{\Gamma_1}]\ne 0$ by
the transitivity of $\g$ and the definition of
$k$ which implies that $e^{\Gamma_1}$ has
zero $\g_{0,1}$-weight, the linear map $\ad e^{\Gamma_1}$
induces an isomorphism between the $\g_0^{[1]}$-modules $P$ and
$\g_{0,1}$. Note that $[f^\Lambda,e^{\Gamma_1}]=e_\beta$, after
rescaling $e_\beta$ possibly. Since $\g$ is transitive, Theorem
\ref{Thm:3.20} yields $\Lambda(h_\Lambda)\ne 0$.  Rescaling
$e^{-\Lambda}$ if necessary we may assume that
$\Lambda(h_\Lambda)=2$.

Suppose $\Phi^{[1]}\cong\mathrm{A}_1$. Then the endomorphism
$$C:=E_\beta E_{-\beta}+E_{-\beta}E_\beta+\frac{1}{2}H_\beta^2$$ of
$\g$ commutes with $\ad \g_0$ and hence acts on $\widetilde{\g}_1$
as $\mu\,\text{id}$ where
$$\mu=\frac{1}{2}\Lambda(h_\beta)^2+\Lambda(h_\beta)=\frac{1}{2}
\beta(h_\beta)^2+\beta(h_\beta)=4.$$ But then
$4e^{\Gamma_1}=C(e^{\Gamma_1})=0$, a contradiction. This shows
that the rank of the root system $\Phi^{[1]}$ is $\ge 2$.

For $\alpha\in (\Phi^{[1]})^+$,  set $e_1=e_{-\alpha}$,
$e_2=f_\Lambda$, $f_1=e_{\alpha}$, $f_2=e_{-\Lambda}$. Clearly,
$[e_i,f_j]=0$ for $i\ne j$ and $[h_i,e_i]=2e_i$, where
$i,j\in\{1,2\}$. Let $\mathfrak k_\alpha$ denote the subalgebra of
$\g$ generated by $e_i,f_i$, $i=1,2$. Set $h_1=-h_\alpha$ and
$h_2=h_\Lambda$. Since $h_1\in\g_{0,1}$ we have $\Lambda(h_1)=
-\beta(h_\alpha)$.   The matrix $A_\alpha$ of the Lie
algebra $\mathfrak{k}_\alpha$ has the form
\begin{equation}\label{eq:4.21}
A_\alpha = \left(\begin{matrix} 2 & \Lambda(h_1)  \\
-\alpha(h_\Lambda) & 2  \\
\end{matrix} \right).
\end{equation}
Since $\alpha(h_\beta)=0$ if and only if $\beta(h_\alpha)=-\Lambda(h_1) = 0$, we
now apply Theorem \ref{Thm:3.10} to deduce that
\begin{equation}\label{eq:4.22}
\alpha(h_\beta)=0\iff \alpha(h_\Lambda)=0.
\end{equation}

Now let $\alpha_s$ be any simple root in $\Phi^{[1]}$ with
$\beta-\alpha_s\in(\Phi^{[1]})^+$ (such a root exists because
$\Phi^{[1]}$ is not of type $\mathrm{A}_1$). Then we have

\begin{eqnarray*}
[e_{-\gamma},[e^{-\Lambda},e^{\Gamma_1}]]\,=\,0\,=\,
[e_\gamma,[f^\Lambda,[f^\Lambda,e_{-\alpha_s}]]]
\end{eqnarray*}

\noindent for all $\gamma\in\Phi^+$. Since
$\Gamma_1-\beta=-\Lambda$ we also have

\begin{gather*}
\big[[f^\Lambda,[f^\Lambda,e_{-\alpha_s}]],[e^{-\Lambda},e^{\Gamma_1}]\big]=
\big[[f^\Lambda,[e^{-\Lambda},e^{\Gamma_1}]],[f^\Lambda,e_{-\alpha_s}]\big]\\
+
\big[f^\Lambda,[[f^\Lambda,[e^{-\Lambda},e^{\Gamma_1}]],e_{-\alpha_s}]\big]\\
 =\Gamma_1(h_\Lambda)[e^{\Gamma_1},[f^\Lambda,e_{-\alpha_s}]]+
\big[[e^{-\Lambda},e_\beta],[f^\Lambda,e_{-\alpha_s}]\big]\\
 +\,\Gamma_1(h_\Lambda)[f^\Lambda,[e^{\Gamma_1},e_{-\alpha_s}]]+
\big[f^\Lambda,[[e^{-\Lambda},e_\beta],e_{-\alpha_s}]\big]\\
=-\Gamma_1(h_\Lambda)[e_\beta,e_{-\alpha_s}]
-[e_{\beta},[e^{-\Lambda},[f^\Lambda,e_{-\alpha_s}]]]+
[f^\Lambda,[e^{-\Lambda},[e_\beta,e_{-\alpha_s}]]]\\
=-\Gamma_1(h_\Lambda)[e_\beta,e_{-\alpha_s}]+[e_\beta,[h_\Lambda,e_{-\alpha_s}]]+
[h_\Lambda,[e_\beta,e_{-\alpha_s}]]\\
=-(\Gamma_1
+\alpha_s-\beta+\alpha_s)(h_\Lambda)[e_\beta,e_{-\alpha_s}]\,=\,
(\Lambda-2\alpha_s)(h_\Lambda)[e_\beta,e_{-\alpha_s}].
\end{gather*}

\noindent Next observe that $2\Lambda-\alpha_s$ and
$-\Lambda+\Gamma_1=-2\Lambda+\beta$ do not vanish on
$\mathfrak{t}^{[i]}$ for $i\ge 2$ because $\alpha_s, \, \beta \in \Phi^{[1]}$ (see Remark \ref{Rem:4.8}). Also,
$(-2\Lambda+\beta)(h_\beta)=-\beta(h_\beta)\ne 0$ and
$$(2\Lambda-\alpha_s)(h_\beta)\,=\,2\beta(h_\beta)-\alpha_s(h_\beta)
\,=\,4-2\frac{(\alpha_s,\beta)}{(\beta,\beta)}\ne 0$$ as
$(\alpha_s,\beta)=(\beta,\alpha_s)>0$. Applying Theorem
\ref{Thm:4.7} to the graded Lie algebra
$\bigoplus_{i\in\Z}\,\g_{2i}$,   we are now able to deduce that
$(\Lambda-2\alpha_s)(h_\Lambda)=0$ since neither
$2\Lambda - \alpha_s$ nor $-2\Lambda + \beta$ is the highest
or lowest weight of a standard representation for either $\hbox{\rm A}_r$
or $\hbox{\rm C}_r$.  It follows that

\begin{equation}\label{eq:4.22a}
\alpha_s(h_\Lambda)=1 \quad\mbox{whenever}\ \ \alpha_s\in\Delta\ \
\mbox{and}\ \ \beta-\alpha_s\in\Phi.
\end{equation}

\noindent Let $\Delta_\beta$ denote the set of all simple roots
$\alpha\in(\Phi^{[1]})^+$ with $(\beta,\alpha)=0$. If $\Phi^{[1]}$
is not of type $\mathrm{A}_r$, $r \geq 2$,  then there is a unique
$\alpha_s\in\Delta$ such $\beta-\alpha_s\in\Phi$ and
$\beta=2\alpha_s+\sum_{\alpha\in\Delta_\beta}\,r_\alpha \alpha$
for some $r_\alpha\in\Z$. If $\Phi^{[1]}$ is of type
$\mathrm{A}_r$, $r\ge 2$, then $\beta-\alpha_s\in\Phi$ for
$s\in\{1,r\}$ and
$\beta=\alpha_1+\alpha_r+\sum_{\alpha\in\Delta_\beta}\,\alpha$. In
conjunction with (\ref{eq:4.22}) and (\ref{eq:4.22a}) this yields
$\beta(h_\Lambda)=2$.   As a result,  if $A_\beta$ denotes the
matrix in \eqref{eq:4.21} with $\beta$ in place of $\alpha$, we have 
\begin{equation}\nonumber
A_\beta = \left(\begin{matrix} 2 & -2  \\
-2 & 2  \\
\end{matrix} \right),
\end{equation}
see (\ref{eq:4.21}).    But then $h_1+h_2$
lies in the center of the subalgebra $\mathfrak k_\beta$
generated by $e_1=e_{-\beta}$,
$e_2=f_\Lambda$, $f_1=e_{\beta}$, $f_2=e_{-\Lambda}$,
and the graded
algebra $\mathfrak{k}_{\beta}/\F\,(h_1+h_2)$ satisfies the
conditions of Proposition \ref{Pro:3.8}.    Since $\g$ is finite-dimensional,
this is impossible in view of Theorem \ref{Thm:3.9}.
Thus, $\Gamma_k \neq -\Lambda$.

If $\Gamma_k(\mathfrak{t}^{[i]}) \neq 0$ for all $i \in \{1,
\dots, \ell\}$, then Theorem \ref{Thm:4.7} would imply that $\Phi
= \Phi^{[1]}$ and either $\Phi =$ A$_m$, $\Lambda \in
\{\varpi_1,\varpi_m\}$, or $\Phi =$ C$_m$, $\Lambda = \varpi_1$.
Since $\ad e^\Gamma \mid_{\mathfrak{g}_{-1}} \in
\Hom(\mathfrak{g}_{-1},\mathfrak{g}_0)^{\mathfrak{t}^{[1]} +
(\mathfrak{n}^{[1]})^-} \cong (\mathfrak{g}_{-1}^* \ot
\mathfrak{g}_0)^{\mathfrak{t}^{[1]} + (\mathfrak{n}^{[1]})^-}$, this
would say that $\mathfrak{g}_{-1}^* \ot \mathfrak{g}_0$ has a
$\mathfrak b^-$-primitive vector of weight zero, which implies that
$\mathfrak{g}_0$ has a minuscule weight (see \eqref{eq:2.8}).    As
we have noted in the proof of Lemma \ref{Thm:4.12}(i),  this is
impossible. Consequently,  we must have $\ell
> 1$ and $\Gamma_k(\mathfrak{t}^{[i]}) = 0$ for some $i \in \{2, \dots,
\ell\}$. We may assume that $i= 2$.  Since by transitivity,
 $\Lambda(\mathfrak{t}^{[l]}) \neq 0$ for all $l \in
\{1, \dots, \ell\}$ (see Remark \ref{Rem:4.8}), we have $[f^\Lambda,e^\Gamma] = \xi e_{\mu}$ and
$[f^\Lambda,e^{\Gamma_k}] = \xi' e_{\nu}$ for some  $\mu \in
\Phi^{[1]}$, $\nu \in \Phi^{[2]}$,  and scalars
$\xi,\xi'$ which are nonzero (again by transitivity).
If $\mu \in \Phi^-$ or $\nu \in \Phi^-$, then $\ell = 1$ by Theorem
\ref{Thm:4.7}. Therefore, it must be that $\mu,\nu \in \Phi^+$. By
adjusting $e^\Gamma$ and $e^{\Gamma_k}$ by scalars if necessary,
we may assume that $\xi = 1 = \xi'$.  Set

$$\begin{array}{cc}
e_1 = [f^\Lambda,e_{-\mu}] &\qquad  e_2 = [f^\Lambda,e_{-\nu}] \\
f_1 = e^\Gamma &\qquad
f_2 = e^{\Gamma_k}\\  h_1 = h_\mu &\qquad h_2 = h_\nu.
\end{array}$$

\noindent It is easy to check that $[e_i,f_j] = \delta_{i,j}h_i$
for $i,j \in \{1,2\}$.  As $\Gamma(\mathfrak{t}^{[1]}) =
\Gamma_k(\mathfrak{t}^{[2]}) = 0$, we have that $\Lambda
\mid_{\mathfrak{t}^{[1]}} = \mu \mid_{\mathfrak{t}^{[1]}}$ and
$\Lambda \mid_{\mathfrak{t}^{[2]}} = \nu
\mid_{\mathfrak{t}^{[2]}}$. This implies that $[h_i,e_i] =
[h_i,f_i] = 0$ (for $i = 1,2$), $[h_1,e_2] = 2e_2, \; [h_2,e_1] =
2 e_1, \; [h_1,f_2] = - 2f_2$ and  $[h_2,f_1] = -2f_1$.   But then
the Lie algebra generated by the $e_i,f_i,h_i$, $i = 1,2$ is
infinite-dimensional by Theorem \ref{Thm:3.9}.  We have reached a
contradiction, so it must be that $\Gamma(\mathfrak{t}^{[i]}) \neq 0$
for all $i$.   \qed

\m

\section {\ Determination of
the local Lie algebra \label{sec:4.6}} 
\m

    Here we determine the Lie algebra  generated by the local part
    $\g_{-1} \oplus \g_0 \oplus \g_1$ and show that it must be one of the
restricted Lie algebras of Cartan type.   If  $\g_{-1}$ has
a $\mathfrak b^+$-primitive vector of weight $\Lambda$ and  $\g_1$ has
a $\mathfrak b^-$-primitive vector 
of weight $\Gamma$,   then a classical Lie algebra is generated if
$\Gamma = -\Lambda$. 
We analyze which classical Lie algebras can be obtained by
adding $-\Lambda$  to the root system of $\g_0$.      Of course, since
we are in the noncontragredient case, there must be more to
$\mathfrak{g}$ than just that larger classical Lie algebra, and we
must determine that part as well.

 Given a rational module $M$ over a reductive group $G$
 we denote by $X_+(M)$ the set of dominant weights of $M$
 relative to a maximal torus of $G$.  \m

\begin{Thm} \label{Thm:4.22}    \ Let $\mathfrak{g}$ be a graded Lie
algebra satisfying conditions (1)-(5).  Suppose $\mathfrak{g}$ is
generated by its local part $\mathfrak{g}_{-1} \oplus
\mathfrak{g}_0 \oplus \mathfrak{g}_1$, and assume there exist a
$\mathfrak b^+$-primitive vector $f^\Lambda \in \mathfrak{g}_{-1}$ of
weight $\Lambda$ and a
$\mathfrak b^-$-primitive vector $e^\Gamma \in \mathfrak{g}_1$ of weight
$\Gamma$ such that
$[f^\Lambda, e^\Gamma] = e_{-{\alpha}}$ where $\alpha \in \Phi^+$.
Then the graded Lie algebra $\mathfrak{g}$ is isomorphic to a
restricted Lie algebra of Cartan type with its natural grading.
More precisely, $\g$ is either $W(m;\underline{1}),\, m\ge 2,$ or
$S(m;\underline{1})^{(1)},\,m\ge 3,$ or $
S(m;\underline{1})^{(1)}\oplus\F{\mathfrak D}_1,\,m\ge 3,$ or
$H(2m,\underline{1})^{(2)},\, m\ge 1,$ or $
H(2m;\underline{1})^{(2)}\oplus\F{\mathfrak D}_1, \,m\ge 3,$ or
$K(2m+1,\underline{1})^{(1)}, m\ge 1$.
\end{Thm}

\pf Combining Proposition \ref{Pro:4.20} with Theorem
\ref{Thm:4.7} we see that the root system $\Phi$ of $\g_0^{(1)}$
is of type A$_m$ or C$_m$; ${\alpha}$ is the highest root of
$\Phi$; and $\g_{-1}$ is a standard $\g_0^{(1)}$-module.   Hence the
weight $\Lambda$ of $\g_{-1}$ is a minuscule weight.
 Replacing the base $\Delta=\{\alpha_1,\dots,\alpha_m\}$ of $\Phi$
 by $\{\alpha_1' = \alpha_m, \alpha_2' =
\alpha_{m-1}, \dots, \alpha_m' =\alpha_1\}$ if necessary, we may
assume that $\Lambda = \varpi_m$ if $\Phi$ is of type A$_m$. Thus,
we may assume  that either $\g_{-1}=V$ and
$\mathfrak{sp}(V)\subseteq \g_0\subseteq \mathfrak{csp}(V)$ or
$\g_{-1}=V^*$ and
$\mathfrak{sl}(V)\subseteq\g_0\subseteq\mathfrak{gl}(V)$, where
$V$ is the standard $\g_0^{(1)}$-module. There is a connected,
reductive subgroup $\widetilde{G}_0$ in $\mathrm{GL}(V)$ with
$\text{Lie}(\widetilde{G}_0)=\g_0$ such that the adjoint action of
$\g_0$ on $\g_{-1}\oplus\g_0$ is induced by the differential of
the natural action of $\widetilde{G}_0$. Moreover,
$\widetilde{G}_0$ contains a maximal torus $\widetilde{T}$ and
maximal unipotent subgroups $N^\pm$ such that $\hbox{\rm
Lie}(T)=\mathfrak{t}\cap \g_0^{(1)}$ and $\hbox{\rm Lie}(N^\pm)=
\mathfrak{n}^\pm$. Let $G_0$ denote the derived subgroup of
$\widetilde{G}_0$ and put $T=\widetilde{T}\cap G_0$. Then $T$ is a
maximal torus of $G_0$ with $\text{Lie}(T)=\tf\cap\g_0^{(1)}$ and
$\text{Lie}(G_0)=\g_0^{(1)}$.

\medskip

\noindent (a) The transitivity of $\g$ allows us to identify
$\g_1$ with a $\g_0$-submodule of $\text{Hom}(\g_{-1},\g_0)$. The
above discussion shows that the adjoint action of $\g_0$ on $\g_1$
is obtained by restricting to $\g_1$ the differential of the
natural rational action of $G_0$ on
$\text{Hom}(\g_{-1},\g_0)\,\cong\,\g_{-1}^*\otimes \g_0$.

Suppose $\g_0^{(1)}$ is of type ${\rm A}_m$, $m\ge 2$. Then
$\text{Hom}(\g_{-1},\g_0)\,
\hookrightarrow\,V\otimes\mathfrak{gl}(V)$ is isomorphic to a
$G_0$-submodule of the rational $G_0$-module $V \ot V \ot V^*$.
Since  the set of $\widetilde{T}$-weights of $V$ equals
$\{\e_1,\ldots,\e_{m+1}\}$ and $p>2$, we have
\begin{eqnarray}\label{eq:A}
X_+\big(\text{Hom}(\g_{-1},\g_0)\big)&\subseteq & X_+(V \ot V \ot V^*)\\
&=&\{2\varpi_1+\varpi_m,\varpi_2+\varpi_m,\varpi_1\}\nonumber\\
&\subseteq& X_1(T),\nonumber
\end{eqnarray}
where $\varpi_2+\varpi_m$ is omitted if $m=2$. Hence any
$\g_0$-submodule of $\g_1$ is $G_0$-stable by Proposition
\ref{Pro:2.805}.

Suppose $\g_0^{(1)}$ is of type $\mathrm{C}_m$, $m\ge 1$. Then
$V\cong V^*$ and $\g_0^{(1)}=\mathfrak{sp}(V)\,\cong\, S^2(V)$ as
$G_0$-modules. Hence $\Hom(\mathfrak{g}_{-1}, \,
\mathfrak{g}_0)\hookrightarrow V\otimes
\big(\mathfrak{sp}(V)\oplus\F\big)$ is isomorphic to a
$G_0$-submodule of $V \oplus \bigl(S^2(V) \otimes V\bigr)$. Since
$X(V)=\{\pm\e_1,\ldots,\pm\e_m\}$ and $p>3,$ we have
\begin{eqnarray}\label{eq:B}
X_+\big(\text{Hom}(\g_{-1},\g_0)\big)
&\subseteq& X_+(V) \cup X_+\big(S^2(V) \otimes V\big)\\
&=&\{3\varpi_1, \varpi_1+\varpi_2,
\varpi_3,\varpi_1\}\nonumber\\&\subseteq& X_1(T),\nonumber
\end{eqnarray}
where $\varpi_1+\varpi_2$ is omitted if $m=1$ and $\varpi_3$ is
omitted if $m=1,2$. Again Proposition \ref{Pro:2.805} applies to
show that any $\g_0$-submodule of $\g_1$ is an $G_0$-stable.

\m

\noindent
(b) Let $\g'_1$ be the $\mathfrak{g}_0$-submodule
of $\mathfrak{g}_1$ generated by $e^\Gamma$.  Let ${\mathfrak
U}(\mathfrak{n}^-)$ denote the universal enveloping algebra of the
Lie algebra $\mathfrak{n}^-$. As $\mathfrak{g}_{-1}$ is
irreducible, $\mathfrak{g}_{-1} = {\mathfrak
U}(\mathfrak{n}^-)f^\Lambda$. Therefore, $[e^\Gamma,
\mathfrak{g}_{-1}] \subseteq {\mathfrak
U}(\mathfrak{n}^-)[e^\Gamma,f^\Lambda] = {\mathfrak
U}(\mathfrak{n}^-)e_{-{\alpha}} = {\mathbb F} e_{-{\alpha}}$.
Moreover, $[\mathfrak{n}^-,\mathfrak{g}_{-1}] \subseteq \ker \ad
e^\Gamma$ as $[\mathfrak{n}^-,e_{-{\alpha}}] = 0$. Since
$\mathfrak{g}_{-1}$ is  standard, we have $\mathfrak{g}_{-1} = {\mathbb F} f^\Lambda
\oplus [\mathfrak{n}^-,\mathfrak{g}_{-1}]$ and
$[e_{-{\alpha}},\mathfrak{g}_{-1}] = {\mathbb F} e^{\Omega}$ where
$e^{\Omega}$ is a $\mathfrak b^-$-primitive vector of $\mathfrak{g}_{-1}$. But
then

\begin{eqnarray*} [e^\Gamma,[\mathfrak{g}_{-1},\mathfrak{g}_{-1}]]  &=&
[e^\Gamma,[({\mathbb F} f^\Lambda +
[\mathfrak{n}^-,\mathfrak{g}_{-1}]),({\mathbb F} f^\Lambda +
[\mathfrak{n}^-,\mathfrak{g}_{-1}])]] \nonumber \\
&=&   [e^\Gamma, [\g_{-1}, [\mathfrak{n}^-,\mathfrak{g}_{-1}]] =
[[e^\Gamma,\g_{-1}], [\mathfrak{n}^-,\mathfrak{g}_{-1}]] \nonumber \\
&=& [e_{-{\alpha}},[\mathfrak{n}^-,\mathfrak{g}_{-1}]] =
[\mathfrak{n}^-,e^\Omega]
= 0.
\end{eqnarray*}

\noindent This implies that
$[\g'_1,\mathfrak{g}_{-2}] = 0$.

\m

\noindent
(c) Let $e^\Xi$ be a $\mathfrak b^-$-primitive vector from
$\mathfrak{g}_1$ of weight $\Xi$.  If $\Xi \neq -\Lambda$, then
$[f^\Lambda, e^\Xi] \in {\mathbb F}^\times e_{-\delta}$ for some $\delta
\in \Phi$ by transitivity \eqref{eq:1.3}.  Since no minuscule weight is a root,
$\Xi(\mathfrak{t}$ $\cap$ $\g_0^{(1)}) \neq 0$.    If $\delta \in
\Phi^-$, then $\delta$ is the lowest root of $\Phi$ and $\Xi$ is the
lowest weight of a standard $\g_0^{(1)}$-module  by
Theorem \ref{Thm:4.7}.  As this situation is easily seen to be
impossible,
$\delta  \in \Phi^+$. Applying Theorem \ref{Thm:4.7} (i), we conclude that
$\delta= \alpha$ and $\Xi = -(\Lambda+{\alpha})$.  Thus, any
$\mathfrak b^-$-primitive vector of $\mathfrak{g}_1$ has weight
$-(\Lambda+{\alpha})$ or $-\Lambda$.

Now if $e$ is any $\mathfrak b^-$-primitive vector from
$\g_1^{-\Lambda}$, then the tuple $(f^\Lambda, \mathfrak{g}_0, e)$
satisfies all the conditions of Theorem \ref{Thm:3.20} (with
$-\Delta$ being the set of simple roots). Hence $[f^\Lambda,[e,
f^\Lambda]]\ne 0$. But $[f^\Lambda, \mathfrak{g}_1^{-\Lambda}]
\subseteq \mathfrak{t}$ and $[\mathfrak{t}, f^\Lambda] \subseteq
{\mathbb F} f^\Lambda$, yielding $[f^\Lambda,[e, f^\Lambda]]  \in
\F f^\Lambda$. From this it is immediate that the  subspace of
$\mathfrak b^-$-primitive vectors from $\mathfrak{g}_1^{-\Lambda}$ has
dimension $\leq 1$.

\m

\noindent (d) Let $e^{-\Lambda}$ be a $\mathfrak b^-$-primitive vector
in $\mathfrak{g}_1^{-\Lambda}$, and let $\widetilde{\g}_1$ be the
$\g_0$-submodule of $\g_1$ generated by $e^{-\Lambda}$. By our
discussion in part~(a), the subspace $\widetilde{\g}_1$ is
$G_0$-stable. Let $\widetilde{\mathfrak{g}}$ denote the graded
subalgebra of $\mathfrak{g}$ generated by $f^\Lambda,
\mathfrak{g}_0$, and $e^{-\Lambda}$. We claim that the graded Lie
algebra $\widetilde{\mathfrak{g}}/\mathcal
M(\widetilde{\mathfrak{g}})$ (where $\mathcal
M(\widetilde{\mathfrak{g}})$ denotes the Weisfeiler radical)
satisfies the hypotheses of Theorem \ref{Thm:3.34}. To prove the
claim we need to show that the $\g_0$-module $\widetilde{\g}_1$ is
irreducible.

Since the torus $T$ acts trivially on $\mathfrak t$, it
preserves the subspace of $\mathfrak b^-$-primitive vectors from
$\widetilde{\g}_1^{-\Lambda}$. By part~(c), the latter is spanned
by $e^{-\Lambda}$, so it must be that $e^{-\Lambda}\in
\g_1^{-\lambda}$ for some $-\lambda\in X(\g_1)$. The equality
$\widetilde{\g}_{1} = {\mathfrak U}(\mathfrak{n}^+)
\,e^{-\Lambda}$ implies that $-\lambda$ is a minimal weight of the
$G_0$-module $\widetilde{\g}_1$. But then $-w_0\lambda$ is a
maximal weight of $\widetilde{\g}_1$.    Also,
$\dim\,(\widetilde{\g}_1)^{-w_0\lambda}=\dim\,(\widetilde{\g}_1)^{-\lambda}=1$,
by part~(c). Since $\mathfrak{t}\cap \g_0^{(1)}=\text{Lie}(T)$,
the restriction of $-w_0\Lambda$ to $\mathfrak{t}\cap\g_0^{(1)}$
coincides with the differential of the rational character
$-w_0\lambda$ at the identity element of $T$.   From this it follows
that the image of $-w_0\lambda$ in $X(T)/pX(T)$ coincides with
$-w_0\Lambda\mid_{\,\tf\,\cap\,\g_0^{(1)}}$. Note that
$-w_0\Lambda(h_\gamma)\in\{0,\pm 1\}$ for all $\gamma\in\Phi$.
Combining this with (\ref{eq:A}) and (\ref{eq:B}) one derives that
$-w_0\lambda=\varpi_1$.    Since it also follows from (\ref{eq:A})
and \ref{eq:B}) that $\nu\ge\varpi_1$ for all $\nu\in X_+(\g_1)$,
one obtains now that $\varpi_1$ is the only dominant weight of the
$G_0$-module $\widetilde{\g}_1$. As a consequence, all weights of
the $G_0$-module $\widetilde{\g}_1$ are conjugate under the Weyl
group of $G_0$. Therefore, the $G_0$-modules $V$ and
$\widetilde{\g}_1$ are isomorphic.   But then $\widetilde{\g}_1$ is
$\g_0$-irreducible, hence the claim.

\m

\noindent (e)  The graded Lie algebra
$\widetilde{\mathfrak{g}}/\mathcal M(\widetilde{\mathfrak{g}})$
(where $\mathcal M(\widetilde{\mathfrak{g}})$ denotes the
Weisfeiler radical) satisfies the hypotheses of Theorem
\ref{Thm:3.34}. Therefore by that theorem, this algebra is
classical. Using \cite{Bou1}, it is easy to check that
$\widetilde{\mathfrak{g}}/\mathcal M(\widetilde{\mathfrak{g}})$
corresponds to one of the following Dynkin diagrams:

\begin{equation}\label{eq:4.23}
{\beginpicture \setcoordinatesystem
units <0.45cm,0.3cm> % sets scale
 \setplotarea x from -4 to 16, y from -2 to 1    % sets plot size up
 \linethickness=0.03pt                          % sets line thickness
  \put{$\circ$} at 0 0
 \put{$\circ$} at 2 0
 \put{$\cdots$} at 5 0
 \put{$\circ$} at  8 0
 \put{$\circ$} at 10  0 \put{$\bullet$} at  12 0
 \plot .15 .1 1.85 .1 /
 \plot 2.15 .1 3.85 .1 /
 \plot 6.15 .1 7.85 .1 /
 \plot 8.15  .1 9.85 .1 /
 \plot 10.15 .1  11.85 .1 /
 \put{$\alpha_1$} at 0 -1.5
  \put{$\alpha_2$} at 2 -1.5
    \put{$\alpha_{m-1}$} at 8 -1.5
      \put{$\alpha_m$} at 10 -1.5
        \put{$-\Lambda$} at 12 -1.5
        \put {$(m \geq 1)$} at 16 0  \endpicture}
\end{equation}

\begin{equation}\label{eq:4.24}
{\beginpicture \setcoordinatesystem
units <0.45cm,0.3cm> % sets scale
 \setplotarea x from -4 to 16, y from -1 to 1    % sets plot size up
 \linethickness=0.03pt                          % sets line thickness
  \put{$\circ$} at 0 0
 \put{$\circ$} at 2 0
 \put{$\cdots$} at 5 0
 \put{$\circ$} at  8 0
 \put{$\circ$} at 10  0 \put{$\bullet$} at  12 0
 \plot .15 .1 1.85 .1 /
 \plot 2.15 .1 3.85 .1 /
 \plot 6.15 .1 7.85 .1 /
 \plot 8.15 .1 9.85 .1 /
  \put {$>$} at 11 0
 \plot 10.15 -.2  11.85 -.2 /
  \plot 10.15 .2 11.85  .2 /
 \put{$\alpha_1$} at 0 -1.5
  \put{$\alpha_2$} at 2 -1.5
    \put{$\alpha_{m-1}$} at 8 -1.5
      \put{$\alpha_m$} at 10 -1.5
        \put{$-\Lambda$} at 12 -1.5
        \put {$(m \geq 2)$} at 16 0  \endpicture}
\end{equation}

\begin{equation}\label{eq:4.25}
{\beginpicture \setcoordinatesystem
units <0.45cm,0.3cm> % sets scale
 \setplotarea x from -6 to 8, y from -1 to 1    % sets plot size up
 \linethickness=0.03pt                          % sets line thickness
  \put{$\circ$} at 0 0
 \put{$\bullet$} at 2 0
 \plot .15 -.3 1.85 -.3 /
  \plot .15 0 1.85  0 /
  \plot .15 .3 1.85 .3 /
  \put {$>$} at 1 0
   \put{$\alpha_1$} at 0 -1.5
  \put{$-\Lambda$} at 2 -1.5 \endpicture}
\end{equation}

\begin{equation}\label{eq:4.26}
{\beginpicture \setcoordinatesystem
units <0.45cm,0.3cm> % sets scale
 \setplotarea x from -6 to 8, y from -1 to 1    % sets plot size up
 \linethickness=0.03pt                          % sets line thickness
  \put{$\circ$} at 0 0
 \put{$\bullet$} at 2 0
 \plot .15 -.2 1.85 -.2 /
  \plot .15 .2 1.85 .2 /
  \put {$>$} at 1 0
  \put{$\alpha_1$} at 0 -1.5
  \put{$-\Lambda$} at 2 -1.5 \endpicture} \end{equation}

\begin{equation}\label{eq:4.27}
 \beginpicture \setcoordinatesystem
units <0.45cm,0.3cm> % sets scale
 \setplotarea x from -4 to 12, y from -1 to 1    % sets plot size up
 \linethickness=0.03pt                          % sets line thickness
  \put{$\bullet$} at 0 0
 \put{$\circ$} at 2 0
  \put{$\circ$} at 4 0
 \put{$\cdots$} at 6 0
 \put{$\circ$} at  8 0
 \put{$\circ$} at 10  0
 \put{$\circ$} at  12 0
 \plot .15 .1 1.85 .1 /
 \plot 2.15 .1 3.85 .1 /
  \plot 8.15 .1 9.85 .1 /
  \put {$<$} at 11 0
 \plot 10.15 -.2  11.85 -.2 /
  \plot 10.15 .2 11.85  .2 /
  \put{$-\Lambda$} at 0 -1.5
  \put{$\alpha_1$} at 2 -1.5
   \put{$\alpha_2$} at 4 -1.5
      \put{$\alpha_{m-1}$} at 10 -1.5
        \put{$\alpha_m$} at 12 -1.5
        \put {$(m \geq 2).$} at 16 0  \endpicture
        \end{equation}

\bi

\noindent (Bear in mind that $\Phi$ is of type A$_m$, $m\ge 2$, or
C$_m$, $m\ge 1$, and $\Lambda = \varpi_m$ if $\Phi = \,$A$_m$  and
$\Lambda = \varpi_1$ if $\Phi =$ C$_m$.)     We may suppose that
the vector $e^{-\Lambda}$ has been scaled if need be, so that
$(f^\Lambda,\,h_\Lambda=[f^\Lambda,e^{-\Lambda}],\, e^{-\Lambda})$
is a standard $\mathfrak{sl}_2$-triple.

If $\widetilde{\mathfrak{g}}/\mathcal M(\widetilde{\mathfrak{g}})$
has type \eqref{eq:4.24},  set

\begin{equation}x_{\lambda} =
[[[e_{-\alpha_m}, \, f^\Lambda], \, f^\Lambda], \,
f^\Lambda] \quad \hbox{\rm and} \quad x_{\gamma} = [[[e_{\alpha_m}, \, e^{-\Lambda}], \,
e^{-\Lambda}], \, e^{\Gamma}],    \nonumber\end{equation}

\noindent and use capital letters to denote the corresponding
adjoint mappings. Then

\begin{eqnarray}[x_\lambda,x_\gamma]  &=&\left(\sum_{j = 0}^3 (-1)^j {3 \choose
j} (F^\Lambda)^j E_{-\alpha_m}(F^{\Lambda})^{3-j}\right)
\bigl([[[e_{\alpha_m},e^{-\Lambda}],e^{-\Lambda}],e^{\Gamma}]\bigr)  \nonumber \\
&&\hspace{-1.7truecm} =  (\Lambda-2 \alpha_m)(h_{\Lambda})
\left( \sum_{j = 0}^2 (-1)^j {3 \choose
j} (F^\Lambda)^j E_{-\alpha_m}(F^{\Lambda})^{2-j} \right) \big([[e_{\alpha_m},e^{-\Lambda}],
e^{\Gamma}]\big)  \nonumber \\
&&\hspace{-1.7truecm} = -2 \Big (3 \Lambda(h_{\alpha_m})\Gamma (h_{\Lambda})
- 3 \alpha_m(h_{\Lambda})\Gamma (h_{\alpha_m}) -\alpha_m(h_\Lambda)\alpha(h_{\alpha_m}) \Big)e_{-\alpha} \nonumber\\
&&\hspace{-1.7truecm} =  -2(-12 +12-2)e_{-\alpha}
= 4 e_{-{\alpha}}.  \nonumber\end{eqnarray}

\noindent As $[x_\lambda,e_{\alpha}] = 0 =
[x_\gamma,e_{-{\alpha}}]$, the elements $x_\lambda, \frac{1}{4}x_\gamma,
e_{-{\alpha}}$ satisfy the conditions of Lemma \ref{Lem:4.5}.  But
then $\lambda(h_{\alpha}) \in \{0,1\}$, which contradicts the fact that
$(3 \Lambda -\alpha_m)(h_{\alpha}) = 2$.  So case \eqref{eq:4.24}
cannot occur.

Now when  $\widetilde{\mathfrak{g}}/\mathcal
M(\widetilde{\mathfrak{g}})$ is of type \eqref{eq:4.25}, set

\begin{eqnarray}
x_\lambda &=& [[[[e_{-\alpha_1},f^{\Lambda}],f^{\Lambda}],
f^{\Lambda}],f^{\Lambda}], \nonumber \\
x_\gamma &=& [[[[e_{\alpha_1},e^{-\Lambda}],e^{-\Lambda}],e^{-\Lambda}],e^{\Gamma}]. \nonumber
\end{eqnarray}

\noindent Then, since here $\alpha_1(h_{\Lambda}) = 3,$
$\Lambda(h_{\alpha_1}) = 1,$ and $[f^{\Lambda}, \, e^{\Gamma}] =
e_{-\alpha_1},$ we have

\begin{eqnarray*}
[x_\lambda,x_\gamma]\hspace{-.3truecm}
&=&\hspace{-.3truecm}\Bigg(\sum_{j = 0}^4 (-1)^j{4 \choose j}(F^{\Lambda})^j
E_{-\alpha_1}(F^{\Lambda})^{4-j}\Bigg)
\,[[[[e_{\alpha_1},e^{-\Lambda}],e^{-\Lambda}],e^{-\Lambda}],e^{\Gamma}] \\
&& \hspace*{-2cm}= 3(\Lambda-\alpha_1)(h_{\Lambda}) \Bigg(\sum_{j = 0}^2
(-1)^j{4 \choose j}(F^{\Lambda})^j
E_{-\alpha_1}(F^{\Lambda})^{3-j}\Bigg)
[[[e_{\alpha_1},e^{-\Lambda}],e^{-\Lambda}],e^{\Gamma}] \\
&& \hspace*{-2cm}=-3(\Lambda-2 \alpha_1)(h_{\Lambda}) \Bigg(\sum_{j = 0}^2
(-1)^j{4 \choose j}(F^{\Lambda})^j
E_{-\alpha_1}(F^{\Lambda})^{2-j}\Bigg)
\,[[e_{\alpha_1},e^{-\Lambda}],e^{\Gamma}] \\
&& \hspace*{-2cm}= 12 \Big( 6\Lambda (h_{\alpha_1})\Gamma (h_{\Lambda}) - 4
\alpha_1(h_{\Lambda})\Gamma (h_{\alpha_1})  \\
&& \hspace*{1.45cm} -\alpha_1(h_\Lambda)\alpha_1(h_{\alpha_1}) -
\Lambda(h_{\alpha_1})\alpha_1(h_\Lambda)\Big)e_{-\alpha_1}\\
&& \hspace*{-2cm}= 12\big(-30+36-6-3\big)e_{-\alpha_1}
= -36e_{-\alpha_1}.  \end{eqnarray*}

\noindent Since $[x_\lambda,e_{\alpha_1}] = 0 =
[x_\gamma,e_{-\alpha_1}]$ and $\lambda(h_{\alpha_1}) = (4 \Lambda
-\alpha_1)(h_{\alpha_1}) = 2$, this case cannot occur by Lemma
\ref{Lem:4.5}.

\m

\noindent (f) Now suppose that $\Phi =$ A$_m$ for $m \geq 2$. Then
$\widetilde{\mathfrak{g}}/\mathcal M(\widetilde{\mathfrak{g}})$
must be of type \eqref{eq:4.23}. As
$\widetilde{\mathfrak{g}}_{-2}=\mathcal
M(\widetilde{\mathfrak{g}})_{-2}$, it follows that
$[\mathfrak{g}_{-2},\widetilde{\mathfrak{g}}_1] = 0$.

According to part~(a), $\g_1$ is a $G_0$-stable subspace of $V \ot
V \ot V^*$. Obviously, $V \ot V \ot V^* \cong M \oplus N$ where $M
\cong S^2(V) \ot V^*$ and $N \cong \wedge^2 \,V \ot V^*$ as
$\widetilde{G}_0$-modules. Since $\alpha=\e_1-\e_{m+1}$, $p>3$,
and the set of weights of $N$ with respect to $\widetilde{T}$
equals $\{\e_i+\e_j-\e_k \mid 1\le i,j,k\le m+1,\,i\ne j\}$, it is
straightforward to see that the $\tf$-weight space $N^{-(\Lambda +
{\alpha})} = N^\Gamma$ is trivial (see the first paragraph of the proof). Therefore $\g'_1$, the $\g_0$-submodule of
$\g_1$ generated by $e^\Gamma$, is contained in $M$. We claim that
$\mathfrak{g}_1$ is isomorphic to a $\mathfrak{g}_0$-submodule of
$M$.

If the natural projection map $\pi: \mathfrak{g}_1 \rightarrow M$
is injective, there is nothing to prove. So we assume that $N \cap
\mathfrak{g}_1 = \ker \pi \neq 0$. If $M$ is an indecomposable
$\g_0^{(1)}$-module, then $\g'_1$ contains a
$\g_0^{(1)}$-submodule isomorphic to $V$ (see Theorem
\ref{Thm:2.82}), hence it has a $\mathfrak b^-$-primitive vector of weight
$-\Lambda$.  Since $N \cap \mathfrak{g}_1$ has no nonzero weight
vectors of weight $-(\Lambda+{\alpha})$, it contains a nonzero
$\mathfrak b^-$-primitive vector of weight $-\Lambda$; see part~(c).  But then
the subspace of $\mathfrak b^-$-primitive vectors from
$\mathfrak{g}_1^{-\Lambda}$ has dimension $\geq 2$, which is
impossible by our final remark in part~(c). Thus,
$\mathfrak{g}'_1$ is an irreducible $\g_0^{(1)}$-module, and
$M=\mathfrak{g}'_1 \oplus \widetilde M$, where $\widetilde M$ is a
$\g_0^{(1)}$-submodule of $M$ isomorphic to $V$ (see Theorem
\ref{Thm:2.82}).

By part~(a), the $\g_0$-module $\ker \pi$ is $G_0$-stable. Assume
$\mu$ is a minimal weight of $\ker \pi$, and let $(\ker \pi)^\mu$
be the corresponding weight space.  As $(\ker \pi)^\mu$ lies in
the subspace of fixed points of $N^-$, we have $[\mathfrak{n}^-,
(\ker \pi)^\mu] = 0$.  Thus, $(\ker \pi)^\mu$ consists of 
$\mathfrak b^-$-primitive vectors. Since $\ker \pi \subseteq N$ contains no weight
vector associated with $-(\Lambda + \alpha)$, the image of $\mu$
in $X(T)/pX(T)$ coincides with the restriction of $-\Lambda$ to
$\tf\cap \g_0^{(1)}$, and $\dim \,(\ker \pi)^\mu = 1$. As $w_0\mu$
is dominant and $p> 3$, it is immediate from (\ref{eq:A}) that
$\mu = -\varpi_m$. Arguing as in part~(d) we now deduce that all
weights of the $G_0$-module $\ker \pi$ are conjugate under the
Weyl group of $G_0$. It follows that the $G_0$-modules $V$ and
$\ker \pi$ are isomorphic.

By the same reasoning one can show that $V\cong\pi^{-1}(\widetilde
M)$ as $G_0$-modules.  As $\ker \pi \subseteq \pi^{-1}(\widetilde
M)$, it must be that $\ker \pi = \pi^{-1}(\widetilde M)$.  Hence,
$\pi(\mathfrak{g}_1) = \mathfrak{g}'_1$.  As $\mathfrak{g}'_1
\subseteq \mathfrak{g}_1$, this yields $\mathfrak{g}_1 =
\mathfrak{g}'_1 \oplus \ker \pi \cong \mathfrak{g}'_1 \oplus
\widetilde M$, proving the claim.

If $\mathfrak{g}'_1 = \mathfrak{g}_1$ or if the
$\g_0^{(1)}$-submodule $S^2(V) \ot V^*$ is completely reducible,
then $\mathfrak{g}_1$ is generated by its $\mathfrak b^-$-primitive vectors.
In other words, $\mathfrak{g}_1 = \mathfrak{g}'_1 \oplus
\widetilde{\mathfrak{g}}_1$. (If $\mathfrak{g}_1$ has no 
$\mathfrak b^-$-primitive vectors of weight $-\Lambda$, we assume that
$\widetilde{\mathfrak{g}}_1 = 0$.)   Since we have proven that
$[\mathfrak{g}_{-2}, \mathfrak{g}'_1 + \widetilde{\mathfrak{g}}_1]
= 0$, we must have $\mathfrak{g}_{-2} = 0$ by \,1-transitivity
\eqref{eq:1.4} in these cases. But then Theorem \ref{Thm:2.81}
applies   showing that in the above cases $\g$ is either $W(m+1;{\un
1})$ or $S(m+1;{\un 1})^{(1)}$ or $ S(m+1;{\un
1})^{(1)}\oplus\F{\mathfrak D}_1$.

\m
\noindent (g)  Now suppose that $S^2(V) \ot V^*$ is not completely
reducible and $\mathfrak{g}_1 \neq \mathfrak{g}'_1$. As
$\widetilde{\mathfrak{g}}_1$ is isomorphic to a maximal submodule
of $S^2(V) \ot V^*$, then $m+2 \equiv 0 \mod p$ and
$\mathfrak{g}_1 \cong S^2(V) \ot V^*$ (see Theorem
\ref{Thm:2.82}). In particular, $\widetilde{\mathfrak{g}}_1 \neq
0$. As $m+1 \not \equiv 0 \mod p$, we have
$\mathfrak{gl}(V)=\mathfrak{sl}(V)\oplus\F\,\text{id}_V$. We
regard $\text{id}_V$ as the degree derivation and embed $\g$ into
the graded Lie algebra $\F\,\text{id}_V+\g$. Thus, we may assume
in this part that $\g_0=\mathfrak{gl}(V)$. Our nearest goal is to
determine the bilinear mapping $V^*\times \big(S^2(V)\otimes
V^*\big)\,\longrightarrow\,\mathfrak{gl}(V)$ given by the Lie
bracket in $\g$.

  First let us show that $\dim \Hom_{\mathfrak{g}_0}(\mathfrak{g}_{-1} \ot
  \mathfrak{g}_1, \mathfrak{g}_0) = 3$.
Clearly,

\begin{eqnarray*} \Hom_{\mathfrak{g}_0}\left (\mathfrak{g}_{-1}
\ot \mathfrak{g}_1, \mathfrak{g}_0\right )
&\cong& \Hom_{\mathfrak{g}_0}(V^* \ot S^2(V) \ot V^*, V \ot V^*)\\
&& \hspace*{-3.2cm} \cong  \Big(V \ot (S^2(V))^* \ot V \ot V^* \ot V\Big)^{\g_0}\\
&& \hspace*{-3.2cm}\cong  \Big(\big(S^2(V) \ot V\big)^* \ot V \ot V \ot V\Big)^{\g_0}\\\
&&\hspace*{-3.2cm}\cong \Hom_{\mathfrak{g}_0}\big(S^2(V) \ot V, V \ot V \ot V\big)  \\
&&\hspace*{-3.2cm}\cong \Hom_{\mathfrak{g}_0}(S^2(V) \ot V, (S^2(V) \oplus \wedge^2(V)) \ot V) \\
&&\hspace*{-3.2cm}\cong \Hom_{\mathfrak{sl}(V)}\left(V(2\varpi_1)
\ot V(\varpi_1), V(2\varpi_1) \ot V(\varpi_1) \oplus V(\varpi_2)
\ot V(\varpi_1)\right),
\end{eqnarray*}
where $( \,\cdot \,)^{\g_0}$ indicates the $\g_0$-invariants.  
Since $m+1 \not \equiv 0 \mod p$, the trace form $\Psi$ of the
linear Lie algebra $\mathfrak{gl}(V)$ is nondegenerate on
$\mathfrak{sl}(V)$. By \cite[Chap.~VIII, Sec.~6]{Bou2}, the
Casimir element of $\mathfrak{U}(\mathfrak{sl}(V))$ associated
with $\Psi$ acts on an $\mathfrak{sl}(V)$-module with highest
weight $\lambda$ as the scalar $(\lambda, \lambda+2\rho)$ where
$\rho$ is the half-sum of the positive roots, and $(\,\, , \,)$ is
the corresponding $W$-invariant form on $\tf^*$. We may assume
that
$$(\mathfrak t \cap \g_0^{(1)})^* = \left \{ \sum_{i = 1}^{m+1} \lambda_i \e_i \,\Big | \,
\lambda_1 + \cdots + \lambda_{m+1} = 0\right \}$$
where the
$\lambda_i \in {\mathbb F}$, and where $(\e_i, \e_j) =
\delta_{i,j}$ for $1 \leq i,j \leq m+1$. Then the fundamental
weights are given by
$$\varpi_i = \e_1 + \cdots + \e_i - \frac{i}{m+1}(\e_1 + \dots + \e_{m+1}),
\quad 1 \leq i \leq m.$$
As $m+2 \equiv 0 \mod p $, $\rho =
\varpi_1 + \cdots + \varpi_m = (m+1)\e_1 + m \e_2 + \cdots + 2
\e_m + \e_{m+1}$, and so $(\varpi_i,\rho) = (m+1) + \cdots +(m-i+2)
- \frac {i}{m+1}(1+2 + \cdots + m+1) =\hf\Big( (m+1)(m+2) -
(m-i+1)(m-i+2) - i(m+2)\Big)$.  Therefore, for $1 \leq i \leq j \leq m$,
we have
\begin{eqnarray*}
(\varpi_i, 2\rho) &=&  -i(i+1) \quad \text {and} \\
(\varpi_i, \varpi_j) &=& i- \frac{ij}{m+1}=i(j+1).
\end{eqnarray*}   Using those values for $i=1,2,3$,  it is easy to see the
following: \vspace{20pt}
$$\vbox{
\offinterlineskip \halign{ \strut \vrule  $\hfil#\hfil$  & \vrule
$\hfil#\hfil$ & \vrule  $\hfil#\hfil$  & \vrule  $\hfil#\hfil$
\vrule\cr \noalign{\hrule} \,\lambda & \,3\varpi_1 & \,\varpi_1 +
\varpi_2 & \,\varpi_3  \cr \noalign{\hrule}\,(\lambda,
\lambda+2\rho) & \,12 & \,6 & \,0 \cr \noalign{\hrule} } }$$
Since
the dominant weights of the $SL(V)$-modules
$V(2\varpi_1) \ot V(\varpi_1)$, $V(\varpi_2) \ot V(\varpi_1)$,
$V(3 \varpi_1)$,  $V(\varpi_1$ $+ \varpi_2)$, and $V(\varpi_3)$
belong to the set $\{  3 \varpi_1, \varpi_1$ $+ \varpi_2,$
$\varpi_3\}$, a standard argument shows that the Weyl modules $V(3
\varpi_1)$, $V(\varpi_1$ $+ \varpi_2)$, and $V(\varpi_3)$ are
irreducible, and the tensor products $V(2\varpi_1) \ot
V(\varpi_1)$ and $V(\varpi_2) \ot V(\varpi_1)$ are completely
reducible with
\begin{eqnarray} & V(2 \varpi_1) \ot V(\varpi_1) \cong V(3 \varpi_1)
\oplus V(\varpi_1 + \varpi_2) \nonumber \\
& V(\varpi_2) \ot V(\varpi_1) \cong V(\varpi_1 + \varpi_2) \oplus
V(\varpi_3).\nonumber
\end{eqnarray}
Thus, it follows from the above
Casimir-element calculations that
\begin{eqnarray*}
&&\hspace*{-.5cm}\Hom_{\mathfrak{g}_0}(\mathfrak{g}_{-1} \ot \mathfrak{g}_1, \mathfrak{g}_0) \\
&& \cong \Hom_{\mathfrak{gl}(V)}\left(V(3 \varpi_1) \oplus V(\varpi_1 + \varpi_2), V(3\varpi_1)
\oplus 2V(\varpi_1 + \varpi_2) \oplus V(\varpi_3)\right)\\
&& \cong {\mathbb F} \oplus {\mathbb F} \oplus {\mathbb F}.
\end{eqnarray*}

Let $v_1, \dots, v_{m+1}$ be a basis for $V$ and let $v_1^*,
\dots, v_{m+1}^*$ be the dual basis in $V^*$.   We use $\{v_iv_j
\mid 1 \leq i,j \leq m+1\}$ as a basis for $S^2(V)$ and for
brevity write $v_iv_i$ as $v_i^2$.   We identify $E_{i,j} \in
\mathfrak{gl}(V)$ with $v_i \ot v_j^* \in V \ot V^*$, so that
$E_{i,j}v_k = v_j^*(v_k)v_i = \delta_{j,k}v_i$ and $E_{i,j}v_k^*= -v_k^*(v_i)v_j^* = -\delta_{i,k}v_j^*$. 
We define three maps
$\zeta, \eta, \theta \in \Hom\left(V^* \ot S^2(V) \ot V^*,
\mathfrak{gl}(V)\right)$ by
\begin{eqnarray}
& \zeta(v_i^* \ot v_j v_k \ot v_\ell^*) = \delta_{i,j} E_{k,\ell}
+ \delta_{i,k} E_{j,\ell} \nonumber \\
& \eta(v_i^* \ot v_j v_k \ot v_\ell^*) =
\delta_{j,\ell}E_{k,i} + \delta_{k,\ell} E_{j,i} \nonumber \\
& \theta(v_i^* \ot v_j v_k \ot v_\ell^*) =
\Big(\delta_{i,j}\delta_{k,\ell} +
\delta_{i,k}\delta_{j,\ell}\Big)\id_V. \nonumber \end{eqnarray}
Direct verification shows that these are $\mathfrak{gl}(V)$-module
homomorphisms.  It is evident that $\eta$ and $\theta$ are
linearly independent, and since $\eta(v_1^* \ot v_1v_1 \ot v_2^*)$
$= 0$ $= \theta(v_1^* \ot v_1v_1 \ot v_2^*)$ and $\zeta(v_1^* \ot
v_1 v_1 \ot v_2^*)$ $= 2E_{1,2}$ $\neq 0$, we have $\zeta$ $\not
\in {\mathbb F} \eta$ $+ {\mathbb F} \theta$.   Therefore,
\begin{equation*}\Hom_{\mathfrak{gl}(V)}\left(V^* \ot S^2(V) \ot V^*,
\mathfrak{gl}(V)\right) = {\mathbb F}\zeta \oplus {\mathbb F} \eta
\oplus {\mathbb F} \theta.\end{equation*}
Let $[\,,\,]$ be the Lie
bracket in $\mathfrak{g}$.  We identify
$\mathfrak{g}_{-1},\mathfrak{g}_0, \mathfrak{g}_1$ with $V^*$,
$\mathfrak{gl}(V)$, and $S^2(V) \ot V^*$ respectively. Then $[
\,,\,]\mid_{\mathfrak{g}_{-1} \times \mathfrak{g}_1} = a \zeta + b
\eta + c \theta$ for some $a,b,c \in {\mathbb F}$. Since
$\mathfrak{g}$ is transitive,  $a \neq 0$, as otherwise $[v_i^*,
v_1^2 \ot v_2^*] = 0$ for all $1 \leq i \leq m+1$. Furthermore,  
\begin{eqnarray*}
[[v_i^*,v_n^*],v_jv_k \ot v_\ell^*] &=&
[v_i^*,[v_n^*,v_jv_k\ot v_\ell^*] - [v_n^*,[v_i^*,v_jv_k \ot v_\ell^*] \\
&& \hspace*{-1.2cm}=  [\,v_i^*, a(\delta_{j,n}E_{k,\ell} + \delta_{k,n}E_{j,\ell}) +
b(\delta_{j,\ell}E_{k,n} + \delta_{k,\ell}E_{j,n})\,]\\
&& \hspace*{-1.2cm} \quad  +[\,v_i^*,c\left(\delta_{j,n} \delta_{k,\ell} +
\delta_{k,n}\delta_{j,\ell}\right)\id\,]  -[\,v_n^*, a(\delta_{i,j}E_{k,\ell}+
\delta_{i,k}E_{j,\ell})\,]  \\
&& \hspace*{-1.2cm}\qquad  -[\,v_n^*, b(\delta_{j,\ell}E_{k,i} + \delta_{k,\ell}E_{j,i}) +
c\left(\delta_{i,j} \delta_{k,\ell} + \delta_{i,k}\delta_{j,\ell}\right)
\id\,]   \\
&& \hspace*{-1.2cm}= a(\delta_{j,n}\delta_{i,k} + \delta_{k,n}\delta_{i,j}) v_\ell^* +
b(\delta_{j,\ell}\delta_{i,k} + \delta_{k,\ell}\delta_{i,j})
v_n^* \\
&& \hspace*{-1.2cm} \quad \; \; + c(\delta_{j,n} \delta_{k,\ell} + \delta_{k,n}
\delta_{j,\ell})v_i^* - a(\delta_{i,j}\delta_{k,n} +
\delta_{i,k}\delta_{j,n})
v_\ell^* \nonumber \\
&& \hspace*{-1.2cm} \qquad \; \; - b(\delta_{j,\ell}\delta_{k,n} +
\delta_{k,\ell}\delta_{j,n}) v_i^* - c(\delta_{i,j} \delta_{k,\ell}
+ \delta_{i,k} \delta_{j,\ell})v_n^* \nonumber \\
&& \hspace*{-1.2cm}= (b-c)\Big( (\delta_{j,\ell}\delta_{i,k} +
\delta_{k,\ell}\delta_{i,j})v_n^* -(\delta_{j,\ell}\delta_{k,n} +
\delta_{k,\ell}\delta_{j,n}) v_i^*\Big) \end{eqnarray*}
for
$1 \leq i,j,k,\ell, n \leq m+1$.   Suppose
$\mathfrak{g}_{-2} \neq 0$. The 1-transitivity \eqref{eq:1.4} of
$\mathfrak{g}$ implies $[\mathfrak{g}_{-2},\mathfrak{g}_1] \neq
0$, which forces $b-c \neq 0$. Thus, $a(b-c) \neq 0$. But then
\begin{eqnarray*} \left[ [v_m^*,v_{m+1}^*], [v_m^2 \ot v_m^*,v_{m+1}^2 \ot
v_1^*]\right] &=& \left[[[v_m^*,v_{m+1}^*],v_m^2 \ot v_m^*\right],
v_{m+1}^2 \ot v_1^*] \\
&=& 2(b-c)[v_{m+1}^*,
v_{m+1}^2 \ot v_1^*] \\
&=& 4a(b-c)E_{m+1,1} \neq 0. \end{eqnarray*}
Set
\begin{equation*}\begin{gathered} \mathfrak{n}^+ = \bigoplus_{1 \leq i < j \leq m+1} {\mathbb F}
E_{i,j}, \qquad  \mathfrak{t} = \bigoplus_{1 \leq i \leq m+1} {\mathbb F}
E_{i,i}, \quad  \hbox{ \rm and }\\ \mathfrak{n}^- = \bigoplus_{1 \leq j < i \leq m+1} {\mathbb F} E_{i,j}, 
\end{gathered} \end{equation*}
and let $\mathfrak b^+ = \mathfrak t \oplus \mathfrak n^+$.  Clearly $u = [v_m^*,v_{m+1}^*]$ is a $\mathfrak b^+$-primitive vector of weight $-(\varepsilon_m+\varepsilon_{m+1})$  in 
$\mathfrak{g}_{-2}$,  and $\mathfrak{g}_{-2}$ $=
\mathfrak{U}(\mathfrak{n}^-)[v_m^*,v_{m+1}^*]$ because $\g_{-2}= [\g_{-1},\g_{-1}]$.   Set $w$ $= [v_m^2 \ot
v_m^*,v_{m+1}^2 \ot v_1^*]$. As $[\mathfrak{n}^-,E_{m+1,1}]$ $=
0$, we have $[E_{i,j}u,w]$ $+ [u,E_{i,j}w]$ $= 0$ for all $i$ $>
j$. If $k$ $< m$, then
\begin{eqnarray*} [[v_k^*,v_{m+1}^*], [v_m^2 \ot v_m^*, v_{m+1}^2 \ot
v_1^*]] &=&
[[[v_k^*,v_{m+1}^*], v_m^2 \ot v_m^*], v_{m+1}^2 \ot v_1^*] \\
&&+ [v_m^2 \ot v_m^*, [[v_k^*,v_{m+1}^*], v_{m+1}^2 \ot v_1^*]] \\
&=& 0. \end{eqnarray*} Similarly,
\begin{eqnarray} [[v_k^*,v_m^*], [v_m^2 \ot
v_m^*, v_{m+1}^2 \ot v_1^*]] &=&
[[[v_k^*,v_m^*], v_m^2 \ot v_m^*], v_{m+1}^2 \ot v_1^*] \nonumber \\
&=& 2(c-b)[v_k^*, v_{m+1}^2 \ot v_1^*] = 0. \nonumber
\end{eqnarray}       Since obviously, $[[v_i^*, \, v_j^*], \,
[v_{m+1}^2 \ot v_m^*, v_{m+1}^2 \ot v_1^*]] = 0$ for all $i,j \leq
m - 1,$ we conclude that $[[\mathfrak{n}^-, \, \mathfrak{g}_{-2}],
\, w] = 0$.   Inasmuch as $[\mathfrak{g}_{-2}, \, w] = {\mathbb
F}E_{m+1,1}$ by the above computations, we have $[\mathfrak{g}_{-2},
[\mathfrak{n}^-, w] ]= [[\mathfrak{g}_{-2}, \mathfrak{n}^-], \, w] =
0.$

Set $\mathfrak{g}_2' \eqdef \{x \in \mathfrak{g}_2 \mid
[\mathfrak{g}_{-2}, \, x] = 0\}.$ We have shown that $w \not \in
\mathfrak{g}_2',$ and that $[\mathfrak{n}^- ,w] \subseteq
\mathfrak{g}_2'.$  Let $\mathfrak{g}^{\{2\}}$ be the graded
subalgebra of $\mathfrak{g}$ generated by $\mathfrak{g}_{-2},$
$\mathfrak{g}_0,$ and $\mathfrak{g}_2$.   Let $M^{\{2\}}$ denote the
maximal graded ideal of $\mathfrak{g}^{\{2\}}$ contained in the
positive part $\bigoplus_{i>0}\mathfrak{g}_i^{\{2\}}$.    Set
$\overline{\mathfrak{g}}^{\{2\}} :=
\mathfrak{g}^{\{2\}}/M^{\{2\}}$. Clearly,
$\overline{\mathfrak{g}}^{\{2\}}$ is graded and transitive. It is
easy to check that $M^{\{2\}}\cap \mathfrak{g}_2 =
\mathfrak{g}_2'.$ Hence, $\overline{\mathfrak{g}}_1^{\{2\}}\cong
\mathfrak{g}_2/\mathfrak{g}_2',$
$\overline{\mathfrak{g}_0}^{\{2\}} \cong \mathfrak{g}_0,$ and
$\overline{\mathfrak{g}}_{-1}^{\{2\}} \cong \mathfrak{g}_{-2}.$

Let $\overline{w}$ denote the image of $w$ under the canonical
epimorphism $\g_2\twoheadrightarrow \g_2/\g'_2$.   Since we know that
$\overline{w}$ is a $\mathfrak b^-$-primitive vector of
$\overline{\mathfrak{g}}_1^{\{2\}}$ and
$$[u, \, \overline{w}] \in {\mathbb F}^\times e_{-\alpha}\nonumber$$
(here $e_{-\alpha} = E_{m+1,1} \in \mathfrak{Z}(\mathfrak{n}^-))$
we identify $u$ with its image in
$\overline{\mathfrak{g}}^{\{2\}}$.    Since $u$ has weight
$-(\varepsilon_m + \varepsilon_{m+1})$ $= \varepsilon_1 + \cdots +
\varepsilon_{m-1}$ $= \varpi_{m-1}$ with respect to the torus
$\mathfrak t\cap \g_0^{(1)}$,  the Lie algebra
$\overline{\mathfrak{g}}^{\{2\}}$ is infinite-dimensional by
Theorem \ref{Thm:4.7}. Therefore, $\g_{-2} = 0,$ yielding
$[\mathfrak{g}_{-1}, \, \mathfrak{g}_{-1}] = 0.$ (Note that in the
case under consideration  $m \geq 3.)$ Thus, $\mathfrak{g} =
\mathfrak{g}_{-1} \oplus \mathfrak{g}_0 \oplus \mathfrak{g}_1
\oplus \cdots \oplus \mathfrak{g}_r,$ where $\mathfrak{g}_0 \cong
\mathfrak{gl}(V),$ and $\mathfrak{g}_{-1} \cong V^*$ as a
$\mathfrak{g}_0$-module. Since $\mathfrak{g}$ is transitive, we
apply Theorem \ref{Thm:2.81} to conclude that $\mathfrak{g}$ is
isomorphic to $W(m+1;\un 1).$

\m

\noindent (h)  It now remains for us to consider cases
\eqref{eq:4.23} for $m=1$, \eqref{eq:4.26} and \eqref{eq:4.27}. In
other words, we may assume that $\g_0^{(1)} = \mathfrak{sp}(V)$
and that $\mathfrak{g}_{-1} \cong V$ as a $\g_0^{(1)}$-module.
Recall from part~(a) that $V \cong V^*,$ and $\mathfrak{sp}(V)
\cong S^2(V)$ as $\mathfrak{sp}(V)$-modules, and
 $\mathfrak{g}_1$ is a $\g_0^{(1)}$-submodule of
$V \oplus \big(V \otimes S^2(V)\big).$ We claim that
$\mathfrak{g}_1$ is a completely reducible $\g_0^{(1)}$-module.

    First we recall that for any natural number $n$ there is an exact sequence,
the Koszul resolution (see \cite[p.~377]{Ja}),

\begin{eqnarray}\label{eq:4.28}
\cdots  \rightarrow S^{n-i}(V)\otimes \wedge^i V &\rphi& S^{n-i+1}(V)\otimes
\wedge^{i-1}V \rightarrow  \\
&\rightarrow& \cdots \rightarrow S^{n-1}(V)\otimes V   \rightarrow
S^n(V) \rightarrow 0. \nonumber\end{eqnarray}

\noindent The map $\varphi_i$ is given by

$$\varphi_i(x \otimes (v_1 \wedge \cdots \wedge v_i)) \eqdef \sum_{j=1}^{i-1}(-1)^ixv_j
\otimes (v_1 \wedge \cdots \wedge \widehat{v_j} \wedge \cdots \wedge v_i)\nonumber$$

\noindent for all $v_j \in V$ and $x \in S^{n-i}(V).$  If $V$ is a
module for a group $H$,   then (\ref{eq:4.28}) is an exact
sequence of $H$-modules.  If $p$ does not divide $n$,
then the exact sequence (\ref{eq:4.28}) splits. More
precisely, define

$$\psi_j: S^j(V) \otimes \wedge^{n-j} V \rightarrow S^{j-1}(V) \otimes
\wedge^{n-j+1}V\nonumber$$

\noindent by

$$\psi_j(w_1 \cdots w_j \otimes w) = \sum_{i=1}^j w_1 \cdots \widehat{w}_i
\cdots w_j \otimes w_i \wedge w\nonumber$$

\noindent for all $w_i \in V$ and $w \in \wedge^{n-j}V.$  Then all
the $\psi_j$ are $H$-module homomorphisms, and

$$\varphi_{i+1} \circ \psi_{n-i} + \psi_{n-i+1} \circ \varphi_i  =
m \cdot \id_{S^{n-i}(V) \otimes \wedge^iV}\nonumber$$

\noindent for all $i,$ $0 \leq i \leq n.$  (We set $\varphi_{n+1}
= \psi_{n+1} = 0.)$   If we set $n = 3$ and $H =G_0= \hbox{\rm
Sp}(V),$ we get that there is a split exact sequence of
$G_0$-modules

$$0 \rightarrow  E\,\, {\buildrel{\eta} \over  \rightarrow}\,\,
  S^2(V) \otimes V \rightarrow  S^3(V) \rightarrow  0.\nonumber$$

\noindent in which $E$ is a homomorphic image of $V \otimes
\wedge^2 V$ (as seen from (\ref{eq:4.28}) with $n=3$).  This implies that $S^2(V) \otimes V \cong S^3(V)
\oplus E$.

    According to part~(a), $\mathfrak{g}_1$ is a $G_0$-stable subspace
    of $V \oplus \bigl(S^2(V) \otimes V\bigr)$.   Thus, we
may assume that $\mathfrak{g}_1$ is a $G_0$-submodule of $V$
$\oplus$ $E$ $\oplus$ $S^3(V).$ Let
$\pi'\colon\,\g_1\,\rightarrow\,V\oplus E$ and
$\pi''\colon\,\mathfrak{g}_1\,\rightarrow\, S^3(V)$ denote the
corresponding projection maps, and put $U' = \ker \pi'$ and
$U''=\ker \pi''$. If $U''=0$,  then $\mathfrak{g}_1$ is isomorphic
to a nonzero submodule of the $G_0$-module $S^3(V)$. Since $p
> 3$,   Proposition  \ref{Pro:2.848} implies  that $S^3(V)\cong L(3\varpi_1)$ is irreducible over
$\g_0^{(1)}$. As a consequence, $\mathfrak{g}_1$ is irreducible if
$U''=0$.

Thus, in proving the claim we may assume that $U''\ne 0$. Note
that

\begin{equation}\label{eq:C}
X_+(V\otimes\wedge^2V)\,=\,\{\varpi_1+\varpi_2,\varpi_3,\varpi_1\},
\end{equation}
where $\varpi_1+\varpi_2$ is omitted if $m=1$ and $\varpi_3$ is
omitted if $m=1,2$. Since $p>3$ and $\alpha=2\e_1$, it is
straightforward to see that the weight subspace $\left(V \otimes
\wedge^2V\right)^{-\Lambda - \alpha}$ with respect to $\mathfrak t
= \hbox{\rm Lie}(\widetilde{T})$ is trivial. Together with our
remarks in part~(c), this yields that $- \lambda = - \varpi_1\in
X(T)$ is the only minimal weight of $U''$ and $ \dim\,
(U'')^{-\lambda}=1$. But then $-w_0\lambda=\varpi_1$ is the only
maximal weight of $U''$. On the other hand, (\ref{eq:B}) says
that $\nu\ge\varpi_1$ for any dominant weight $\nu$ of the
$G_0$-module $\g_1$. It follows that $\varpi_1$ is the only
dominant weight of $U''$. Hence, $U''\cong L(\varpi_1)\cong V$ as
$G_0$-modules.

Now since $e^\Gamma\in\g_1$ and $(U'')^{\Gamma} = (U'')^{-\Lambda - \alpha} = 0$, we
have $\pi''\neq 0$. Hence, $\pi''(\mathfrak{g}_1) = S^3(V)$,
thanks to Proposition \ref{Pro:2.848}. Therefore $\g_1^{3\varpi_1} \ne
0$. On the other hand, (\ref{eq:C}) shows  that
$\g_1^{3\varpi_1}\subseteq U'$. Since $U'\subseteq S^3(V)$, Proposition
\ref{Pro:2.848} implies  that $U'\cong \pi''(\g_1)$. But then
$\mathfrak{g}_1 \cong V \oplus S^3(V)$ in the present case, to
prove the claim.

\m
\noindent (i) It follows from our discussion in part~(h) that
$\g_1=\g'_1\,\oplus\, \widetilde{\g}_1$, where we assume that
$\widetilde{\g}_1=0$ if $\g_1$ has no $\mathfrak b^-$-primitive vectors of
weight $-\Lambda$. If $\mathfrak{g}_1$ $= \mathfrak{g}'_1$, then
we have $[\mathfrak{g}_{-2},$ $\mathfrak{g}_1]$ by part~(b), so
that $\mathfrak{g}_{-2}$ $= 0$ by the 1-transitivity
 of $\mathfrak{g}.$ Therefore, Theorem \ref{Thm:2.81}
applies, and we can conclude that $\mathfrak{g}$ is either
$H(2m;{\un 1})^{(2)}$ or $H(2m;{\un 1})^{(2)}\oplus\F{\mathfrak
D}_1.$

Thus, in what follows, we may assume that $\mathfrak{g}_1=
\mathfrak{g}'_1\,\oplus\, \widetilde{\mathfrak{g}}_1$,  where
$\widetilde{\mathfrak{g}}_1\cong V$ is generated by a nonzero
$\mathfrak b^-$-primitive vector $e^{-\Lambda}$ $\in \mathfrak{g}_1$.   We have
already established that $[\mathfrak{g}_{-2},$ $\mathfrak{g}'_1]=
0$.   As before, denote by $\widetilde{\mathfrak{g}}$ the subalgebra
generated by $f^{\Lambda},$ $\mathfrak{g}_0,$ and $e^{-\Lambda}.$
We know from part~(e) that the Lie algebra
$\widetilde{\mathfrak{g}}/\mathcal M(\widetilde{\mathfrak{g}})$ is
either of type \eqref{eq:4.26}  or \eqref{eq:4.27} or of type
\eqref{eq:4.23} for $m=1$. This implies that $\dim
\mathfrak{g}_{-2}/\mathcal M(\widetilde{\mathfrak{g}})_{-2}\le 1$
and $\mathfrak{g}_{-3}= \mathcal
M(\widetilde{\mathfrak{g}})_{-3}$. Since $[\mathcal
M(\widetilde{\mathfrak{g}})_{-2},\,\widetilde{\mathfrak{g}}_1]= 0
= [\mathfrak{g}_{-2},\,\mathfrak{g}'_1]$ and $\mathfrak{g}_1=
\mathfrak{g}'_1 \, \oplus \, \widetilde{\mathfrak{g}}_1$,  we must
have $\mathcal M(\widetilde{\mathfrak{g}})_{-2}= 0$ by the
1-transitivity of $\mathfrak{g}.$  Hence, $\dim
\mathfrak{g}_{-2}\le 1.$  If $\mathfrak{g}_{-3}$ $=
[\mathfrak{g}_{-2},$ $\mathfrak{g}_{-1}]$ $\neq 0,$ then the
$\mathfrak{sp}(V)$-modules $\mathfrak{g}_{-3}$ and
$\mathfrak{g}_{-1}$ are isomorphic.   Moreover, because
$[\mathfrak{g}_{-3},\,\widetilde{\mathfrak{g}}_1]\subseteq
\mathcal M(\widetilde{\mathfrak{g}})_{-2} \, = 0$,  it follows that
$[\mathfrak{g}_{-3},\,\mathfrak{g}'_1]= \mathfrak{g}_{-2}$ by the
$1$-transitivity  of $\mathfrak{g}.$ But then $\mathfrak{g}'_1
\cong \mathfrak{g}_{-3}^*$ as $\g_0^{(1)}$-modules, for both
$\mathfrak{g}_{-3}$ and $\mathfrak{g}'_1$ are irreducible. This,
however, contradicts the fact that $\dim V$ $< \dim S^2(V).$ Thus 
$\mathcal M(\widetilde{\mathfrak{g}})= 0$;  that is,  $\dim
\mathfrak{g}_{-2}\le 1,$ and $\mathfrak{g}_{-3}$ $= 0.$ One can
now apply Theorem \ref{Thm:2.81} to conclude that either
$\mathfrak{g}\cong K(2m+1;{\un 1})^{(1)}$ or $\mathfrak{g}\cong
W(2;{\un 1})$.

    The proof of Theorem \ref{Thm:4.22}   is now finished. \qed

 \m
\section {\ The irreducibility of \,$\boldsymbol{\g_1}$ \label{sec:4.7}}
\m

    Here we investigate the case where
$[f^\Lambda,e^\Gamma] = e_{\alpha}$ for some positive root
$\alpha\in\Phi$. This case occurs for Lie algebras of type $S$, $CS$, $H$, or
$CH$  when the natural grading is reversed, and our ultimate goal
is to exclude all other possibilities. As a first step, we will
show that   $\g_1$  is an irreducible $\g_0$-module.
Essential to the proof is the similarity between the structure of an irreducible restricted
$\g_0^{(1)}$-module and that of a rational module over the simply connected algebraic
group corresponding to $\g_0^{(1)}$.

 \m  Under the assumptions of the next three lemmas,  $\g_0^{(1)}$ is of type A$_m$ or
C$_m$, and so $\g_0^{(1)} \cong \mathfrak{sl}_{m+1}$, $\mathfrak{psl}_{m+1}$,  or
$\mathfrak{sp}_{2m}$.     Given $\psi\in\tf^*$,   we denote by $\bar{\psi}$
 the restriction of $\psi$ to the subspace
$\tf\cap\g_0^{(1)}$.      We adopt the notation of Section
\ref{sec:2.3}.\bi

\begin{Lem} \label{Lem:4.30}   \ Let $\mathfrak{g}$ be a
graded Lie algebra satisfying conditions
(1)-(4) and (6).   Suppose there exist a ${\mathfrak b}^+$-primitive
vector $f^\Lambda \in \mathfrak{g}_{-1}$ of weight $\Lambda$ and 
a $\mathfrak b^-$-primitive
vector $e^\Gamma \in \mathfrak{g}_1$ of weight $\Gamma$ such that
$[f^\Lambda,e^\Gamma] = e_{\alpha}$  for some root $\alpha \in \Phi^+$.
Then $\g_0^{(1)}$ consists of only one summand, which is of type
\hbox{\rm A}$_m$ or  \hbox{\rm C}$_m$, and $\mathfrak{g}_1$ is an
irreducible standard $\g_0^{(1)}$-module.  \end{Lem}

\pf  Replacing $\mathfrak{g}$ by its subalgebra
$\widehat{\mathfrak{g}}$ generated by $\mathfrak{g}_{-1} \oplus
\mathfrak{g}_0 \oplus \mathfrak{g}_1$ and passing to the graded
quotient $\widehat{\g}/{\mathcal M}(\widehat{\g})$, we can assume
that $\mathfrak{g}$ is generated by its local part
$\mathfrak{g}_{-1} \oplus \mathfrak{g}_0 \oplus \mathfrak{g}_1$.
By Proposition \ref{Pro:4.20},  $\Gamma(\mathfrak{t}^{[i]}) \neq
0$ for any $i = 1, \dots, \ell$. Then Theorem \ref{Thm:4.7}
applies to show  $\g_0^{(1)}$ has a unique summand which is of
type A$_m$ or C$_m$;  \  $\alpha$ is the highest root;  \ and
$\Gamma$ is the lowest weight of an irreducible standard
$\g_0^{(1)}$-submodule of $\mathfrak{g}_1$.   Thus, we have
\begin{equation*}
\begin{array}{ccc}
&\alpha=\varpi_1+\varpi_m; \ \
-\bar{\Gamma}\in\{\varpi_1,\varpi_m\};  \ \  \hbox{\rm and} \\
& \bar{\Lambda}=2\varpi_1+\varpi_m \  \hbox{\rm or} \ 
\varpi_1+2\varpi_m& \quad \  (\hbox{\rm A}_m) \\
&\alpha=2\varpi_1; \ \ -\bar{\Gamma}=\varpi_1;  \ \
\hbox{\rm and}\ \ \bar{\Lambda}=3\varpi_1 
& \quad \ \, (\hbox{\rm C}_m).  \end{array}  \end{equation*}

Let $V$ be an irreducible $\g_0$-submodule of $\g_1$. Since $\g_1$
embeds into $\text{Hom}(\g_{-1},\g_0)$ by transitivity
\eqref{eq:1.3}, the $\g_0^{(1)}$-module $V$ is restricted. Let
$e=e^\Theta$ be any $\mathfrak b^-$-primitive vector from $V$ of
weight $\Theta\in\tf^*$. Then $V=\mathfrak{U}({\mathfrak
n}_+)\,e$. Since $\g$ is irreducible and transitive,
$[f^\Lambda,e]\ne 0$.  Let $\widetilde{\g}$ denote the subalgebra
of $\g$ generated by $f^\Lambda, \mathfrak{g}_0$, and $e$, and let
${\mathcal M}(\widetilde{\g})$ be the Weisfeiler radical of
$\widetilde{\g}$.

If $\Theta \neq -\Lambda$, then Theorem \ref{Thm:4.7} shows that
$\Theta = \Gamma$.   We suppose that $\Theta=-\Lambda$. Then
the Lie algebra $\widetilde{\g}/{\mathcal M}(\widetilde{\g})$ is
classical, and its grading is standard according to  Theorem
\ref{Thm:3.34}. If $\g_0^{(1)}$ is of type $\mathrm{A}_m$ or
$\mathrm{C}_m$, where $m\ge 2$, then using \cite{Bou1} it is not
difficult to verify that $\bar{\Lambda}$ is in the following list:

\begin{eqnarray} & \{2 \varpi_1, \varpi_1,
\varpi_2, \varpi_{m-1},
\varpi_m, 2 \varpi_m\}, \qquad \quad \  m \geq 8, \nonumber \\
& \{2 \varpi_1, \varpi_1, \varpi_2, \varpi_3, \varpi_5,\varpi_6,
\varpi_7, 2 \varpi_7\}, \qquad m = 7, \nonumber \\
& \{\varpi_i \mid 1 \leq i \leq m \} \cup \{2 \varpi_1, 2
\varpi_m\}, \qquad
2 \leq m \leq 6, \nonumber \\
\nonumber\end{eqnarray}

\noindent for type A$_m$, and
\begin{eqnarray}
& \{\varpi_1\}, \qquad \qquad \quad  m \geq 4, \nonumber \\
& \{\varpi_1,\varpi_m\}, \qquad m = 2,3, \nonumber \\
\nonumber\end{eqnarray}
\noindent for type C$_m$.   But we have
shown that $\bar{\Lambda} = 2\varpi_1 + \varpi_m$ or
$\bar{\Lambda} = \varpi_1 + 2 \varpi_m$ if $\g_0^{(1)}$ is of type
A$_m$, and $\bar{\Lambda} = 3 \varpi_1$ if $\g_0^{(1)}$ is of type
C$_m$.  These linear functions do not appear on the list.  Hence,
$\Theta \neq -\Lambda$ for $m \geq 2$.

If $\g_0^{(1)} = \mathfrak{sl}_2$, then $\bar{\Lambda} = 3
\varpi_1$. Suppose $\Theta = -\Lambda$.  Then the classical Lie
algebra $\widetilde{\g}/{\mathcal M}(\widetilde{\g})$ has the
following Dynkin diagram of type G$_2$:

\begin{equation}\label{eq:4.620}
{\beginpicture \setcoordinatesystem
units <0.45cm,0.3cm> % sets scale
 \setplotarea x from -6 to 8, y from -1 to 1    % sets plot size up
 \linethickness=0.03pt                          % sets line thickness
  \put{$\circ$} at 0 0
 \put{$\bullet$} at 2 0
 \plot .15 -.3 1.85 -.3 /
  \plot .15 0 1.85  0 /
  \plot .15 .3 1.85 .3 /
  \put {$<$} at 1 0
   \put{$\alpha_1$} at 0 -1.5
  \put{$-\Lambda$} at 2 -1.5 \endpicture}\end{equation}

\bi

 Set $f' = [f^\Lambda,[f^\Lambda,e_{-\alpha_1}]]$ and
$e' = [e^\Gamma,e^{-\Lambda}]$.   Then

\begin{eqnarray}
[f',e'] &=& \Big( (F^{\Lambda})^2 E_{-\alpha_1} - 2 F^{\Lambda}
E_{-\alpha_1}F^{\Lambda}
+ E_{-\alpha_1}(F^{\Lambda})^2 \Big)\left ([e^\Gamma,e^{-\Lambda}]\right) \nonumber \\
&=& \Big(- 2 F^{\Lambda} E_{-\alpha_1}
+ E_{-\alpha_1}F^{\Lambda}  \Big)\left([e_{\alpha_1},e^{-\Lambda}] +
[h_{-\Lambda},e^\Gamma]\right) \nonumber \\
&=& 2 \Lambda(h_{\alpha_1}) h_{-\Lambda} - \Big(
\alpha_1(h_{-\Lambda}) + 2
+ \alpha_1(h_{-\Lambda})\Big) h_{\alpha_1} \nonumber \\
&=& 6 h_{-\Lambda}. \nonumber \\
\nonumber\end{eqnarray} \noindent Here $F^{\Lambda}$ $= \ad
f^{\Lambda},$ $E_{-\alpha_1}$ $= \ad e_{-\alpha_1}$, and we assume
that

$$[e^{-\Lambda}, f^\Lambda] = h_{-\Lambda}, \quad
[h_{-\Lambda},e^{-\Lambda}] = 2 e^{-\Lambda}, \quad  \hbox{\rm and} \quad
[h_{-\Lambda}, f^\Lambda] = -2f^\Lambda.$$  Moreover,
$(2\Lambda - \alpha_1)(h_{\alpha_1})= 4$ and $[\mathfrak{n}^+,$
$f']$ $= 0$ $= [\mathfrak{n}^-, e']$. Setting

\begin{equation*}\begin{array}{lll}  &e_1 = e_{\alpha_1},  &\qquad\qquad f_1 = e_{-\alpha_1},
 \\
&e_2 = \frac {1}{3}[e^\Gamma,e^{-\Lambda}],  &\qquad\qquad f_2 =
-\frac{1}{3}[f^\Lambda,[f^\Lambda,e_{-\alpha_1}]],\\
&h_1 = h_{\alpha_1},  &\qquad\qquad h_2 = \frac {2}{3} h_{-\Lambda}, \\
\end{array}\end{equation*}

\noindent we have from the calculations in the previous
paragraph,

\begin{equation*}\begin{array}{llll}&[e_i,f_j] = \delta_{i,j}h_i,& \qquad &[h_i,h_j] = 0,  \\
&[h_i,e_j] = A_{i,j}e_j, & \qquad &[h_i,f_j] = -A_{i,j}f_j
\qquad\quad (1 \leq i,j \leq 2), \end{array}\end{equation*} where
$A = (A_{i,j})$ is the matrix
$$ A=\left(\begin{matrix} 2 & -4 \nonumber \\
-\frac {2}{3} & 2 \end{matrix} \right).$$

\m

\noindent
Now put $e_1'\eqdef \,[e_1,[e_1,\, e_2]],\,$
$f_1'\eqdef\, [f_1,[f_1,\, f_2]]$, $u_1\eqdef \,[e_1,[e_1,[e_1,\,
e_2]]],$  and $v_1\eqdef\, [f_1,[f_1,[f_1,\, f_2]]]$. Then
\begin{eqnarray*}[e_1',f_1']&=&[[e_1,[e_1, \, e_2]],  \, [f_1,[f_1, \,
f_2]]]\\
&=&\Big(E_1^2E_2 - 2E_1E_2E_1 + E_2E_1^2\Big)\left([f_1, \,[f_1,
\,
 f_2]]\right) \\
&=&\left (-2E_1E_2 + E_2E_1\right)\left(6 [f_1, \, f_2]\right)\\
&=& 6\left(-2E_1\left (-\textstyle{\frac{2}{3}}f_1\right) +
E_2\left (4f_2\right)\right) \\
&=&8(h_1 + 3h_2).
\end{eqnarray*}
Also, $[[e_1,e_2],
\,[f_1,f_2]]\,=\,[e_1,\,[f_1,h_2]]-[e_2,\,[h_1,f_2]]\,=\,-\frac{2}{3}h_1-4h_2$,
and
\begin{eqnarray*}[u_1,v_1]&=&\big[[e_1,[e_1,[e_1 \, e_2]]],\, [f_1,[f_1,[f_1 \,
f_2]]]\big]\\
&=&\Big(E_1^3E_2 - 3E_1^2E_2E_1
+3E_1E_2E_1^2-E_2E_1^3\Big)\left([f_1,[f_1, [f_1, \,
 f_2]]]\right) \\
&=&\left(-3E_1^2E_2 + 3E_1E_2E_1-E_2E_1^2\right)\left(6 [f_1,[f_1 \, f_2]]\right)\\
&=&\left(3E_1E_2-E_2E_1\right)\left(36 [f_1 \, f_2]\right)\\
&=& 36\left(3E_1\left (-\textstyle{\frac{2}{3}}f_1\right) -
E_2\left (4f_2\right)\right) \\
&=&-72(h_1 + 2h_2).
\end{eqnarray*}
Note that $[-\frac{2}{3}h_1-4h_2,e_1]=\frac{4}{3}e_1,\,$
$[-\frac{2}{3}h_1-4h_2,e_2]=-\frac{16}{3}e_2,\,$
$[h_1+2h_2,e_1]=\frac{2}{3}e_1,$ and $[h_1+2h_2,e_2]=0$. Now put
$e_1''\eqdef\,[[e_1,e_2],\,u_1]$ and
$f_1''\eqdef\,[[f_1,f_2],\,v_1]$. Since
$[[e_1,e_2],\,v_1]=[[f_1,f_2],\,u_1]=0$, we have that
\begin{eqnarray*}[e_1'',f_1'']&=&\big[[[e_1,e_2],\,u_1],\,
[[f_1,f_2],\,v_1]\big]\\
&=&\big[[e_1,e_2],\,[[f_1,f_2],\,[u_1,v_1]]\big]-
\big[u_1,\,[[[e_1,e_2],[f_1,f_2]],\,v_1]\big]\\
&=&\big[[e_1,e_2],\,[[f_1,f_2],\,-72(h_1+2h_2)]\big]-
\big[u_1,\,[-\textstyle{\frac{2}{3}}h_1-4h_2,\,v_1]\big]\\
&=&72\big[[e_1,e_2],\,-\textstyle{\frac{2}{3}}[f_1,f_2]\big]-
[u_1,\textstyle{\frac{4}{3}}v_1]\\
&=&-48(-\textstyle{\frac{2}{3}}h_1-4h_2)-\textstyle{\frac{4}{3}}\big(-72(h_1+2h_2)\big)\\
&=&128(h_1 + 3h_2).
\end{eqnarray*}
Next observe that
\begin{eqnarray*}[e_1',f_1'']&=&\big[[e_1,[e_1,\,e_2]],\,
[[f_1,f_2],\,v_1]\big]\\
&=&\big[\big[[e_1,[e_1,\,e_2]],\,
[f_1,f_2]\big],v_1\big]\\ && \hspace{.4truein}+\big[[f_1,f_2],\,\big[[e_1,[e_1,\,e_2]],\,
[f_1,[f_1,[f_1,\,f_2]]]\big]\big]\\
&=&\big[[[[e_1,[e_1,\,e_2]],\,
f_1],f_2],\,v_1\big]+\Big[[f_1,f_2],\,[[e_1',f_1],
f_1']+[f_1,[e_1',f_1']]\Big] \end{eqnarray*}
\noindent so that 

\begin{eqnarray*} [e_1',f_1'']&=&-6\big[[[e_1,\,e_2],f_2],\,v_1\big]\\
&& \hspace{.4truein}+\Big[[f_1,f_2],\,\big[[e_1,[e_1,\,e_2]],f_1],
f_1'\big]+8[f_1,h_1+3h_2]\Big]\\
&=&-6[[e_1,h_2],\,v_1]-6\big[[f_1,f_2],\,[[e_1,\,e_2],\,f_1']\big]\\
&=&-4\big[e_1,\,[f_1,[f_1,[f_1,\,f_2]]]\big]+6\big[[f_1,f_2],\,[e_2,[e_1,[f_1,[f_1,\,f_2]]]]\big]\\
&=&-24[f_1,[f_1,\,f_2]]+36\big[[f_1,f_2],\,[e_2,[f_1,\,f_2]]\big]\\
&=&-24[f_1,[f_1,\,f_2]]+36\big[[f_1,f_2],\,[f_1,h_2]\big]\\
&=&-24[f_1,[f_1,\,f_2]]-24[[f_1,f_2],\,f_1]\,=\,0.
\end{eqnarray*}
Similarly, $[e_1'',f_1']=0$. As $[h_1+3h_2, e_1'] = 4e_1'\ne 0$
and $[h_1+3h_2, e_1''] = 8e_1''\ne 0,$ Theorem \ref{Thm:3.10} now
implies that $e_1',$ $f_1',$ $e_1'',$ and $f_1''$ generate an
infinite-dimensional Lie algebra. Since $\g$ is
finite-dimensional, we conclude that $\Theta=\Gamma$.

We let  $-\Delta$ be the base of the root system $\Phi$ and
reverse the grading of
$\widetilde{\g}/\mathcal{M}(\widetilde{\g})$, so that
$\big(\widetilde{\g}/\mathcal{M}(\widetilde{\g})\big)_i\,=\,\widetilde{\g}_{-i}/
\mathcal{M}(\widetilde{\g})_{-i}$ for all $i$. Note that
$\big(\widetilde{\g}/\mathcal{M}(\widetilde{\g})\big)_{-1}\cong V$
and $\,\big(\widetilde{\g}/\mathcal{M}(\widetilde{\g})\big)_1\cong
\g_{-1}$ as $\g_0$-modules. Then it is straightforward to see that
the Lie algebra $\widetilde{\g}/\mathcal{M}(\widetilde{\g})$
(viewed with its new grading) satisfies all the conditions of
Theorem \ref{Thm:4.22}. Since
$\big(\widetilde{\g}/\mathcal{M}(\widetilde{\g})\big)_1$ is an
irreducible $\g_0$-module, Theorem \ref{Thm:2.82} now shows that
if $\g_0^{(1)}\cong \mathfrak{sl}_{m+1}$, $m\ge 1$, then
$m+2\,\not\equiv 0 \mod p$.

     As $\g_0^{(1)}$ is almost simple and $\mathfrak
{Z}(\mathfrak{g}_0)$ acts faithfully on $\mathfrak{g}_1$,  we may
assume that $\mathfrak{g}_0 \subseteq \mathfrak{gl}(V)$. When
$\g_0^{(1)}$ is of type A$_m$, then $\mathfrak{sl}(V) \subseteq
\mathfrak{g}_0\subseteq \mathfrak{gl}(V)$, and when $\g_0^{(1)}$
is of type C$_m$, then $\mathfrak{sp}(V) \subseteq \mathfrak{g}_0
\subseteq \mathfrak{csp}(V)$. Let $G_0,\,$ $T\,$, and $N^{\pm}$ be
the algebraic groups introduced at the beginning of the proof of
Theorem \ref{Thm:4.22}. To simplify the notation, we will identify the
weights $\nu\in X_1(T)$ with their differentials
$(\mathrm{d}\nu)_e\in \big(\tf\cap\g_0^{(1)}\big)^*\cong
\,X(T)\otimes_{\Z}\F$. Since $G_0$ is a simply connected group and
$\g_{-1}$ is an irreducible restricted $\g_0^{(1)}$-module, there
is a unique $\lambda\in X_1(T)$ such that $\g_{-1}\cong
L(\lambda)$ as $\g_0^{(1)}$-modules, where the action of
$\g_{0}^{(1)}$ on $L(\lambda)$ is induced by the differential of
the rational action of $G_0$ on $L(\lambda)$; see Proposition
\ref{Pro:2.803}. The image of $\lambda$ in $X(T)/p\,
X(T)\hookrightarrow \big(\tf\cap \g_0^{(1)}\big)^*$ equals
$\bar{\Lambda}$.

{F}rom now on we will identify $\g_{-1}$ with $L(\lambda)$. The
group $G_0$ acts on $\g_0$ by conjugation, and hence it acts
rationally on the vector space $\text{Hom}(\g_{-1},\g_0)\cong
\g_{-1}^*\otimes \g_0$. By the choice of $L(\lambda)$, the
differential of this action coincides with the natural action of
$\g_0^{(1)}$ on  $\text{Hom}(\g_{-1},\g_0)$.

First we suppose that $\g_0^{(1)}$ is of type A$_m$ for $m \geq 2$.
Renumbering the simple roots if necessary,   we may assume that
$\bar{\Gamma} = -\varpi_m$ (see the proof of Theorem
\ref{Thm:4.22} for a similar argument). Then
$\g_{-1}=L(\varpi_1+2\varpi_m)$. Since $w_0 \varpi_i =
-\varpi_{m+1-i}$, $1 \leq i \leq m$, we have  that $V\cong
L(\varpi_1)$, and $\g_{-1}^* \cong L(2\varpi_1 + \varpi_m)$ as
$G_0$-modules; see Proposition \ref{Pro:2.802}. From this it
follows that  $\g_{-1}^*$ is isomorphic to a composition
factor of the $G_0$-module $S^2(V)\ot V^*$.  Because  the
$G_0$-module $\g_0$ is isomorphic to a $G_0$-submodule of $V\ot
V^*$, any weight of $\text{Hom}(\g_{-1},\g_0)$ relative $T$
belongs to the set
$$R \eqdef \{ \e_{i_1} + \e_{i_2} + \e_{i_3} - \e_{j_1}
- \e_{j_2} \mid 1 \leq i_1, i_2, i_3, \,j_1, j_2 \leq m+1\}.$$

\noindent Since $-\e_{m+1}=\varpi_m$ and
$-\e_{m}-\e_{m+1}=\varpi_{m-1}$ as rational characters of $T$ and
since the Weyl group $W$ of $G_0$ permutes the $\e_i$'s, we have
that

\begin{eqnarray*} R \cap X(T)_+ &=&
\big \{3 \varpi_1 + 2 \varpi_m, \, 3\varpi_1+\varpi_{m-1},\,2
\varpi_1 + \varpi_m,\\
&& \hskip 0.17 truein
 \varpi_1 + \varpi_2 + 2 \varpi_m, \, \varpi_3 + 2 \varpi_m, \,
 \varpi_1 + \varpi_2 + \varpi_{m-1},\\
 && \hskip 0.17 truein \varpi_3+\varpi_{m-1},\,\varpi_2+\varpi_m, \, \varpi_1\big
\},
\end{eqnarray*}

\m

\noindent where $\varpi_3+\varpi_{m-1}$ is omitted if $m\le 3$,
and $\varpi_1+\varpi_2+\varpi_{m-1},\,$ $\varpi_3+2\varpi_m$ are
omitted if $m=2$.   Clearly, $R \cap X(T)_+ \subset X_1(T)$. But then
$\g_1$ is a $G_0$-stable subspace of $\text{Hom}(\g_{-1}, \g_0)$
by Proposition \ref{Pro:2.805}.

Let $\mu$ be a minimal weight of the $G_0$-module
$\mathfrak{g}_1$, and let $\mathfrak{g}_1^\mu$ denote the
corresponding weight space. Since $\mathfrak{g}_1^\mu$ is
contained in the subspace of fixed points of the subgroup $N^-$,
we have $[\mathfrak{n}^-,\mathfrak{g}_1^\mu] = 0$. As $T$ acts
as the identity  on $\tf$, the weight space $\g_1^\mu$ is $\tf$-stable.
Let $e_1=e^{\Theta_1}$ be any nonzero weight vector for $\tf$
contained in $\g_1^\mu$.    It follows that $e_1$ is a $\mathfrak b^-$-primitive vector, and its weight $\bar{\Theta}_1$ coincides with the
image of $\mu$ in $X(T)/p\,X(T)\hookrightarrow
\big(\tf\cap\g_0^{(1)}\big)^*$. If $\Theta_1 \neq \Gamma$, then
$\Theta_1 =-\Lambda $ by Theorem \ref{Thm:4.7} (as in part (c) of the
proof of Theorem \ref{Thm:4.22}).

Suppose $\Theta_1=-\Lambda$. Since $w_0\mu$ is a maximal weight of
$\g_1$, it belongs to $R \,\cap X(T)_+$. On the other hand, the
image of $w_0\mu$ in $X(T)/p\,X(T)$ equals $w_0\bar{\Lambda}$.
As $p>3$, this implies that $w_0\mu=2\varpi_1+\varpi_m$. Let
$M$ denote the $G_0$-submodule of $\g_1$ generated by $e_1$. The
orbit $N_{G_0}(T)\cdot e_1$ contains an eigenvector for $B^+=TN^+$
of weight $w_0\mu$. Therefore, $M$ is a homomorphic image of the
Weyl module $V(2\varpi_1+\varpi_m)$; see Proposition
\ref{Pro:2.804}\,(a). Combining Theorem \ref{Thm:2.82} with Weyl's
dimension formula \eqref{eq:Weyldim},  it is now easy to observe that
$$\dim V(2\varpi_1+\varpi_m)=\,\dim \left(S^2(V)\otimes V^*\right) -\dim
V^*=\,\dim W(m+1)_1^\dagger,$$
where $W(m+1)_1^\dagger$ is as in \eqref{eq:2.82a}.   Since $m+2\not\equiv 0\mod p$ by
our earlier remarks, Theorem \ref{Thm:2.82} says that the
$\mathfrak{sl}_{m+1}$-module $W(m+1)_1^\dagger$ is irreducible and
contains a $\mathfrak b^+$-primitive vector of weight $2\varpi_1+\varpi_m$.
But then Proposition \ref{Pro:2.803} yields $\dim
L(2\varpi_1+\varpi_m)=\dim W(m+1)_1^\dagger$. It follows that
$V(2\varpi_1+\varpi_m)\cong L(2\varpi_1+\varpi_m)\cong M$. As a
result, $M$ is an irreducible $\g_0$-module. Hence, we can replace
the triple $(f^\Lambda,\g_0,e)$ by the triple
$(f^\Lambda,\g_0,e_1)$ and argue as before  to conclude that
$\Theta_1\ne -\Lambda$. As a consequence, $\g_1^\mu\subseteq
\g_1^{\Gamma}$;  that is,   the image of $\mu$ in $X(T)/p\,X(T)$
coincides with $\Gamma$.

{F}rom  the description of $R$ above and the fact that $p>3$, 
we see that $\mu = -\varpi_m$.  There is a
linear function $\delta \in (\mathfrak{g}_1^\mu)^*$ such that
$[f^\Lambda, x] = \delta(x)e_{\alpha}$ for any $x \in
\mathfrak{g}_1^\mu$. Since $\g$ is irreducible and transitive, and
$\g_1^\mu$ consists of $\mathfrak b^-$-primitive vectors, we have
$[f^\Lambda,x] \neq 0$ whenever $0 \neq x \in \mathfrak{g}_1^\mu$.
This yields $\ker \delta=0$, i.e., $\dim \mathfrak{g}_1^\mu = 1$.
If $\nu$ is a weight of the $G_0$-module $\mathfrak{g}_{-1}^* \ot
\mathfrak{g}_0$, then $w(\nu) \in \Delta \cap X(T)_+$ for some $w
\in W$. Using the list of weights in $R \cap X(T)_+$ given above,
it is easy to see that $w(\nu) \geq \varpi_1$. But then $w_0
w(\nu) \leq w_0 (\varpi_1) = \mu$.  By the minimality of $\mu$ we
have $\mu = w_0 w (\nu)$.  Hence, all weights of the $G_0$-module
$\mathfrak{g}_1$ are conjugate under $W$. As $\dim
\mathfrak{g}_1^\mu = \dim \mathfrak{g}_1^{w \mu}$ for any $w \in
W$, each weight occurs with multiplicity one. Now it is
straightforward to see that $\mathfrak{g}_1=V$ is an irreducible
standard $\g_0^{(1)}$-module.

Suppose then that $\g_0^{(1)}$ is of type C$_m$ where $m \geq 1,$
so that $G_0=\mathrm{Sp}(V)$ where $\dim V=2m.$
 Then  $\Gamma =
-\varpi_1$, $\Lambda = 3 \varpi_1$ and $\mathfrak{g}_0 \subseteq
\g_0^{(1)} \oplus {\mathbb F}$, where ${\mathbb F}$ is a trivial
$G_0$-module and $\g_0^{(1)} \cong S^2(V)$ as $\g_0^{(1)}$-modules
(see \eqref{eq:spiso}).      As $\mathfrak{g}_{-1}$ is
irreducible, the $\g_0^{(1)}$-modules $\mathfrak{g}_{-1}$ and
$S^3(V)$ are isomorphic (see Proposition \ref {Pro:2.848} or  Theorem \ref{Thm:2.24}).
Identifying
$\mathfrak{g}_{-1}$ with $S^3(V)$ and taking into account the fact
that $S^2(V) \cong S^2(V)^*$, we can embed $\mathfrak{g}_1$ in the
$\g_0^{(1)}$-module
$${\mathcal N} := S^3(V)^* + \left (S^3(V)^* \ot
S^2(V)^* \right),$$ since
\begin{eqnarray*}
&&\mathfrak{g}_1  \subseteq   \Hom(\mathfrak{g}_{-1},\mathfrak{g}_0)
\cong  \mathfrak{g}_{-1}^* \ot \mathfrak{g}_0
\subseteq  \mathfrak{g}_{-1}^* \ot (\g_0^{(1)} \oplus
{\mathbb F}) \quad\,\, \hbox{\rm where}  \\
&&\hspace*{.5cm}\mathfrak{g}_{-1}^* \ot (\g_0^{(1)} \oplus
{\mathbb F}) \cong   S^3(V)^* \ot (S^2(V)^* \oplus {\mathbb F})  \\
&&\hspace*{3.3cm}\cong  S^3(V)^* + \left (S^3(V)^* \ot S^2(V)^* \right) = \mathcal N.
\end{eqnarray*}

The group \, $G_0$ \, acts by automorphisms on the symmetric algebra
$S(V) = \bigoplus_{j \geq 0}\, S^j(V)$, and the subspace
${\mathcal N} \subset S(V)^*$ is invariant under the dual action
of $G_0$ on $S(V)^*$. The action of $\g_0^{(1)}$ is the
differential of this action of $G_0$. Now the multiplication in
$S(V)$ induces a surjective $G_0$-module homomorphism $S^2(V) \ot
S^3(V) \rightarrow S^5(V)$. Hence, there is an embedding of
$G_0$-modules

$$\eta: S^5(V)^* \rightarrow S^2(V)^* \ot S^3(V)^*.\nonumber$$

\noindent Let
$${\mathcal N}_0 := \eta(S^5(V)^*)$$ be the image and
define

$$R' = \{\pm \e_{i_1} \pm \e_{i_2}\pm \e_{i_3}\pm \e_{i_4}\pm \e_{i_5}
\mid 1 \leq i_1,i_2,i_3,i_4,i_5 \leq m \}.\nonumber$$

\smallskip

\noindent Since $X(V)=\{\pm \e_1,\ldots, \pm \e_m\}$, it is
straightforward to see that $X({\mathcal N}) \subseteq R'$ and
$\dim {\mathcal N}^{\,5\varpi_1} = \dim {\mathcal
N}_0^{\,5\varpi_1}=1$. From this it follows that
\begin{eqnarray*} X_+({\mathcal N}/{\mathcal N}_0) &\subseteq& \big
\{3 \varpi_1 + \varpi_2, \, 3\varpi_1,\,2\varpi_1+\varpi_3,\\
&& \hskip 0.32 truein \varpi_1+\varpi_4,\, \varpi_1 + \varpi_2, \,
\varpi_5,\, \varpi_3,\,\varpi_1\big \},
\end{eqnarray*}

\smallskip

\noindent where $\varpi_5$ is omitted if $m\le 4,\,$
$\varpi_1+\varpi_4$ is omitted if $m\le 3,\,$ $2\varpi_1+\varpi_3$
and $\varpi_3$ are omitted if $m=1,2$, and $3\varpi_1+\varpi_2,\,$
$\varpi_1+\varpi_2$ are omitted if $m=1$. As a consequence,
$X_+(\mathcal{N}/\mathcal{N}_0) \subset X_1(T)$, and any weight in
$X_+(\mathcal{N}/\mathcal{N}_0)$ dominates $\varpi_1$.

We claim that $\mathfrak{g}_1 \cap \mathcal{N}_0 = 0$. Indeed,
suppose the contrary. Then Theorem \ref{Thm:4.7} shows that the
$\g_0^{(1)}$-module $S^5(V)^*\cong\mathcal{N}_0$ contains a
$\mathfrak b^-$-primitive vector $u$ of weight
$\bar{\Gamma}=-\varpi_1$ or $-\bar{\Lambda}=-3\varpi_1$. If $p
> 5$, then $S^5(V)$ is an irreducible $\g_0^{(1)}$-module (see
Lemma \ref{Lem:2.87}). Therefore, any $\mathfrak b^-$-primitive vector
from $S^5(V)^*$ has weight $-5\varpi_1
\not\in\{-\varpi_1,-3\varpi_1\}$, which is impossible.    

Now assume $p = 5$.    By Lemma  \ref{Lem:2.87},
 there exists a trivial $\g_0^{(1)}$-submodule
$Y$ of $S^5(V)$,  and the
quotient $Y': = S^5(V)/Y$ has  a $\mathfrak b^+$-primitive vector of weight 
$3 \varpi_1 + \varpi_2$ if $m\geq 2$
($3\varpi_1$ if $m=1$).     Set $Y^\perp := \{\phi \in S^5(V)^* \mid
\phi(Y) = 0\}$.   The $\g_0^{(1)}$-modules $Y^\perp$ and $(Y')^*$
(resp. $S^5(V)^*/Y^\perp$ and $Y^*$) are isomorphic.   If $m \geq 2$,
then any $\mathfrak b^-$-primitive vector from $(Y')^*$ has weight
$-3\varpi_1 -\varpi_2 \not \in \{-\varpi_1, -3\varpi_1\}$.
But then  $u \not \in
Y^\perp$, and so $S^5(V)^*/Y^\perp \cong Y^*$  contains a nonzero vector of weight $\bar \Gamma =
-\varpi_1$ or $-\bar \Lambda =-3\varpi_1$.   As the
$\g_0^{(1)}$-module $Y^*$ is trivial, this is impossible.       If $m=1$, this argument shows that $u$ has
$\tf$-weight $-\bar \Lambda = -3\varpi_1$; \
$\,Y^\perp \cong \mathfrak{U}(\mathfrak{n}^+)\,u$; \  and $\eta( \mathfrak{U}(\mathfrak{n}^+)\,u)$
is an irreducible submodule of the $\g_0$-module $\g_1$ (one
should also take into account that
$\g_0=\g_0^{(1)}\oplus{\mathfrak Z}(\g_0)$).   But then one can
replace the triple $(f^\Lambda,\g_0,e)$ by the triple
$(f^\Lambda,\g_0,\eta(u))$ and argue as before, to conclude that
that $\eta(u)$ has weight $\Gamma$. This is a contradiction, and
it proves the claim.

Thus, the canonical epimorphism $\mathcal{N} \rightarrow
\mathcal{N}/\mathcal{N}_0$ induces an injection $
\mathfrak{g}_{-1} \rightarrow \mathcal{N}/\mathcal{N}_0$. As
$X_+(\mathcal{N}/\mathcal{N}_0) \subset X_1(T)$, Proposition
\ref{Pro:2.805} shows that the adjoint action of $\g_0^{(1)}$ on
$\g_{-1}$ is induced by the differential of a rational action of
$G_0$. The weights of the $G_0$-module $\g_1$ belong to
$X(\mathcal{N}/\mathcal{N}_0)$.

Again let  $\mu$ be a minimal weight
 of the $G_0$-module $\mathfrak{g}_1$. Arguing as
 before we see that the weight space $\mathfrak{g}_1^\mu$ is $\tf$-stable and
consists of $\mathfrak b^-$-primitive vectors relative to $\g_0$.
 Let $e_2=e^{\Theta_2}$
 be a nonzero weight vector for $\tf$ contained in $\g_1^\mu$.
By Theorem \ref{Thm:4.7} (again, as in part (c) of the proof of Theorem \ref{Thm:4.22}), the image of $\mu$ in $X(T)/p\,X(T)$
coincides with $-3\varpi_1$ or $-\varpi_1$. Since $p>3$, the
description of $X_+(\mathcal{N}/\mathcal{N}_0)$ given above yields
that either $w_0\mu=3\varpi_1$ or $w_0\mu=\varpi_1$. Let $M'$
denote the $G_0$-module generated by $e_2$. The orbit
$N_{G_0}(T)\cdot e_2$ contains an eigenvector for $B^+=TN^+$ of
weight $w_0\mu$. Thus, $M'$ is a homomorphic image of the Weyl
module $V(w_0\mu)$; see Proposition \ref{Pro:2.804}\,(a).

Assume that $\mu=-3\varpi_1$. Then $\Theta_2=-\Lambda$.  Now
 $S^3(V)$ is an irreducible
$\g_0^{(1)}$-module generated by a $\mathfrak b^+$-primitive vector of weight
$3\varpi_1$.   (This can be seen directly or by identifying $\g_0^{(1)}$
with $H(2m)_0$ and $S^3(V)$
with $\mathcal O(2m)_3$, which has a $\mathfrak b^+$-primitive vector $x_1^{(3)}$, and then by appealing to
\cite[Lem.~5.2.2]{St} or Proposition \ref{Pro:2.848}\,(a)).
By Weyl's dimension formula \eqref{eq:Weyldim}, \ $\dim V(3\varpi_1)=\dim
S^3(V)$.    In conjunction with Proposition \ref{Pro:2.803} this
implies that $S^3(V)\cong L(3\varpi_1)$ as $\g_0^{(1)}$-modules.
But then $V(3\varpi_1)\cong L(3\varpi_1)\cong M'$ as
$G_0$-modules; see Proposition~\ref{Pro:2.804}(a). Hence, $M'$ is
an irreducible $\g_0$-module, and we can, as before, replace the
triple $(f^\Lambda,\g_0,e)$ by the triple $(f^\Lambda,\g_0,e_2)$
to deduce that $\Theta_2\ne -\Lambda$. This contradiction shows
that $-\varpi_1$ is the only minimal weight of the $G_0$-module
$\g_1$.

Since $\g$ is irreducible and transitive, it must be that $\dim
\g_1^\mu=1$. We have already established that $\nu \geq \varpi_1$
for any $\nu\in X_+(\g_1)$. From this it is immediate that all
weights of the $G_0$-module $\mathfrak{g}_1$ are conjugate under
the Weyl group of $G_0$.  Using these properties and reasoning as
above, we see that $\mathfrak{g}_1 = V$.      We have considered all
the cases, so the lemma is proved.  \qed \m

\section  {\ Determining the negative part  when 
$\boldsymbol{\g_1}$ is irreducible \label{sec:4.8} } 
\m

In Lemma \ref{Lem:4.33}  below,  we show that if the Lie algebra
$\widehat \g$ generated by the local part of $\g$ is isomorphic to
one of the graded Lie algebras of Cartan type $S$, $CS$, $H$, or
$CH$ (so that $\g_1$ is an irreducible $\g_0$-module), then $\g$
itself must be depth-one. Except for a couple of special cases
which must be treated separately, we are able to verify the lemma
by applying previous results to the Lie algebra generated by
$\g_0$,  a homogeneous space  $\widehat \g_k$ of $\widehat \g$ for
$k \geq 1$, and a certain $\g_0$-submodule $V_{-k}$ of
$\g_{-k}$.\bi

\begin{Lem} \label{Lem:4.33}  \  Let $\mathfrak{g}$ be a
finite-dimensional graded Lie algebra satisfying conditions
(1)-(5). Suppose that the subalgebra $\widehat{\mathfrak{g}}$
generated by the local part $\mathfrak{g}_{-1} \oplus
\mathfrak{g}_0 \oplus \mathfrak{g}_1$ is isomorphic one of the
restricted Lie algebras of Cartan type
$S(m+1;\underline{1})^{(1)}, m\ge 2,\,$ $S(m+1;
\underline{1})^{(1)}\oplus\F{\mathfrak D}_1,\, m\ge 2,\,$
$H(2m;\underline{1})^{(2)}, m\ge 1,$ or $
H(2m;\underline{1})^{(2)}\oplus\F{\mathfrak D}_1, m\ge 1,$ with
its natural grading. Then $\mathfrak{g}_{-2} = 0$.
\end{Lem}

\pf (a) Set $\widehat{\mathfrak{g}}_j = \widehat{\mathfrak{g}}
\cap \mathfrak{g}_j$ for $ -q \leq j \leq r$. By our hypothesis on
$\widehat{\mathfrak{g}},$ the commutator $\g_0^{(1)}$ is a
classical simple Lie algebra of type A$_m$ or C$_m$.  Let
$\mathfrak t$ denote the toral subalgebra of $\g_0$, and let
$f^\Lambda$ be a $\mathfrak b^+$-primitive vector of weight $\Lambda$ from
$\mathfrak{g}_{-1}$.    In the A$_m$-case, we assume as before
that $\Lambda = \varpi_m$.   In this case, we have by Lemma
\ref{Lem:2.85}\,(a)   that when $m \geq 2$, the $\g_0^{(1)}$-module
$\widehat{\mathfrak{g}}_k$ is generated by a $\mathfrak b^-$-primitive vector
$e^{-\varpi_1 - (k+1)\varpi_m}$, $k = 2,3$.   If $\Phi$ is of
type C$_m$, $m \geq 1$, and $p>5$, the same lemma shows that
$\widehat{\mathfrak{g}}_k$ is generated by a $\mathfrak b^-$-primitive vector
$e^{-(k+2)\varpi_1}$, $k = 2,3$.  If $p = 5$ and $\Phi$ is of
type C$_m$, $m \geq 2$, then $\widehat{\mathfrak{g}}_2$ is
generated by a $\mathfrak b^-$-primitive vector $e^{-4\varpi_1}$ and
$\widehat{\mathfrak{g}}_3$ is generated by a $\mathfrak b^-$-primitive vector
$e^{-3\varpi_1 - \varpi_2}$; see Lemma \ref{Lem:2.85}\,(b).   Finally, in the A$_1$, $p = 5$ case,
$\widehat{\mathfrak{g}}_2$ is generated by a $\mathfrak b^-$-primitive vector
$e^{-4\varpi_1}$,  and $\widehat{\mathfrak{g}}_3$ is generated by
a $\mathfrak b^-$-primitive vector $e^{-3\varpi_1}$, again by Lemma \ref{Lem:2.85}\,(b).
\m

\noindent (b) By $1$-transitivity \eqref{eq:1.4}, the
$\mathfrak{g}_0$-submodule $\mathfrak{g}_{-i}$ is isomorphic to a
submodule of $\Hom(\mathfrak{g}_1,\mathfrak{g}_{-(i-1)})$ ($i
> 0$). This can be used to show that each $\mathfrak{g}_{-i}$ with $i>0$ is a
restricted $\g_{0}^{(1)}$-module.  Indeed for $i=1$, this is 
Remark \ref{Rem:4.10}.   By transitivity \eqref{eq:1.3}, $\mathfrak{g}_{1}$
is isomorphic to a submodule of $\Hom(\mathfrak{g}_{-1},
\mathfrak{g}_0);$ hence it is a restricted $\g_0^{(1)}$-module.
(Note that $\mathfrak{g}_0$ is a restricted $\g_0^{(1)}$-module
too, as the restriction map on $\g_0^{(1)} $ is induced by that on
$\mathfrak{g}_0.)$    By our induction assumption, the
$\g_0^{(1)}$-module $\mathfrak{g}_{-(i-1)}$ is restricted. But
then so is the $\g_0^{(1)}$-module
$\Hom(\mathfrak{g}_1,\mathfrak{g}_{-(i-1)})$ along with its
submodule $\g_{-i}$. This completes the induction step. \m

\noindent (c) Suppose that $\Phi \neq$ A$_1$ if $p = 5$. We will
show that $[\mathfrak{g}_{-k},\widehat{\mathfrak{g}}_k] = 0$ for
$k = 2,3$. Suppose the contrary, and let $\mathfrak{a}_{-k} := \{x
\in \mathfrak{g}_{-k} \mid [x,\widehat{\mathfrak{g}}_k] = 0\}$ for
$k = 2,3.$ Let $V_{-k}$ denote a $\mathfrak{g}_0$-submodule of
$\mathfrak{g}_{-k}$ containing $\mathfrak{a}_{-k}$ such that
$U_{-k} \eqdef V_{-k}/\mathfrak{a}_{-k}$ is irreducible.  Both $V_{-k}$
and ${\mathfrak a}_{-k}$ are $\g_0^{(1)}$-submodules of
$\mathfrak{g}_{-k},$ and hence are restricted over $\g_0^{(1)}$,
by (b). Therefore, so is the quotient module $U_{-k}:=
V_{-k}/{\mathfrak a}_{-k}.$\m

\noindent (c1) By Engel's theorem, the subalgebra ${\mathfrak
n}^{+}$ of $\g_0^{(1)}$ annihilates a nonzero vector on any
$\mathfrak{g}_0$-module $E$ endowed with a $p$-character which
vanishes on $\g_0^{(1)}$. This means that $E_0=\{v\in E\mid
e_\alpha.v = 0\ \hbox{\rm for all }\ \alpha\in \Phi^+\}$ is a
nonzero $\mathfrak t$-module. But $\mathfrak t$ is diagonalizable
for any $\g_0$-module affording a $p$-character. So there is
$\lambda \in \mathfrak t^*$ such that $E_{0}^\lambda = \{v\in
E_0\mid t.v= \lambda(t)v \hbox{ for all } t\in\mathfrak t\}$ is
nonzero. \m

\noindent (c2) Let $S^{\{k\}} = \bigoplus_{i \in k\Z} S_i^{\{k\}}
$ be the graded subalgebra of $\mathfrak{g}$ generated by
$V_{-k}$, $\mathfrak{g}_0$, and $\widehat{\mathfrak{g}}_k$, and
let $\mathcal M^{\{k\}} $ be the Weisfeiler radical $\mathcal M(S^{\{k\}})$
of $S^{\{k\}}$. Obviously, $\mathcal M^{\{k\}} \cap S_{-k}^{\{k\}}  =
\mathfrak{a}_{-k}$, and $S_k^{\{k\}}  = \mathfrak{g}_k$. Let
$\overline{S}^{\{k\}}  := S^{\{k\}}/\mathcal M^{\{k\}} .$ Then
$\overline{S}_{-k}^{\{k\}}  \cong U_{-k}$ as $\g_0$-modules, and
$\overline{S}_k^{\{k\}} \cong \mathfrak{g}_k$. Part~(c1) shows
that we can find a $\mathfrak b^+$-primitive vector in $U_{-k}$, say
$f^{\mu_k}$. Because $[\overline{S}_{-k}^{\{k\}},
\overline{S}_k^{\{k\}}] \neq 0$, we must have
$[f^{\mu_k},e^{\lambda_k}] \neq 0$ where $e^{\lambda_k}$ is a
$\mathfrak b^-$-primitive vector of $\mathfrak{g}_k$ generating
$\mathfrak{g}_k$ as a $\g_0^{(1)}$-module. \m

\noindent (c3) If $\mu_k \neq -\lambda_k$ and $\mu_k(\mathfrak t
\cap \g_0^{(1)}) \neq 0,$ then the graded Lie algebra
$\overline{S}^{\{k\}}$ satisfies all the conditions of Theorem
\ref{Thm:4.7}.  But then $-\lambda_k \in \{ \varpi_1, \varpi_m,
2\varpi_1 + \varpi_m, \varpi_1 + 2\varpi_m \}$ if $\Phi \cong
$A$_m, \, m \geq 2,$ and $-\lambda_k \in \{\varpi_1, 3\varpi_1 \}$
if $\Phi \cong$ C$_m, \, m \geq 1.$  Using our description of the
$\lambda_k$'s in part~(a), it is easy to check that this is not
the case. Hence, either $\mu_k(\mathfrak t \cap \g_0^{(1)}) = 0,$
or $\mu_k = -\lambda_k.$

If $\mu_k(\mathfrak t \cap \g_0^{(1)}) = 0,$ then it follows from
the classification of the irreducible restricted representations
of a classical Lie algebra that $U_{-k}$ is
a trivial $\g_0^{(1)}$-module. But then $\ad f^{\mu_k}$ induces a
nonzero $\g_0^{(1)}$-module homomorphism $\phi_k: \mathfrak{g}_k
\rightarrow \mathfrak{g}_0.$ Since $\mathfrak{g}_k$ is generated
by $e^{\lambda_k}$ as a $\g_0^{(1)}$-module, it must be that
$\phi_k(e^{\lambda_k})$ is a $\mathfrak b^-$-primitive vector in
$\mathfrak{g}_0.$  But $-\lambda_k$ is neither zero nor  the
highest root of $\Phi$; see part~(a) and our assumption in
part~(c). Hence, $\mu_k(\mathfrak t \cap \g_0^{(1)}) \neq 0.$ \m

\noindent (c4) Now suppose $\mu_k = -\lambda_k$. We first assume
that $\g_k$ is an irreducible $\g_0^{(1)}$-module and let $
\overline{\mathcal N}^{\{k\}}$ denote the maximal graded ideal of
$\overline{S}^{\{k\}}$ contained in $\bigoplus_{i\ge
1}\,\overline{S}_i^{\{k\}}$. Theorem \ref{Thm:3.34} shows that
$\overline{S}^{\{k\}}/\overline{\mathcal N}^{\{k\}}$ is a classical Lie
algebra with one of its standard gradings (note that
$\overline{S}^{\{k\}}/\overline{\mathcal N}^{\{k\}}$ cannot be a Melikyan
algebra in view of our assumption in part (c)). As in the proof of
Lemma \ref{Lem:4.30},  one then observes that the 
the $\g_0^{(1)}$-module $\overline{S}_{-k}^{\{k\}}\cong U_{-k}$  must
have a $\mathfrak b^+$-primitive vector whose weight  is
in the list
\begin{eqnarray}\label{eq:sasha0}
 & \{2 \varpi_1, \varpi_1,
\varpi_2, \varpi_{m-1},
\varpi_m, 2 \varpi_m\}, \qquad \quad \quad  \  m \geq 8, \nonumber \\
& \{2 \varpi_1, \varpi_1, \varpi_2, \varpi_3, \varpi_5,\varpi_6,
\varpi_7, 2 \varpi_7\}, \qquad \quad m = 7,  \\
 &\{\varpi_i \mid 1 \leq i \leq m \} \cup \{2 \varpi_1, 2
\varpi_m\}, \quad \quad 2\, \leq m \leq 6, \nonumber
\end{eqnarray} for $\g_0^{(1)}$ of type A$_m$, and in
the list
\begin{eqnarray}
& \{\varpi_1\}, \qquad \qquad \quad \qquad m \geq 4, \nonumber \\
& \{\varpi_1, \varpi_m\}, \qquad \qquad m = 2,3, \nonumber\\
& \{\varpi_1, 2\varpi_1, 3\varpi_1\}, \qquad\ \, m = 1, \nonumber
\end{eqnarray} for $\g_0^{(1)}$ of type C$_m$. On the other hand, it follows
from our discussion in part (a) that $\mu_k= \varpi_1 +
(k+1)\varpi_m,\,$ $k=2,3,\,$ if $\g_0^{(1)}$ is of type A$_m$,
$m\ge 2$, and $\mu_k=(k+2) \varpi_1,\,$ $k=2,3,$ for $p>5$ and
$\mu_k\in \{4\varpi_1,\,3\varpi_1+\varpi_2\}$ for $p=5$ and $m\ge
2,$ if $\g_0^{(1)}$ is of type C$_m,$ $m\ge 1$. These linear
functions do not appear on the list. Therefore, the case we are
considering cannot occur.  \m

\noindent (c5) Next assume that $\mu_k=-\lambda_k$ and $\g_k$ is a
reducible $\g_0^{(1)}$ module. Proposition \ref{Pro:2.848}\,(b)
together with Theorem \ref{Thm:2.82}\,(i) show that $m+k\equiv 0\mod
p$ and $\widehat{\g}$ is isomorphic to one of
$S(m+1;\underline{1})^{(1)}$ or $
S(m+1;\underline{1})^{(1)}\oplus\F{\mathfrak D}_1$. Theorem
\ref{Thm:2.82} then says that $\g_k\cong\overline{S}_k^{\{k\}}$
contains an irreducible $\g_0^{(1)}$-submodule $\g_k^\sharp$  
with a $\mathfrak b^+$-primitive vector of weight $k\varpi_1$ such that the quotient module
$\g_k/\g_k^\sharp$ is irreducible. Let $e^{-k\varpi_m}$ be a
$\mathfrak b^-$-primitive vector of the $\g_0^{(1)}$-module
$\g_k^\sharp\subset\overline{S}_k^{\{k\}}$. Note that
$[f^{\mu_k},e^{k\varpi_m}]$ is a multiple of a highest root vector
in $\overline{S}_0^{\{k\}}\cong\g_0$. Since $k\in\{2,3\}$, Theorem
\ref{Thm:4.7} implies $[f^{\mu_k},e^{k\varpi_m}]=0$. But then
$[f^{\mu_k},\g_k^\sharp]=0$, forcing
$[\overline{S}_k^{\{k\}},\g_k^\sharp]=0$.

As in part (c4),  we let $\overline{\mathcal N}^{\{k\}}$ denote the maximal
graded ideal of $\overline{S}^{\{k\}}$ contained in
$\bigoplus_{i\ge 1}\,\overline{S}_i^{\{k\}}$. The above discussion
then shows that $\overline{\mathcal N}_k^{\{k\}}=\g_k^\sharp$ and
$$\overline{S}_k^{\{k\}}/\overline{\mathcal N}_k^{\{k\}}\,\cong\,\g_k/\g_k^\sharp\,\cong\,
U_{-k}^*\,=\,\big(\overline{S}_{-k}^{\{k\}}\big)^*$$ as
$\g_0^{(1)}$-modules. But then Theorem \ref{Thm:3.34} applies to
the graded Lie algebra
$\overline{S}^{\{k\}}/\overline{\mathcal N}^{\{k\}}$. Repeating the
argument from part (c4) we now deduce that the present case cannot
occur. This shows that under the main assumption of part (c) we
have $[\g_k,\widehat{\g}_{k}]=0$ for $k=2,3$. \m

\noindent (c6) As $\widehat{\mathfrak{g}}_2 \neq 0,$ it follows
that $[\mathfrak{g}_{-1}, \, [\mathfrak{g}_{-1}, \,
\widehat{\mathfrak{g}}_2] \neq 0$ by transitivity. As
$[\mathfrak{g}_{-2}, \, \widehat{\mathfrak{g}}_2] = 0,$ it follows
that $[\mathfrak{g}_{-2}, \, [\mathfrak{g}_{-1}, \,
[\mathfrak{g}_1, \, \widehat{\mathfrak{g}}_2]]] = 0.$ By
assumption, $\mathfrak{g}_{-2} \neq 0,$ so that
$[\mathfrak{g}_{-2}, \, \mathfrak{g}_1] = \mathfrak{g}_{-1},$ by
irreducibility  and  1-transitivity. But then
\begin{equation*}\begin{array}{ccc}
[\mathfrak{g}_{-3}, \, \widehat{\mathfrak{g}}_3]&\supseteq &
[[\mathfrak{g}_{-1}, \, \mathfrak{g}_{-2}], \, [\mathfrak{g}_1, \,
\widehat{\mathfrak{g}}_2]]
 \ =\ [\mathfrak{g}_{-1}, \,
[\mathfrak{g}_{-2}, \, [\mathfrak{g}_1, \,
\widehat{\mathfrak{g}}_2]]]\\
& =&\ [\mathfrak{g}_{-1}, \, [[\mathfrak{g}_{-2}, \,
\mathfrak{g}_1], \, \widehat{\mathfrak{g}}_2]] \ = \
[\mathfrak{g}_{-1}, \, [\mathfrak{g}_{-1}, \,
\widehat{\mathfrak{g}}_2]]\, \neq \, 0.
\end{array}
\end{equation*}
This contradiction shows that $\mathfrak{g}_{-2} = 0$ if $p > 5$,
or if $p = 5$ and $\Phi$ is not of type $\mathrm{A}_1.$ \m

\noindent (d) Now suppose that $p = 5$ and $\Phi$ is of type
$\mathrm{A}_1$.  By our initial assumption, the graded Lie algebra
$\widehat{\mathfrak{g}}$ is then isomorphic to either
$H(2;\un{1})^{(2)}$ or $H(2;\un{1})^{(2)}\oplus\,\F{\mathfrak
D}_1$. By Proposition \ref{Pro:2.848}(a) and Lemma
\ref{Lem:2.85}(b), the $\g_0^{(1)}$-module $\g_2$ is irreducible
and generated by a $\mathfrak b^-$-primitive vector $e^{-4\varpi_1}$. If we
reason as above, it is not difficult to observe that
$[\mathfrak{g}_{-2}, \, \widehat{\mathfrak{g}}_2] = 0.$  This
forces $[\mathfrak{g}_{-3}, \, \widehat{\mathfrak{g}}_3] \neq 0;$
see the above calculation. Since $\widehat{\g}$ is generated by
its local part, the commutator subalgebra
$\widehat{\mathfrak{g}}^{(1)}$ is simple. Moreover,
$\widehat{\mathfrak{g}}_i \subset \widehat{\mathfrak{g}}^{(1)}$
for any $i \neq 0$ and $\widehat{\mathfrak{g}}^{(1)} \cap
\mathfrak{g}_{2(p-1)-2} = 0$; see Lemma \ref{Lem:2.88}.  As $p =
5$, this gives $\widehat{\mathfrak{g}}_6 = 0.$  By Lemma
\ref{Lem:2.88} again, $\widehat{\mathfrak{g}}_3$ is an irreducible
$\g_0^{(1)}$-module generated by a $\mathfrak b^-$-primitive vector
$e^{-3\varpi_1}$. Consequently, $[\mathfrak{g}_{-3}, \,
e^{-3\varpi_1}] \neq 0.$

As before, let $V_{-3}$ be a $\mathfrak{g}_0$-submodule of
$\mathfrak{g}_{-3}$ containing $\mathfrak{a}_{-3} \eqdef \{x \in
\mathfrak{g}_{-3} \mid [x, \, \widehat{\mathfrak{g}}_3] = 0 \}$
and such that the quotient module $U_{-3} =
V_{-3}/\mathfrak{a}_{-3}$ is irreducible.  Let $S^{\{3\}}$ be the
graded subalgebra generated by $V_{-3},$ $\mathfrak{g}_0$,  and
$\widehat{\mathfrak{g}}_3.$ Let $\mathcal M^{\{3\}}$ denote the maximal
ideal of $S^{\{3\}}$ contained in the negative part
$\bigoplus_{i<0}S_i^{\{3\}}$.  Set $\overline {S}^{\{3\}} \eqdef
S^{\{3\}}/\mathcal M^{\{3\}}$.    Clearly, $\overline{S}^{\{3\}}$ is
graded:
$$\overline{S}^{\{3\}} = \bigoplus_{i \in 3{\mathbb Z}}(\overline{S}^{\{3\}})_i\nonumber$$
and $\overline{S}^{\{3\}}_{-3} \cong U_{-3},$
$\overline{S}_0^{\{3\}}\cong \mathfrak{g}_0,$ and $S_3^{\{3\}}
\cong \overline{\mathfrak{g}}_3.$  By similar reasoning, we find a
$\mathfrak b^+$-primitive vector $f^{\mu_3} \in U_{-3}$ and show that
$[f^{\mu_3}, \overline{e^{-3\varpi_1}}] \neq 0;$ here
$\overline{e^{-3\varpi_1}}$ denotes the image of $e^{-3\varpi_1}$
in $\overline{S}^{\{3\}}.$

If $\mu_3(\mathfrak t \cap \g_0^{(1)}) = 0,$ then, as above,
$U_{-3} = {\mathbb F}f^{\mu_3},$ and $\ad f^{\mu_3}$ induces a
nonzero $\g_0^{(1)}$-module homomorphism from
$\widehat{\mathfrak{g}}_3$ into $\mathfrak{g}_0.$   Since this
situation is easily seen to be impossible, we have
$\mu_3(\mathfrak t \cap \g_0^{(1)}) \neq 0.$  If $\mu_3 \neq
3\varpi_1,$ Theorem \ref{Thm:4.7} applies to yield that
$[f^{\mu_3}, \, \overline{e^{-3\varpi_1}}] = f_1,$ a $\mathfrak b^-$-primitive
vector of $\g_0^{(1)}\cong \mathfrak{sl}_2$.  Hence,
$\mu_3= \varpi_1.$ Let $f^{3\varpi_1}$ be a $\mathfrak b^+$-primitive vector of $\widehat{\mathfrak{g}}_3,$ and let
$\overline{f^{3\varpi_1}}$ denote its image in
$\overline{S}^{\{3\}}.$  Because $[f^{\mu_3}, \,
\overline{f^{3\varpi_1}}]$ $\in \mathfrak{g}_0^{4\varpi_1}$ we
have $[f^{\mu_3}, \, \overline{f^{3\varpi_1}}]$ $= 0.$  Therefore,
$$[f^{\mu_3}, \, [\overline{f^{3\varpi_1}}, \, \overline{e^{-3\varpi_1}}]]\, =
- \,[f_1, \, \overline{f^{3\varpi_1}}] \neq 0,$$ forcing
$(\overline{S}^{\{3\}})_{6}\neq 0.$  Since this contradicts
$\widehat{\mathfrak{g}}_6 = 0$,   we conclude that $\mu_3 =
3\varpi_1.$

    Let $\mathcal J$ be the maximal graded ideal of $\overline {S}^{\{3\}}$
contained in the subspace $\bigoplus_{\mid i\mid
>3}\overline{S}_i^{\{3\}}$.   By construction, the graded Lie algebra
$\overline {S}^{\{3\}}/{\mathcal J}$ satisfies all the conditions of Theorem
\ref{Thm:3.34}. Since $\overline {S}^{\{3\}}/{\mathcal J}$ is generated by
its local part, it cannot be isomorphic to a Melykian algebra.
Hence, $\overline {S}^{\{3\}}/{\mathcal J}$ is isomorphic to a graded Lie
algebra of type G$_2$ associated with the Dynkin diagram
\begin{equation}\label{eq:4.34}{\beginpicture \setcoordinatesystem
units <0.45cm,0.3cm> % sets scale
 \setplotarea x from -6 to 7, y from -1 to 1    % sets plot size up
 \linethickness=0.15pt                          % sets line thickness
  \put{$\circ$} at 0 0
 \put{$\bullet$} at 2 0
 \plot .15 -.3 1.85 -.3 /
  \plot .15 0 1.85  0 /
  \plot .15 .3 1.85 .3 /
  \put {$>$} at 1 0
   \put{$\alpha_1$} at 0 -1.5
  \put{$-\mu_3$} at 2 -1.5 \endpicture}\end{equation}
(see Theorem \ref{Thm:3.34} and \cite{Bou1}).  But then
$(\overline{S}^{\{3\}})_6 \neq \left ((\overline{S}^{\{3\}})_6
\cap {\mathcal J} \right)$  to force $\widehat{\mathfrak{g}}_6 \neq 0.$  This
contradiction shows that $\mathfrak{g}_{-2}$ must be zero, as
desired. \qed \m

\section {\ Determining the negative part 
 when $\boldsymbol{\g_1}$ is reducible \label{sec:4.9}}  
\m

Here we show that if the local Lie algebra $\widehat \g$ is a
restricted Lie algebra of Cartan type $W$ or $K$, then the
negative part of $\g$ must coincide with the negative part of
$\widehat \g$.   Under these hypotheses, $\g_1$ is the direct sum
of two irreducible $\g_0$-modules in most cases, so we can look at
the Lie algebras generated by $\g_{-1}$, $\g_0$, and  various
submodules of $\g_1$ and apply previous results to (quotients of)
them.\bi 

\begin{Lem} \label{Lem:4.36}  \  Let $\mathfrak{g} =
\bigoplus_{i=-q}^r \mathfrak{g}_i$ be a finite-dimensional graded
Lie algebra which satisfies conditions (1)-(5). Suppose that the
subalgebra $\widehat{\mathfrak{g}}$ generated by the local part
$\mathfrak{g}_{-1} \oplus \mathfrak{g}_0 \oplus \mathfrak{g}_1$ is
isomorphic to a graded Lie algebra of Cartan type $W(m;\un 1), \,
m \geq 2$, or $K(2m+1;\un 1)^{(1)}.$  Then $\mathfrak{g}_i =
\widehat{\mathfrak{g}}_i$ for $i < 0.$
\end{Lem}

\pf (a) We first assume that $\widehat{\mathfrak{g}}\not\cong
W(m;\un 1)$ if $m+1\equiv 0\mod p$. In this case the
$\mathfrak{g}_0$-module $\mathfrak{g}_1$ decomposes into the
direct sum of two $\mathfrak{g}_0$-submodules
$\mathfrak{g}_1^\sharp$ and $\mathfrak{g}_1^\dagger$. For
$\widehat{\mathfrak g}\cong W(m;\un{1})$  this
follows from Theorem \ref{Thm:2.82}.   For
$\widehat{\g}\cong K(2m+1,\un{1})^{(1)}$ we let
$\mathfrak{g}_1^\dagger$ denote the span of all
$D_K(x_i^{(1)}x_j^{(1)}x_k^{(1)})$ with $1\le i,j,k\le 2m$, and we
put $\mathfrak{g}_1^\sharp\eqdef
\{D_K(x_i^{(1)}x_{2m+1}^{(1)}) \mid 1\le i\le 2m\}$. Recall that
here we have $\g_0^{(1)}\cong \mathfrak{sp}(V)$ where $V$ is the
span of $x_i^{(1)}$ with $1\le i\le 2m$. Using Proposition
\ref{Pro:2.40} it is straightforward to see that
$\g_1^\dagger\cong S^3(V)$ and $\g_1^\sharp\cong V\cong\g_{-1}$ as
$\g_0^{(1)}$-modules.

Let $\widehat{\g}^\dagger$ and $\widehat{\g}^\sharp$ denote the
graded subalgebras of $\g$ generated by $\mathfrak{g}_{-1} \oplus
\mathfrak{g}_0 \oplus \mathfrak{g}_1^\dagger$ and
$\mathfrak{g}_{-1} \oplus \mathfrak{g}_0 \oplus
\mathfrak{g}_1^\sharp$, respectively. It follows from our general
assumption and Proposition \ref{Pro:2.40}(ii) that
$\widehat{\g}^\sharp$ is $1$-transitive and
$\widehat{\g}_{-2}={\mathfrak Z}(\widehat{\g}^\dagger)$. Applying
Theorem \ref{Thm:2.81} we now deduce that the graded Lie algebra
$\widehat{\mathfrak{g}}^\dagger/\widehat{\g}_{-2}$ is isomorphic
to one of Cartan type Lie algebras $S(m;{\un
1})^{(1)}\oplus\F{\mathfrak D}_1,$ $m\ge 3,$ or $H(2m;{\un
1})^{(2)}\oplus\F{\mathfrak D}_1,$ $m\ge 1$, while
$\widehat{\mathfrak{g}}^\sharp$ is isomorphic to a classical Lie
algebra of type $\mathrm{A}_m$ or $\mathrm{C}_{m+1}$.

\m

\noindent (b) Let us denote by $\mathfrak{g}^\sharp$ the graded
subalgebra of $\mathfrak{g}$ generated by the subspaces
$\mathfrak{g}_1^\sharp$ and $\mathfrak{g}_i$ with $i \leq 0$.   Let
$\mathcal M(\mathfrak{g}^\sharp)$ be the maximal ideal of
$\mathfrak{g}^\sharp$ contained in the negative part $\bigoplus_{i
< 0}\mathfrak{g}_i$ of $\g^\sharp$. Clearly, $\mathcal
M(\mathfrak{g}^\sharp) \subseteq \bigoplus_{i \leq
-2}\mathfrak{g}_i.$ Since $\g_1^\sharp\cong\g_{-1}^*$ as
$\g_0$-modules, the quotient algebra
$\overline{\mathfrak{g}}^\sharp \eqdef
\mathfrak{g}^\sharp/\mathcal M(\mathfrak{g}^\sharp)$ satisfies all
the conditions of Theorem \ref{Thm:3.34}. Our general assumption
implies $[[\mathfrak{g}_{-1},\mathfrak{g}_{-1}],\mathfrak{g}_{-1}]
= 0$. Since the subalgebra generated by the local part of a
Melikyan algebra is isomorphic to a Lie algebra of type
$\mathrm{G}_2$ graded according to the {\it long} simple root,
Theorem \ref{Thm:3.34} yields that
 $\overline{\mathfrak{g}}^\sharp$ is classical.  As a corollary, if
$\widehat{\mathfrak{g}}\cong W(m;\un 1)$,  then $\mathfrak{g}_{-2}
\subseteq \mathcal M(\mathfrak{g}^\sharp)$, and if
$\widehat{\mathfrak{g}}\cong K(2m+1;\un 1)^{(1)},$ then
$\mathfrak{g}_{-3} \subseteq {\mathcal M}(\mathfrak{g}^\sharp)$
and ${\mathcal M}(\mathfrak{g}^\sharp)_{-2}$ has codimension $1$
in $\g_{-2}$.

Clearly, $\mathcal M(\mathfrak{g}^\sharp)_{-2} = \{x \in
\mathfrak{g}_{-2}  \mid  [x, \mathfrak{g}_1^\sharp] = 0 \}$.   It
follows easily from this and our general assumption that
$\widehat{\mathfrak{g}}_{-2} \cap \mathcal
M(\mathfrak{g}^\sharp)_{-2} = 0$  and that
$\widehat{\mathfrak{g}}_{-2} \neq 0$ only if
$\widehat{\mathfrak{g}}\cong K(2m+1;\un 1)^{(1)}$.  Therefore,
$$\mathfrak{g}_{-2} =\, \widehat{\mathfrak{g}}_{-2} \oplus
\mathcal M(\mathfrak{g}^\sharp)_{-2}.$$

\m

\noindent (c) Let $\mathfrak{g}^\dagger$ denote the subalgebra of
$\mathfrak{g}$ generated by the subspaces $\mathfrak{g}_1^\dagger$
and $\mathfrak{g}_i, \, i \leq 0.$  Let $\mathcal
M(\mathfrak{g}^\dagger)$ be the maximal ideal of
$\mathfrak{g}^\dagger$ contained in the negative part
$\bigoplus_{i<0}\mathfrak{g}_i$ of $\g^\dagger$. Our discussion in
part~(a) shows that $\widehat{\g}_{-2}\subseteq \mathcal
M(\mathfrak{g}^\dagger)$ and the quotient algebra $\overline
{\mathfrak{g}^\dagger} \eqdef \mathfrak{g}^\dagger/\mathcal
M(\mathfrak{g}^\dagger)$ satisfies all of the conditions of Lemma
\ref{Lem:4.33}. Applying this lemma, we obtain that $\mathcal
M(\mathfrak{g}^\dagger)_{-2} = \mathfrak{g}_{-2}.$ As $\mathcal
M(\mathfrak{g}^\sharp)_{-2} = \{x \in \mathfrak{g}_{-2} \mid
[x, \mathfrak{g}_1^\sharp] = 0 \},$ we have $[\mathcal
M(\mathfrak{g}^\dagger)_{-2},\, \mathfrak{g}_1^\sharp+
\mathfrak{g}_1^\dagger] = 0.$ Since $\g$ is $1$-transitive and
$\mathfrak{g}_1^\sharp + \mathfrak{g}_1^\dagger = \mathfrak{g}_1$
by part~(a), we get $\mathcal M(\mathfrak{g}^\sharp)_{-2} = 0$.
Our final remark in part~(b) then yields $\mathfrak{g}_{-2} =
\widehat{\mathfrak{g}}_{-2}.$ If $\widehat{\mathfrak{g}}\cong
W(m;\un 1),$ this gives the result. If $\widehat
{\mathfrak{g}}\cong K(2m+1;{\un 1})^{(1)},$ then $\mathcal
M(\mathfrak{g}^\dagger)_{-3} = \mathcal
M(\mathfrak{g}^\sharp)_{-3} = \mathfrak{g}_{-3}$ by our remarks
above.

\m

\noindent (d) Suppose $\mathfrak{g}_{-3} \neq 0$.   We have
$[\mathfrak{g}_{-3},\mathfrak{g}_1^\sharp] = 0$ by part~(b), and
$[{\g}_{-3},\g_1]=[\mathfrak{g}_{-3}, \mathfrak{g}_1^\dagger] =
\widehat{\mathfrak{g}}_{-2}$ by part~(c) and 1-transitivity. Since $\dim
\widehat{\g}_{-2}=1$ and $\g$ is $1$-transitive, the
$\g_0^{(1)}$-module $\mathfrak{g}_{-3}$ embeds into
$\Hom(\mathfrak{g}^\dagger_1, \, \mathfrak{g}_{-2}) \cong
(\mathfrak{g}_1^\dagger)^*$.   Moreover, because $\widehat{\mathfrak{g}}\cong
K(2m+1;{\un 1})^{(1)}$ in the present case, we have that
$\g_0^{(1)}=\mathfrak{sp}(V)$ and $\mathfrak{g}_1^\dagger \cong
S^3(V)$ as $\g_0^{(1)}$-modules; see part~(a). Since $p > 3$,
Proposition \ref{Pro:2.848} shows that $\g_1^\dagger$ is an
irreducible (and self-dual) $\mathfrak{g}_0^{(1)}$-module. It is
generated by $D_K(x_1^{(3)})$, a $\mathfrak b^+$-primitive vector of weight
$3\varpi_1.$ Let $f^{3\varpi_1}$ be a  $\mathfrak b^+$-primitive vector of
$\mathfrak{g}_{-3},$ and let $e^{-3\varpi_1}$ and $e^{-\varpi_1}$
be $\mathfrak b^-$-primitive vectors of $\mathfrak{g}_1^\dagger$ and
$\mathfrak{g}_1^\sharp$, respectively.   Set $w \eqdef
[f^{3\varpi_1},  e^{-3\varpi_1}]$.  By the above, both
$\mathfrak{g}_{-3}$ and $\mathfrak{g}_1^\sharp$ are irreducible
over $\g_0^{(1)}$. Since $[\g_{-3},\g_1^\sharp]=0$, we have
$[f^{3\varpi_1},  e^{-\varpi_1}] = 0$, while
 the $1$-transitivity of $\g$
implies that $w$ spans $\mathfrak{g}_{-2}=\widehat{\g}_{-2}$. By
our discussion in part~(a), $\widehat{\g}^\sharp$ is a Lie algebra
of type $\mathrm{C}_{m+1}$. Because $[[e_{-2\varepsilon_1},
e_{\varepsilon_1 - \varepsilon_2}],, e_{\varepsilon_1 - \varepsilon_2}]\neq 0$
in any Lie algebra of type C$_{m+1}$ over $\F$
in which $\widehat \g_{-2} = \F e_{-2 \varepsilon_1}$, it must be that
\begin{equation}[f^{3\varpi_1},  [[e^{-3\varpi_1},
e^{-\varpi_1}],  e^{-\varpi_1}]] = \,[[w, e^{-\varpi_1}],
e^{-\varpi_1}] \in {\mathbb F}^\times e_{-\alpha},\nonumber
\end{equation}
where $\alpha$ is the highest root of $\Phi.$  It is now immediate
that the elements $f^{3\varpi_1}$ and $[[e^{-3\varpi_1},$ $\,
e^{-\varpi_1}],$ $\, e^{-\varpi_1}]$ satisfy all of the conditions
of Theorem \ref{Thm:4.7}, contrary to the fact
that $S^3(V)\not\cong V$ as $\mathfrak{sp}(V)$-modules. Therefore,
$\mathfrak{g}_{-3} = 0.$

\m

\noindent (e) We turn now to the case in which
$\widehat{\mathfrak{g}}\cong W(m;\un 1),$ where $m+1 \equiv 0 \mod
p.$ It is {\it much} more complicated. As above, we let
$\widehat{\g}^\dagger$ denote the graded subalgebra of $\g$
generated by $\mathfrak{g}_{-1} \oplus \mathfrak{g}_0 \oplus
\g_1^\dagger$. Note that in the present case $m\ge 4$. Theorem
\ref{Thm:2.82} shows that the $\g_0$-module $\g_1$ has a unique
composition series $\g_1\supset\mathfrak{g}_1^\dagger \supset
\mathfrak{g}_1^\sharp\supset 0$. Furthermore,
$\mathfrak{g}_1/\mathfrak{g}_1^\dagger \cong \mathfrak{g}_{-1}^*$,
and $\widehat{\g}^\dagger$ is isomorphic to a graded Lie algebra
$S(m;{\un 1})^{(1)}\oplus \,\F{\mathfrak D}_1;$ see Theorem
\ref{Thm:2.81}. Setting $\mathfrak{g}^\dagger$ as above and using
Lemma \ref{Lem:4.33}, we see now that $\mathfrak{g}_{-2} =
\mathcal M(\mathfrak{g}^\dagger)_{-2}.$  In other words,
$[\mathfrak{g}_{-2}, \mathfrak{g}_1^\dagger] = 0.$ By
$1$-transitivity, $\mathfrak{g}_{-2}$ is then isomorphic to a
submodule of the $\mathfrak{g}_0$-module
$\Hom(\mathfrak{g}_1/\g_1^\dagger,\,\mathfrak{g}_{-1}) \cong
\mathfrak{g}_{-1} \otimes \mathfrak{g}_{-1}$.   Since
$\mathfrak{g}_0 \cong \mathfrak{gl}(V)$ and $\mathfrak{g}_{-1}
\cong V^*,$ it follows that $\mathfrak{g}_{-2}$ is isomorphic to a
$\mathfrak{gl}(V)$-submodule of $V^* \otimes\, V^* \cong S^2(V^*)
\oplus \,\wedge^2V^*.$ In particular, $\g_{-2}$ is a completely
reducible $\g_0$-module.

\m

\noindent (f) Suppose $\g_{-2}\ne 0$ and let $\mathfrak{g}_{-2}'$
be an irreducible $\mathfrak{g}_0$-submodule of
$\mathfrak{g}_{-2}$. Clearly, either $\mathfrak{g}_{-2}' \cong
S^2(V^*)$ or $\mathfrak{g}_{-2}' \cong \wedge^2V^*.$  Due to
the irreducibility and 1-transitivity of $\g$,  we have
$[\mathfrak{g}_{-2}', \, \mathfrak{g}_1] = \mathfrak{g}_{-1}.$
Therefore, $[\mathfrak{g}_{-2}', \, [\mathfrak{g}_1, \,
\mathfrak{g}_1^\dagger]] = [\mathfrak{g}_{-1}, \,
\mathfrak{g}_1^\dagger] \neq 0$ by transitivity. As a consequence,
$[\mathfrak{g}_{-2}', \widehat{\mathfrak{g}}_2] \neq 0.$

Since $m+2 \not \equiv 0 \mod p,$ Theorem \ref{Thm:2.82} shows
that $\widehat{\mathfrak{g}}_2 = (\widehat{\mathfrak{g}}_2)^\sharp
\oplus (\widehat{\mathfrak{g}}_2)^\dagger,$ where both
$(\widehat{\mathfrak{g}}_2)^\sharp$ and
$(\widehat{\mathfrak{g}}_2)^\dagger$ are irreducible
$\mathfrak{g}_0^{(1)}$-modules. Moreover,
$(\widehat{\mathfrak{g}}_2)^\sharp$ and
$(\widehat{\mathfrak{g}}_2)^\dagger$ are generated by 
$\mathfrak b^-$-primitive vectors $e^{-2\Lambda}$ and $e^{-\mu},$ respectively, where
$\mu = 2\Lambda + \alpha$,  $\Lambda = \varpi_{m-1}$, and $\alpha
= \alpha_1 + \cdots + \alpha_{m-1} = \varpi_1 + \varpi_{m-1}$; see
Theorem \ref{Thm:2.82}, Lemma \ref{Lem:2.85}, and \cite{Bou1}. Let
$f$ be a $\mathfrak b^+$-primitive vector of $\mathfrak{g}_{-2}'.$ Since
$\mathfrak{g}_{-2}'$ is isomorphic to a submodule of $S^2(V^*)
\oplus \wedge^2V^*,$ the weight of $f$ is equal to either
$2\Lambda$ or $\varpi_{m-2}.$ Since $[\mathfrak{g}_{-2}', \,
\widehat{\mathfrak{g}}_2] \neq 0$, it follows that either $[f,
e^{-2\Lambda}] \neq 0$ or else $[f, e^{-\mu}] \neq 0$.
 Applying Theorem \ref{Thm:4.7} shows that the case $[f,
e^{-\mu}] \neq 0$ is impossible. Thus $f$ has weight $2\Lambda,$
and $[f, e^{-\mu}] = 0.$ This, in turn, yields
$\g_{-2}=\,\g_{-2}'\cong S^2(V^*)$.

At this point, it will be convenient to identify $\widehat{\g}$
with the graded Lie algebra $W(m;\un{1})$, adopt the notation of
Section~\ref{sec:2.13}, and choose $\mathfrak t$ to be the
Cartan subalgebra $\h$ from the proof of Theorem \ref{Thm:2.82}.
The above shows that $[f, e^{-2\Lambda}]$ is a nonzero element in
$\h$. Let $S^{\{2\}}$ denote the graded Lie subalgebra of $\g$
generated by $\g_{-2}$, $\mathfrak{g}_0,$ and
$(\widehat{\g}_2)^\sharp$. Let $I$ be the maximal graded ideal of
$S^{\{2\}}$ contained in the subspace $\bigoplus_{\mid i\,\mid\ge
4}\, S^{\{2\}}_i$, and set $\overline{S}^{\{2\}}\eqdef\,
S^{\{2\}}/I$. Then $\overline{S}^{\{2\}}_0\cong \g_0$ and both
$\overline{S}^{\{2\}}_{-2}$ and $\overline{S}^{\{2\}}_2$ are
irreducible $\overline{S}^{\{2\}}_0$-modules. Applying Theorem
\ref{Thm:3.24} we derive that $\overline{S}^{\{2\}}$ is isomorphic
to a classical Lie algebra with a standard grading. Our earlier
remarks show that the Dynkin diagram of $\overline{S}^{\{2\}}$
must have the following form:
\begin{equation}\label{eq:4.37}{\beginpicture \setcoordinatesystem
units <0.45cm,0.3cm> % sets scale
 \setplotarea x from -5 to 16, y from -2 to 1    % sets plot size up
 \linethickness=0.15pt                          % sets line thickness
  \put{$\circ$} at 0 0
 \put{$\circ$} at 2 0
 \put{$\cdots$} at 5 0
 \put{$\circ$} at  8 0
 \put{$\bullet$} at 10  0
 \plot .15 .1 1.85 .1 /
 \plot 2.15 .1 3.85 .1 /
 \plot 6.15 .1 7.85 .1 /
 \plot 8.15 -.2  9.85 -.2 /
 \plot 8.15 .2 9.85  .2 /
 \put {$<$} at 9 0
  \put{$\alpha_1$} at 0 -1.5
   \put{$\alpha_2$} at 2 -1.5
    \put{$\alpha_{m-1}$} at 8 -1.5
     \put{$-2\Lambda$} at 10 -1.4
      \endpicture}
\end{equation}
\noindent Therefore, it can be assumed,  after rescaling $f$ if
necessary,  that $e^{-2\Lambda}=x_{m}^{(2)}{\mathfrak D}_1$ and $[f,
e^{-2\Lambda}] = h_{2\Lambda},$ where $\Lambda(h_{2\Lambda}) = 1,$
$\alpha_{m-1}(h_{2\Lambda}) = -1,$ and $\alpha_i(h_{2\Lambda}) =
0$ for $i \leq m - 2$.

\m

\noindent (g) By Theorem \ref{Thm:2.82}  the $\g_0^{(1)}$-module
$\g_1^\sharp$ is generated by a $\mathfrak b^-$-primitive vector
$e^{-\Lambda}=x_m^{(1)}{\mathfrak D}_1$. Since
$\g_1^\sharp\subset\g_1^\dagger$ and $[\g_{-2},\g_1^\dagger]=0$ by
part~(e), we have
\begin{equation}\label{eq:sasha1}
[f,e^{-\Lambda}]=0.
\end{equation}
 Since $m + 3 \not \equiv\, 0 \mod p,$ Theorem
\ref{Thm:2.82} shows that the $\g_0^{(1)}$-module
$\widehat{\mathfrak{g}}_3$ is a direct sum of its irreducible
submodules $(\widehat{\g}_3)^\sharp$ and
$(\widehat{\mathfrak{g}}_3)^\dagger$.    Moreover, the
$\g_0^{(1)}$-module $(\widehat{\mathfrak{g}}_3)^\sharp$ is
generated by a $\mathfrak b^-$-primitive vector $e^{-3\Lambda}\eqdef
x_m^{(3)}{\mathfrak D}_1.$    In view of  \eqref{eq:degdermult}, we have
\begin{equation}\label{eq:sasha2}
[x_m^{(i)}{\mathfrak D}_1,\,x_m^{(j)}{\mathfrak D}_1]
\,=\,(j-i){i+j\choose i}
 x_m^{(i+j)}{\mathfrak D}_1\qquad\ \
\big(\forall\, i,j\in\mathbb N\big).
\end{equation}
 Since
$\di\big(x_m^{(3)}D_m-x_m^{(2)}{\mathfrak
D}_1\big)=-(m+1)x_m^{(2)}=0$ and $(\widehat{\g}_2)^\dagger$
consists of all divergence-free derivations in
$\widehat{\g}_2=W(m;\un{1})_2$, it must be that
$-x_m^{(3)}D_m+x_m^{(2)}{\mathfrak
D}_1\in(\widehat{\g}_2)^\dagger.$      Since the irreducible
$\g_0^{(1)}$-module $(\widehat{\mathfrak{g}}_3)^\dagger$ contains
$x_m^{(4)}D_1$ and
$$[x_m^{(1)}(x_1^{(1)}D_{1} - x_2^{(1)}D_{2}), \,
[x_m^{(1)}(x_1^{(1)}D_{1} - x_2^{(1)}D_{2}), \, x_m^{(2)}D_{1}]]
\,\in\, {\mathbb F}^\times x_m^{(4)}D_{1},$$ we have the inclusion
$(\widehat{\g}_3)^\dagger \subset \mathfrak{g}^\dagger$ which
shows that
\begin{equation}\label{eq:sasha3}
[f,\,x_m^{(3)}D_m-x_m^{(2)}{\mathfrak D}_1]=0,
\end{equation}
since, as we noted
in part (e),  $[\g_{-2}, \, \g_1^\dagger] = 0$.

Set $f^{\Lambda}:=D_m$, a $\mathfrak b^+$-primitive vector of the
$\g_0^{(1)}$-module $\g_{-1}$. Combining (\ref{eq:sasha1}),
(\ref{eq:sasha2}), (\ref{eq:sasha3}) with the equality
$\Lambda(h_{2\Lambda}) = 1$ we now obtain
\begin{eqnarray*}
\big[[f,f^\Lambda],\,[e^{-2\Lambda},e^{-\Lambda}]\big]&=&\big[[f,[e^{-2\Lambda},e^{-\Lambda}]],
\,f^\Lambda\big]+
\big[f,\,[f^\Lambda,\,[e^{-2\Lambda},e^{-\Lambda}]]\big]\\
&=&\big[[h_{2\Lambda},e^{-\Lambda}],f^\Lambda\big]+\big[f,\,[D_m,[x_m^{(2)}{\mathfrak
D}_1,x_m^{(1)}{\mathfrak D}_1]]\big]\\
&=&-\Lambda(h_{2\Lambda})[e^{-\Lambda},f^\Lambda]-3\big[f,\,[D_m,\,x_m^{(3)}{\mathfrak
D}_1]\big]\\
&=&[D_m,\,x_m^{(1)}{\mathfrak D}_1]-3[f,\,x_m^{(2)}{\mathfrak
D}_1+x_m^{(3)}D_m]\\
&=&{\mathfrak D}_1+x^{(1)}_mD_m- 3[f,\,2x_m^{(2)}{\mathfrak
D}_1+ x_m^{(3)}D_m-x_m^{(2)}{\mathfrak D}_1]\\
&=&{\mathfrak D}_1+x^{(1)}_mD_m-6h_{2\Lambda}.
\end{eqnarray*}
Since $\Lambda(h_{2\Lambda})=1$, we have $[{\mathfrak
D}_1+x^{(1)}_mD_m-6h_{2\Lambda},\,D_m]=-8D_m\ne 0$. Consequently,
$[f,f^\Lambda]\ne 0$ is a $\mathfrak b^+$-primitive vector of weight
$3\Lambda=3\varpi_{m-1}$ in $\g_{-3}$.

\m

\noindent (h) Let $\g_{-3}'$ be the $\g_0$-submodule of $\g_{-3}$
generated by $[f,f^\Lambda]$, and denote by $S^{\{3\}}$ the Lie
subalgebra of $\g$ generated by $\g_{-3}'$, $\mathfrak{g}_0,$ and
$(\widehat{\g}_3)^\sharp$. We give $S^{\{3\}}$ a $\mathbb
Z$-grading by setting $S_i^{\{3\}}=\,\g_{-3i}\cap S^{\{3\}}$ for
all $i\in \mathbb Z$. Let $J$ be the maximal graded ideal of
$S^{\{3\}}$ contained in the subspace $\bigoplus_{\mid i\,\mid\ge
1}\, S^{\{3\}}_i$, and set $\overline{S}^{\{3\}}\eqdef\,
S^{\{3\}}/J$. Our computation in part~(g) implies that
$[f,f^\lambda]\not\in J$. From this it is easy to deduce that the
graded Lie algebra $\overline{S}^{\{3\}}$ is irreducible,
transitive, and $1$-transitive. Obviously,
$\overline{S}^{\{3\}}_0\cong \g_0$.

Let $\Lambda'=-w_0(\Lambda)=\varpi_1$ and extend $\Lambda'$ to a
linear function on $\h$ by setting $\Lambda({\mathfrak D}_1)=1$.
Let $u^{3\Lambda'}$ denote a $\mathfrak b^+$-primitive vector of
$\overline{S}^{\{3\}}_{-1}\cong S^3(V)$, and let $v^\lambda$ be
any $\mathfrak b^-$-primitive vector of $\overline{S}^{\{3\}}_1$. Then
$[u^{3\Lambda'},v^\lambda]\ne 0$ by transitivity. If $\lambda\ne
-3\Lambda'$, then Proposition \ref{Pro:4.20} says that $\lambda$
does not vanish on $\h\cap\g_0^{(1)}$. Since $3\varpi_1 \not\in \{
\varpi_1, \varpi_{m-1}, 2\varpi_1 + \varpi_{m-1}, \varpi_1 +
2\varpi_{m-1}\},$ this contradicts Theorem \ref{Thm:4.7}. Thus,
the $\mathfrak b^-$-primitive vectors of $\overline{S}^{\{3\}}_1$
must have weight $-3\Lambda'$.
Let $M$ be an irreducible $\g_0$-submodule of
$\overline{S}^{\{3\}}_1$. Since $\overline{S}^{\{3\}}_{-1}$ is a
restricted $\g_0^{(1)}$-module, so is the $\g_0^{(1)}$-module
$M\hookrightarrow
\text{Hom}\big(\overline{S}^{\{3\}}_{-1},\,\g_0\big)$. The
preceding remark then implies that $M$ is generated by a 
$\mathfrak b^-$-primitive vector of weight $-3\Lambda'$. In this situation,  Theorem
\ref{Thm:4.12}(ii) shows that $\overline{S}^{\{3\}}_1$ is an
irreducible $\g_0^{(1)}$-module. But then $\overline{S}^{\{3\}}$
satisfies all of the conditions of Theorem \ref{Thm:3.34}. It
follows that $\overline{S}^{\{3\}}$ is a classical graded Lie
algebra whose null component is of type $\mathrm{A}_{m-1}$, $m\ge
4$. Yet $3\varpi_1$ is not listed in (\ref{eq:sasha0}),  which
exhibits all possible weights of $\mathfrak b^+$-primitive vectors in $\overline{S}_{-1}^{\{3\}}$ that
can occur in this situation. This implies that in the present case
$\g_{-2}=0$, completing the proof. \qed

\m

\section {The case that
$\boldsymbol{\g_0}$  is abelian \label{sec:4.10}} 
\m

  Here we verify that $\g_0$ can be abelian  only
in a Zassenhaus Lie algebra $W(1;\un n)$ or in  $\mathfrak{sl}_2$.  We use the fact that by transitivity, $\g_0$
must be one-dimensional, so that by irreducibility,
$\g_{-1}$ must be one-dimensional  as well.  \bi 

\begin{Lem} \label{Lem:4.39}  \  Let $\mathfrak{g} = \bigoplus_{i=-q}^r \mathfrak{g}_i$ be a
finite-dimensional graded Lie algebra which satisfies conditions
(1)-(5), and suppose further that ${\mathfrak{g}}_0^{(1)} = 0$ .
If $r
> 1,$ then $q = 1,$ and $\mathfrak{g}$ is isomorphic to a graded Lie algebra
$W(1;\un n)$ for some $n \geq 1$.    If $q > 1,$ then $r = 1,$ and $\mathfrak{g} =
\bigoplus_{i=-1}^{q}\mathfrak{g}_{-i}$ is isomorphic to a graded
Lie algebra $W(1;\un n'))$  for some $n' \geq 1$.    If $r = q = 1,$ then $\mathfrak{g}
\cong \mathfrak{sl}_2.$
\end{Lem}

\pf Since $\g_0$ is abelian, conditions (3) and (4) show that
$\g_0$ and $\g_{-1}$ are one-dimensional, say $\g_0 = \F z_0$, and
$\g_{-1} = \F z_{-1}$, and $z_0$ acts as a nonzero scalar on
$\g_{-1}$. Since the pairing $\mathfrak{g}_1 \times \g_{-1}
\rightarrow \g_0$ given by the product is nondegenerate by
transitivity,  $\g_1 \cong \mathfrak{g}_{-1}^*$.  Thus $\g_1$ is
one-dimensional also, say $\g_1 = \F z_1$.

Denote by $\mathfrak{g}^{\geq -1}$ and $\mathfrak{g}^{\leq 1}$ the
graded subalgebras of $\mathfrak{g}$ generated by all
$\mathfrak{g}_i$ with $i \ge -1$ and $i \le 1$, respectively. Both
$\mathfrak{g}^{\geq -1}$ and $\mathfrak{g}^{\leq 1}$ satisfy
conditions (1)-(5). If $\mathfrak{g} = \mathfrak{g}^{\geq -1}$ or
$\mathfrak{g} = \mathfrak{g}^{\leq 1}$, then our statement follows
from Theorem \ref{Thm:2.81}. Thus, in what follows we will assume
that $\min(r,q)> 1.$

As $[\mathfrak{g}_{-1},  \mathfrak{g}_{-1}] = [\mathfrak{g}_1,
\mathfrak{g}_1] = 0,$ Theorem \ref{Thm:2.81} shows that
$\mathfrak{g}^{\geq -1} \cong W(1;\un n)$ and $\mathfrak{g}^{\leq
1} \cong W(1;\un n')$ for some positive integers $n$ and $n'$.  In
particular, $\dim \mathfrak{g}_i = 1$ for $-q \leq i \leq r$, and
$[\mathfrak{g}_1,  \mathfrak{g}_{p-2}] = 0 = [\mathfrak{g}_{-1},
\mathfrak{g}_{-(p-2)}].$ For $-p+2 \leq i \leq p-2,$ choose $z_i
\in \mathfrak{g}_i \backslash \{0\}$. By transitivity,  $(\ad
z_{-1})^{p-3}(z_{p-2}) \neq 0$, while
$$(\ad
z_{-1})^{p-3}([\mathfrak{g}_{-(p-2)},\mathfrak{g}_{p-2}])\subseteq
{\mathbb F}(\ad z_{-1})^{p-3}(z_0) = 0,$$ since $p\ge 5$. But then
\begin{eqnarray*}
0 &=& (\ad z_{-1})^{p-3}([\mathfrak{g}_{-(p-2)},\,z_{p-2}]) \\
&=& [\mathfrak{g}_{-(p-2)},(\ad z_{-1})^{p-3}(z_{p-2})]\\
&=&
 [\mathfrak{g}_{-(p-2)}, \, \mathfrak{g}_1],
\nonumber\end{eqnarray*}

\noindent violating condition (5).  Therefore, $\min(r,q) = 1,$
proving the lemma. \qed

\m
\section {Completion of the proof
of the Main Theorem \label{sec:4.11}} 
\m

\m Suppose $\g = \bigoplus_{j=-q}^r \g_j$ is a graded Lie algebra
satisfying the hypotheses of the Main Theorem, and let
$\widehat{\g}$ denote the Lie subalgebra of $\g$ generated by
$\g_{-1}\oplus\g_0\oplus \g_1$. Note that $\widehat{\g}$ satisfies
the hypotheses of the Main Theorem. If $\g_0^{(1)} = 0$, then the
Main Theorem follows from Lemma \ref{Lem:4.39}. So assume from now
that $\g_0^{(1)} \neq 0$. Since $[[\g_{-1},\g_1],\g_1]\ne 0$ by
$1$-transitivity, Theorem \ref{Thm:1.48} shows $\mathfrak{g}_{-1}$
is a restricted $\g_0^{(1)}$-module. Since $\g_1$ embeds into
$\text{Hom}(\g_{-1},\g_0)$ by transitivity, it is a restricted
$\g_0^{(1)}$-module too.

Let $f^\Lambda$ be a $\mathfrak b^+$-primitive vector of
weight $\Lambda$ in the $\g_0$-module
$\mathfrak{g}_{-1}$, and let $e^\Gamma$ be a $\mathfrak b^-$-primitive vector of
weight $\Gamma$
in the $\g_0$-module $\g_1$. Then $[f^\Lambda,e^\Gamma]\ne 0$ by
transitivity. If $[f^\Lambda,e^\Gamma]\not\in\mathfrak t$, then it
can be assumed (after rescaling if necessary) that
$[f^\Lambda,e^\Gamma]=e_{-\alpha}$ for some $\alpha\in\Phi$. If
$\alpha$ is a positive root, then Theorem \ref{Thm:4.22} shows
that $\widehat{\g}$ is a restricted Lie algebra of Cartan type
with its natural grading.  More precisely, $\widehat{\g}$ is either
$W(m;\un 1),$ $m\ge 2$, or $S(m;\un 1)^{(1)},$ $m\ge 3$, or
$S(m;\un 1)^{(1)}\oplus \F{\mathfrak D}_1,$ $m\ge 3$, or
$H(2m;\underline{1})^{(2)},$ $m\ge 1$, or $H(2m;\un 1)^{(2)}\oplus
\F{\mathfrak D}_1,$ $m\ge 1$, or $K(2m+1;\un 1)^{(1)},$ $m\ge 1$.
In this situation,  Lemmas \ref{Lem:4.33} and \ref{Lem:4.36} yield
that $\g$ is isomorphic to a Lie algebra of Cartan type with its
natural grading.

If $\alpha$ is a negative root, then Lemma \ref{Lem:4.30} says
that $\g_1$ is an irreducible $\g_0$-module.  In this case,  we
can regard $\g$ with its reverse grading without violating the
hypotheses of the Main Theorem. Interchanging the roles of
positive and negative roots and arguing as above, we are now able
to conclude that $\g$ itself is isomorphic to a Lie algebra of
Cartan type with the reverse grading.

Thus, from now on we can assume that $[f^\Lambda,
e^\Gamma]\in\mathfrak t$ for {\it any} $\mathfrak b^-$-primitive vector
$e^\Gamma\in\g_1$. Clearly, this implies that any irreducible
$\g_0$-submodule of $\g_1$ is generated by a $\mathfrak b^-$-primitive vector
of weight $-\Lambda$. In this situation Theorem \ref{Thm:4.12}(ii)
applies showing that $\g_1$ is an irreducible $\g_0$-module.  But
then Theorem \ref{Thm:3.34} says that either $\g$ is a classical Lie
algebra with one of its standard gradings  or else $p=5$ and $\g$
is a Melikyan algebra with its natural grading or the
reverse of its natural grading. \qed

\backmatter
%-----------------------------------------------------------------------------
% Beginning of biblio.tex 8-19-05 version
%-----------------------------------------------------------------------------
%
% AMS-LaTeX 1.2 sample file for a monograph, based on amsbook.cls.
% This is a data file input by chapter.tex.
%%%%%%%%%%%%%%%%%%%%%%%%%%%%%%%%%%%%%%%%%%%%%%%%%%%%%%%%%%%%%%%%%%%%%%%%
%\setcounter{page}{124}
\bibliographystyle{amsalpha}

%-----------------------------------------------------------------------------
% End of biblio.tex
%-----------------------------------------------------------------------------
%\bibitem[J1]{J1}N.~Jacobson, {\it Lectures in Abstract Algebra II},  Van Nostrand, Toronto,
%1953.  Reprinted,  Graduate Texts in Math. \textbf{31} Springer-Verlag, New York
%1975.

%\printindex
 
\end{document}